\documentclass[11pt]{amsart}

\usepackage{amssymb,enumitem}
\usepackage[mathscr]{eucal}
\usepackage[pdftex,pagebackref,bookmarksnumbered,bookmarksopen,colorlinks,citecolor=blue,linkcolor=blue]{hyperref}

\newtheorem{theorem}{Theorem}[section]
\newtheorem{lemma}[theorem]{Lemma}
\newtheorem{corollary}[theorem]{Corollary}
\theoremstyle{definition}
\newtheorem{proposition}[theorem]{Proposition}
\theoremstyle{remark}
\newtheorem{remark}[theorem]{Remark}
\numberwithin{equation}{section}

\newlist{theoremitems}{enumerate}{1}
\setlist[theoremitems]{label=\textup{(\roman*)},ref=(\roman*),leftmargin=2.3em,topsep=.3em,itemsep=.4em,parsep=0pt}

\newcommand{\R}{\mathbb R}
\newcommand{\T}{\mathbb T}
\newcommand{\M}{\mathcal M}
\newcommand{\N}{\mathcal N}
\hypersetup{pdftitle={Two-Dimensional Local Smoothing in Noncommutative Spaces and Its Applications},pdfsubject={Harmonic analysis of the planar wave equation with semifinite tracial coefficients}}

\begin{document}

\title[Two-Dimensional Local Smoothing]{Two-Dimensional Local Smoothing in Noncommutative Spaces and Its Applications}
\author{Xudong Lai and Siyu Liu}

\address{Xudong Lai:
Institute for Advanced Study in Mathematics\\
Harbin Institute of Technology\\
Harbin
150001\\
China;
Zhengzhou Advanced Research Institute\\
Harbin Institute of Technology\\
Zhengzhou
450000\\
China}
\email{xudonglai@hit.edu.cn\ xudonglai@mail.bnu.edu.cn}

\address{Siyu Liu:
Institute for Advanced Study in Mathematics\\
Harbin Institute of Technology\\
Harbin
150001\\
China}
\email{syliu@stu.hit.edu.cn\ siyulliu@icloud.com}


\begin{abstract}
We prove local smoothing for planar half-waves with coefficients in a semifinite tracial von Neumann algebra \(\M\). Let \(U(t)=e^{it|D|}\) for \(t\in\R\). For \(N\ge2\) and \(P_N=\psi(|D|/N)\), where \(\psi\in C_c^\infty((1/2,2))\), we obtain the estimate
\begin{equation*}
 \|U(t)P_N f\|_{L^4([1,2]\times\R^2;L^4(\M))}
 \le C_\psi\sqrt{\log(2+N)}\,\|f\|_{L^4(\R^2;L^4(\M))},
\end{equation*}
uniformly in \(\M\) and its trace. Summing over frequency annuli gives local smoothing at every Sobolev order \(s>0\) for \(p=4\). The difficulty in the energy comparison is a determinant contribution that can have either sign for operator-valued waves. We cancel it by pairing the Dirichlet energies for \(\M\) and its opposite algebra. Interpolation before summation gives spacetime \(L^p\) bounds for \(2\le p\le\infty\) and \(s>\sigma(p)=\max\{0,1/2-2/p\}\). For \(2<p<\infty\) we also obtain logarithmically weighted bounds at \(s=\sigma(p)\). The finite-\(p\) estimates transfer to quantum tori and generalized Moyal planes, answering the question of Hong, the first author and Wang about the Sobolev range \(s>\sigma(p)\) on generalized Moyal planes for \(2<p<\infty\).

For \(1<p<\infty\) and complex orders \(\alpha\) with \(\operatorname{Re}\alpha>\max\{2/p-1,\,1/p-1/2,\,-1/p\}\), we prove a maximal inequality over all radii for generalized circular means on the plane, quantum tori and generalized Moyal planes. For an action \(\beta\) of \(\R^2\) on \(\M\) by normal trace-preserving \(*\)-automorphisms with norm-continuous \(L^p\) orbits, the orbit means \(A_{t,\beta}^\alpha x\) converge bilaterally almost uniformly to \(\pi x/\Gamma(1+\alpha)\) as \(t\downarrow0\) for every \(x\in L^p(\M)\). We also obtain estimates for truncated cone multipliers on three-dimensional Euclidean and quantum Euclidean spaces, and full-radius row and column bounds for Stein's square function on the plane, quantum tori and generalized Moyal planes.
\end{abstract}

\maketitle
\newpage
\tableofcontents

\section{Introduction}

A high-frequency wave on the plane can focus near a point at a prescribed time.
For \(2\le p<\infty\), the sharp fixed-time \(L^p\) estimate for the scalar half-wave at nonzero time requires \(s_p=1/2-1/p\) derivatives~\cite{sss}.
Sogge's local smoothing problem asks how much of this loss can be recovered by integrating in time~\cite{sogge}.
We ask the same question for operator-valued waves, with a bound uniform in the coefficient algebra.

Let \(I=[1,2]\). For \(t\in\R\) write \(U(t)=e^{it|D|}\) for
the unitary half-wave group on \(L^2(\R^2)\). The multipliers \(|D|\) and
\(J=(1-\Delta)^{1/2}\) map \(H^1(\R^2)\) to \(L^2(\R^2)\) and have
symbols \(|\xi|\) and \(\langle\xi\rangle=(1+|\xi|^2)^{1/2}\),
respectively, for \(\xi\in\R^2\). In the scalar case the cone square
function estimate of Guth, Wang and Zhang gives the planar local
smoothing conjecture \cite{gwz}. For \(2\le p<\infty\) and every real
\(s>\sigma(p)\), it gives
\[
 \|U(t)f\|_{L^p(I\times\R^2)}
    \lesssim_{p,s}\|J^s f\|_{L^p(\R^2)}
 \qquad\bigl(J^sf\in L^p(\R^2)\bigr).
\]
The critical order and the gain over the fixed-time estimate are
\begin{equation}\label{intro:exponents}
 \sigma(p)=\max\{0,\tfrac12-\tfrac2p\},\qquad
 g(p)=s_p-\sigma(p)=\min\{\tfrac12-\tfrac1p,\tfrac1p\}.
\end{equation}
Here \(g(p)\) is the supremum of the gains \(s_p-s\) over admissible Sobolev orders \(s\). It reaches \(1/4\) at \(p=4\), where the two branches of \(\sigma(p)\) meet.

For matrices of size \(d\in\mathbb N\), we ask in particular whether the scalar range at \(p=4\) holds in the Schatten space \(S_d^4\) with a constant independent of \(d\). More generally, let \(\M\) be a von Neumann algebra with a normal faithful semifinite trace \(\operatorname{tr}_{\M}\). For a measure space \((Y,dy)\), write \(\N_Y=L^\infty(Y)\bar\otimes\M\) with the product trace. We use Lebesgue measure on Euclidean domains and probability Haar measure on tori. In particular, for \(1\le p<\infty\) the spacetime norm is given by
\begin{equation*}
 \|u\|_{L^p(\N_{I\times\R^2})}^p
 =\int_I\int_{\R^2}\operatorname{tr}_{\M}(|u(t,x)|^p)\,dx\,dt .
\end{equation*}

Hong, the first author and Wang obtained spacetime bounds for operator-valued waves on \(\R^2\) and for waves on generalized Moyal planes. For \(2<p<\infty\) their estimates require Sobolev orders \(s>s_p-g(p)/2\) for the initial displacement and \(s-1\) for the initial velocity~\cite[Theorems~1.2--1.3]{hlw}. They asked whether every \(s>\sigma(p)=s_p-g(p)\) suffices on generalized Moyal planes~\cite[discussion after Theorem~1.2]{hlw}. We prove that it does. The proof starts with the Euclidean half-wave. On a frequency annulus \(|\xi|\asymp N\), their \(L^4\) estimate has a factor of \(C_\varepsilon N^{1/8+\varepsilon}\) for each \(\varepsilon>0\)~\cite[Proposition~3.5]{hlw}. We obtain \(C\sqrt{\log(2+N)}\) by pairing the energies for \(\M\) and its opposite algebra. This gives every \(s>0\) at \(p=4\). Interpolation and summation over annuli give every \(s>\sigma(p)\) for \(2<p<\infty\), and transference gives the same Sobolev range on quantum tori and generalized Moyal planes.

\subsection{Main results and applications}
We first state the Euclidean local smoothing estimates in (A) and transfer them to quantum waves in (B). The circular maximal estimates in (C) give the ergodic bounds and convergence in (D). Local smoothing also yields bounds for truncated cone multipliers in (E) and Stein's square functions in (F).

\medskip
\paragraph{\textbf{(A) Euclidean local smoothing.}}
The multipliers below act on the scalar factor of elementary tensors \(f\otimes a\), where \(f\in\mathscr S(\R^2)\) and \(a\in\M\) is bounded with \(\operatorname{tr}_{\M}(s(|a|))<\infty\). Here \(s(|a|)=1_{(0,\infty)}(|a|)\) is the right support projection of \(a\). The estimates determine the extensions of these multipliers.

\begin{theorem}\label{thm:four}
For every \(\psi\in C_c^\infty((1/2,2))\) there is \(C_\psi\) such that for every semifinite tracial coefficient algebra and every real \(N\ge2\), the spatial cutoff \(P_N=\psi(|D|/N)\) satisfies
\begin{equation}\label{intro:l4}
 \|U(t)P_NF\|_{L^4(\N_{I\times\R^2})}
 \le C_\psi\sqrt{\log(2+N)}\,\|F\|_{L^4(\N_{\R^2})}
\end{equation}
for all \(F\in L^4(\N_{\R^2})\).
\end{theorem}

We interpolate~\eqref{intro:l4} with Plancherel's identity at \(p=2\) and with the spatial kernel bound in \(L^1(\R^2)\) at \(p=\infty\). On one frequency annulus the resulting powers of \(N\) and \(\log(2+N)\) are \(\sigma(p)\) and \(\beta(p)=2g(p)\), respectively, with \(1/\infty=0\) in~\eqref{intro:exponents}. For \(s>\sigma(p)\) the annular estimates are summable for \(U(t)J^{-s}\). For the estimate at \(s=\sigma(p)\), define the positive weights on \([0,\infty)\) by
\[
 \ell_1(r)=\log(2+r),\qquad \ell_2(r)=\log(2+\ell_1(r)),\qquad
 L_{B,C}(r)=\ell_1(r)^B\ell_2(r)^C \qquad (B,C\in\R)\,.
\]

For \(p<\infty\), each Sobolev space below is the completion of the test tensor span for the model in question in the stated norm.

\begin{theorem}\label{thm:lp}
For any semifinite tracial coefficient algebra
\(\M\), the following estimates hold with constants independent of \(\M\)
and its trace.
\begin{theoremitems}
\item\label{thm:lp-annular}
For \(2\le p\le\infty\), \(N\ge2\), the cutoff \(P_N\) of
Theorem~\ref{thm:four}, and \(F\in L^p(\N_{\R^2})\),
\begin{equation}\label{intro:annular}
 \|U(t)P_NF\|_{L^p(\N_{I\times\R^2})}
 \le C_{p,\psi}N^{\sigma(p)}[\log(2+N)]^{\beta(p)}
                 \|F\|_{L^p(\N_{\R^2})}.
\end{equation}
\item\label{thm:lp-sob}
Let \(2\le p\le\infty\) and \(s>\sigma(p)\), or let \(p=2\) and \(s=0\).
For every \(F\) with \(J^sF\in L^p(\N_{\R^2})\),
\begin{equation}\label{intro:sob}
 \|U(t)F\|_{L^p(\N_{I\times\R^2})}
        \le C_{p,s}\|J^sF\|_{L^p(\N_{\R^2})}.
\end{equation}
For \(p=\infty\), \(G\in\N_{\R^2}\), and \(t\in I\), we define
\begin{equation*}
 U(t)(J^{-s}G)=K_t^{(s)}*G,\qquad
 K_t^{(s)}=\bigl(e^{it|\xi|}\langle\xi\rangle^{-s}\bigr)^\vee\,.
\end{equation*}
The convolution is defined by duality with the integrable kernel constructed in Section~\ref{lp:sum-annuli}.
\item\label{thm:lp-weight}
Let \(2<p<\infty\). Suppose that \(B,C\in\R\) satisfy
\begin{equation}\label{intro:logrange}
  B>\beta(p)+s_p \qquad\text{or}\qquad
  B=\beta(p)+s_p\ \ \text{and}\ \ C>s_p .
\end{equation}
Then
\begin{equation}\label{intro:weighted}
 \|U(t)F\|_{L^p(\N_{I\times\R^2})}\le C_{p,B,C}
 \|J^{\sigma(p)}L_{B,C}(J)F\|_{L^p(\N_{\R^2})}\,.
\end{equation}
\end{theoremitems}
\end{theorem}

The estimates of Theorems~\ref{thm:four} and~\ref{thm:lp} also hold for \(U(-t)\) on \(I\). By~\eqref{intro:half-wave-adjoint}, taking adjoints sends the symbol \(e^{it|\xi|}\psi(|\xi|/N)\) to \(e^{-it|\xi|}\overline{\psi(|\xi|/N)}\) and preserves the spatial and spacetime \(L^p\) norms.

On one annulus, the optimal logarithmic exponent for the class \(\psi\in C_c^\infty((1/2,2))\) at \(p=4\) lies between \(1/4\) and \(1/2\). The upper bound is~\eqref{intro:l4}. For the lower bound, take \(\psi=1\) on \([5/4,7/4]\). The scalar Schwartz functions of~\cite[Proposition~8.2]{openai} have Fourier support in \(5N/4<|\xi|<7N/4\) and give an annular operator norm of at least \(c_p(\log N)^{1/2-1/p}\) for \(2<p\le4\).

Proposition~\ref{lp:sob-necessary} shows that the Sobolev orders in Theorem~\ref{thm:lp}\ref{thm:lp-sob} are sharp for \(2\le p\le4\) and for \(p=\infty\). For \(4<p<\infty\), the scalar focusing example of~\cite[Proposition~8.1]{openai} forces \(s\ge\sigma(p)=1/2-2/p\). Whether the unweighted estimate~\eqref{intro:sob} holds at equality for general coefficient algebras remains open. At \(p=4\), where~\eqref{intro:sob} requires \(s>0\), the weighted estimate~\eqref{intro:weighted} holds at \(s=0\) with \(B=3/4\) and \(C>1/4\). For all \(2<p<\infty\), the weights in~\eqref{intro:logrange} agree with the scalar estimate of~\cite[Corollary~7.1]{openai}.

\medskip
\paragraph{\textbf{(B) Wave estimates on quantum tori and generalized Moyal planes.}}
For finite \(p\), Proposition~\ref{q:spacetime} transfers the Euclidean spacetime estimates through the gauge action on quantum tori~\cite{cxy} and the translation action on generalized Moyal planes~\cite{hlw}. At \(p=\infty\) the annular estimate follows from the spatial kernel bound in \(L^1(\R^2)\). Related Fourier estimates on quantum spaces include the first author's Bochner--Riesz bounds~\cite{lai} and the restriction estimates of Hong, the first author and Wang~\cite{restriction}.

Let \(\mathcal Q_\Theta\) be the two-dimensional quantum torus of period \(2\pi\) or the generalized Moyal plane, for a real skew-symmetric \(2\times2\) matrix \(\Theta\). For any semifinite tracial algebra \(\mathcal C\), put \(\mathcal D=\mathcal C\bar\otimes\mathcal Q_\Theta\) and \(\mathcal D_I=L^\infty(I)\bar\otimes\mathcal D\), with product traces. The quantum multipliers \(|D_\Theta|\) and \(J_\Theta\) have symbols \(|\xi|\) and \(\langle\xi\rangle\), with \(\xi\in\mathbb Z^2\) on the torus and \(\xi\in\R^2\) on the plane. For \(t\in\R\) write \(U_\Theta(t)=e^{it|D_\Theta|}\). We first take \(x,y\) in the linear span of elementary tensors \(c\otimes q\), where \(c\in\mathcal C\) is bounded and its support projections have finite trace. The factor \(q\) is a quantum trigonometric polynomial on the torus or a quantum Schwartz element on the plane (Section~\ref{q:models}). The solution with initial displacement \(x\) and initial velocity \(y\) is
\[
 u_\Theta(t;x,y)=\cos(t|D_\Theta|)x+
                \frac{\sin(t|D_\Theta|)}{|D_\Theta|}y,
 \qquad t\in\R,
\]
where the sine quotient takes the value \(t\) at zero frequency.

\begin{theorem}\label{thm:quantum}
For both quantum models, the following bounds hold uniformly in \(\Theta\), \(\mathcal C\) and its trace.
\begin{theoremitems}
\item\label{thm:q-wave} If \(2\le p<\infty\) and \(s>\sigma(p)\), or if \(p=2,\ s=0\), then
\begin{equation}\label{intro:quantum}
 \|u_\Theta(t;x,y)\|_{L^p(\mathcal D_I)}
 \le C_{p,s}\big(\|J_\Theta^s x\|_{L^p(\mathcal D)}
                 +\|J_\Theta^{s-1}y\|_{L^p(\mathcal D)}\big).
\end{equation}
\item\label{thm:q-annular} For \(2\le p\le\infty\), every real \(N\ge2\), every \(\psi\in C_c^\infty((1/2,2))\), and every \(x\in L^p(\mathcal D)\),
\begin{equation*}
 \|U_\Theta(t)\psi(|D_\Theta|/N)x\|_{L^p(\mathcal D_I)}
 \le C_{p,\psi}N^{\sigma(p)}[\log(2+N)]^{\beta(p)}\|x\|_{L^p(\mathcal D)}\,.
\end{equation*}
\item\label{thm:q-weight} Let \(2<p<\infty\) and let \(B,C\in\R\) satisfy \eqref{intro:logrange}. Then the half-wave satisfies
\begin{equation*}
 \|U_\Theta(t)x\|_{L^p(\mathcal D_I)}
 \le C_{p,B,C}\|J_\Theta^{\sigma(p)}L_{B,C}(J_\Theta)x\|_{L^p(\mathcal D)}\,.
\end{equation*}
The wave with initial displacement \(x\) and initial velocity \(y\) satisfies
\begin{equation*}
 \|u_\Theta(t;x,y)\|_{L^p(\mathcal D_I)}
 \le C_{p,B,C}\bigl(
     \|J_\Theta^{\sigma(p)}L_{B,C}(J_\Theta)x\|_{L^p(\mathcal D)}
     +\|J_\Theta^{\sigma(p)-1}L_{B,C}(J_\Theta)y\|_{L^p(\mathcal D)}
     \bigr)\,.
\end{equation*}
\end{theoremitems}
\end{theorem}

The annular and weighted half-wave estimates also hold for \(U_\Theta(-t)\) on \(I\) by the adjoint identity for orbit multipliers in Section~\ref{q:secspacetime}.

\medskip
\paragraph{\textbf{(C) Maximal bounds for generalized circular means.}}
The connection between local smoothing and circular maximal estimates goes back to Mockenhaupt, Seeger and Sogge \cite{mss}. For \(2\le p<\infty\), the gain \(g(p)\) in local smoothing gives a maximal estimate for every complex order \(\alpha\) with \(\Re\alpha>-g(p)\). For \(t>0\) and \(\alpha\in\mathbb C\), let \(A_t^\alpha:\mathscr S(\R^2)\to\mathscr S(\R^2)\) have Fourier symbol \(m_\alpha(t|\xi|)\), acting on the scalar factor of test tensors, where
\[
 m_\alpha(r)=\pi 2^\alpha r^{-\alpha}J_\alpha(r)\quad(r>0),
 \qquad m_\alpha(0)=\pi/\Gamma(1+\alpha).
\]
Here \(J_\alpha\) is the Bessel function of the first kind and \(r^{-\alpha}=e^{-\alpha\log r}\). Thus \(A_t^0\) is \(\pi\) times the probability average on the circle of radius \(t\). We measure the family in the noncommutative maximal norm \(L^p(\N_{\R^2};\ell^\infty)\), whose factorization definition is given in \eqref{max:norm}. Put
\begin{equation}\label{intro:max-threshold}
 a(p)=\begin{cases}
  2/p-1,&1<p\le2,\\
  1/p-1/2,&2<p\le4,\\
  -1/p,&4<p<\infty .
 \end{cases}
\end{equation}

In particular, \(a(p)=-g(p)\) for \(2\le p<\infty\).

\begin{theorem}\label{thm:maximal}
Let \(1<p<\infty\) and \(\alpha\in\mathbb C\) with \(\operatorname{Re}\alpha>a(p)\).
\begin{theoremitems}
\item For every \(F\in L^p(\N_{\R^2})\),
\begin{equation}\label{intro:maximal}
 \|(A_t^\alpha F)_{t>0}\|_{L^p(\N_{\R^2};\ell^\infty)}
 \le C_{p,\alpha}\|F\|_{L^p(\N_{\R^2})},
\end{equation}
uniformly over \(\M\) and its trace.
\item\label{thm:q-max} For either quantum model \(\mathcal D=\mathcal C\bar\otimes\mathcal Q_\Theta\) of Theorem~\ref{thm:quantum} and every \(x\in L^p(\mathcal D)\),
\begin{equation*}
 \big\|(m_\alpha(t|D_\Theta|)x)_{t>0}\big\|_{L^p(\mathcal D;\ell^\infty)}
 \le C_{p,\alpha}\|x\|_{L^p(\mathcal D)}\,,
\end{equation*}
uniformly in \(\Theta\), \(\mathcal C\) and its trace.
\end{theoremitems}
\end{theorem}

For \(2<p<\infty\), Theorem~\ref{thm:maximal} improves the range \(\Re\alpha>-g(p)/2\) of Hong, the first author and Wang \cite[Theorem 1.4]{hlw} to \(\Re\alpha>-g(p)\). This agrees with the strict scalar range of Liu, Shen, Song and Yan, who also proved the necessary condition \(\Re\alpha\ge-g(p)\) \cite[Theorem 1.1]{lssy}. In the Bessel frequency decomposition of~\cite{hlw}, Theorem~\ref{thm:lp} and a bound for the radius derivative give positive majorants for each band. The band norms are summable in the improved range. Interpolation with positive ball averages gives \(1<p<2\).

For every \(2<p\le4\) and every fixed \(b\in\R\), the scalar maximal operator at \(\alpha=1/p-1/2+ib\) fails to map \(L^p(\R^2)\) to weak \(L^p(\R^2)\), already for radii in \([1,2]\) (Proposition~\ref{max:endpoint}). For \(4<p<\infty\), the full-radius endpoint \(\Re\alpha=-1/p\) on the whole of \(L^p(\R^2)\) remains unresolved here. Quantitative bounds for the growth of the maximal norm near the critical line are given in Corollaries~\ref{max:growth} and~\ref{max:lowergrowth}.

\medskip
\paragraph{\textbf{(D) Ergodic averages and bilateral almost uniform convergence.}}
The circular maximal bound~\eqref{intro:maximal} also transfers to actions of \(\R^2\) by normal trace-preserving \(*\)-automorphisms of a semifinite tracial algebra \(\mathcal D\). Fix \(1<p<\infty\), assume that the action \(\beta\) has norm-continuous \(L^p(\mathcal D)\) orbits, and let \(\alpha\in\mathbb C\) satisfy \(\Re\alpha>a(p)\). The orbit averages \(A_{t,\beta}^\alpha\), defined by the compactly supported kernels of \(A_t^\alpha\) on smooth orbits and extended to \(L^p(\mathcal D)\) in Theorem~\ref{q:action}, satisfy the same maximal bound uniformly in \(\mathcal D\) and \(\beta\). For every \(x\in L^p(\mathcal D)\),
\[
 A_{t,\beta}^\alpha x\longrightarrow\frac{\pi}{\Gamma(1+\alpha)}x
 \qquad\text{bilaterally almost uniformly as }t\downarrow0
\]
by Corollary~\ref{q:bau}. The distributional definition of the orbit averages and the projection formulation of this convergence are given in Theorem~\ref{q:action} and Section~\ref{q:secconvergence}, respectively.

\medskip
\paragraph{\textbf{(E) Truncated cone multipliers.}}
Let \((\xi,\rho)\in\R^2\times\R\), let \(\chi\in C_c^\infty((0,\infty))\), and let \(1<p<\infty\). For each real \(\lambda>\max\{0,2|1/p-1/2|-1/2\}\), the truncated cone multiplier with symbol
\begin{equation*}
 \chi(\rho)(1-|\xi|^2/\rho^2)_+^\lambda\qquad(\rho>0)
\end{equation*}
and value zero for \(\rho\le0\) is bounded on \(L^p(\N_{\R^3})\). The same bound holds on three-dimensional quantum Euclidean spaces, uniformly in the deformation and the tracial coefficients (Theorem~\ref{cone:theorem} and Corollary~\ref{cone:quantum}). To obtain it, we decompose the symbol near \(\rho=|\xi|\) into shells of thickness \(0<\delta\le1/2\). Averaging the half-wave against a time kernel of width \(\delta^{-1}\) gives an \(L^4\) bound of \(C[1+\log(1/\delta)]^{1/2}\) for each normalized smooth shell in Proposition~\ref{cone:shell}. Interpolation gives the bounds needed to sum the shells.

\medskip
\paragraph{\textbf{(F) Stein's row and column square functions.}}
For real \(\nu>1/2\) and \(r>0\), let \(B_r^\nu:L^2(\N_{\R^2})\to L^2(\N_{\R^2})\) be the Bochner--Riesz multiplier with symbol \((1-|\xi|^2/r^2)_+^\nu\), and define \(D_r^\nu=r\partial_r B_r^\nu\) on test tensors, with the derivative interpreted distributionally in \(r\) (Section~\ref{sec:stein}).

Local smoothing controls a single boundary shell on one dyadic interval of radii. At \(p=4\), summing the square integrals over all such intervals introduces traces of products from different frequency scales. For separated scales, Lemma~\ref{stein:four} controls these products by estimating a spatial increment of the lower-frequency square integral. The decay in the separation of the scales makes the resulting errors summable.

\begin{theorem}\label{thm:stein}
Let \(2\le p<\infty\) and let \(\nu>\max\{1/2,1-2/p\}\) be real. For every \(F\in L^p(\N_{\R^2})\),
\begin{multline*}
 \max\!\left\{
 \left\|\left(\int_0^\infty (D_r^\nu F)^*D_r^\nu F\,\frac{dr}{r}\right)^{1/2}\right\|_{L^p(\N_{\R^2})},\,
 \left\|\left(\int_0^\infty D_r^\nu F(D_r^\nu F)^*\,\frac{dr}{r}\right)^{1/2}\right\|_{L^p(\N_{\R^2})}
 \right\}\\
 \le C_{p,\nu}\|F\|_{L^p(\N_{\R^2})}\,.
\end{multline*}
The family \((D_r^\nu F)_{r>0}\), initially defined on \(L^2(\N_{\R^2})\cap L^p(\N_{\R^2})\) as in Section~\ref{sec:stein}, extends uniquely by completion in both displayed square-function norms (Theorem~\ref{stein:theorem}). The constant is independent of \(\M\) and its trace.
\end{theorem}

These are the column and row versions of Stein's square function in the planar scalar range of Carbery \cite{carbery,gow}. Corollary~\ref{stein:quantum} gives the same estimates on quantum tori and generalized Moyal planes, uniformly in the deformation and the tracial coefficient algebra.

\subsection{The local smoothing argument}
The main difficulty in proving~\eqref{intro:l4} is a sign in the scalar energy argument of~\cite[Section~2, Appendix~D and Lemma~6.1]{openai}. For a wave supported in one frequency cap, the time density of a conserved current on an enlarged spacetime supplies the source for an auxiliary Dirichlet problem. Green's formula for its energy contains a tangent-plane determinant. In the scalar argument the spatially averaged current has nonnegative determinant, which gives the inequality needed to compare energies at different radii.

With coefficients in \(\M\), the fourth power of the pointwise norm of the wave \(u\) is \(\operatorname{tr}_{\M}((u^*u)^2)\). To retain this order of multiplication, we polarize the unaveraged current and take trace coordinates of the quadratic products of Fourier coefficients. The resulting currents still satisfy the conservation laws. Their tangent-plane determinants, however, can have a negative sum. To see where the sign enters, consider a finite Fourier sum supported in one cap. Let \(X_0\in L^4(\M)\) denote its Gaussian-localized sum and \(X_1\in L^4(\M)\) the sum with each summand further weighted by its real frequency slope, as in Section~\ref{cur:sec}. By Proposition~\ref{cur:cancel}, the sum of the determinant contributions over the trace coordinates is a nonnegative geometric factor times
\begin{equation*}
 \mathfrak D_{\M}
 :=\|X_0X_1^*\|_{L^2(\M)}^2-\|X_1^*X_0\|_{L^2(\M)}^2 .
\end{equation*}

In a commutative coefficient algebra, the two norms agree. For \(2\times2\) matrix coefficients, even a wave with just two frequencies can have \(\mathfrak D_{\M}<0\), as the example following Proposition~\ref{cur:cancel} shows.

The opposite algebra \(\M^{\mathrm{op}}\), with multiplication \(a^{\mathrm{op}}b^{\mathrm{op}}=(ba)^{\mathrm{op}}\), reverses the products while preserving the Fourier phases and the fourth moment. For the same wave viewed in \(\M^{\mathrm{op}}\), it gives
\begin{equation*}
 \mathfrak D_{\M^{\mathrm{op}}}=-\mathfrak D_{\M} .
\end{equation*}
We first sum over the trace coordinates within each algebra, then add the two Dirichlet energies. The determinant contributions cancel pointwise on the lifted spacetime. The remaining interior squares enter with a nonpositive sign in the energy Laplacian, restoring the comparison inequality (Lemma~\ref{en:green}).

Gaussian localization relates the paired energy to the wave. Near the center of the auxiliary disk, the localized time density controls the fourth moment (Lemma~\ref{en:recovery}). At the boundary, the flux decays exponentially. For a sufficiently narrow frequency cap, this decay compensates for the loss in recovering the fourth moment. Comparing energies first outward to a fixed radius and then inward to a circle near the inner boundary costs one factor of \(\log N\) in each direction. Taking the fourth root of the resulting \((\log N)^2\) bound gives~\eqref{intro:l4}.

\medskip

\paragraph{\textbf{Organization.}}
Section~\ref{sec:prelim} reduces~\eqref{intro:l4} to the periodic estimate for one frequency cap in Proposition~\ref{en:cap}. Section~\ref{cur:sec} proves the determinant cancellation, and Section~\ref{en:section} uses it in the two radial energy comparisons that complete the proof of Theorem~\ref{thm:four}. Section~\ref{sec:lp} derives the Sobolev estimates and transfers them to quantum waves. Section~\ref{max:sec} uses local smoothing to control circular means and their limits under trace-preserving actions. The long-time wave estimate of Lemma~\ref{cone:long-wave} is used in Section~\ref{sec:cone} to estimate thin cone shells and in Section~\ref{sec:stein} to bound Stein square functions on a single frequency annulus.

\section{Preliminaries and localization}\label{sec:prelim}

We record the trace identities, approximation arguments and scalar kernel bounds used in the proof of Theorem~\ref{thm:four}. Section~\ref{pre:periodic-cap} states the periodic fourth-moment estimate for one frequency cap and derives the Euclidean estimate~\eqref{intro:l4} from it.

\subsection{Tracial \texorpdfstring{\(L^p\)}{Lp} spaces and approximation}\label{pre:tracial}

Throughout, \(\operatorname{tr}_{\M}:\M_+\to[0,\infty]\) denotes the normal faithful semifinite trace on \(\M\). We use the same notation for its extension by spectral integration to positive measurable operators affiliated with \(\M\) and for the resulting linear trace \(\operatorname{tr}_{\M}:L^1(\M)\to\mathbb C\).
For \(1\le p<\infty\), \(L^p(\M)\) consists of the affiliated measurable
operators \(a\) for which \(\|a\|_{L^p(\M)}^p=\operatorname{tr}_{\M}(|a|^p)<\infty\).
For measurable operators we take the closure of the product on its natural domain. We write \(|a|=(a^*a)^{1/2}\), and put \(L^\infty(\M)=\M\) with its operator norm. The trace Hilbert space \(L^2(\M)\) has inner product \(\langle a,b\rangle=\operatorname{tr}_{\M}(a^*b)\), linear in \(b\). On \(M_d(\mathbb C)\), \(d\in\mathbb N\), with ordinary trace, the \(L^2\) norm is the Hilbert--Schmidt norm. The fourth power of the \(L^4(\M)\) norm, which enters the wave estimate, is \(\operatorname{tr}_{\M}((a^*a)^2)\). 

If \(Y\) is a Euclidean domain or a flat torus with measure \(dy\),
write \(\N_Y=L^\infty(Y)\,\bar\otimes\,\M\) with the product trace.
For \(1\le p<\infty\) we identify \(F\in L^p(\N_Y)\) with a Bochner measurable function \(F:Y\to L^p(\M)\), up to equality almost everywhere, by
\begin{equation}\label{pre:fubini}
 L^p(\N_Y)=L^p(Y;L^p(\M)),\qquad
 \|F\|_{L^p(\N_Y)}^p=\int_Y\operatorname{tr}_{\M}(|F(y)|^p)\,dy.
\end{equation}
For the Euclidean and circle cases of~\eqref{pre:fubini}, see~\cite[Introduction]{parcet} and~\cite[p.~1483]{px}, respectively.
Thus the fourth power of the norm on the left of \eqref{intro:l4} integrates
\(\|(U(t)P_NF)(x)\|_{L^4(\M)}^4\) against \(dx\,dt\). In \eqref{intro:maximal}
the radii instead index a family in a single spatial algebra.

Let \((\mathcal A_i,\tau_i)\), \(i\in\{1,2\}\), be semifinite tracial von Neumann algebras and let \(1\le p<\infty\). Le Merdy, Ricard and Roydor give the norm density~\cite[Section~5, p.~864]{lrr}
\begin{equation*}
 \overline{L^p(\mathcal A_1,\tau_1)\odot L^p(\mathcal A_2,\tau_2)}^{\,\|\cdot\|_{L^p(\mathcal A_1\bar\otimes\mathcal A_2,\tau_1\bar\otimes\tau_2)}}
 = L^p(\mathcal A_1\bar\otimes\mathcal A_2,\tau_1\bar\otimes\tau_2)\,,
\end{equation*}
where \(\odot\) denotes the algebraic tensor product. This applies to the tensor products in Section~\ref{q:models}.

For \(a\in L^p(\M)\) and \(n\in\mathbb N\), put
\(a_n=a\,1_{[1/n,n]}(|a|)\). Normality of the trace, applied to the
spectral integral of \(|a|^p\), gives \(\|a-a_n\|_{L^p(\M)}\to0\).
Moreover \(\|a_n\|_{L^\infty(\M)}\le n\), and its right support has trace
at most \(\operatorname{tr}_{\M}(1_{[1/n,\infty)}(|a|))
\le n^p\|a\|_{L^p(\M)}^p\). Its left support has the same trace by polar
decomposition.

The product trace gives \(\|x_1\otimes x_2\|_{L^p(\mathcal A_1\bar\otimes\mathcal A_2,\tau_1\bar\otimes\tau_2)}=\|x_1\|_{L^p(\mathcal A_1,\tau_1)}\|x_2\|_{L^p(\mathcal A_2,\tau_2)}\) for \(x_i\in L^p(\mathcal A_i,\tau_i)\), \(i=1,2\), so we may approximate the factors of each finite tensor sum in their respective \(L^p\) norms. In \(L^p(\N_Y)\) we therefore use the dense subspace of finite sums \(\sum_{i=1}^m f_i\otimes a_i\), \(m\in\mathbb N\), with \(f_i\in C_c^\infty(Y)\) for a Euclidean domain or \(f_i\in C^\infty(Y)\) for a flat torus, and \(a_i\in\M\) satisfying \(\operatorname{tr}_{\M}(1_{(0,\infty)}(|a_i|))<\infty\). The scalar factors are dense in \(L^p(Y)\) and the coefficient factors are dense by the spectral truncation above. For \(Y=\R^d\), \(d\in\mathbb N\), we also use Schwartz scalar factors, since \(\mathscr S(\R^d)\) contains \(C_c^\infty(\R^d)\).

For \(1\le q,r,s\le\infty\) with \(1/q=1/r+1/s\), multiplication maps \(L^r(\M)\times L^s(\M)\) into \(L^q(\M)\) with
\(\|ab\|_{L^q(\M)}\le\|a\|_{L^r(\M)}\|b\|_{L^s(\M)}\) for \(a\in L^r(\M)\) and \(b\in L^s(\M)\)
by H\"older's inequality for the trace when \(r,s<\infty\) \cite[Theorem 4.2]{fk}.
When \(r=\infty\) or \(s=\infty\), the stated inequality is the
bound for left or right multiplication by a bounded operator; if both
are infinite, it is submultiplicativity of the operator norm.

We also use cyclicity of the trace for products \(a_1\cdots a_m\), where \(m\ge2\) is an integer, \(a_j\in L^{p_j}(\M)\) with \(1\le p_j<\infty\) for \(1\le j\le m\), and \(\sum_{j=1}^m 1/p_j=1\). Repeated applications of H\"older's inequality put both \(a_1(a_2\cdots a_m)\) and \((a_2\cdots a_m)a_1\) in \(L^1(\M)\). Applying Brown and Kosaki's trace identity~\cite[Theorem~17]{bk} to \(a_1\) and \(a_2\cdots a_m\) gives
\begin{equation*}
 \operatorname{tr}_{\M}(a_1\cdots a_m)=\operatorname{tr}_{\M}(a_2\cdots a_m a_1)\,.
\end{equation*}
Iterating this identity gives invariance under every cyclic permutation.

\subsection{Scalar multipliers and Fourier localization}\label{pre:fourier}
For \(d\in\mathbb N\) and scalar Schwartz functions \(f,a:\R^d\to\mathbb C\), we use
\begin{equation*}
 \widehat f(\xi)=\int_{\R^d}e^{-ix\cdot\xi}f(x)\,dx,
 \qquad
 \check a(x)=(2\pi)^{-d}\int_{\R^d}e^{ix\cdot\xi}a(\xi)\,d\xi
 \qquad(x,\xi\in\R^d).
\end{equation*}
We also write \(a^\vee=\check a\). The Fourier transform acts on the scalar factor of finite tensors. For \(1\le p\le\infty\), the multiplier \(a(D)\) maps \(\mathscr S(\R^d)\odot L^p(\M)\) into itself by \(a(D)F=\check a*F\). For a multi-index \(\nu=(\nu_1,\ldots,\nu_d)\in\{0,1,\ldots\}^d\), write \(|\nu|=\nu_1+\cdots+\nu_d\) and \(\partial^\nu=\partial_1^{\nu_1}\cdots\partial_d^{\nu_d}\). We write \(A\lesssim B\) when \(A\le CB\) with the constant dependencies specified in the argument, and \(A\asymp B\) when \(A\lesssim B\) and \(B\lesssim A\).

For a bounded multiplier symbol \(m:\R^2\to\mathbb C\), set \(m^\sharp(\xi)=\overline{m(-\xi)}\) for \(\xi\in\R^2\). Taking the pointwise adjoint \(F^*(x)=F(x)^*\) gives \(\widehat{F^*}(\xi)=\widehat F(-\xi)^*\) on test tensors. Consequently,
\begin{equation}\label{intro:half-wave-adjoint}
 [m(D)F]^*=m^\sharp(D)F^*\,.
\end{equation}

For linearly independent \(v_1,v_2\in\R^2\), the associated \emph{full-rank period lattice} \(\Gamma\) and \emph{fundamental cell} \(P\) are
\begin{equation*}
 \Gamma=\mathbb Z v_1+\mathbb Z v_2,\qquad
 P=\{t_1v_1+t_2v_2:0\le t_1,t_2<1\},
\end{equation*}
where \(t_1,t_2\) are real. The translates \(P+\gamma\), \(\gamma\in\Gamma\), are pairwise disjoint and cover \(\R^2\), by uniqueness of the integer and fractional parts of the coordinates in this basis. Write \(|P|\) for the Lebesgue area of \(P\). Probability Haar measure on \(\R^2/\Gamma\) is \(dx/|P|\) on \(P\). The dual lattice is
\begin{equation*}
 \Gamma^*=\{\xi\in\R^2:\xi\cdot\gamma\in2\pi\mathbb Z
                  \text{ for every }\gamma\in\Gamma\}.
\end{equation*}
For a scalar symbol \(a:\Gamma^*\to\mathbb C\), the periodic multiplier \(a(D)\) maps finite trigonometric polynomials with coefficients in \(L^p(\M)\) into themselves by multiplying the coefficient of \(e^{ix\cdot\xi}\) by \(a(\xi)\), for each \(\xi\in\Gamma^*\) and \(1\le p\le\infty\).

\begin{samepage}
Let \(Y=\R^d\), \(d\in\mathbb N\), with Lebesgue measure, or \(Y=\R^2/\Gamma\) with probability Haar measure, and denote either measure by \(d\mu\). For a scalar kernel \(K\in L^1(Y)\), define convolution on finite test tensors by
\begin{equation*}
 (K*F)(x)=\int_Y K(y)F(x-y)\,d\mu(y),\qquad
 \widetilde K(y)=\overline{K(-y)}\qquad(x,y\in Y).
\end{equation*}
For duality between \(L^1(\N_Y)\) and \(\N_Y\), we use the sesquilinear pairing \(\langle H,F\rangle_Y=\operatorname{tr}_{\N_Y}(H^*F)\), where \(H\in L^1(\N_Y)\), \(F\in\N_Y\), and \(\operatorname{tr}_{\N_Y}\) is the product trace.
\end{samepage}

\begin{lemma}
\MakeLinkTarget{theorem}\label{pre:scalar-kernels}
\begin{theoremitems}
\item\label{pre:kernel-bounds} For every compact set \(B\subset\R^2\), there is a constant \(C_B>0\) such that every \(a\in C_c^\infty(\R^2)\) supported in \(B\) satisfies
\begin{equation}\label{pre:kern}
 \|\check a\|_{L^1(\R^2)}
 \le C_B\sum_{|\nu|\le4}\|\partial^\nu a\|_{L^\infty(\R^2)}.
\end{equation}
For every \(d\in\mathbb N\), \(a\in\mathscr S(\R^d)\), and real \(h>0\),
\begin{equation*}
 \|(a(h\,\cdot))^\vee\|_{L^1(\R^d)}=\|\check a\|_{L^1(\R^d)}.
\end{equation*}
\item\label{pre:periodic-kernel} For every full-rank lattice \(\Gamma\subset\R^2\) with fundamental cell \(P\) as above and every \(a\in C_c^\infty(\R^2)\), the multiplier with symbol \(a|_{\Gamma^*}\) has convolution kernel \(K_{\mathrm{per}}:\R^2/\Gamma\to\mathbb C\) satisfying
\begin{equation*}
 K_{\mathrm{per}}(x)=|P|\sum_{\gamma\in\Gamma}\check a(x+\gamma),
 \qquad
 \|K_{\mathrm{per}}\|_{L^1(d\mu)}\le\|\check a\|_{L^1(\R^2)}.
\end{equation*}
\item\label{pre:convolution} For \(Y=\R^d\) or \(Y=\R^2/\Gamma\) with the measures above, every scalar \(K\in L^1(Y)\) defines a bounded map \(F\mapsto K*F\) from \(L^p(\N_Y)\) to itself for \(1\le p\le\infty\), with
\begin{equation*}
 \|K*F\|_{L^p(\N_Y)}\le\|K\|_{L^1(Y)}\|F\|_{L^p(\N_Y)}.
\end{equation*}
At \(p=\infty\), this map is the adjoint of \(H\mapsto\widetilde K*H\) on \(L^1(\N_Y)\) under the pairing above. All bounds are independent of \(\M\) and its trace.
\end{theoremitems}
\end{lemma}

\begin{proof}
To prove~\eqref{pre:kern}, integrate by parts with \((1-\Delta_\xi)^2\), where \(\Delta_\xi=\partial_{\xi_1}^2+\partial_{\xi_2}^2\). For every \(x\in\R^2\),
\begin{equation*}
 |\check a(x)|
 \le (2\pi)^{-2}(1+|x|^2)^{-2}\|(1-\Delta_\xi)^2a\|_{L^1(\R^2)}
 \le C_B(1+|x|^2)^{-2}\sum_{|\nu|\le4}\|\partial^\nu a\|_{L^\infty(\R^2)}.
\end{equation*}
Integration in \(x\) gives~\eqref{pre:kern}. For \(a\in\mathscr S(\R^d)\) and \(h>0\), the substitution \(\eta=h\xi\) in Fourier inversion gives
\begin{equation*}
 (a(h\,\cdot))^\vee(x)=h^{-d}\check a(x/h)\qquad(x\in\R^d).
\end{equation*}
Substituting \(z=x/h\) in the kernel integral proves that its \(L^1(\R^d)\) norm is unchanged.

Since \(\check a\) is Schwartz, its periodization converges absolutely and uniformly on \(P\). The factor \(|P|\) cancels the factor \(|P|^{-1}\) in \(d\mu\): for each \(\xi\in\Gamma^*\),
\begin{equation*}
 \int_P K_{\mathrm{per}}(x)e^{-ix\cdot\xi}\,\frac{dx}{|P|}
 =\sum_{\gamma\in\Gamma}\int_{P+\gamma}\check a(z)e^{-iz\cdot\xi}\,dz
 =a(\xi).
\end{equation*}
This identifies the periodic multiplier, and the disjoint covering by \(P+\gamma\) gives
\begin{equation*}
 \|K_{\mathrm{per}}\|_{L^1(d\mu)}
 \le\sum_{\gamma\in\Gamma}\int_P|\check a(x+\gamma)|\,dx
 =\int_{\R^2}|\check a(z)|\,dz.
\end{equation*}

For \(1\le p<\infty\), the Bochner triangle inequality and translation invariance give
\begin{equation*}
 \|K*F\|_{L^p(\N_Y)}
 \le\int_Y|K(y)|\,\|F(\cdot-y)\|_{L^p(\N_Y)}\,d\mu(y)
 =\|K\|_{L^1(Y)}\|F\|_{L^p(\N_Y)}.
\end{equation*}
Density extends convolution from test tensors to \(L^p(\N_Y)\). At \(p=\infty\), the \(L^1\) bound for \(H\mapsto\widetilde K*H\) defines its adjoint by
\begin{equation*}
 \langle H,K*F\rangle_Y=\langle\widetilde K*H,F\rangle_Y
 \qquad(H\in L^1(\N_Y),\ F\in\N_Y).
\end{equation*}
Its norm is at most \(\|\widetilde K\|_{L^1(Y)}=\|K\|_{L^1(Y)}\), and Fubini's theorem shows that it agrees with scalar convolution on finite test tensors.
\end{proof}

For \(d\in\mathbb N\) and \(1\le p<\infty\), an element \(F\in L^p(\R^d;L^p(\M))\) has an \(L^p(\M)\)-valued tempered Fourier transform defined by
\begin{equation*}
 \langle\widehat F,\varphi\rangle
 =\int_{\R^d}F(x)\widehat\varphi(x)\,dx
 \qquad(\varphi\in\mathscr S(\R^d));
\end{equation*}
H\"older's inequality makes this a continuous map from \(\mathscr S(\R^d)\) to \(L^p(\M)\). For a closed set \(K\subset\R^d\), the condition \(\operatorname{supp}\widehat F\subset K\) means that \(\langle\widehat F,\varphi\rangle=0\) for every \(\varphi\in C_c^\infty(\R^d\setminus K)\).

\begin{lemma}\label{pre:localized-approximation}
Let \(d\in\mathbb N\), \(1\le p<\infty\), and \(F\in L^p(\R^d;L^p(\M))\). Suppose that \(\operatorname{supp}\widehat F\subset K\) for a compact set \(K\subset\R^d\), and that \(\chi\in C_c^\infty(\R^d)\) equals one on a neighborhood of \(K\). Then \(\chi(D)F=F\). There is a sequence of finite sums
\begin{equation*}
 G_n=\sum_{j=1}^{m_n}f_{n,j}\otimes b_{n,j},
 \qquad m_n\in\mathbb N,\quad n\in\mathbb N,
\end{equation*}
converging to \(F\) in \(L^p(\R^d;L^p(\M))\), with \(f_{n,j}\in\mathscr S(\R^d)\), \(b_{n,j}\in\M\), and \(\operatorname{tr}_{\M}(1_{(0,\infty)}(|b_{n,j}|))<\infty\) for \(1\le j\le m_n\). Each \(\chi(D)G_n\) is a finite sum of the same type, with the same coefficients \(b_{n,j}\), and
\begin{equation}\label{pre:localized-core}
 \begin{aligned}
 \operatorname{supp}\widehat{\chi(D)G_n}&\subset\operatorname{supp}\chi,\\
 \|\chi(D)G_n-F\|_{L^p(\R^d;L^p(\M))}
 &\le\|\check\chi\|_{L^1(\R^d)}\|G_n-F\|_{L^p(\R^d;L^p(\M))}\longrightarrow0
 \qquad(n\to\infty).
 \end{aligned}
\end{equation}
\end{lemma}

\begin{proof}
We first prove \(\chi(D)F=F\), then apply Lemma~\ref{pre:scalar-kernels} to \(\check\chi*(G_n-F)\). For \(\varphi\in\mathscr S(\R^d)\), choose \(\theta\in C_c^\infty(\R^d)\) equal to one near zero. For every real \(R\ge1\), the test function \(\theta(\cdot/R)(1-\chi)\varphi\) is compactly supported away from \(K\). Hence its pairing with \(\widehat F\) vanishes. As \(R\to\infty\), these tests converge to \((1-\chi)\varphi\) in the Schwartz topology, so continuity gives \(\chi\widehat F=\widehat F\). Testing convolution against Schwartz functions and applying Fubini's theorem gives
\begin{equation*}
 \widehat{\chi(D)F}=\widehat{\check\chi*F}
 =\chi\widehat F=\widehat F.
\end{equation*}
Injectivity of the Fourier transform on tempered distributions yields \(\chi(D)F=F\) in \(L^p(\R^d;L^p(\M))\).

Choose \(G_n\to F\) from the dense tensor class in Section~\ref{pre:tracial}. Since \(\chi\widehat f_{n,j}\in C_c^\infty(\R^d)\),
\begin{equation*}
 \chi(D)G_n=\sum_{j=1}^{m_n}(\chi\widehat f_{n,j})^\vee\otimes b_{n,j}
\end{equation*}
is again a finite sum with Schwartz scalar factors and the same coefficients. Its Fourier support lies in \(\operatorname{supp}\chi\). Finally,
\begin{equation*}
 \chi(D)G_n-F=\check\chi*(G_n-F),
\end{equation*}
and Lemma~\ref{pre:scalar-kernels}\ref{pre:convolution} gives the norm estimate and convergence in~\eqref{pre:localized-core}.
\end{proof}

If \(K\) is a compact annulus away from the origin, \(\chi\) may be supported in any slightly larger open annulus containing \(K\), and the functions \(\chi_k(\xi)=\chi(2^{-k}\xi)\), \(k\in\mathbb Z\), equal one near \(2^kK\) and satisfy \(\|\check\chi_k\|_{L^1(\R^d)}=\|\check\chi\|_{L^1(\R^d)}\) by Lemma~\ref{pre:scalar-kernels}\ref{pre:kernel-bounds}.

\subsection{Reduction to a periodic single-cap estimate}\label{pre:periodic-cap}\label{sec:period}

To isolate the estimate that will be proved by the energy argument, let $\Gamma\subset\R^2$ be a full-rank period lattice and put $\T^2=\R^2/\Gamma$. We retain the dual lattice $\Gamma^*$ and probability Haar measure from Section~\ref{pre:fourier}, writing $dX$ for this measure. For a finite set $\Lambda\subset\Gamma^*$ and coefficients $A_\xi\in L^4(\M)$, $\xi\in\Lambda$, define the periodic half-wave $u:\R\times\T^2\to L^4(\M)$ and its fourth moment $q_u:\R\to[0,\infty)$ by
\begin{align}
 u(T,X)&=\sum_{\xi\in\Lambda}A_\xi e^{i(X\cdot\xi-T|\xi|)},
       \qquad T\in\R,\quad X\in\T^2,\label{pre:wave}\\
 q_u(T)&=\|u(T,\cdot)\|_{L^4(\T^2;L^4(\M))}^4
       =\int_{\T^2}\|u(T,X)\|_{L^4(\M)}^4\,dX\,.
       \label{pre:qu}
\end{align}
To specify a cap, choose vectors $e_x,e_{b_0}\in\R^2$ with $|e_x|=|e_{b_0}|=1$ and $e_x\cdot e_{b_0}=0$. The vector $e_{b_0}$ gives the cap axis. For a point $X$ on the spatial cover $\R^2$ and a frequency $\xi\in\R^2$, define their real coordinates by
\begin{equation*}
 x=X\cdot e_x,\qquad b_0=X\cdot e_{b_0},\qquad
 \xi_x=\xi\cdot e_x,\qquad \xi_{b_0}=\xi\cdot e_{b_0}\,.
\end{equation*}
Thus $X=xe_x+b_0e_{b_0}$ and $\xi=\xi_xe_x+\xi_{b_0}e_{b_0}$. We write $X=(x,b_0)$ and $\xi=(\xi_x,\xi_{b_0})$ for these coordinates. For $\xi_{b_0}>0$, put
\begin{equation}\label{pre:slope}
 \kappa_\xi=\frac{|\xi|+\xi_{b_0}}{\sqrt2}>0,
 \qquad v_\xi=\frac{\xi_x}{\kappa_\xi}\in\R\,.
\end{equation}
The condition $|v_\xi|\le\gamma$ specifies the cap aperture about the positive $e_{b_0}$ direction. The orthonormal frame $(e_x,e_{b_0})$ can be chosen independently of $\Gamma$. Changing this frame gives orthogonal coordinate descriptions of both $\Gamma$ and $\Gamma^*$.

\begin{proposition}\label{en:cap}
Fix real numbers $0<c_0<C_0$. There are an aperture $\gamma>0$ and a constant $C>0$ such that, for every real $N\ge2$, every full-rank period lattice $\Gamma\subset\R^2$, every orthonormal frame $(e_x,e_{b_0})$ of $\R^2$, and every semifinite tracial coefficient algebra $\M$, the following holds. Let $\Lambda\subset\Gamma^*$ be a finite set such that
\begin{equation*}
 c_0\le |\xi|/N\le C_0,\qquad \xi_{b_0}>0,\qquad
 |v_\xi|\le\gamma\qquad(\xi\in\Lambda)\,.
\end{equation*}
For coefficients $A_\xi\in L^4(\M)$, $\xi\in\Lambda$, define $u$ and $q_u$ by~\eqref{pre:wave} and~\eqref{pre:qu}. Then
\begin{equation}\label{en:capestimate}
 \int_1^2 q_u(T)\,dT\le C[\log(2+N)]^2q_u(0)\,.
\end{equation}
\end{proposition}

\begin{proof}[Proof of Theorem~\ref{thm:four} assuming Proposition~\ref{en:cap}]
Fix $\psi\in C_c^\infty((1/2,2))$, a real $N\ge2$, and $h=N^{-1}$. Recall that $I=[1,2]$. We first estimate $U(-T)P_NF$ on $I$. Taking pointwise adjoints will then give~\eqref{intro:l4}.

A finite cap partition gives the annular estimate on the torus. To apply it on each Euclidean unit square, we restrict \(F\) to the concentric square of side ten and extend this restriction periodically. Decay of the Euclidean kernel away from the propagation region controls the contributions from outside the larger square and from the extra periodic copies.

\medskip
For the periodic estimate, take $\T^2=\R^2/(100\mathbb Z^2)$, so its dual lattice is $(2\pi/100)\mathbb Z^2$. Apply Proposition~\ref{en:cap} with $c_0=1/2$ and $C_0=2$, and choose a smooth finite partition of the annulus into caps of its aperture $\gamma$. In the scaled frequency variable $\eta\in\R^2$ this gives
\begin{equation*}
 \psi(|\eta|)=\sum_{d=1}^{d_{\max}}\psi_d(\eta),\qquad
 d_{\max}\in\mathbb N,\quad \psi_d\in C_c^\infty(\R^2)\,.
\end{equation*}
For each $1\le d\le d_{\max}$, the symbol $\psi_d$ is supported in $1/2\le|\eta|\le2$ and in a cap satisfying the slope condition of Proposition~\ref{en:cap} in an orthonormal frame $(e_x^{(d)},e_{b_0}^{(d)})$, with $e_{b_0}^{(d)}$ pointing along that cap. The frames $(e_x^{(d)},e_{b_0}^{(d)})$ and symbols $\psi_d$ are fixed independently of $N$ and $\M$. In each frame, we express both the spatial variable and the period lattice in the new coordinates; normalized Haar integrals are unchanged. Thus the proposition applies to $U(-T)\psi_d(hD)F$ for every $F\in L^4(\N_{\T^2})$. The frequency support is finite, since the dual lattice meets the compact set $\operatorname{supp}(\psi_d(h\,\cdot))$ in finitely many points. Its Fourier coefficients belong to $L^4(\M)$ by the Bochner triangle inequality.

By Lemma~\ref{pre:scalar-kernels}, $\|\psi_d(hD)F\|_{L^4(\N_{\T^2})}\le\|\check\psi_d\|_{L^1(\R^2)}\|F\|_{L^4(\N_{\T^2})}$. Taking fourth roots in~\eqref{en:capestimate} changes the logarithmic factor to $([\log(2+N)]^2)^{1/4}=[\log(2+N)]^{1/2}$. The triangle inequality in $L^4(I\times\T^2;L^4(\M))$ then gives
\begin{equation}\label{pre:torus-annulus}
 \begin{aligned}
 \|U(-T)P_NF\|_{L^4(I\times\T^2;L^4(\M))}
 &\le\sum_{d=1}^{d_{\max}}
   \|U(-T)\psi_d(hD)F\|_{L^4(I\times\T^2;L^4(\M))}\\
 &\le C_\psi[\log(2+N)]^{1/2}\|F\|_{L^4(\N_{\T^2})}\,.
 \end{aligned}
\end{equation}
Changing probability Haar measure to Lebesgue measure on a fundamental square multiplies both fourth powers by its area, so~\eqref{pre:torus-annulus} also holds with that measure.

\medskip
To pass from~\eqref{pre:torus-annulus} to the plane, we split the Euclidean convolution into a local part and a part controlled by kernel decay. Let $F\in L^4(\N_{\R^2})$. For each $T\in I$, the spatial convolution kernel $K_{h,T}:\R^2\to\mathbb C$ of $U(-T)P_N$ is
\begin{equation*}
 K_{h,T}(x)=(2\pi h)^{-2}\int_{\R^2}
      e^{i(x\cdot\eta-T|\eta|)/h}\psi(|\eta|)\,d\eta,
 \qquad x\in\R^2\,.
\end{equation*}
Fix the partition $\mathscr D=\{k+[0,1)^2:k\in\mathbb Z^2\}$ of $\R^2$. For $D\in\mathscr D$, let $10D$ denote the concentric square of side ten. Define $F_{\mathrm{near}},F_{\mathrm{far}}:I\times\R^2\to L^4(\M)$ by
\begin{equation*}
 \begin{aligned}
 F_{\mathrm{near}}(T,x)&=\sum_{D\in\mathscr D}1_D(x)\,[K_{h,T}*(F1_{10D})](x),\\
 F_{\mathrm{far}}(T,x)&=\sum_{D\in\mathscr D}1_D(x)\,[K_{h,T}*(F1_{(10D)^c})](x)\,.
 \end{aligned}
\end{equation*}
Each sum has at most one nonzero term at any $x\in\R^2$.
Thus $U(-T)P_NF=F_{\mathrm{near}}+F_{\mathrm{far}}$, and the triangle inequality gives
\begin{equation}\label{pre:physical-split}
 \|U(-T)P_NF\|_{L^4(I\times\R^2;L^4(\M))}
 \le \|F_{\mathrm{far}}\|_{L^4(I\times\R^2;L^4(\M))}
 +\|F_{\mathrm{near}}\|_{L^4(I\times\R^2;L^4(\M))}\,.
\end{equation}
We first bound $F_{\mathrm{far}}$ by kernel decay. We then compare $F_{\mathrm{near}}$ on each square with a periodic wave to which~\eqref{pre:torus-annulus} applies.

For the decay estimate, fix an integer $m>2$. When $|x|\ge3$ and $1\le T\le2$, the phase $\phi(\eta)=x\cdot\eta-T|\eta|$ satisfies $|\nabla_\eta\phi|\ge|x|-2\ge|x|/3$ on the fixed annulus supporting $\psi(|\eta|)$. Apply the formal adjoint of $h\nabla\phi\cdot\nabla_\eta/(i|\nabla\phi|^2)$ to the amplitude $m+2$ times. Each derivative of $\nabla\phi/|\nabla\phi|^2$ on this annulus is $O_m(|x|^{-1})$. The $m+2$ integrations by parts, together with the prefactor $h^{-2}$, therefore give
\begin{equation}\label{pre:kernel-tail}
 |K_{h,T}(x)|\le C_{m,\psi}h^m(1+|x|)^{-m}
 \qquad(|x|\ge3,\ T\in I)\,.
\end{equation}
For $x\in D$ and $y\notin10D$ we have $|x-y|\ge9/2$. Consequently the Bochner triangle inequality bounds $\|F_{\mathrm{far}}(T,x)\|_{L^4(\M)}$ by the convolution of $y\mapsto\|F(y)\|_{L^4(\M)}$ with $C_{m,\psi}h^m(1+|\cdot|)^{-m}$. This scalar majorant is independent of $D$ and $T$. It is integrable on $\R^2$ because $m>2$. Young's inequality and $|I|=1$ give
\begin{equation}\label{pre:far-bound}
 \|F_{\mathrm{far}}\|_{L^4(I\times\R^2;L^4(\M))}
 \le C_{m,\psi}h^m\|F\|_{L^4(\N_{\R^2})}\,.
\end{equation}

\medskip
To estimate $F_{\mathrm{near}}$ by~\eqref{pre:torus-annulus}, fix $D\in\mathscr D$. Use the square of side $100$ concentric with $D$ as a fundamental square, and extend $g_D=F1_{10D}$ with period $100$ in each coordinate. Denote this extension by $g_D^{\mathrm{per}}:\T^2\to L^4(\M)$. With Lebesgue measure on the fundamental square, its fourth power norm is exactly $\|F\|_{L^4(10D;L^4(\M))}^4$. Periodizing the kernel gives, for $(T,x)\in I\times D$,
\begin{equation*}
 \begin{aligned}
 (U(-T)P_Ng_D^{\mathrm{per}})(x)
 &=K_{h,T}*g_D(x)+R_D(T,x),\\
 R_D(T,x)
 &=\sum_{k\in\mathbb Z^2\setminus\{0\}}
       \int_{10D}K_{h,T}(x-y+100k)F(y)\,dy\,.
 \end{aligned}
\end{equation*}
Here $R_D:I\times D\to L^4(\M)$ is the contribution from the additional periodic copies. For $x\in D$, $y\in10D$ and $k\in\mathbb Z^2\setminus\{0\}$, the bound $|x-y|\le11\sqrt2/2$ gives $|x-y+100k|\gtrsim1+|k|$. In particular these distances exceed $3$, so~\eqref{pre:kernel-tail} and the Bochner triangle inequality yield
\begin{equation*}
 \begin{aligned}
 \|R_D(T,x)\|_{L^4(\M)}
 &\le C_{m,\psi}h^m
     \sum_{k\in\mathbb Z^2\setminus\{0\}}(1+|k|)^{-m}
           \int_{10D}\|F(y)\|_{L^4(\M)}\,dy\\
 &\le C'_{m,\psi}h^m\|F\|_{L^4(10D;L^4(\M))}\,.
 \end{aligned}
\end{equation*}
The lattice sum $\sum_{k\in\mathbb Z^2\setminus\{0\}}(1+|k|)^{-m}$ converges because $m>2$, and the last inequality follows from H\"older's inequality on a square of area $100$. The same majorant gives absolute convergence in $L^4(\M)$ for the series defining $R_D(T,x)$. Since $|I\times D|=1$, the same bound holds in $L^4(I\times D;L^4(\M))$.

Since $D$ lies in the fundamental square, the $L^4(I\times D;L^4(\M))$ norm of $U(-T)P_Ng_D^{\mathrm{per}}$ is bounded by its norm on $I$ times the full fundamental square. We estimate the latter by~\eqref{pre:torus-annulus} applied to $g_D^{\mathrm{per}}$. Substituting $K_{h,T}*g_D=U(-T)P_Ng_D^{\mathrm{per}}-R_D$ on $D$ and using $h\le1/2$ now gives
\begin{equation}\label{pre:near-square}
 \|K_{h,T}*(F1_{10D})\|_{L^4(I\times D;L^4(\M))}
 \le C_\psi[\log(2+N)]^{1/2}
        \|F\|_{L^4(10D;L^4(\M))}\,.
\end{equation}
The sets $I\times D$, $D\in\mathscr D$, partition $I\times\R^2$. Taking fourth powers in~\eqref{pre:near-square} and summing over $D\in\mathscr D$ therefore gives
\begin{equation*}
 \begin{aligned}
 \|F_{\mathrm{near}}\|_{L^4(I\times\R^2;L^4(\M))}^4
 &\le C_\psi[\log(2+N)]^2
             \sum_D\|F\|_{L^4(10D;L^4(\M))}^4,\\
 \sum_D\|F\|_{L^4(10D;L^4(\M))}^4
 &=\int_{\R^2}\|F(y)\|_{L^4(\M)}^4
                        \sum_D1_{10D}(y)\,dy
 \le11^2\|F\|_{L^4(\N_{\R^2})}^4\,.
 \end{aligned}
\end{equation*}
There are at most eleven admissible square centers in each coordinate, which gives the last inequality. Tonelli's theorem applies because every summand is nonnegative. Taking fourth roots and combining the resulting estimate with~\eqref{pre:far-bound} in~\eqref{pre:physical-split} gives the required bound for $U(-T)P_NF$.

For each $T\in I$, apply~\eqref{intro:half-wave-adjoint} to the symbol $m_T:\R^2\to\mathbb C$ given by $m_T(\xi)=e^{-iT|\xi|}\overline\psi(|\xi|/N)$. Since $m_T$ is even, $m_T^\sharp(\xi)=e^{iT|\xi|}\psi(|\xi|/N)$, and hence
\begin{equation*}
 [U(-T)\overline\psi(|D|/N)(F^*)]^*
       =U(T)\psi(|D|/N)F\,.
\end{equation*}
Taking adjoints preserves the $L^4(\N_{\R^2})$ and $L^4(\N_{I\times\R^2})$ norms. Applying the bound for $U(-T)\overline\psi(|D|/N)$ to $F^*$ therefore proves~\eqref{intro:l4} with constant $C_{\overline\psi}$. The maximum of $C_\psi$ and $C_{\overline\psi}$ works for both $U(T)P_N$ and $U(-T)P_N$, independently of $N$, $\M$, and its trace.
\end{proof}

It remains to prove Proposition~\ref{en:cap}. Section~\ref{cur:sec} encodes the localized quadratic density of $u$ in conserved currents and pairs the currents for $\M$ and its opposite algebra to cancel their determinant contributions. Section~\ref{en:section} uses the time components as sources for an auxiliary Dirichlet energy. Two radial comparisons of that energy, one outward to a fixed radius and one back to the frequency scale, produce the two logarithms in~\eqref{en:capestimate}.

\section{Conserved currents and cancellation via the opposite algebra}\label{cur:sec}

We now begin the proof of Proposition~\ref{en:cap}. The first task is to encode the localized quadratic density of the periodic wave by conserved currents. Their time components will be the sources of the Dirichlet equations in Section~\ref{en:section}, and their conservation laws will give the energy identity~\eqref{en:balance}. The scalar construction of~\cite[Lemma~2.1 and Proposition~2.2]{openai} produces currents whose real columns have zero divergence and whose ranges lie in a quaternionic plane. In the scalar energy argument, the determinant on this plane is nonnegative. After polarization with operator coefficients, the sum of these determinants can have either sign. Proposition~\ref{cur:cancel} will show that passing from $\M$ to the opposite algebra $\M^{\mathrm{op}}$ reverses this sum pointwise. Adding the energy identities for the two algebras will therefore cancel the determinant contributions.

\subsection{The lifted scalar current}\label{cur:lift}

The auxiliary variables allow us to localize the wave in its frequency slope while retaining an exact conservation law. Attach a complex variable to each of the physical coordinates $(T,x,b_0)$, and denote these variables by $(\tau,P,Z)$. We take them in the plane $\mathbb C_k=\mathbb R+k\mathbb R$ inside the quaternion algebra $\mathbb H$, where $k^2=-1$. At fixed $(T,\tau)\in\R\times\mathbb C_k$, the variables left on the fiber are the physical position $X=(x,b_0)\in\T^2$ and the auxiliary pair $(P,Z)\in\mathbb C_k^2$. The Dirichlet equation in Section~\ref{en:section} will be solved on a bounded part of this fiber; its integration domain and measure will be specified there.

Recall that $X=(x,b_0)$ denotes coordinates in the orthonormal frame $(e_x,e_{b_0})$ of Section~\ref{pre:periodic-cap}. Fix real numbers $0<c_0<C_0$ and $N\ge2$. Let $\Lambda$ be a finite subset of the dual period lattice $\Gamma^*$ such that
\begin{equation*}
 \xi_{b_0}>0,\qquad c_0\le |\xi|/N\le C_0
 \qquad(\xi\in\Lambda)\,.
\end{equation*}
For $\xi=(\xi_x,\xi_{b_0})\in\Lambda$, recall the frequency parameters $\kappa_\xi,v_\xi$ from~\eqref{pre:slope} and define the phase $\theta_\xi$ by
\begin{equation*}
 \kappa_\xi=(|\xi|+\xi_{b_0})/\sqrt2,\qquad
 v_\xi=\xi_x/\kappa_\xi,\qquad
 \theta_\xi=X\cdot\xi-T|\xi|\,.
\end{equation*}
Fix the disk radius $R=1/10$. For now the aperture $0<\gamma<R/8$ in $|v_\xi|\le\gamma$ is arbitrary; Section~\ref{en:section} will choose it independently of $N$, $\Lambda$, the cap axis, and the coefficient algebra. An orthogonal change of $(e_x,e_{b_0})$ changes the components used to describe the dual lattice.

We distinguish the auxiliary plane $\mathbb C_k$ from the plane $\mathbb C_i=\mathbb R+i\mathbb R$ containing the scalar Fourier amplitudes and the phases $e^{i\theta_\xi}$. The quaternionic units satisfy $i^2=j^2=k^2=-1$, $ij=k$, and $ki=j$. We use quaternionic conjugation $q\mapsto\bar q$, conjugate transpose $A^*=(\overline{A_{bd}})_{d,b}$, and right quaternionic column spaces.

To express the phase as a quadratic polynomial in the slope $v_\xi$, introduce null coordinates by the same orthogonal change in the physical and auxiliary variables:
\begin{equation}\label{cur:frames}
\begin{gathered}
 \begin{array}{c|c|c}
  &\text{physical frame}&\text{null frame}\\ \hline
  \text{real spacetime}& (T,x,b_0)&(t,x,s_0)\\
  \text{auxiliary variables}&(\tau,P,Z)&(Y,P,S)
 \end{array}\\[.5em]
 \begin{aligned}
  t&=(T+b_0)/\sqrt2, &s_0&=(b_0-T)/\sqrt2,\\
  Y&=(\tau+Z)/\sqrt2, &S&=(Z-\tau)/\sqrt2\,.
 \end{aligned}
\end{gathered}
\end{equation}
Indeed,
\begin{equation*}
 |\xi|-\xi_{b_0}=\frac{\xi_x^2}{|\xi|+\xi_{b_0}}
     =\frac{\kappa_\xi v_\xi^2}{\sqrt2},\qquad
 |\xi|=\frac{\kappa_\xi}{\sqrt2}(1+v_\xi^2/2)\,,
\end{equation*}
and hence
\begin{equation*}
 \theta_\xi=\kappa_\xi(s_0+v_\xi x-v_\xi^2t/2)\,.
\end{equation*}
In particular $\kappa_\xi/N$ has positive upper and lower bounds depending only on the fixed annulus. The quadratic dependence on $v_\xi$ is what permits the Gaussian localization below.

On the auxiliary chart $Y\ne0$, put $W=P/Y\in\mathbb C_k$ and $e=Y/|Y|$. We use $W$ to distinguish the real slopes $v_\xi$: the factor $\exp(-\kappa_\xi|Y|\,|W-v_\xi|^2/2)$ is concentrated near $W=v_\xi$, with width $(\kappa_\xi|Y|)^{-1/2}$. To combine this Gaussian with the complex Fourier amplitudes, define $\rho_e:\mathbb C_i\to\mathbb H$ by $\rho_e(\alpha+i\beta)=\alpha-ei\beta$ for $\alpha,\beta\in\R$. Since $ic=\bar c i$ for $c\in\mathbb C_k$, the unit $ei$ squares to $-1$; thus $\rho_e$ is a real algebra embedding preserving conjugation.

Fix scalar Fourier coefficients $c=(c_\xi)_{\xi\in\Lambda}\in\mathbb C_i^\Lambda$. For each $v\in\R$, let $V^v\in\mathbb C_k^3$ be the column below, written in the null frame. Define the real Gaussian weights $m_\xi$ and the lifted wave $\chi(c)$, with values in $\mathbb H^3$, by
\begin{align}
 m_\xi&=(\kappa_\xi/|Y|)^{1/2}
       \exp\{-\kappa_\xi|Y|\,|W-v_\xi|^2/2\},\qquad
 V^v=(1,(v+W)/2,-Wv/2)^{\mathsf t},\notag\\
 \chi(c)&=\sum_{\xi\in\Lambda}V^{v_\xi}m_\xi
                         \rho_e(c_\xi e^{i\theta_\xi})\,.
 \label{cur:scalar}
\end{align}
These fields are defined on the lifted spacetime chart $Y\ne0$. The column $V^{v_\xi}$ multiplies the embedded Fourier amplitude on the left. The conjugation in the scalar construction~\cite[(2.14)]{openai} determines the order of the trace products after polarization.

We also need to restrict the other auxiliary direction. For the resulting current to remain conserved, this restriction must be constant in the directions of $\chi(c)$. The appropriate holomorphic function is $\Phi=S+P^2/(2Y)$. Lemma~\ref{cur:conserve} will show that each column $V^v$ is tangent to its levels.

Fix a real number $K\ge1$ and put $a=K/N$. For fixed $\tau\in\mathbb C_k$, recall the null coordinates $Y=(\tau+Z)/\sqrt2$ and $S=(Z-\tau)/\sqrt2$ from~\eqref{cur:frames}. The forward map sends $(P,Z)\in\mathbb C_k^2$ with $\tau+Z\ne0$ to $(W,Q)\in\mathbb C_k^2$ by
\begin{subequations}\label{cur:chart}
\begin{equation}\label{cur:chart-forward}
 \begin{aligned}
 W&=\frac{P}{Y}=\frac{\sqrt2\,P}{\tau+Z},\\
 Q&=\frac{\Phi}{a}
   =\frac{1}{a\sqrt2}\left(Z-\tau+\frac{P^2}{\tau+Z}\right)\,.
 \end{aligned}
\end{equation}
For $|\tau|>a$, $|W|<R=1/10$ and $|Q|<1/10$, the inverse map recovers $(P,Z)$ from $(W,Q)$ by
\begin{equation}\label{cur:chart-inverse}
 Y=\frac{\sqrt2\,\tau+aQ}{1+W^2/2},\qquad
 P=WY,\qquad Z=\sqrt2\,Y-\tau\,.
\end{equation}
\end{subequations}
On these ranges, the numerator and denominator in~\eqref{cur:chart-inverse} are nonzero, and $|Y|\asymp|\tau|$ with constants depending only on the fixed coordinate ranges. Both maps leave the physical position $X$ fixed. We use these ranges to construct the currents; the full fiber domain, its volume Jacobian, and the boundary notation are given in Section~\ref{en:section}.

Choose a scalar function $\chi_{\mathrm{fib}}\in C_c^\infty(\mathbb C_k;[0,1])$ equal to one on $|Q|\le1/100$ and zero on $|Q|\ge1/20$. Define the quaternionic Hermitian current by
\begin{equation}\label{cur:current}
 A(c)=\chi_{\mathrm{fib}}(Q)\chi(c)\chi(c)^*\,.
\end{equation}
The next lemma proves that multiplying by this function of $Q$ preserves the conservation law. We first specify the differential conventions used to state that law and to form the Hessians in the energy identity.

Number the three coordinates \(d=0,1,2\) in the order of either frame in the array in~\eqref{cur:frames}. Write \(z_{d0}\) for the real spacetime coordinate, \(\omega_d=z_{d1}+kz_{d2}\) for the attached auxiliary coordinate, and \(\partial_{d\lambda}=\partial/\partial z_{d\lambda}\) for \(\lambda\in\{0,1,2\}\). In \(\mathbb C_k\) our Wirtinger convention is
\begin{equation}\label{cur:wirtinger}
 \partial_{\omega_d}=\tfrac12(\partial_{d1}-k\partial_{d2}),\qquad
 \partial_{\bar\omega_d}=\tfrac12(\partial_{d1}+k\partial_{d2}).
\end{equation}
For a holomorphic function \(F\) from an open subset of \(\mathbb C_k^3\) to \(\mathbb C_k\), \(\partial F=(\partial_{\omega_d}F)_{d=0}^2\) is a row. In particular, this convention applies to \(\Phi\), \(W\), and \(Q\) in~\eqref{cur:chart-forward}. Set \(e_0=1,e_1=i,e_2=j\). For a smooth real function \(f\) on an open coordinate chart in \(\R^9\), define \(G(f)\) and \(\mathsf H(f)\) on that chart with values in \(\mathbb H^3\) and the space of \(3\times3\) quaternionic Hermitian matrices, respectively, by
\begin{equation}\label{cur:hess}
 G_d(f)=\sum_{\lambda=0}^2e_\lambda\partial_{d\lambda}f,\qquad
 \mathsf H_{db}(f)=\sum_{\lambda,\mu=0}^2e_\lambda\bar e_\mu
                         \partial_{d\lambda}\partial_{b\mu}f,
 \qquad d,b\in\{0,1,2\}.
\end{equation}
For \(f\) depending only on auxiliary coordinates, \(ki=j\) gives \(G_d(f)=2(\partial_{\bar\omega_d}f)i\). In the Hessian the two mixed coefficients are \(i\bar j=-k\) and \(j\bar i=k\), so that
\[
 \mathsf H_{db}(f)
 =(\partial_{d1}+k\partial_{d2})(\partial_{b1}-k\partial_{b2})f
 =4\partial_{\bar\omega_d}\partial_{\omega_b}f .
\]
In particular the two-dimensional Euclidean Laplacian in a single auxiliary coordinate is \(4\partial_{\bar\omega}\partial_\omega\).

For a square matrix \(A\) over \(\mathbb H\) write \(\mathrm{tr}_{\mathbb H} A=\sum_d A_{dd}\). Real parts of quaternionic products satisfy \(\operatorname{Re}(uv)=\operatorname{Re}(vu)\), so \(\operatorname{Re}\mathrm{tr}_{\mathbb H}\) is cyclic in matrix products. For a \(3\times3\) quaternionic Hermitian matrix \(A\) define its real \(9\times9\) tensor by
\(A^{\mathbb R}_{d\lambda,b\mu}=\operatorname{Re}(\bar e_\lambda A_{db}e_\mu)\).
Conjugating the latter scalar shows that \(A^{\mathbb R}\) is symmetric. For $m\in\mathbb N$ and real matrices $B,C\in\R^{m\times m}$, define their Frobenius inner product by
\begin{equation*}
 B:C=\sum_{r,s=1}^{m}B_{rs}C_{rs}
     =\operatorname{tr}(B^{\mathsf t}C)\,.
\end{equation*}
Here $\mathsf t$ denotes transpose and $\operatorname{tr}$ is the ordinary matrix trace. Expanding the real part of the quaternionic trace gives
\[
 \operatorname{Re}\mathrm{tr}_{\mathbb H}(A\mathsf H(f))
 =\sum_{d,b,\lambda,\mu}\operatorname{Re}(A_{db}e_\mu\bar e_\lambda)
                \partial_{b\mu}\partial_{d\lambda}f
 =A^{\mathbb R}:D^2 f .
\]
For \(A=\chi\chi^*\), \(A^{\mathbb R}_{d\lambda,b\mu}\) is the real inner product of \(\bar e_\lambda\chi_d\) and \(\bar e_\mu\chi_b\), with \((u,v)_{\mathbb R}=\operatorname{Re}(u\bar v)\).

\begin{lemma}\label{cur:conserve}
For every $c\in\mathbb C_i^\Lambda$, each real column of $A(c)^{\mathbb R}$ has zero divergence in the nine real coordinates. The range of $A(c)$ is contained in $E=\ker(\partial\Phi:\mathbb H^3\to\mathbb H)$, where $\partial\Phi=(-W^2/2,W,1)$ is the holomorphic derivative in the auxiliary null frame. In the physical frame $V_T^v=(1+Wv/2)/\sqrt2$.
\end{lemma}
\begin{proof}
To prove conservation of the real columns of $A(c)^{\mathbb R}$, we first construct a tensor $B_0$ with divergence-free columns in coordinates where the wave satisfies a symmetric first-order system. Passing to the variables $(Y,P,S)$ introduces the Frobenius inner product of $B_0$ with the Hessian of each real component of the coordinate map. The two diagonal blocks of $B_0$ in the auxiliary state variables $w_1,w_2\in\R^2$ introduced below agree, and the mixed $w_1,w_2$ block is zero. The Cauchy--Riemann equations then imply that these inner products vanish. We identify the transformed tensor with $(\chi(c)\chi(c)^*)^{\mathbb R}$ and then use $\partial\Phi\,V^v=0$ to justify multiplication by $\chi_{\mathrm{fib}}(Q)$.

\medskip
To construct $B_0$, write $\mathsf X=(t,x,s_0)\in\R^3$. Use the real symmetric matrices
$C^t=2E_{00}$, $C^x=E_{01}+E_{10}$ and $C^{s_0}=-E_{11}$, with matrix units indexed by $0,1$. For each frequency $\xi\in\Lambda$, put $(\kappa,v)=(\kappa_\xi,v_\xi)$ and $p=d\theta_\xi$. In the coordinates $(t,x,s_0)$, the components of $p$ are
\[
 (p_t,p_x,p_{s_0})
 =(\partial_t\theta_\xi,\partial_x\theta_\xi,\partial_{s_0}\theta_\xi)
 =\kappa(-v^2/2,v,1).
\]
The symbol $C(p)$ of the first-order system satisfies
\[
 C(p)=\sum_{\nu=t,x,s_0}p_\nu C^\nu
   =\kappa\begin{pmatrix}-v^2&v\\v&-1\end{pmatrix},
 \qquad C(p)(1,v)^{\mathsf t}=0 .
\]
Since $\partial_\nu e^{i\theta_\xi}=ip_\nu e^{i\theta_\xi}$ for $\nu\in\{t,x,s_0\}$, the Fourier sum $\mathbf f=\sum_{\xi\in\Lambda}\sqrt{\kappa_\xi}(1,v_\xi)^{\mathsf t}c_\xi e^{i\theta_\xi}$ solves $\sum_\nu C^\nu\partial_{\mathsf X^\nu}\mathbf f=0$. The finite Fourier sum extends to an entire function $\mathbf f:\mathbb C_i^3\to\mathbb C_i^2$ satisfying the same identity.

For two real state vectors $w_1,w_2\in\mathbb R^2$ define the displacement $q=q(w_1,w_2)\in\R^3$ and the field $G_0:\R^3\times\R^2\times\R^2\to\mathbb C_i^2$ by
\[
 q^\nu=\tfrac12\sum_{\ell=1}^2w_\ell^{\mathsf t}C^\nu w_\ell,\qquad
 G_0=\overline{\mathbf f(\mathsf X-iq)} .
\]
This displacement will produce the Gaussian in $m_\xi$. On $\R^3_{\mathsf X}\times\R^2_{w_1}\times\R^2_{w_2}$ define, for $\ell=1,2$ and $b=0,1$,
\[
 V_{\ell b}=\sum_\nu (C^\nu w_\ell)_b\partial_{\mathsf X^\nu},
 \qquad D_{\ell b}=\partial_{w_{\ell b}}+iV_{\ell b},\qquad
 Z_\ell=\sum_b(G_0)_bD_{\ell b}.
\]
Since $\partial_{w_{\ell b}}q^\nu=(C^\nu w_\ell)_b$, the field $D_{\ell b}$ annihilates each component of $\mathsf X-iq$. Holomorphy of $\mathbf f$ gives $D_{\ell b}\bar G_0=0$ and, by conjugation, $\bar D_{\ell b}G_0=0$. Each $D_{\ell b}$ has divergence zero, and $D_{\ell b}-\bar D_{\ell b}=2iV_{\ell b}$. Thus, for $\ell=1,2$,
\begin{align*}
 \operatorname{div}Z_\ell
 &=\sum_{b=0}^1 D_{\ell b}(G_0)_b
  =\sum_{b=0}^1 (D_{\ell b}-\bar D_{\ell b})(G_0)_b\\
 &=2i\sum_{b=0}^1 V_{\ell b}(G_0)_b
  =2i w_\ell^{\mathsf t}\sum_{\nu\in\{t,x,s_0\}}
        C^\nu\partial_{\mathsf X^\nu}G_0=0 ,
\end{align*}
where the last equality uses the real coefficients of the system. When $Z_\ell$ differentiates $\bar Z_\ell$, it annihilates the coefficients $\bar G_0$. The derivative of $\bar D_{\ell b}$ in direction $D_{\ell a}$ is $-i\sum_\nu C^\nu_{ba}\partial_{\mathsf X^\nu}$. Hence
\[
 \nabla_{Z_\ell}\bar Z_\ell
 =-i\sum_\nu(G_0^{\mathsf t}C^\nu\bar G_0)\partial_{\mathsf X^\nu}.
\]
Each scalar in parentheses is real by symmetry of $C^\nu$. For
$B_0=\mathrm{Re}\sum_{\ell=1}^2 Z_\ell\otimes\bar Z_\ell$ we therefore have
\[
 \operatorname{div} B_0
 =\mathrm{Re}\sum_{\ell=1}^2\big((\operatorname{div} Z_\ell)\bar Z_\ell+
                                      \nabla_{Z_\ell}\bar Z_\ell\big)=0 .
\]
Here the divergence is taken in the first tensor index, so this identity holds for each real column.

\medskip
We transport $B_0$ to the variables $(Y,P,S)$ and check that the real columns of the transported tensor density still have zero divergence. For $b\in\{0,1\}$, put $\zeta_b=w_{1b}+kw_{2b}$. Adjoin two real coordinates $\phi_1,\phi_2$ and put $\Phi=\phi_1+k\phi_2$. Extend $Z_\ell$, $\ell\in\{1,2\}$, by zero in the $\phi_1,\phi_2$ directions, so the corresponding rows and columns of $B_0$ are zero. Define $\Psi:\R^3\times\mathbb C_k^3\to\R^3\times\mathbb C_k^3$ by
\[
 (\mathsf X,\zeta_0,\zeta_1,\Phi)\longmapsto(\mathsf X,Y,P,S),
 \qquad (Y,P,S)=(\zeta_0^2,\zeta_0\zeta_1,\Phi-\zeta_1^2/2).
\]
Write $Y=Y_1+kY_2$, $P=P_1+kP_2$ and $S=S_1+kS_2$ with real components. The real coordinate lists for $\mathsf z=\Psi(\mathsf y)$ are
\begin{align*}
 \mathsf y&=(t,x,s_0,w_{10},w_{11},w_{20},w_{21},\phi_1,\phi_2),\\
 \mathsf z&=(t,Y_1,Y_2,\ x,P_1,P_2,\ s_0,S_1,S_2).
\end{align*}
For $\nu\in\{t,x,s_0\}$, write $\omega_\nu=\omega_{\nu,1}+k\omega_{\nu,2}$, where $(\omega_t,\omega_x,\omega_{s_0})=(Y,P,S)$. Thus $\mathsf z$ is ordered by the triples $(\mathsf X^\nu,\omega_{\nu,1},\omega_{\nu,2})$ in the order $t,x,s_0$, as in~\eqref{cur:hess}.
Use $A,B\in\{1,\ldots,9\}$ to index $\mathsf y$ and $I,J\in\{1,\ldots,9\}$ to index $\mathsf z$. In particular,
\[
 \Psi^I=\mathsf z^I\circ\Psi,\qquad (D\Psi)_{IA}=\partial_{\mathsf y^A}\Psi^I.
\]
Restrict $\Psi$ to the open set
\[
 \bigl\{(\mathsf X,\zeta_0,\zeta_1,\Phi)\in\R^3\times\mathbb C_k^3:
         \zeta_0\ne0\bigr\}.
\]
On this set $Y\circ\Psi=\zeta_0^2\ne0$. With $W=P/Y$ for $Y\ne0$, the coordinate formulas give
\[
 \left(S+\frac{P^2}{2Y}\right)\circ\Psi=\Phi,\qquad
 W\circ\Psi=\frac{\zeta_1}{\zeta_0}.
\]
The complex differential on the auxiliary coordinates is
\[
 \frac{\partial(Y,P,S)}{\partial(\zeta_0,\zeta_1,\Phi)}
 =\begin{pmatrix}2\zeta_0&0&0\\ \zeta_1&\zeta_0&0\\0&-\zeta_1&1\end{pmatrix}.
\]
Its determinant is $2\zeta_0^2$, so $\Psi$ is locally invertible at every point with $\zeta_0\ne0$. The real volume Jacobian is $J_\Psi=|\det D\Psi|=4|\zeta_0|^4$. Choose an open set $\mathcal O\subset\{\mathsf y:\zeta_0\ne0\}$ on which $\Psi$ is a diffeomorphism onto $\Psi(\mathcal O)$. Define $\widetilde B$ on $\Psi(\mathcal O)$ by
\[
 \widetilde B\circ\Psi=J_\Psi^{-1}D\Psi\,B_0\,D\Psi^{\mathsf t}.
\]
The factor $J_\Psi^{-1}$ accounts for the change of volume $d\mathsf z=J_\Psi\,d\mathsf y$ under $\mathsf z=\Psi(\mathsf y)$.

Fix $J\in\{1,\ldots,9\}$ and $\varphi\in C_c^\infty(\Psi(\mathcal O);\R)$. To compute the divergence of the column $(\widetilde B_{IJ})_{I=1}^9$ in the coordinates $\mathsf z$, we change variables and use $\operatorname{div}B_0=0$ in the coordinates $\mathsf y$ to obtain
\begin{align*}
 \int_{\Psi(\mathcal O)}\sum_{I=1}^9
     \widetilde B_{IJ}\partial_{\mathsf z^I}\varphi\,d\mathsf z
  &=\int_{\mathcal O}\sum_{A,B=1}^9(B_0)_{AB}
      \,\partial_{\mathsf y^A}(\varphi\circ\Psi)
      \,\partial_{\mathsf y^B}\Psi^J\,d\mathsf y\\
  &=-\int_{\mathcal O}(\varphi\circ\Psi)\sum_{A,B=1}^9(B_0)_{AB}
       \,\partial_{\mathsf y^A}\partial_{\mathsf y^B}\Psi^J\,d\mathsf y.
\end{align*}
Only second derivatives in $w_1,w_2$ can contribute, since $\Psi$ is linear in $\mathsf X,\Phi$ and its nonlinear part depends on $\zeta_0,\zeta_1$ alone. For $a,b\in\{0,1\}$, the $(a,b)$ entry of both the $w_1$ and the $w_2$ diagonal blocks of $B_0$ is $\mathrm{Re}((G_0)_a\overline{(G_0)_b})$. The mixed $w_1,w_2$ block is zero. If $\Psi^J$ is one of $t,x,s_0$, its second derivatives in $w_1,w_2$ vanish. For the real and imaginary components of $(Y,P,S)\circ\Psi$, we apply the Cauchy--Riemann equations twice. In both cases,
\[
 \partial_{w_{1a}}\partial_{w_{1b}}\Psi^J
 +\partial_{w_{2a}}\partial_{w_{2b}}\Psi^J=0,
 \qquad a,b\in\{0,1\},\quad J\in\{1,\ldots,9\}.
\]
Thus for each $J\in\{1,\ldots,9\}$,
\[
 B_0:D_{\mathsf y}^2\Psi^J
 =\sum_{a,b=0}^1\operatorname{Re}\bigl((G_0)_a\overline{(G_0)_b}\bigr)
 \bigl(\partial_{w_{1a}}\partial_{w_{1b}}
       +\partial_{w_{2a}}\partial_{w_{2b}}\bigr)\Psi^J
 =0.
\]
The integral against $\varphi$ therefore vanishes, so $\sum_{I=1}^9\partial_{\mathsf z^I}\widetilde B_{IJ}=0$ on $\Psi(\mathcal O)$ for every $J\in\{1,\ldots,9\}$.

\medskip
It remains to identify $\widetilde B$ with $(\chi(c)\chi(c)^*)^{\mathbb R}$. In the coordinates $\mathsf y$, write the components of $Z_\ell$ for $\ell\in\{1,2\}$ as
\[
 Z_\ell=\sum_{A=1}^9 Z_\ell^A\partial_{\mathsf y^A},\qquad
 Z_\ell^A=Z_\ell(\mathsf y^A)\in\mathbb C_i
 \quad (A\in\{1,\ldots,9\}).
\]
For $A\in\{1,\ldots,9\}$, put $h_A=Z_1^A+kZ_2^A\in\mathbb H$. The real planes $\mathbb C_i$ and $k\mathbb C_i$ are orthogonal, and left multiplication by $k$ preserves the real inner product. Thus for all $A,B\in\{1,\ldots,9\}$,
\[
 \operatorname{Re}(h_A\overline{h_B})
 =\operatorname{Re}\bigl(Z_1^A\overline{Z_1^B}
                      +Z_2^A\overline{Z_2^B}\bigr)
 =(B_0)_{AB}.
\]
For $I\in\{1,\ldots,9\}$, set
\[
 g_I=J_\Psi^{-1/2}\sum_{A=1}^9
               (\partial_{\mathsf y^A}\Psi^I)h_A .
\]
Substituting $\widetilde B\circ\Psi=J_\Psi^{-1}D\Psi\,B_0\,D\Psi^{\mathsf t}$ gives
$(\widetilde B_{IJ})\circ\Psi=\operatorname{Re}(g_I\overline{g_J})$ for all $I,J\in\{1,\ldots,9\}$.

For $\nu\in\{t,x,s_0\}$ and $b\in\{0,1\}$, differentiating $\omega_\nu\circ\Psi$ with respect to $\zeta_b$ gives
\[
 L_{\nu b}
 =\frac{\partial(\omega_\nu\circ\Psi)}{\partial\zeta_b}
 =\bigl(C^\nu(w_1+kw_2)\bigr)_b.
\]
In the $\mathsf y$ coordinates, the values of $h_A$ for $\mathsf y^A=w_{1b},w_{2b}$ are $(G_0)_b$ and $k(G_0)_b$, respectively. For the coordinate $\mathsf y^A=\mathsf X^\nu$, expanding $Z_1$ and $Z_2$ gives
\[
 U_\nu=Z_1(\mathsf X^\nu)+kZ_2(\mathsf X^\nu)
      =\sum_{b=0}^1 L_{\nu b}\,i(G_0)_b,
 \qquad U=(U_t,U_x,U_{s_0})^{\mathsf t}\in\mathbb H^3.
\]
Write $L_{\nu b}=\ell_1+k\ell_2$ with $\ell_1,\ell_2\in\R$. The block of $D\Psi$ with columns $w_{1b},w_{2b}$ and rows $\omega_{\nu,1},\omega_{\nu,2}$ acts on the corresponding $h_A$ by
\[
 \begin{pmatrix}\ell_1&-\ell_2\\ \ell_2&\ell_1\end{pmatrix}
 \begin{pmatrix}(G_0)_b\\ k(G_0)_b\end{pmatrix}
 =
 \begin{pmatrix}\bar L_{\nu b}(G_0)_b\\ k\bar L_{\nu b}(G_0)_b\end{pmatrix}.
\]
For each $\nu\in\{t,x,s_0\}$, use $\mathsf X\circ\Psi=\mathsf X$ for the $\mathsf X^\nu$ coordinate and $\bar L_{\nu b}=-iL_{\nu b}i$, $ki=j$ for the two $\omega_\nu$ components. This gives
\[
 {\renewcommand{\arraystretch}{1.18}
 \begin{array}{c|ccc}
 \partial_{\mathsf z^I}
 &\partial_{\mathsf X^\nu}&\partial_{\omega_{\nu,1}}&
                              \partial_{\omega_{\nu,2}}\\[.25em]
 g_I & J_\Psi^{-1/2}U_\nu & -iJ_\Psi^{-1/2}U_\nu &
                             -jJ_\Psi^{-1/2}U_\nu
 \end{array}} .
\]
By the real representation following~\eqref{cur:hess}, these values give
\[
 \widetilde B\circ\Psi=\bigl(J_\Psi^{-1}UU^*\bigr)^{\mathbb R}.
\]

To compute $J_\Psi^{-1/2}U$, fix $\xi\in\Lambda$ and recall that $(\kappa,v)=(\kappa_\xi,v_\xi)$ and $p=d\theta_\xi$. The differential $L$ and the shift $\mathsf X\mapsto\mathsf X-iq$ give
\begin{align*}
 L(1,v)^{\mathsf t}
    &=(2\zeta_0,\zeta_1+\zeta_0v,-\zeta_1v)^{\mathsf t}
       =2\zeta_0(V^v\circ\Psi),\\
 p\cdot q
    &=\tfrac12\sum_{\ell=1}^2 w_\ell^{\mathsf t}C(p)w_\ell
     =-\tfrac{\kappa}{2}|\zeta_1-v\zeta_0|^2
     =-\tfrac{\kappa|\zeta_0|^2}{2}|W\circ\Psi-v|^2 .
\end{align*}
At $\mathsf X-iq$, the exponential for $\xi$ is $e^{i\theta_\xi(\mathsf X-iq)}=e^{i\theta_\xi(\mathsf X)}e^{p\cdot q}$. Conjugation in $G_0$ preserves the real factor $e^{p\cdot q}$ and conjugates $c_\xi e^{i\theta_\xi(\mathsf X)}$. Put $d=\zeta_0/|\zeta_0|$. Since $2\zeta_0/\sqrt{J_\Psi}=d/\sqrt{|Y\circ\Psi|}$, summing over $\xi\in\Lambda$ gives
\[
 J_\Psi^{-1/2}U=\sum_{\xi\in\Lambda}
   \bigl((V^{v_\xi}m_\xi)\circ\Psi\bigr)\,di\,
   \overline{c_\xi e^{i\theta_\xi(\mathsf X)}}.
\]
Right multiplication of each $g_I$ by the unit $(di)^{-1}$ preserves $\operatorname{Re}(g_I\overline{g_J})$. Since $Y\circ\Psi=\zeta_0^2$, we have $d^2=e\circ\Psi$. Also $d i d^{-1}=d^2i$ because $d\in\mathbb C_k$. For $u=\alpha+i\beta\in\mathbb C_i$ with $\alpha,\beta\in\R$,
\[
 (di)\bar u(di)^{-1}=\alpha-(d i d^{-1})\beta
                     =\alpha-(e\circ\Psi)i\beta=\rho_{e\circ\Psi}(u).
\]
By~\eqref{cur:scalar}, $J_\Psi^{-1/2}U(di)^{-1}=\chi(c)\circ\Psi$. Hence $\widetilde B=(\chi(c)\chi(c)^*)^{\mathbb R}$ on $\Psi(\mathcal O)$. For each $\mathsf z$ with $Y\ne0$, the two preimages share $\mathsf X$ and are given by
\[
 \zeta_0^2=Y,\qquad
 \zeta_1=P/\zeta_0,\qquad
 \Phi=S+P^2/(2Y).
\]
Both choices of $\zeta_0$ give the same tensor $(\chi(c)\chi(c)^*)^{\mathbb R}$ at $\mathsf z$. Thus the real columns of $(\chi(c)\chi(c)^*)^{\mathbb R}$ have zero divergence throughout $\{\mathsf z:Y\ne0\}$.

Finally, on $\{\mathsf z:Y\ne0\}$, recall that $\Phi=S+P^2/(2Y)$ and $Q=\Phi/a$. Since $\partial\Phi\,V^v=-W^2/2+W(v+W)/2-Wv/2=0$, the column $\chi(c)$ lies in $\ker\partial\Phi$. For a real function of $(Y,P,S)$ alone, $G_d(f)=2(\partial_{\bar\omega_d} f)i$. Applying this identity and the chain rule gives
\[
 G(\chi_{\mathrm{fib}}(Q))=(\partial\Phi)^*\tfrac{2}{a}(\partial_{\bar Q}\chi_{\mathrm{fib}})(Q)i .
\]
Thus $\chi(c)^*G(\chi_{\mathrm{fib}}(Q))=0$. To see why this
suffices in the real conservation law, put $f=\chi_{\mathrm{fib}}(Q)$
and $A=\chi(c)\chi(c)^*$. Use the coordinates $\mathsf z$ with the three groups indexed by $d=0,1,2$ in the order $t,x,s_0$, as in~\eqref{cur:hess}. Then
\[
 (A^{\mathbb R}\nabla f)_{d\lambda}
 =\sum_{b,\mu}\operatorname{Re}(\bar e_\lambda A_{db}e_\mu)
                    \partial_{b\mu}f
 =\operatorname{Re}\big(\bar e_\lambda \chi_d(c)\,\chi(c)^*G(f)\big)=0 .
\]
The product rule proves conservation after multiplication by $\chi_{\mathrm{fib}}(Q)$ and the range remains in $E$. The constant orthogonal change from $(t,x,s_0;Y,P,S)$ to $(T,x,b_0;\tau,P,Z)$ in~\eqref{cur:frames} preserves column divergence. In particular $V_T^v=(V_t^v-V_{s_0}^v)/\sqrt2=(1+Wv/2)/\sqrt2$.
\end{proof}

\subsection{Tracial polarization and determinant cancellation}

To prove Proposition~\ref{cur:cancel}, we polarize the scalar currents of Lemma~\ref{cur:conserve} and show that the determinant sums for $\M$ and $\M^{\mathrm{op}}$ cancel pointwise.

To polarize the currents, return to the coefficients $A_\xi\in L^4(\M)$, $\xi\in\Lambda$, in the cap-localized wave of Section~\ref{pre:periodic-cap}. We carry out the construction for each choice $\mathcal A\in\{\M,\M^{\mathrm{op}}\}$. Define the coefficients $a_\xi\in L^4(\mathcal A)$ for $\xi\in\Lambda$ and the wave $u_{\mathcal A}$ by
\begin{equation*}
 a_\xi=
 \begin{cases}
 A_\xi, & \mathcal A=\M,\\
 A_\xi^{\mathrm{op}}, & \mathcal A=\M^{\mathrm{op}},
 \end{cases}
 \qquad
 u_{\mathcal A}=\sum_{\xi\in\Lambda}e^{i\theta_\xi}a_\xi\,.
\end{equation*}
We construct scalar currents from the trace coordinates of the products $a_\eta^*a_\xi$, $\xi,\eta\in\Lambda$.

The opposite algebra has the same underlying complex vector space and adjoint as \(\M\), with \(a^{\mathrm{op}}b^{\mathrm{op}}=(ba)^{\mathrm{op}}\). Thus $(e^{i\theta_\xi}A_\xi)^{\mathrm{op}}=e^{i\theta_\xi}A_\xi^{\mathrm{op}}$ for each $\xi\in\Lambda$. Give $\M^{\mathrm{op}}$ the transported trace $\mathrm{tr}_{\M^{\mathrm{op}}}(a^{\mathrm{op}})=\mathrm{tr}_{\M}(a)$ for $a\ge0$, and write $\mathrm{tr}_{\mathcal A}$ for the trace in either algebra.

The polar decomposition gives $\mathrm{tr}_{\mathcal M}(|a^*|^p)=\mathrm{tr}_{\mathcal M}(|a|^p)$, so for $1\le p<\infty$ the map $a\mapsto a^{\mathrm{op}}$ is an isometry from $L^p(\M)$ onto $L^p(\M^{\mathrm{op}})$ since $|a^{\mathrm{op}}|=|a^*|^{\mathrm{op}}$. The same identity preserves the operator norm for $p=\infty$. One may first apply these identities to bounded elements in finite-trace corners, then use the \(L^p\) truncations of Section~\ref{sec:prelim}. By H\"older's inequality and the $L^4$ isometry, products of $L^4$ truncations converge in $L^2$ in both algebras. Passing to the limit using the $L^2$ isometry gives, for every $a,b\in L^4(\M)$,
\[
 (ab)^{\mathrm{op}}=b^{\mathrm{op}}a^{\mathrm{op}}
 \qquad\text{in }L^2(\M^{\mathrm{op}}).
\]
Since $u_{\M^{\mathrm{op}}}=(u_\M)^{\mathrm{op}}$, the isometry at $p=4$ gives $\|u_{\M^{\mathrm{op}}}\|_{L^4(\M^{\mathrm{op}})}=\|u_\M\|_{L^4(\M)}$ at each spacetime point.

For $\xi,\eta\in\Lambda$, H\"older's inequality for the trace gives $a_\eta^*a_\xi\in L^2(\mathcal A)$. For every $H\in L^2(\mathcal A)$, both $H a_\eta^*a_\xi$ and $a_\eta^*a_\xi H$ belong to $L^1(\mathcal A)$. Since $a_\eta^*,a_\xi\in L^4(\mathcal A)$ and $1/2+1/4+1/4=1$, the cyclic trace identity in Section~\ref{pre:tracial} gives
\[
 \mathrm{tr}_{\mathcal A}(H a_\eta^*a_\xi)
 =\mathrm{tr}_{\mathcal A}(a_\eta^*a_\xi H).
\]

For each choice of $\mathcal A$, all the required trace coordinates lie in the Hilbert space
\[
 Z_{\mathcal A}=\operatorname{span}_{\mathbb C}
 \{a_\eta^*a_\xi:\xi,\eta\in\Lambda\}\subset L^2(\mathcal A),\qquad
 \langle B,C\rangle=\mathrm{tr}_{\mathcal A}(B^*C).
\]
The space $Z_{\mathcal A}$ is finite dimensional and closed under adjoints, since $(a_\eta^*a_\xi)^*=a_\xi^*a_\eta$. We use its Hilbert space structure to take coordinates. Multiplication of two elements of $Z_{\mathcal A}$ takes place in the ambient algebra of measurable operators and gives an $L^1$ trace pairing.

Choose a real orthonormal basis $\mathcal H_{\mathcal A}$ of the self-adjoint part of $Z_{\mathcal A}$ for the real inner product $\mathrm{tr}_{\mathcal A}(HP)$. Every element of $Z_{\mathcal A}$ is the sum of its self-adjoint real and imaginary parts with coefficients $1$ and $i$. Thus $\mathcal H_{\mathcal A}$ also forms a complex orthonormal basis of $Z_{\mathcal A}$. For each $H\in\mathcal H_{\mathcal A}$ form the matrix $k_H:\mathbb C^\Lambda\to\mathbb C^\Lambda$ with entries
\[
 k_H(\xi,\eta)=\mathrm{tr}_{\mathcal A}(H a_\eta^*a_\xi),\qquad
 \xi,\eta\in\Lambda.
\]
Since $H=H^*$, we have for $\xi,\eta\in\Lambda$,
\[
 \overline{k_H(\xi,\eta)}
 =\mathrm{tr}_{\mathcal A}(a_\xi^*a_\eta H)
 =\mathrm{tr}_{\mathcal A}(H a_\xi^*a_\eta)
 =k_H(\eta,\xi).
\]
Take a spectral decomposition $k_H=\sum_l\lambda_l c_lc_l^*$, where $\lambda_l\in\mathbb R$ and $c_l\in\mathbb C^\Lambda$ for $1\le l\le|\Lambda|$. Define the quaternionic Hermitian field $A_H=\sum_l\lambda_l A(c_l)$ on the chart~\eqref{cur:chart}. The vectors $H,c_l$ and the numbers $\lambda_l$ are independent of the lifted spacetime variables, so Lemma~\ref{cur:conserve} gives
\[
 \operatorname{div}(A_H^{\mathbb R})=0,\qquad
 (\partial\Phi)A_H=0.
\]
Here the divergence is taken columnwise, and $(\partial\Phi)A_H=0$ means that the range of $A_H$ lies in $E=\ker\partial\Phi$. Section~\ref{en:section} solves one scalar Dirichlet equation for each $H\in\mathcal H_{\mathcal A}$, with source $(A_H)_{TT}$. Since the basis is independent of spacetime, the solutions $v_H$ of those scalar equations can be combined as $\sum_{H\in\mathcal H_{\mathcal A}}v_H H$.

We give the entries of $A_H$ to identify both its time component and its determinant on $E$. At each point of the chart, define a column $D\in\mathbb C_k^3$, a $3\times2$ matrix $J_0$ and two elements $X_0,X_1\in L^4(\mathcal A)$ by
\begin{equation*}
 D=(0,1/2,-W/2)^{\mathsf t},\quad J_0=[\,V^0\ D\,],\qquad
 X_m=\sum_{\xi\in\Lambda}v_\xi^m m_\xi e^{i\theta_\xi}a_\xi
 \quad(m=0,1)\,.
\end{equation*}
Since $V^v=V^0+vD$, the expansion of the outer products in~\eqref{cur:current} gives
\begin{equation}\label{cur:polar}
 A_H=\chi_{\mathrm{fib}}(Q)J_0\rho_e\!\left(
       \big(\mathrm{tr}_{\mathcal A}(H X_n^*X_m)\big)_{m,n=0,1}
                       \right)J_0^*\,.
\end{equation}
To verify the order of the factors, expansion of the scalar currents in the spectral sum for \(A_H\) gives
\[
 A_H=\chi_{\mathrm{fib}}(Q)
 \sum_{\xi,\eta\in\Lambda}m_\xi m_\eta
 V^{v_\xi}\rho_e\big(e^{i(\theta_\xi-\theta_\eta)}k_H(\xi,\eta)\big)(V^{v_\eta})^* .
\]
Only the scalar amplitudes are combined inside \(\rho_e\). Substituting \(V^{v_\xi}=V^0+v_\xi D\) gives \eqref{cur:polar}, since
\begin{equation*}
 \mathrm{tr}_{\mathcal A}(H X_n^*X_m)
  =\sum_{\xi,\eta\in\Lambda}v_\xi^m v_\eta^n m_\xi m_\eta
                   e^{i(\theta_\xi-\theta_\eta)}k_H(\xi,\eta)\,.
\end{equation*} Here \(\rho_e\) acts entry by entry on the \(2\times2\) numerical matrix. In particular, $A_H$ depends only on $k_H$, so any of its spectral decompositions gives the same current.

For real $c,d$, write $\mathrm{Re}_k(c+kd)=c$ and $\mathrm{Im}_k(c+kd)=d$. Fix the auxiliary variables. On spatial trigonometric polynomials with frequencies in $\Lambda$ and scalar or $L^4(\mathcal A)$ coefficients, let $S_0,S_1$ act into the same polynomial space by the respective scalar symbols
\begin{equation*}
 m_\xi\mathrm{Re}_k V_T^{v_\xi}
      =m_\xi(1+v_\xi\mathrm{Re}_k W/2)/\sqrt 2,\qquad
 m_\xi\mathrm{Im}_k V_T^{v_\xi}
      =m_\xi v_\xi\mathrm{Im}_k W/(2\sqrt 2)\,.
\end{equation*} The time components $\mu_H\in\R$ and the density $\mu_{\mathcal A}\in L^2(\mathcal A)_+$ satisfy
\begin{equation}\label{cur:density}
 \mu_H=(A_H)_{TT}=\mathrm{tr}_{\mathcal A}(H\mu_{\mathcal A}),\qquad
 \mu_{\mathcal A}=\chi_{\mathrm{fib}}(Q)\sum_{j=0}^1(S_j u_{\mathcal A})^*S_j u_{\mathcal A}\,.
\end{equation}
For scalar amplitudes $c\in\mathbb C_i^\Lambda$, put $u_c=\sum_{\xi\in\Lambda}c_\xi e^{i\theta_\xi}$. The time component is $\chi(c)_T=\rho_e(S_0u_c)+k\rho_e(S_1u_c)$. The two summands lie in orthogonal real planes, obtained from $\mathbb C_i,k\mathbb C_i$ by conjugation by a square root of $e$ in $\mathbb C_k$. Since $\rho_e$ preserves the modulus,
\[
 |\chi(c)_T|^2=|S_0u_c|^2+|S_1u_c|^2 .
\]
To compute $\mu_H$, take the $TT$ entry of $A_H=\sum_l\lambda_l A(c_l)$. Expanding $|S_j u_{c_l}|^2$ in the frequency pairs $\xi,\eta\in\Lambda$ and substituting $k_H(\xi,\eta)=\sum_l\lambda_l(c_l)_\xi\overline{(c_l)_\eta}$ gives
\begin{equation*}
 \mu_H
  =\chi_{\mathrm{fib}}(Q)\sum_{j=0}^1\sum_{l=1}^{|\Lambda|}
       \lambda_l|S_j u_{c_l}|^2
 =\chi_{\mathrm{fib}}(Q)\sum_{j=0}^1
      \mathrm{tr}_{\mathcal A}\bigl(H(S_j u_{\mathcal A})^*S_j u_{\mathcal A}\bigr).
\end{equation*}
This proves~\eqref{cur:density}.

\medskip
With~\eqref{cur:density} established, we use~\eqref{cur:polar} to compute the sum of the determinants of $A_H$ on $E=\ker\partial\Phi$. The geometric weight in this sum comes from the Hessian of $|\Phi|^2$. By holomorphy, \(\partial_{\omega_b}|\Phi|^2=\bar\Phi\,\partial_{\omega_b}\Phi\) and \(\partial_{\bar\omega_d}\partial_{\omega_b}\Phi=0\). Hence \eqref{cur:wirtinger}--\eqref{cur:hess} give
\[
 B=\mathsf H(|\Phi|^2)=4(\partial\Phi)^*\partial\Phi,\qquad
 \beta=\mathrm{tr}_{\mathbb H}B
 =4\bigl(1+|W|^2+|W|^4/4\bigr)=4(1+|W|^2/2)^2 .
\]
The kernel of \(B\) is the right quaternionic plane \(E=\ker\partial\Phi\).
In this section and in Section~\ref{en:section}, $\beta$ without an argument denotes the numerical trace of $B$ at the point under consideration. We reserve $\beta(p)$ for the logarithmic exponent in~\eqref{intro:annular}. Write $\alpha_H=A_H|_E$ for the compression to \(E\) in an orthonormal quaternionic basis. In general \(L|_E\) means the matrix of \(P_E L P_E\) on \(E\), where \(P_E:\mathbb H^3\to E\) is the orthogonal projection.

For a two-by-two Hermitian quaternionic matrix with diagonal entries $r,s\in\R$ and off-diagonal entry $q\in\mathbb H$, the Moore determinant is $rs-|q|^2$. For a quaternionic Hermitian \(n\times n\) matrix \(L\), \(n\in\mathbb N\), we denote the product of its real eigenvalues by \(\det_{\mathbb H}L\). To see that \(L\) has an orthonormal eigenbasis, maximize the real Rayleigh quotient \(v^*Lv\) on the unit sphere in \(\mathbb H^n\). At a maximizer \(v\) with value \(\lambda\in\R\), real differentiation gives \(\operatorname{Re}(w^*(Lv-v\lambda))=0\) for every \(w\in\mathbb H^n\), so \(Lv=v\lambda\). The right quaternionic complement \(v^\perp\) is invariant by self-adjointness, and induction diagonalizes \(L\). In dimension two the equality of the product of the eigenvalues with \(rs-|q|^2\) also follows from \(((\mathrm{tr}_{\mathbb H}L)^2-\mathrm{tr}_{\mathbb H}(L^2))/2\).

\begin{proposition}\label{cur:cancel}
At every point of the fiber chart \eqref{cur:chart}, each of the waves $u_{\mathcal A}$ in $\mathcal A=\M$ and $\mathcal A=\M^{\mathrm{op}}$ satisfies
\begin{equation}\label{cur:ordered}
 \sum_{H\in\mathcal H_{\mathcal A}}\frac{\det_{\mathbb H}\alpha_H}{\beta}
  =\frac{\chi_{\mathrm{fib}}(Q)^2}{16}
       \left(\|X_0X_1^*\|_{L^2(\mathcal A)}^2
                     -\|X_1^*X_0\|_{L^2(\mathcal A)}^2\right)\,.
\end{equation}
The sum of~\eqref{cur:ordered} over $\mathcal A=\M,\M^{\mathrm{op}}$ is zero.
\end{proposition}
\begin{proof}
To prove~\eqref{cur:ordered}, we separate the geometry of $E$ from the dependence on the Fourier coefficients. In~\eqref{cur:polar}, $J_0$ contributes a Gram determinant common to all $H$ and both algebras. The remaining two-by-two determinant involves only trace coordinates of $X_n^*X_m$, which we sum by Parseval's identity. We then reverse the products to compare the two algebras.

The columns \(V^0,D\) of \(J_0\) are independent and lie in the two-dimensional plane \(E\). In an orthonormal frame of \(E\), denote their square coordinate matrix by \(J_E\) and the middle matrix in \eqref{cur:polar} by \(L_H=\rho_e((\mathrm{tr}_{\mathcal A}(H X_n^*X_m))_{m,n})\). Then \(\alpha_H=\chi_{\mathrm{fib}}(Q)J_EL_HJ_E^*\). We will use the congruence formula
\[
 \det_{\mathbb H}(J_E L_H J_E^*)=\det_{\mathbb H}(J_E^*J_E)\det_{\mathbb H}L_H .
\]
Since \(J_E^*J_E=J_0^*J_0\) and \((V^0)^*D=\bar W/4\), the Gram determinant of \(V^0,D\) is
\begin{equation*}
 \det_{\mathbb H}(J_E^*J_E)=|V^0|^2|D|^2-|(V^0)^*D|^2
 =\left(1+\frac{|W|^2}{4}\right)\frac{1+|W|^2}{4}
                       -\frac{|W|^2}{16}=\frac{\beta}{16}\,.
\end{equation*}
To verify the congruence formula even for indefinite $L_H$, diagonalize $L_H$ by a quaternionic unitary and multiply $J_E$ on the right by that unitary. Let $\lambda_1,\lambda_2\in\R$ be the eigenvalues in this diagonalization. For the transformed columns $l_1,l_2\in\mathbb H^2$, choose the first basis vector in their span along $l_1$. Writing $l_2=(c,d)^{\mathsf t}$ in this basis, the determinant of $\lambda_1l_1l_1^*+\lambda_2l_2l_2^*$ is \((\lambda_1|l_1|^2+\lambda_2|c|^2)\lambda_2|d|^2-|\lambda_2 c\bar d|^2=\lambda_1\lambda_2|l_1|^2|d|^2\). The quantity $|l_1|^2|d|^2$ is the determinant of the column Gram matrix, invariant under the preceding unitary conjugation.

Substituting $\alpha_H=\chi_{\mathrm{fib}}(Q)J_EL_HJ_E^*$ into the congruence formula gives, for $H\in\mathcal H_{\mathcal A}$,
\[
 \frac{\det_{\mathbb H}\alpha_H}{\beta}
 =\frac{\chi_{\mathrm{fib}}(Q)^2}{16}\det_{\mathbb H}L_H.
\]
The factor $\chi_{\mathrm{fib}}(Q)^2$ comes from scaling a $2\times2$ matrix. It remains to sum $\det_{\mathbb H} L_H$ in the basis $\mathcal H_{\mathcal A}$. The embedding $\rho_e$ preserves real parts, conjugation and modulus. The fixed basis is real orthonormal on the self-adjoint part and complex orthonormal on all of $Z_{\mathcal A}$. Applying Parseval's identity to $X_m^*X_m$, $m=0,1$, and to $X_1^*X_0$ gives
\begin{align*}
 \sum_{H\in\mathcal H_{\mathcal A}}
 \mathrm{tr}_{\mathcal A}(HX_0^*X_0)\mathrm{tr}_{\mathcal A}(HX_1^*X_1)
 &=\mathrm{tr}_{\mathcal A}(X_0^*X_0X_1^*X_1)\,,\\
 \sum_{H\in\mathcal H_{\mathcal A}}
 |\mathrm{tr}_{\mathcal A}(HX_1^*X_0)|^2
 &=\|X_1^*X_0\|_{L^2(\mathcal A)}^2\,.
\end{align*}
Summing the identity for $\det_{\mathbb H}\alpha_H/\beta$ over $H\in\mathcal H_{\mathcal A}$ and using the two Parseval identities gives
\begin{align*}
 \frac{16}{\beta}\sum_{H\in\mathcal H_{\mathcal A}}\det_{\mathbb H}\alpha_H
  &=\chi_{\mathrm{fib}}(Q)^2\left(\mathrm{tr}_{\mathcal A}(X_0^*X_0X_1^*X_1)
                                   -\|X_1^*X_0\|_{L^2(\mathcal A)}^2\right).
\end{align*}
Here $X_n^*X_m\in Z_{\mathcal A}$, so the sums in Parseval's identity include all their trace coordinates. By cyclicity of the trace,
\[
 \mathrm{tr}_{\mathcal A}(X_0^*X_0X_1^*X_1)
 =\mathrm{tr}_{\mathcal A}(X_1X_0^*X_0X_1^*)
 =\|X_0X_1^*\|_{L^2(\mathcal A)}^2.
\]
For $m\in\{0,1\}$, write $X_m\in L^4(\M)$ for the value obtained with $a_\xi=A_\xi$. With $a_\xi=A_\xi^{\mathrm{op}}$ the corresponding value in $L^4(\M^{\mathrm{op}})$ is $X_m^{\mathrm{op}}$, since $a\mapsto a^{\mathrm{op}}$ is complex linear and the scalar coefficients $v_\xi^m m_\xi e^{i\theta_\xi}$ agree in both algebras. Applying $(ab)^{\mathrm{op}}=b^{\mathrm{op}}a^{\mathrm{op}}$ to the products $X_1^*X_0$ and $X_0X_1^*$ gives
\[
 X_0^{\mathrm{op}}(X_1^{\mathrm{op}})^*=(X_1^*X_0)^{\mathrm{op}},
 \qquad
 (X_1^{\mathrm{op}})^*X_0^{\mathrm{op}}=(X_0X_1^*)^{\mathrm{op}} .
\]
The isometry $L^2(\M)\to L^2(\M^{\mathrm{op}})$ then makes the right-hand sides of~\eqref{cur:ordered} for the two algebras negatives of each other. Thus the sum of~\eqref{cur:ordered} over $\mathcal A=\M,\M^{\mathrm{op}}$ vanishes pointwise.
\end{proof}

\begin{remark}
Section~\ref{en:section} adds the energy identities for $\M$ and $\M^{\mathrm{op}}$ to cancel the weighted determinant sum in~\eqref{cur:ordered}, which can be strictly negative in $\M$.

Take $\M=M_2$ with its ordinary trace and matrix units $E_{ij}$, $i,j\in\{1,2\}$. For $X_0=E_{11}$ and $X_1=E_{12}$ we have $X_0X_1^*=0$ and $X_1^*X_0=E_{21}$, so $\|X_0X_1^*\|_{L^2(M_2)}^2-\|X_1^*X_0\|_{L^2(M_2)}^2=-1$. To realize these values, fix a point with $\chi_{\mathrm{fib}}(Q)>0$ and take $\Lambda=\{\xi_1,\xi_2\}$ with $v_{\xi_1}\ne v_{\xi_2}$. At the chosen point, the amplitudes $a_{\xi_1},a_{\xi_2}\in M_2$ determine $X_0,X_1$ by
\[
 \begin{pmatrix}X_0\\ X_1\end{pmatrix}
 =
 \begin{pmatrix}
 m_{\xi_1}e^{i\theta_{\xi_1}} & m_{\xi_2}e^{i\theta_{\xi_2}}\\
 v_{\xi_1}m_{\xi_1}e^{i\theta_{\xi_1}} & v_{\xi_2}m_{\xi_2}e^{i\theta_{\xi_2}}
 \end{pmatrix}
 \begin{pmatrix}a_{\xi_1}\\a_{\xi_2}\end{pmatrix}.
\]
The determinant of this scalar matrix is $(v_{\xi_2}-v_{\xi_1})m_{\xi_1}m_{\xi_2}e^{i(\theta_{\xi_1}+\theta_{\xi_2})}\ne0$, so we can choose the two amplitudes to give $X_0=E_{11}$ and $X_1=E_{12}$. For instance, frequencies $(2\pi/100)(1,n)$ and $(2\pi/100)(-1,n)$ on the square torus $\R^2/(100\mathbb Z^2)$ lie in a common fixed cap for all sufficiently large $n\in\mathbb N$.
They have equal length, so choosing $r\in(c_0,C_0)$ and $N=|\xi|/r$
for either frequency places both on the prescribed annulus with $N\ge2$. At the chosen point, the weighted determinant sum in~\eqref{cur:ordered} is $-\chi_{\mathrm{fib}}(Q)^2/16$ in $M_2$ and $+\chi_{\mathrm{fib}}(Q)^2/16$ in $M_2^{\mathrm{op}}$ by Proposition~\ref{cur:cancel}.
\end{remark}

\section{The auxiliary Dirichlet energy: proof of Theorem~\ref{thm:four}}\label{en:section}

To prove the periodic single-cap estimate~\eqref{en:capestimate} of Proposition~\ref{en:cap}, we compare the fourth moment $q_u$ from~\eqref{pre:qu} with an auxiliary Dirichlet energy $\mathcal E$. We solve the Dirichlet equations with sources $\mu_H=(A_H)_{TT}$ from Section~\ref{cur:sec} and sum the energies for $\M$ and $\M^{\mathrm{op}}$. The determinant contributions cancel by Proposition~\ref{cur:cancel}. Green's formula then gives $\Delta_z\mathcal E=\mathcal B-\mathcal Q$, where $\Delta_z$ is the Laplacian in the three real coordinates of $z=(T,\tau)\in\R\times\mathbb C_k$, $\mathcal B$ is the boundary flux, and $\mathcal Q\ge0$ is the interior contribution.

To bound $\mathcal E(T,\tau)$ from below in terms of $q_u(T)$, we use the Gaussian factor in $m_\xi$ on $|W|\le\gamma$, where $|W-v_\xi|\le2\gamma$ for $\xi\in\Lambda$. On $|W|=R$, we have $|W-v_\xi|\ge R-\gamma$. Since $\kappa_\xi\asymp N$ and $|Y|\asymp|\tau|$, the Gaussian factor in~\eqref{cur:scalar} is bounded by $e^{-cN|\tau|}$ on this boundary for some $c>0$. With $\gamma$ sufficiently small, this decay gives a bound for $|\mathcal B|$ in terms of $\mathcal E$ and hence a differential inequality for $\mathcal E$. Two radial comparisons give the two logarithmic factors in~\eqref{en:capestimate}. We use the disk geometry and the annular comparison of~\cite[Appendix D and Lemma 6.1]{openai}.

\subsection{The fiberwise Dirichlet problem and its energy}\label{en:fiber}

We construct the smooth, nonnegative paired Dirichlet energy $\mathcal E$ in~\eqref{en:positive}. Green's formula~\eqref{en:balance} uses the expression for $\mathcal E$ as an integral over each fiber. The upper and lower bounds~\eqref{en:energyupper} and~\eqref{en:lower} use its expression as a quadratic form of the Dirichlet inverses. We first write the differential operator $L_B$ in disk coordinates, then solve the Dirichlet problem and pair its sources and solutions for $\M$ and $\M^{\mathrm{op}}$.

Recall $a=K/N$ from Section~\ref{cur:lift} and let $z=(T,\tau)$ range over $\mathcal U_a=\R\times\{\tau\in\mathbb C_k:a<|\tau|<4\}$. Set $R=1/10$ and define
\[
 D=\{W\in\mathbb C_k:|W|<R\},\qquad
 D_Q=\{Q\in\mathbb C_k:|Q|<1/10\}\,.
\]
Write $\T^2_X$ for the spatial torus $\R^2/\Gamma$ in the coordinates $(X_1,X_2)=X=(x,b_0)$. With $(W,Q)$ given by~\eqref{cur:chart-forward}, the fiber in the original flat coordinates is
\begin{equation*}
 \Omega_z=\bigl\{(X,P,Z)\in\T^2_X\times\mathbb C_k^2 : \tau+Z\ne0,\ (W,Q)\in D\times D_Q\bigr\}\,.
\end{equation*}
Its measure is $d\zeta=dX\,d^2P\,d^2Z$, where $dX$ is probability Haar measure and $d^2$ denotes planar Lebesgue measure.

To define the differential operator $L_B$ on each fiber $\Omega_z$, recall the quaternionic Hessian $\mathsf H$ from~\eqref{cur:hess} and the holomorphic function $\Phi=aQ$ from~\eqref{cur:chart-forward}. In the physical frame let
\[
 B=\mathsf H(|\Phi|^2)=4(\partial\Phi)^*\partial\Phi\,,
\]
where $\partial\Phi$ is the row of holomorphic derivatives in $(\tau,P,Z)$. A subscript $f$ denotes the principal block with indices $1,2$, corresponding to $X_1,X_2$. For integers $m\ge1$, let $\operatorname{Herm}_m(\mathbb H)$ be the real vector space of $m\times m$ quaternionic Hermitian matrices. Write $\mathrm{tr}_{\mathbb H}$ for the sum of the diagonal entries of a numerical quaternionic matrix. Define $C_2:\operatorname{Herm}_2(\mathbb H)\to\operatorname{Herm}_2(\mathbb H)$ by $C_2(L)=(\mathrm{tr}_{\mathbb H}L)I_2-L$. On each fiber the operator $L_B:C^2(\Omega_z;\R)\to C(\Omega_z;\R)$ is given by
\begin{equation}
 L_B v=C_2(B_f):\mathsf H_f(v),\qquad
 S:L=\operatorname{Re}\mathrm{tr}_{\mathbb H}(SL).
 \label{en:lb}
\end{equation}
Here $S:L$ is the real inner product on $\operatorname{Herm}_m(\mathbb H)$. Taking the real part of $\mathrm{tr}_{\mathbb H}$ gives a cyclic trace, by $\operatorname{Re}(uv)=\operatorname{Re}(vu)$ for quaternions $u,v$. We use $\mathrm{tr}_{\mathcal A}$ for the trace of a coefficient algebra $\mathcal A$.

At fixed $z$, write $\partial_{P,Z}=(\partial_P,\partial_Z)$, $c=(c_1,c_2)=\partial_{P,Z}\Phi$ and $d=(c_2,-c_1)^{\mathsf t}\in\mathbb C_k^2$. Recall $Y=(\tau+Z)/\sqrt2$ from~\eqref{cur:frames}.

\begin{lemma}\label{en:disk-chart}
For every $z\in\mathcal U_a$, the chart~\eqref{cur:chart-forward} identifies $\Omega_z$ with $\T^2_X\times D\times D_Q$. The complex Jacobian of $(P,Z)\mapsto(W,\Phi)$ at fixed $z$ is
\begin{equation}
 \mathscr J=\det_{\mathbb C_k}\frac{\partial(W,\Phi)}{\partial(P,Z)}
 =(\partial_{P,Z}W)d=\frac{1+W^2/2}{\sqrt2Y}.
 \label{en:jac}
\end{equation}
The volume form and coefficient bounds in this chart are
\begin{equation}
 d\zeta=a^2|\mathscr J|^{-2}\,dX\,d^2W\,d^2Q,\qquad
 |Y|\asymp |\tau|,\quad |d|\asymp1,\quad |\mathscr J|\asymp |\tau|^{-1}.
 \label{en:measure}
\end{equation}
The comparison constants depend only on $D$ and $D_Q$. For $v\in C^2(\Omega_z;\R)$, write $\widetilde v$ for its expression in $(X,W,Q)$ coordinates. In these coordinates,
\begin{equation}
 L_B v=4|\mathscr J|^2\Delta_W\widetilde v+
 4\sum_{i,j=1}^2\operatorname{Re}(d_i\bar d_j)\partial_{X_iX_j}\widetilde v\,.
 \label{en:disk}
\end{equation}
All derivatives are taken at fixed $z=(T,\tau)$ and $Q$. The real Laplacian $\Delta_W$ acts on the two real coordinates of $W$ with $X$ fixed, and $\partial_{X_iX_j}$ differentiates in $X_i$ and $X_j$ with $W$ fixed.
\end{lemma}

\begin{proof}
We first prove the Jacobian and measure formulas~\eqref{en:jac} and~\eqref{en:measure}, including the coefficient bounds. The inverse~\eqref{cur:chart-inverse} is defined for $(W,Q)\in D\times D_Q$ when $|\tau|>a$ and leaves $X$ fixed. At fixed $z$, the formulas $W=P/Y$ and $\Phi=(Z-\tau)/\sqrt2+P^2/(2Y)$ give
\[
 \partial_P W=Y^{-1},\qquad \partial_Z W=-W/(\sqrt2Y),\qquad
 c=\big(W,(1-W^2/2)/\sqrt2\big)\,.
\]
Expanding the complex determinant along the first row gives $(\partial_{P,Z}W)d$. Substituting the derivatives above gives the numerator $1-W^2/2+W^2=1+W^2/2$ in~\eqref{en:jac}. Since $\Phi=aQ$, the complex Jacobian for $(W,Q)$ is $\mathscr J/a$. Its squared modulus is the real Jacobian, which gives the volume form in~\eqref{en:measure}. In the inverse chart, $|\sqrt2\tau+aQ|$ lies between $(\sqrt2-1/10)|\tau|$ and $(\sqrt2+1/10)|\tau|$. On $D$, $|1\pm W^2/2|$ lies between $1-R^2/2$ and $1+R^2/2$. These bounds give the comparisons for $|Y|$, $|d|$, and $|\mathscr J|$ in~\eqref{en:measure}.

For~\eqref{en:disk}, we separate the terms in $\mathsf H_f(v)$ according to whether the derivatives are in $X$, in $(P,Z)$, or in both. Since $B_f=4c^*c$ and $c$ is a row over the commutative plane $\mathbb C_k$, the identity $|c|^2I_2-c^*c=dd^*$ gives $C_2(B_f)=4dd^*$. Thus $L_Bv=4\operatorname{Re}(d^*\mathsf H_f(v)d)$. The terms with one physical spatial and one auxiliary derivative have coefficients in $\mathbb C_k i$. Multiplication on either side by elements of $\mathbb C_k$ leaves that plane invariant, so their real parts vanish. The second derivatives in $(P,Z)$ contribute $16\,d^*(\partial_{\bar\omega_i}\partial_{\omega_j}v)_{i,j=1,2}d$, where $(\omega_1,\omega_2)=(P,Z)$. The coordinates $W,Q$ are holomorphic in $(P,Z)$ and satisfy
\begin{equation*}
 (\partial_{P,Z}W)d=\mathscr J,\qquad
 (\partial_{P,Z}Q)d=cd/a=0\,.
\end{equation*}
The mixed chain rule, contracted pointwise with $\bar d_i d_j$ and summed over $i,j\in\{1,2\}$, therefore reads
\begin{equation*}
 \sum_{i,j=1}^2\bar d_i
       (\partial_{\bar\omega_i}\partial_{\omega_j}v)d_j
 =|\mathscr J|^2\partial_{\bar W}\partial_W\widetilde v\,.
\end{equation*}
Indeed the mixed derivatives of the holomorphic coordinate functions vanish, and every term containing a $Q$ derivative has a factor $(\partial_{P,Z}Q)d$ or its conjugate. Thus the auxiliary derivatives contribute $16|\mathscr J|^2\partial_{\bar W}\partial_W\widetilde v=4|\mathscr J|^2\Delta_W\widetilde v$. The pure $X$ terms contribute $4\sum_{i,j=1}^2\operatorname{Re}(d_i\bar d_j)\partial_{X_iX_j}\widetilde v$, since the chart leaves $X$ fixed. Adding these contributions proves~\eqref{en:disk}.
\end{proof}

By Lemma~\ref{en:disk-chart}, we can solve the fiber equation at fixed $(z,Q)$ by taking Fourier coefficients in $X$. For each $\mathcal A\in\{\M,\M^{\mathrm{op}}\}$, recall from Section~\ref{cur:sec} the finite-dimensional Hilbert space $Z_{\mathcal A}$ and the real orthonormal basis $\mathcal H_{\mathcal A}$ of its self-adjoint part. For $H\in\mathcal H_{\mathcal A}$, let $\mu_H=(A_H)_{TT}$ be the real source from~\eqref{cur:density}, written in disk coordinates, and set $\mu_{\mathcal A}=\sum_{H\in\mathcal H_{\mathcal A}}\mu_H H$. A hat denotes a Fourier coefficient in $X$ with respect to $dX$. For $\ell=(\ell_1,\ell_2)\in\Gamma^*$, write $\ell\cdot d=\sum_{j=1}^2\ell_jd_j$.

\begin{lemma}\label{en:fiber-solve}
For each $\mathcal A\in\{\M,\M^{\mathrm{op}}\}$ and $H\in\mathcal H_{\mathcal A}$, there is a jointly smooth real function $v_H:\mathcal U_a\times\T^2_X\times\overline D\times D_Q\to\R$. For each fixed $(z,Q)\in\mathcal U_a\times D_Q$, the slice $v_H(z,\cdot,\cdot,Q)$ is the unique smooth solution of
\[
 L_Bv_H=\mu_H\quad\text{in }\T^2_X\times D,\qquad
 v_H|_{\T^2_X\times\partial D}=0\,,
\]
where $L_B$ is written in disk coordinates. The function $v_H$ vanishes wherever $\chi_{\mathrm{fib}}(Q)=0$.

Set $V_{\mathcal A}=\sum_{H\in\mathcal H_{\mathcal A}}v_H H$, with values in the self-adjoint part of $Z_{\mathcal A}$. For each $(z,Q)\in\mathcal U_a\times D_Q$ and $\ell\in\Gamma^*$,
\begin{equation}
 P_\ell=-\Delta_W+\frac{|\ell\cdot d|^2}{|\mathscr J|^2},\qquad
 f_\ell=\frac{\widehat\mu_{\mathcal A}(\ell)}{4|\mathscr J|^2},\qquad
 -\widehat V_{\mathcal A}(\ell)=P_\ell^{-1}f_\ell\,.
 \label{en:fourier}
\end{equation}
Here $P_\ell$ acts on complex $L^2(D)$ with domain $H^2(D)\cap H^1_0(D)$. Its inverse acts on each coordinate of $f_\ell:D\to Z_{\mathcal A}$. We suppress the dependence of $P_\ell$ on $(\tau,Q)$ and of $f_\ell$ on $(z,Q,\mathcal A)$ in~\eqref{en:fourier}. The bound $\|P_\ell^{-1}\|_{L^2(D)\to L^2(D)}\le\lambda_1(D)^{-1}$ is uniform in $(z,Q,\ell)$, where $\lambda_1(D)$ is the first eigenvalue of the scalar Dirichlet Laplacian on $D$. Moreover, $\widehat V_{\mathcal A}(\ell)=0$ for $\ell\in\Gamma^*\setminus(\Lambda-\Lambda)$.
\end{lemma}

\begin{proof}
We solve the Dirichlet problem one spatial Fourier mode at a time. Each product in the formula for $\mu_{\mathcal A}$ in~\eqref{cur:density} has spatial phases $e^{iX\cdot(\xi-\eta)}$ with $\xi,\eta\in\Lambda$. Thus $\widehat\mu_{\mathcal A}(\ell)=0$ for $\ell\in\Gamma^*\setminus(\Lambda-\Lambda)$. On the character $e^{iX\cdot\ell}$, the $X$ derivatives in~\eqref{en:disk} give the multiplier $-4|\ell\cdot d|^2$. After dividing by $4|\mathscr J|^2$, the equation to solve is $P_\ell(-\widehat V_{\mathcal A}(\ell))=f_\ell$ with the definitions in~\eqref{en:fourier}. At fixed $(z,Q,\ell)$ the potential of $P_\ell$ is real, smooth on $\overline D$ and nonnegative. The source $f_\ell$ is smooth up to $\partial D$ by~\eqref{cur:density} and the inverse chart~\eqref{cur:chart-inverse}.

The first eigenvalue of the scalar Dirichlet Laplacian on $D$ is
\begin{equation}
 \lambda_1(D)=\inf_{0\ne w\in H^1_0(D)}
     \frac{\int_D|\nabla_Ww|^2\,d^2W}{\int_D|w|^2\,d^2W}>0 .
 \label{en:gap}
\end{equation}
The Poincar\'e inequality on the bounded disk gives positivity,
and the compact embedding of \(H^1_0(D)\) in \(L^2(D)\) gives a minimizer
of the Rayleigh quotient. For $w\in H^1_0(D)$, the quadratic form of $P_\ell$ satisfies
\begin{equation*}
 \int_D\left(|\nabla_W w|^2+\frac{|\ell\cdot d|^2}{|\mathscr J|^2}|w|^2\right)d^2W
 \ge \int_D|\nabla_W w|^2\,d^2W
 \ge \lambda_1(D)\|w\|_{L^2(D)}^2\,.
\end{equation*}
For each scalar source $f\in L^2(D)$, the Lax--Milgram theorem on $H^1_0(D)$ gives a unique weak Dirichlet solution $w$. Elliptic regularity gives $w\in H^2(D)\cap H^1_0(D)$, so $w=P_\ell^{-1}f$. Integration by parts identifies the quadratic form at $w$ with $\operatorname{Re}\langle P_\ell w,w\rangle_{L^2(D)}$. Its lower bound and the Cauchy--Schwarz inequality give
\begin{equation*}
 \lambda_1(D)\|w\|_{L^2(D)}^2
 \le \operatorname{Re}\langle P_\ell w,w\rangle_{L^2(D)}
 = \operatorname{Re}\langle f,w\rangle_{L^2(D)}
 \le \|f\|_{L^2(D)}\|w\|_{L^2(D)}\,.
\end{equation*}
For $w\ne0$, dividing by $\lambda_1(D)\|w\|_{L^2(D)}$ gives $\|w\|_{L^2(D)}\le\lambda_1(D)^{-1}\|f\|_{L^2(D)}$. This estimate also holds for $w=0$, so $\|P_\ell^{-1}\|_{L^2(D)\to L^2(D)}\le\lambda_1(D)^{-1}$.

For $F\in L^2(D;Z_{\mathcal A})$, write $F=\sum_{H\in\mathcal H_{\mathcal A}}F_HH$ in the complex orthonormal basis $\mathcal H_{\mathcal A}$ of $Z_{\mathcal A}$, with $F_H\in L^2(D)$. Since $P_\ell^{-1}F=\sum_{H\in\mathcal H_{\mathcal A}}(P_\ell^{-1}F_H)H$, orthonormality and the scalar bound give
\begin{equation*}
\begin{split}
      \|P_\ell^{-1}F\|_{L^2(D;Z_{\mathcal A})}^2
  &=\sum_{H\in\mathcal H_{\mathcal A}}\|P_\ell^{-1}F_H\|_{L^2(D)}^2
 \le\lambda_1(D)^{-2}\sum_{H\in\mathcal H_{\mathcal A}}\|F_H\|_{L^2(D)}^2 \\
 &=\lambda_1(D)^{-2}\|F\|_{L^2(D;Z_{\mathcal A})}^2\,.
\end{split}
\end{equation*}
Taking square roots gives the bound $\lambda_1(D)^{-1}$ for the norm of $P_\ell^{-1}$ on $L^2(D;Z_{\mathcal A})$.

For the sources in~\eqref{en:fourier}, smoothness up to $\partial D$ follows from the classical Dirichlet existence and boundary regularity theorems~\cite[Sections 6.3--6.4]{gt}. Indeed, the second-order part of $-P_\ell$ is the Laplacian, its zeroth-order coefficient is nonpositive, and the potential and source are smooth, the disk has smooth boundary, and the Dirichlet condition is homogeneous. The classical solution agrees with the weak solution by uniqueness.

Since $P_\ell^{-1}f_\ell=0$ for $\ell\in\Gamma^*\setminus(\Lambda-\Lambda)$, at fixed $(z,Q)\in\mathcal U_a\times D_Q$ we define $V_{\mathcal A}$ by
\begin{equation*}
 V_{\mathcal A}(z,X,W,Q)
 =-\sum_{\ell\in\Lambda-\Lambda}e^{iX\cdot\ell}(P_\ell^{-1}f_\ell)(W),
 \qquad (X,W)\in\T^2_X\times\overline D\,.
\end{equation*}
For $H\in\mathcal H_{\mathcal A}$, let $v_H$ be the coefficient of $H$ in $V_{\mathcal A}$. By the coordinatewise definition of $P_\ell^{-1}$, the Fourier coefficient at $\ell\in\Gamma^*$ is $\widehat v_H(\ell)=-P_\ell^{-1}\bigl(\widehat\mu_H(\ell)/(4|\mathscr J|^2)\bigr)$. The real source $\mu_H$ has conjugate Fourier coefficients at $\ell$ and $-\ell$. Since $P_\ell=P_{-\ell}$ has real coefficients, the inverses preserve this symmetry, so $v_H$ is real.

Applying~\eqref{en:disk} to each mode of $v_H$ gives
\begin{equation*}
 L_Bv_H
  =-4|\mathscr J|^2\sum_{\ell\in\Lambda-\Lambda}e^{iX\cdot\ell}
       P_\ell\widehat v_H(\ell)
 =\sum_{\ell\in\Lambda-\Lambda}e^{iX\cdot\ell}\widehat\mu_H(\ell)
 =\mu_H\,.
\end{equation*}
Each $\widehat v_H(\ell)$ is smooth up to $\partial D$ and has zero Dirichlet trace because it lies in $H^1_0(D)$. Thus every summand of $v_H$ vanishes at $W\in\partial D$, and $v_H|_{\T^2_X\times\partial D}=0$. Conversely, Fourier transforming any smooth solution vanishing on \(\T^2_X\times\partial D\) gives the same Dirichlet equation for each $\ell\in\Gamma^*$. Invertibility of $P_\ell$ gives uniqueness of each $v_H$. By~\eqref{cur:density}, the source vanishes wherever $\chi_{\mathrm{fib}}(Q)=0$. At each such $Q$ uniqueness gives $v_H=0$, in particular when $|Q|\ge1/20$.

\medskip
We now check that the solutions depend smoothly on $(z,Q)$. For fixed $\ell\in\Gamma^*$, varying $(z,Q)$ changes only the smooth multiplication potential and the smooth source on the fixed disk $D$. To justify differentiation in one real parameter $\upsilon$ among the coordinates of $(z,Q)$, write $w_\upsilon=P_\ell(\upsilon)^{-1}f_\ell(\upsilon)$. For real $\varepsilon\ne0$ with both parameter values in a fixed compact neighborhood, write $\Delta_\varepsilon h=h(\upsilon+\varepsilon)-h(\upsilon)$ for the increment of a parameter-dependent quantity $h$. Subtraction of the two Dirichlet equations gives
\[
 P_\ell(\upsilon)\frac{w_{\upsilon+\varepsilon}-w_\upsilon}{\varepsilon}
 = \frac{f_\ell(\upsilon+\varepsilon)-f_\ell(\upsilon)}{\varepsilon}
 -\frac{P_\ell(\upsilon+\varepsilon)-P_\ell(\upsilon)}{\varepsilon}
       w_{\upsilon+\varepsilon}.
\]
Before division by $\varepsilon$, the two equations give
\begin{equation*}
 P_\ell(\upsilon)(w_{\upsilon+\varepsilon}-w_\upsilon)
 =f_\ell(\upsilon+\varepsilon)-f_\ell(\upsilon)
 -\bigl(P_\ell(\upsilon+\varepsilon)-P_\ell(\upsilon)\bigr)w_{\upsilon+\varepsilon}\,.
\end{equation*}
The difference $P_\ell(\upsilon+\varepsilon)-P_\ell(\upsilon)$ is multiplication by $\Delta_\varepsilon(|\ell\cdot d|^2/|\mathscr J|^2)$ by~\eqref{en:fourier}. Applying $P_\ell(\upsilon)^{-1}$ and using its uniform bound on $L^2(D;Z_{\mathcal A})$ gives
\begin{equation*}
 \begin{aligned}
 \|w_{\upsilon+\varepsilon}-w_\upsilon\|_{L^2(D;Z_{\mathcal A})}
 &\le \lambda_1(D)^{-1}
 \|f_\ell(\upsilon+\varepsilon)-f_\ell(\upsilon)\|_{L^2(D;Z_{\mathcal A})}\\
 &\quad+\lambda_1(D)^{-1}
 \|P_\ell(\upsilon+\varepsilon)-P_\ell(\upsilon)\|_{L^2(D;Z_{\mathcal A})\to L^2(D;Z_{\mathcal A})}
 \|w_{\upsilon+\varepsilon}\|_{L^2(D;Z_{\mathcal A})}\,.
 \end{aligned}
\end{equation*}
The inverse bound at $\upsilon+\varepsilon$ gives $\|w_{\upsilon+\varepsilon}\|_{L^2(D;Z_{\mathcal A})}\le\lambda_1(D)^{-1}\|f_\ell(\upsilon+\varepsilon)\|_{L^2(D;Z_{\mathcal A})}$, which remains bounded as $\varepsilon\to0$. By smoothness, $f_\ell(\upsilon+\varepsilon)\to f_\ell(\upsilon)$ in $L^2(D;Z_{\mathcal A})$ and the potentials converge uniformly on $\overline D$. Thus both terms in the estimate tend to zero, and $w_{\upsilon+\varepsilon}\to w_\upsilon$ in $L^2(D;Z_{\mathcal A})$.

As $\varepsilon\to0$, $(\Delta_\varepsilon f_\ell)/\varepsilon\to\partial_\upsilon f_\ell$ in $L^2(D;Z_{\mathcal A})$, and $(\Delta_\varepsilon P_\ell)/\varepsilon\to\partial_\upsilon P_\ell$ in the norm of multiplication operators on that space. Using $w_{\upsilon+\varepsilon}\to w_\upsilon$ in the divided equation and applying $P_\ell(\upsilon)^{-1}$ gives the parameter derivative in $L^2(D;Z_{\mathcal A})$ and its Dirichlet equation, with source $\partial_\upsilon f_\ell-(\partial_\upsilon P_\ell)w_\upsilon$.

To obtain joint smoothness up to the disk boundary, fix a compact parameter neighborhood $\mathscr K\Subset\mathcal U_a\times D_Q$. For each fixed $\ell$, the potential of $P_\ell$ and all its derivatives are bounded on $\overline D\times\mathscr K$. The Dirichlet elliptic estimates on the fixed disk therefore give, for every integer $m\ge0$,
\begin{equation*}
 \|w\|_{H^{m+2}(D;Z_{\mathcal A})}
 \le C_{m,\ell,\mathscr K}
 \bigl(\|P_\ell(z,Q)w\|_{H^m(D;Z_{\mathcal A})}
                   +\|w\|_{L^2(D;Z_{\mathcal A})}\bigr)\,,
\end{equation*}
uniformly for $(z,Q)\in\mathscr K$ and smooth $w:\overline D\to Z_{\mathcal A}$ vanishing on $\partial D$. Applied to $P_\ell(z,Q)w=f_\ell(z,Q)$, this gives a uniform $H^{m+2}$ bound for the solutions on $\mathscr K$, since $f_\ell$ is bounded there in $H^m$ and $P_\ell^{-1}$ is bounded on $L^2(D;Z_{\mathcal A})$. In the equation for $w_{\upsilon+\varepsilon}-w_\upsilon$, $\Delta_\varepsilon f_\ell\to0$ in $H^m$, and $\Delta_\varepsilon P_\ell$ tends to zero in the norm of multiplication operators on $H^m$. The elliptic estimate therefore gives $w_{\upsilon+\varepsilon}\to w_\upsilon$ in $H^{m+2}$.

It now follows that the right-hand side of the divided difference equation converges in $H^m$ to $\partial_\upsilon f_\ell-(\partial_\upsilon P_\ell)w_\upsilon$. Applying the elliptic estimate to $(\Delta_\varepsilon w)/\varepsilon-\partial_\upsilon w_\upsilon$ upgrades the convergence of $(\Delta_\varepsilon w)/\varepsilon$ to $H^{m+2}$. For a multi-index $\alpha$ in the real parameters $(z,Q)$, the equation for $w_\alpha=\partial^\alpha w$ has the form
\begin{equation*}
 P_\ell w_\alpha=f_{\ell,\alpha},\qquad
 f_{\ell,\alpha}=\partial^\alpha f_\ell-
   \sum_{0<\beta\le\alpha}\binom{\alpha}{\beta}
       (\partial^\beta P_\ell)\,\partial^{\alpha-\beta}w\,.
\end{equation*}
Here $\beta\le\alpha$ is understood componentwise and $w_\alpha$ has zero Dirichlet trace. If derivatives of $w$ through order $|\alpha|$ are continuous in $H^{m+2}$, then $f_{\ell,\alpha}$ is continuously differentiable in $H^m$, since $|\alpha-\beta|<|\alpha|$ in the sum. In the difference quotient equation, replace $w,f_\ell$ by $w_\alpha,f_{\ell,\alpha}$ and vary any one real parameter. The $L^2$ bound for $P_\ell^{-1}$ gives convergence of the quotient in $L^2$, and the Dirichlet estimate upgrades it to $H^{m+2}$. Sobolev embedding and the finite Fourier sum give joint smoothness of $v_H$ in $(z,X,W,Q)$ up to $\partial D$. The local regularity constants may depend on the frequency $\ell$ and the parameter neighborhood $\mathscr K$. The bound $\lambda_1(D)^{-1}$ for $P_\ell^{-1}$ on $L^2(D;Z_{\mathcal A})$ is uniform in $\ell$ and $(z,Q)$.
\end{proof}

\begin{samepage}
We use the same symbols $v_H,V_{\mathcal A},\mu_H,\mu_{\mathcal A}$ for the fields in disk coordinates and for their pullbacks to the flat coordinates. The solutions from Lemma~\ref{en:fiber-solve} define the paired energy $\mathcal E:\mathcal U_a\to[0,\infty)$ by
\begin{align}
 \mathcal E(z)
 &=\tfrac12\sum_{\mathcal A}\int_{\Omega_z}
     \mathrm{tr}_{\mathcal A}((-V_{\mathcal A})\mu_{\mathcal A})\,d\zeta
  =\tfrac12\sum_{\mathcal A,H}\int_{\Omega_z}(-v_H)\mu_H\,d\zeta\notag\\
 &=2a^2\sum_{\mathcal A,\ell}\int_{D_Q}
  \langle f_\ell,P_\ell^{-1}f_\ell\rangle_{L^2(D;Z_{\mathcal A})}\,d^2Q\ge0 .
 \label{en:positive}
\end{align}
\end{samepage}
Here and below $\sum_{\mathcal A,H}$ sums over $\mathcal A=\M,\M^{\mathrm{op}}$ and $H\in\mathcal H_{\mathcal A}$, and the sum over $\ell\in\Gamma^*$ has finite support. Expanding the trace pairing in the orthonormal basis $\mathcal H_{\mathcal A}$ gives the sum of scalar integrals. Plancherel's identity in $X$, with real part taken in the Hilbert inner product, pairs $-\widehat V_{\mathcal A}(\ell)=P_\ell^{-1}f_\ell$ with $\widehat\mu_{\mathcal A}(\ell)=4|\mathscr J|^2f_\ell$. Combining the factor $4|\mathscr J|^2$ with $\tfrac12 a^2|\mathscr J|^{-2}$ from~\eqref{en:measure} gives $2a^2$ in~\eqref{en:positive}. The positive Dirichlet inverse gives $\mathcal E\ge0$.

In the fixed chart, the integrand in the first line of~\eqref{en:positive} is smooth in $z$ up to $|W|=R$ and vanishes for $|Q|\ge1/20$. We may therefore differentiate under the integral on compact subsets of $\mathcal U_a$ to obtain smoothness of $\mathcal E$.

We can now state the three estimates that will close the proof of Proposition~\ref{en:cap}. Define the circle mean $H:\R\times(a,4)\to[0,\infty)$ by
\begin{equation}\label{en:circlemean}
 H(T,s)=\frac{1}{2\pi}\int_0^{2\pi}\mathcal E(T,se^{k\theta})\,d\theta.
\end{equation}
The aperture $\gamma$ will be chosen in Section~\ref{en:comparison}, independently of $N$ and the coefficient algebra; we then fix $K=10$. Recall $a=K/N$ and suppose for the moment that $2a<1/4$.

\begin{samepage}
The estimates to be proved below are
\begin{equation}\label{en:closure}
 \begin{aligned}
 a^2\int_1^2 q_u(T)\,dT
 &\lesssim \int_1^2 H(T,2a)\,dT
 &&\text{(Lemma~\ref{en:recovery})}\\
 &\lesssim \log^2(3/(2a))\,H(0,4a)
 &&\text{(Lemma~\ref{en:annulus})}\\
 &\lesssim a^2\log^2(3/(2a))\,q_u(0)
 &&\text{(Lemma~\ref{en:upper}).}
 \end{aligned}
\end{equation}
\end{samepage}
Section~\ref{en:comparison} combines these estimates to prove~\eqref{en:capestimate} and also treats the bounded range $2a\ge1/4$.

The middle inequality in~\eqref{en:closure} depends on
\begin{equation*}
 \Delta_z\mathcal E(T,\tau)\le V_N(|\tau|)\mathcal E(T,\tau),
 \qquad 2a\le|\tau|\le3,
\end{equation*}
where $\Delta_z=\partial_T^2+\partial_{\operatorname{Re}_k\tau}^2+\partial_{\operatorname{Im}_k\tau}^2$ and the nonnegative radial function $V_N:[2a,3]\to[0,\infty)$ will be specified in~\eqref{en:sub}. To obtain this differential inequality, Section~\ref{en:green-section} derives $\Delta_z\mathcal E+\mathcal Q=\mathcal B$, with $\mathcal Q\ge0$ and boundary flux $\mathcal B$ defined in Lemma~\ref{en:green}. Lemmas~\ref{en:upper} and~\ref{en:recovery} then show that $|\mathcal B(T,\tau)|\le V_N(|\tau|)\mathcal E(T,\tau)$: the Gaussian decay on $|W|=R$ dominates the loss in recovering $q_u(T)$ near $W=0$.

\subsection{Green’s formula for the paired energy}\label{en:green-section}

The goal of this subsection is the identity $\Delta_z\mathcal E+\mathcal Q=\mathcal B$ in Lemma~\ref{en:green}. We first specify the boundary measure and the normal derivative that enter the formula for $\mathcal B$.

By Lemma~\ref{en:fiber-solve}, $v_H$ vanishes for $|Q|\ge1/20$. This support property and periodicity in $X$ leave $|W|=R$ as the only boundary contributing to Green's formula. Set $r=|W|^2/R^2-1$. For $W\in\partial D$, let $n(W)=W/R$ be the outward unit normal, let $dS_W$ denote arc length, and write $\partial_n$ for differentiation in the direction $n(W)$ at fixed $(z,X,Q)$. Let $\delta_0$ be the Dirac mass at zero in the real variable $r$. Since $\partial_n r=2/R$ on $\partial D$, the volume formula~\eqref{en:measure} gives the boundary measure
\begin{equation*}
 \delta_0(r)\,d\zeta
 =\frac{R}{2}a^2|\mathscr J|^{-2}\,dX\,dS_W\,d^2Q\,.
\end{equation*}
For any real function $v$ smooth up to $|W|=R$ and vanishing there, $v/r$ extends smoothly to the boundary. Indeed, the fundamental theorem of calculus in the transverse coordinate $r$ writes $v$ as $r$ times a smooth function. Differentiating in the normal direction gives the boundary value $(v/r)|_{\partial D}=(R/2)\partial_n v$. The vector $G(r)$ from~\eqref{cur:hess} encodes all nine flat derivatives of $r$, including the derivatives in $\tau$ that record the motion of the boundary as $z$ varies.

To compute $\Delta_z\mathcal E$, we use the mixed cofactor $M$ of the Hessians $\mathsf H(v_H)$ and $B=\mathsf H(|\Phi|^2)$. We first prove that the columns of its real tensor have zero divergence in flat coordinates, which lets us integrate by parts. The fiber equation gives $M_{TT}=\mu_H$. The smooth extension of $v_H/r$ and the identity $\operatorname{ran}\mathsf H(r)=G(r)\mathbb H$ on $r=0$ show that the compression of $\mathsf H(v_H)$ to $(G(r)\mathbb H)^\perp$ vanishes there. The cofactor formula then gives $G(r)^*MG(r)=0$ on $r=0$.

\begin{lemma}\label{en:cofactor-divergence}
Let $C:\operatorname{Herm}_3(\mathbb H)\to\operatorname{Herm}_3(\mathbb H)$ be the gradient of the Moore determinant for the pairing in~\eqref{en:lb}. For $U,V\in\operatorname{Herm}_3(\mathbb H)$ let $M(U,V)\in\operatorname{Herm}_3(\mathbb H)$ be the coefficient of the real parameter $t$ in $C(U+tV)$. If $u,v:O\to\R$ are smooth on an open set $O\subset\R^9$ in the flat coordinates of~\eqref{cur:hess}, every real column of $M(\mathsf H(u),\mathsf H(v))^{\mathbb R}$ has zero divergence.
\end{lemma}
\begin{proof}
The column divergence is a bilinear differential expression with constant coefficients in $u$ and $v$. Its values on pairs of real exponentials determine all of its coefficients. The Hessians of the exponentials have rank at most one, which reduces the calculation to an algebraic identity for the mixed cofactor.

For $L\in\operatorname{Herm}_3(\mathbb H)$, the product of its three real eigenvalues is
\[
 \det_{\mathbb H} L=\tfrac16\big(p_1^3-3p_1p_2+2p_3\big),
 \qquad p_m=\operatorname{Re}\mathrm{tr}_{\mathbb H}(L^m),\quad m=1,2,3.
\]
Indeed unitary diagonalization identifies \(p_m\) with the sum of the \(m\)-th powers of the eigenvalues, and expansion of the cubic leaves six times their product. For a variation $D_0\in\operatorname{Herm}_3(\mathbb H)$, differentiating \((L+tD_0)^m\) at zero inserts \(D_0\) in each of the \(m\) possible positions. Cyclicity of the real trace gives \(dp_m(L)[D_0]=m\operatorname{Re}\mathrm{tr}_{\mathbb H}(L^{m-1}D_0)\). Thus the cubic is differentiable also at repeated or zero eigenvalues, and its gradient is
\[
 C(L)=L^2-(\mathrm{tr}_{\mathbb H}L)L+
 \tfrac12\big((\mathrm{tr}_{\mathbb H}L)^2-\mathrm{tr}_{\mathbb H}(L^2)\big)I_3 .
\]
Expanding $C(U+tV)$ and using $\operatorname{Re}\mathrm{tr}_{\mathbb H}(UV)=\operatorname{Re}\mathrm{tr}_{\mathbb H}(VU)$ gives
\begin{equation}
 M(U,V)= UV+VU-(\mathrm{tr}_{\mathbb H}U)V-(\mathrm{tr}_{\mathbb H}V)U
 +\big((\mathrm{tr}_{\mathbb H}U)(\mathrm{tr}_{\mathbb H}V)
                   -\operatorname{Re}\mathrm{tr}_{\mathbb H}(UV)\big)I_3 .
 \label{en:cofactorpoly}
\end{equation}

First evaluate~\eqref{en:cofactorpoly} on
$U=gg^*$ and $V=ll^*$, where $g,l\in\mathbb H^3$ are columns. Set $a_g=g^*g$, $a_l=l^*l$ (both real), and $c=g^*l$, so $l^*g=\bar c$. The products in their given orders are $UV=gc\,l^*$ and $VU=l\bar c\,g^*$, and $\operatorname{Re}\mathrm{tr}_{\mathbb H}(UV)=\operatorname{Re}(l^*g\,c)=|c|^2$. Consequently, for $M_0=M(gg^*,ll^*)$,
\[
 M_0=gc\,l^*+l\bar c\,g^*-a_gll^*-a_lgg^*
                         +(a_ga_l-|c|^2)I_3 .
\]
Multiplying by $g$ and by $l$ on the right gives
\begin{align*}
 M_0g&=gc\bar c+l\bar c\,a_g-a_g l\bar c-a_lg a_g
                           +(a_ga_l-|c|^2)g=0,\\
 M_0l&=gc\,a_l+l\bar c c-a_g l a_l-a_lgc
                           +(a_ga_l-|c|^2)l=0 .
\end{align*}
In the first line the two multiples of $l\bar c$ cancel because $a_g$ is real. The remaining terms are $g|c|^2-ga_ga_l+g(a_ga_l-|c|^2)$. In the second line the two multiples of $gc$ cancel because $a_l$ is real, and the rest is $l|c|^2-la_ga_l+l(a_ga_l-|c|^2)$. Both cancellations use the real identities $c\bar c=\bar c c=|c|^2$.

To pass to real column divergences, let $\mathbf y=(y_{d\kappa})_{d,\kappa=0}^2$ denote the nine flat real coordinates, and let $\eta,\xi\in\R^9$ be arbitrary. For the test functions $u=e^{\eta\cdot\mathbf y}$, $v=e^{\xi\cdot\mathbf y}$ set
$g_d=\sum_\kappa e_\kappa\eta_{d\kappa}$ and $l_d=\sum_\kappa e_\kappa\xi_{d\kappa}$.
Equation~\eqref{cur:hess} gives $\mathsf H(u)=u\,gg^*$ and $\mathsf H(v)=v\,ll^*$. The real component of the cofactor in row $(d,\kappa)$ and column $(b,\lambda)$ is therefore $uv\,\operatorname{Re}(\bar e_\kappa(M_0)_{db}e_\lambda)$. For every fixed $(b,\lambda)$ its column divergence is
\begin{equation*}
 uv\sum_{d,\kappa}(\eta_{d\kappa}+\xi_{d\kappa})
                       \operatorname{Re}(\bar e_\kappa(M_0)_{db}e_\lambda)
 =uv\,\operatorname{Re}\!\left(
           \big[(g+l)^*M_0\big]_b e_\lambda\right)=0 .
\end{equation*}
Indeed, $M_0$ is Hermitian and annihilates $g$ and $l$ on the right. For arbitrary $u,v$ each such divergence is a finite constant-coefficient expression
$\sum_{\mu,\nu}c_{\mu\nu}(\partial^\mu u)(\partial^\nu v)$ with multi-indices $\mu,\nu\in\{0,1,\ldots\}^9$. Its value on these exponentials divided by $uv$ is the polynomial $\sum_{\mu,\nu}c_{\mu\nu}\eta^\mu\xi^\nu$. Since this polynomial vanishes for every $(\eta,\xi)\in\R^9\times\R^9$, every coefficient is zero after terms with the same derivative indices have been collected. This proves the identity locally for all smooth real-valued functions $u,v$ on $O$, including where the original spatial variables are periodic.
\end{proof}

Let $t_0=(1,0,0)^{\mathsf t}$ be the physical time vector, let $E=\ker B$, and write $\beta=\mathrm{tr}_{\mathbb H}B>0$. Since $\mathrm{tr}_{\mathbb H}B_f=4|d|^2>0$ and $B$ has rank one, orthogonal projection onto $E$ is $P_E=I_3-B/\beta$. The vectors $P_Et_0/|P_Et_0|$ and $(0,d/|d|)$ form an orthonormal frame of $E$. In fact $cd=0$, the second vector is spatial, and $|P_Et_0|^2=1-B_{TT}/\beta=4|d|^2/\beta$. For each $\mathcal A\in\{\M,\M^{\mathrm{op}}\}$ and $H\in\mathcal H_{\mathcal A}$, set $D_H=\mathsf H(v_H)$ and let $b_H=D_H|_E$ and $\alpha_H=A_H|_E$ denote the two-dimensional compressions in this frame, with indices $0,1$. Here a compression to $E$ means $P_E D_H P_E$ or $P_E A_H P_E$ viewed as an operator on $E$.

\begin{lemma}\label{en:green}
For every $z=(T,\tau)\in\mathcal U_a$, the energy $\mathcal E$ satisfies
\begin{equation}\label{en:balance}
 \Delta_z\mathcal E(z)+\mathcal Q(z)=\mathcal B(z),
\end{equation}
where $\mathcal Q:\mathcal U_a\to[0,\infty)$ and $\mathcal B:\mathcal U_a\to\R$ are defined by
\begin{align*}
 \mathcal Q(z)&=\sum_{\mathcal A,H}\int_{\Omega_z}\beta
 \left|(b_H)_{01}+\frac{(\alpha_H)_{01}}{\beta}\right|^2\,d\zeta\ge0,\\
 \mathcal B(z)&=\sum_{\mathcal A,H}\int\delta_0(r)\,
        \frac{v_H}{r}\,G(r)^*A_HG(r)\,d\zeta.
\end{align*}
On $|W|=R$, use $v_H/r=(R/2)\partial_n v_H$ and the boundary measure $\delta_0(r)\,d\zeta$ fixed at the start of Section~\ref{en:green-section}.
\end{lemma}

\begin{proof}
Write $D_H=\mathsf H(v_H)$ and $M=M(D_H,B)$. We shall apply integration by parts on the moving fiber to
\[
 S=A_H-\tfrac12 M.
\]
This choice gives $S_{TT}=\mu_H/2$, the density factor in~\eqref{en:positive}. To prove~\eqref{en:balance} we evaluate $G(r)^*SG(r)$ on $r=0$ and write the interior pairing $S:D_H$ as the sum of a square and the determinant term that cancels by Proposition~\ref{cur:cancel}.

\medskip
We first derive the integration-by-parts identity for any $\operatorname{Herm}_3(\mathbb H)$-valued tensor $S$ on the moving fibers that is smooth up to the wall and has divergence-free real columns. For a test function $\varphi$ of $z$ alone, $S^{\mathbb R}:D^2\varphi=S_{TT}\Delta_z\varphi$ since $\operatorname{Re}(\bar e_\kappa S_{TT}e_\lambda)=S_{TT}\delta_{\kappa\lambda}$. Suppose $v$ is real and smooth up to the wall, vanishes there and has support away from $\partial D_Q$. Let $q=v/r$, with its smooth extension given by $q=\int_0^1(\partial_r v)(tr)\,dt$ in local boundary coordinates, keeping the tangential coordinates fixed. Extend $v S^{\mathbb R}$ by zero across $r=0$. For the real indices $I,J$ of the total space the distributional double divergence is
\begin{equation*}
 \sum_{I,J}\partial_I\partial_J
    (\mathbf1_{r<0}\,v S^{\mathbb R}_{IJ})
  =\sum_{I,J}\partial_I(\mathbf1_{r<0}\, S^{\mathbb R}_{IJ}\partial_J v)
 =\mathbf1_{r<0}\,S:\mathsf H(v)
       -\delta_0(r)\,q\,G(r)^*SG(r).
\end{equation*}
Indeed the boundary contribution in the first derivative contains $v|_{r=0}$, and the derivatives on $S^{\mathbb R}$ vanish after summing, by symmetry and column conservation. In the second derivative $\partial_I\mathbf1_{r<0}=-\delta_0(r)\partial_I r$ and $\partial_Jv=q\partial_Jr$ on $r=0$. With $I,J$ ranging over all nine real coordinates, including the auxiliary time coordinates, the definitions of $G$ in~\eqref{cur:hess} and of $S^{\mathbb R}$ give
\[
 \sum_{I,J} S^{\mathbb R}_{IJ}\,\partial_I r\,\partial_J r
 =G(r)^*SG(r).
\]
Let $F_{IJ}=\mathbf1_{r<0}\,v S^{\mathbb R}_{IJ}$. For $\varphi\in C_c^\infty(\mathcal U_a)$ regarded as a function of $z$ alone in the full coordinates, the distributional pairing is
\begin{equation*}
 \left\langle\sum_{I,J}\partial_I\partial_J F_{IJ},\varphi\right\rangle
 =\int_{\mathcal U_a}\Delta_z\varphi(z)
                 \left(\int_{\Omega_z}v S_{TT}\,d\zeta\right)dz\,.
\end{equation*}
Here $dz=dT\,d^2\tau$, and we used $S^{\mathbb R}:D^2\varphi=S_{TT}\Delta_z\varphi$ in the original flat coordinates. Periodicity in $X$ and the support away from $\partial D_Q$ ensure that the double divergence has no further boundary terms. Pairing its interior and boundary expressions with $\varphi$ and changing signs gives
\begin{equation}
 \Delta_z\int_{\Omega_z}(-v)S_{TT}\,d\zeta
 =-\int_{\Omega_z}S:\mathsf H(v)\,d\zeta
   +\int\delta_0(r)\,q\,G(r)^*SG(r)\,d\zeta .
 \label{en:move}
\end{equation}
In the fixed chart $(W,Q)$ all three integrals are smooth in $z$ (the wall integral has $|W|=R$). Thus the distributional identity also holds pointwise.

\medskip
To apply~\eqref{en:move}, take $M=M(D_H,B)$ from~\eqref{en:cofactorpoly}. Its real columns have zero divergence by Lemma~\ref{en:cofactor-divergence} applied with $u=v_H$ and $v=|\Phi|^2$. Since $\operatorname{div}(A_H^{\mathbb R})=0$ by the construction of $A_H$ in Section~\ref{cur:sec}, the real columns of $S=A_H-M/2$ also have zero divergence. We shall use
\begin{equation}
 M:D_H=2\beta\det_{\mathbb H} b_H,\qquad
 M_{TT}=C_2(B_f):\mathsf H_f(v_H)=\mu_H .
 \label{en:minors}
\end{equation}
To compute $M:D_H$, use an orthonormal frame adapted to $E\oplus E^\perp$ so that $B=\operatorname{diag}(0,0,\beta)$ and write
$D_H=\left(\begin{smallmatrix}b_H&u\\u^*&d_0\end{smallmatrix}\right)$.
Then $D_HB+BD_H$ has upper-left block zero, off-diagonal blocks $\beta u,\beta u^*$ and lower-right entry $2\beta d_0$. In~\eqref{en:cofactorpoly} the term $-\beta D_H$ cancels those off-diagonal blocks. The scalar identity term has coefficient $\beta\,\mathrm{tr}_{\mathbb H} b_H$, since $\mathrm{tr}_{\mathbb H}(D_H B)=\beta d_0$. Including $-(\mathrm{tr}_{\mathbb H}D_H)B$ leaves $M=\operatorname{diag}(\beta C_2(b_H),0)$. For $b_H=\left(\begin{smallmatrix}b_{00}&b_{01}\\\bar b_{01}&b_{11}\end{smallmatrix}\right)$ the identity $C_2(b_H):b_H=2(b_{00}b_{11}-|b_{01}|^2)$ gives the formula for $M:D_H$ in~\eqref{en:minors}.

To compute $M_{TT}$ and then $G(r)^*MG(r)$ on the wall, we use the following consequence of the cofactor formula: for any $L\in\operatorname{Herm}_3(\mathbb H)$ and any unit column $w\in\mathbb H^3$,
\[
 w^*C(L)w=\det_{\mathbb H}(L|_{w^\perp}).
\]
Indeed in a frame starting with $w$ write
$L=\left(\begin{smallmatrix}d_0&u^*\\u&B_2\end{smallmatrix}\right)$. The first diagonal entry of $L^2$ is $d_0^2+|u|^2$, and $\mathrm{tr}_{\mathbb H}L^2=d_0^2+2|u|^2+\mathrm{tr}_{\mathbb H}B_2^2$. Substitution in $C(L)$ gives
$\tfrac12((\mathrm{tr}_{\mathbb H}B_2)^2-\mathrm{tr}_{\mathbb H}B_2^2)=\det_{\mathbb H}B_2$.
Apply the cofactor identity to $L=D_H+tB$ and take the coefficient of $t$. Expanding the two-dimensional determinant gives, for every unit $w\in\mathbb H^3$,
\begin{equation*}
 w^*M(D_H,B)w=C_2(B|_{w^\perp}):(D_H|_{w^\perp})\,.
\end{equation*}
With $w=t_0$ this gives $M_{TT}=C_2(B_f):\mathsf H_f(v_H)=\mu_H$, where the last equality is the Dirichlet equation. This proves~\eqref{en:minors}.

On $r=0$ put $\nu=G(r)/|G(r)|$, the full unit normal in its quaternionic encoding. Holomorphy of $W$ in the three auxiliary variables gives
\[
 \mathsf H(r)=\frac{4}{R^2}(\partial W)^*\partial W,\qquad
 G(r)=\frac{2}{R^2}(\partial W)^*W i .
\]
Since $|W|=R>0$ and $\partial_P W=1/Y\ne0$, the range of $\mathsf H(r)$ is the right quaternionic line $G(r)\mathbb H$. Let $P_{\nu^\perp}$ be the orthogonal projection onto the right quaternionic complement of this line. The product rule for $v_H=rq$ gives at the wall
\[
 D_H=q\,\mathsf H(r)+G(q)G(r)^*+G(r)G(q)^*.
\]
We have $P_{\nu^\perp}\mathsf H(r)=0$, $P_{\nu^\perp}G(r)=0$, and $G(r)^*P_{\nu^\perp}=0$. Thus $P_{\nu^\perp}D_HP_{\nu^\perp}=0$. Taking $w=\nu$ in the mixed cofactor formula gives $\nu^*M\nu=C_2(B|_{\nu^\perp}):(D_H|_{\nu^\perp})=0$. Since $G(r)=|G(r)|\nu$ on $r=0$, this gives $G(r)^*MG(r)=0$ there. Thus the boundary integral in~\eqref{en:move} uses $G(r)^*SG(r)=G(r)^*A_HG(r)$ for $S=A_H-M/2$.

\medskip
We have evaluated the boundary pairing $G(r)^*SG(r)$. We now compute the interior pairing $S:D_H$. In~\eqref{en:move} take $v=v_H$ and $S=A_H-M/2$. By~\eqref{en:minors}, $S_{TT}=\mu_H/2$ and $S:D_H=A_H:D_H-\beta\det_{\mathbb H}b_H$. Write $\alpha=\alpha_H$, $b=b_H$. Since $A_H=P_EA_HP_E$,
\begin{equation*}
 \mu_H=t_0^*A_Ht_0=|P_Et_0|^2\alpha_{00}
       =4|d|^2\frac{\alpha_{00}}{\beta}\,.
\end{equation*}
On the other hand $C_2(B_f)=4dd^*$, and $(0,d/|d|)$ is the second vector of the chosen frame of $E$. The fiber equation gives
\begin{equation*}
 \mu_H=C_2(B_f):\mathsf H_f(v_H)=4|d|^2 b_{11}\,.
\end{equation*}
Cancelling the common positive factor $4|d|^2$ yields $b_{11}=\alpha_{00}/\beta$. In the same frame $A_H:D_H=\operatorname{Re}\mathrm{tr}_{\mathbb H}(\alpha b)$, since $A_H$ vanishes on $E^\perp$. Expanding in this two-plane gives
\begin{align*}
 A_H:D_H-\beta\det_{\mathbb H} b
 &=\alpha_{00}b_{00}+\alpha_{11}b_{11}
    +2\operatorname{Re}(\alpha_{01}\bar b_{01})
    -\beta(b_{00}b_{11}-|b_{01}|^2)\\
 &=\frac{\alpha_{00}\alpha_{11}}{\beta}
       +2\operatorname{Re}(\alpha_{01}\bar b_{01})+\beta |b_{01}|^2\\
 &=\beta|b_{01}+\alpha_{01}/\beta|^2+
              \frac{\det_{\mathbb H}\alpha}{\beta}.
\end{align*}
Summing~\eqref{cur:ordered} over $\mathcal A\in\{\M,\M^{\mathrm{op}}\}$ gives $\sum_{\mathcal A,H}\det_{\mathbb H}(\alpha_H)/\beta=0$ pointwise on $\Omega_z$. Summing~\eqref{en:move} and using~\eqref{en:positive} on the left gives $\Delta_z\mathcal E=-\mathcal Q+\mathcal B$, as claimed in~\eqref{en:balance}.
\end{proof}

\subsection{Upper bounds for the energy and boundary flux}\label{en:upper-section}

At $|\tau|=4a$, an upper bound for $\mathcal E(0,\tau)$ in terms of $q_u(0)$ will give the last inequality in~\eqref{en:closure}. We estimate $\mathcal E(T,\tau)$ throughout $2a\le|\tau|\le3$ and also bound the flux $\mathcal B(T,\tau)$ in~\eqref{en:balance}. The flux bound has an exponential factor $e^{-c_wN|\tau|}$, with $c_w>0$ independent of the aperture $\gamma$. This independence will let us choose $\gamma$ so that the decay dominates the loss in Lemma~\ref{en:recovery}.

Fix $0<c_0<C_0$, $R=1/10$ and the fiber cutoff $\chi_{\mathrm{fib}}$ of Section~\ref{cur:sec}. For scaled frequencies $\eta=\xi/N=(\eta_x,\eta_{b_0})\in\R^2$ with $\eta_{b_0}>0$, use $\kappa_\eta,v_\eta$ from~\eqref{pre:slope}. Their homogeneity gives $\kappa_{\xi/N}=\kappa_\xi/N$ and $v_{\xi/N}=v_\xi$. Choose a radial profile $\phi\in C_c^\infty((c_0/2,2C_0))$ equal to one on $[c_0,C_0]$ and a real smooth angular cutoff equal to one on $|v_\eta|\le R/8$ in this cap, supported where $\eta_{b_0}>0$ and $|v_\eta|<R/3$. Write $\vartheta_{\rm wide}(\eta)$ for the product of $\phi(|\eta|)$ and this angular cutoff. On its support, $\kappa_\eta$ has a positive lower bound, and $\kappa_\eta,v_\eta$ have bounded derivatives through order four. For each orthonormal frame $(e_x,e_{b_0})$, we use the same cutoff profiles in its frequency coordinates.

For now take real parameters $K\ge1$, $0<\gamma<R/8$ and $N\ge2$, with $a=K/N$. The wave $u$ has finite spatial Fourier support in $c_0\le |\xi|/N\le C_0$, $|v_\xi|\le\gamma$ in the chosen cap. In the inequalities below $z=(T,\tau)$ is fixed and $2a\le|\tau|\le3$. In particular $N|\tau|\ge2K\ge2$. We obtain the exponential rate $c_w$ at the wall and the coefficient $C_*$ in the recovery bound before selecting $\gamma$ and $K$.

\begin{lemma}\label{en:upper}
For every $T\in\R$, every $\tau\in\mathbb C_k$ with $2a\le |\tau|\le3$, and every finite Fourier sum $u$ as above,
\begin{align}
 \mathcal E(T,\tau)&\le C a^2(N|\tau|)^{18}q_u(T)\,,\label{en:energyupper}\\
 |\mathcal B(T,\tau)|&\le C(a/|\tau|)^2(N|\tau|)^{18}e^{-c_wN|\tau|}q_u(T)\,.
 \label{en:upperbounds}
\end{align}
The constants $C,c_w>0$ depend only on $c_0,C_0,R,\chi_{\mathrm{fib}}$ and $\vartheta_{\rm wide}$.
\end{lemma}

\begin{proof}
We prove~\eqref{en:energyupper} by bounding the Dirichlet source in~\eqref{en:positive} and using the uniform $L^2(D)$ bound for $P_\ell^{-1}$. For~\eqref{en:upperbounds}, we combine an estimate for $\partial_n v_H$ with the exponential bound for $G(r)^*A_HG(r)$ on $r=0$. We first express $G(r)^*A_HG(r)$ in terms of spatial multipliers, as in~\eqref{cur:density} for the time density $\mu_H$.

In the notation of~\eqref{cur:density}, $S_0,S_1$ have symbols $m_\xi$ times the two real components in $\mathbb C_k$ of $V_T^{v_\xi}=(1+Wv_\xi/2)/\sqrt2$. To represent $G(r)^*A_HG(r)$, we first compute $G(r)^*V^v$. In null coordinates $\partial W=(-W,1,0)/Y$, so that $\partial W\,V^v=(v-W)/(2Y)$. Taking the quaternionic adjoint of the formula for $G(r)$ gives
\begin{equation}
 G(r)^* V^v=-i\,\frac{\bar W(v-W)}{R^2Y}
       =\frac{W(v-\bar W)}{R^2 |Y|}\,\rho_e(i)\,.
 \label{en:fullnormal}
\end{equation}
Indeed $i\zeta=\bar\zeta i$ for $\zeta\in\mathbb C_k$ and $1/\bar Y=e/|Y|$. Define scalar spatial multipliers $T_0,T_1$ on the trigonometric polynomials with frequencies in $\Lambda$, mapping that space into itself, by the respective symbols
\begin{equation*}
 \frac{m_\xi}{R^2|Y|}(v_\xi\operatorname{Re}_k W-|W|^2)
 \quad\text{and}\quad
 \frac{m_\xi}{R^2|Y|}v_\xi\operatorname{Im}_k W\,.
\end{equation*}
They act on scalar or $L^4(\mathcal A)$ coefficients. For the scalar wave $u_c=\sum_{\xi\in\Lambda}c_\xi e^{i\theta_\xi}$ with $c\in\mathbb C_i^\Lambda$, formula~\eqref{en:fullnormal} gives
\begin{equation*}
 G(r)^*\chi(c)=\big[\rho_e(T_0u_c)+k\rho_e(T_1u_c)\big]\rho_e(i)\,.
\end{equation*}
Here $\rho_e(i)$ commutes with the embedded amplitudes in $\rho_e(\mathbb C_i)$. Since $|\rho_e(i)|=1$ and $\rho_e(\mathbb C_i)$ and $k\rho_e(\mathbb C_i)$ are orthogonal, we have $G(r)^*A(c)G(r)=\chi_{\mathrm{fib}}(Q)\sum_{j=0}^1|T_j u_c|^2$. For $H\in\mathcal H_{\mathcal A}$, insert this formula into $A_H=\sum_l\lambda_l A(c_l)$. In the expansion of the squared moduli, the coefficient at each pair $\xi,\eta\in\Lambda$ is then $\sum_l\lambda_l(c_l)_\xi\overline{(c_l)_\eta}=\mathrm{tr}_{\mathcal A}(Ha_\eta^*a_\xi)$. Thus $G(r)^*A_HG(r)=\mathrm{tr}_{\mathcal A}(HK_{\mathcal A})$, where
\begin{equation}
 K_{\mathcal A}=\chi_{\mathrm{fib}}(Q)\sum_{j=0}^1(T_j u_{\mathcal A})^*T_j u_{\mathcal A}\,.
 \label{en:wallK}
\end{equation}
The term $-|W|^2$ in the symbol of $T_0$ comes from expanding $W(v_\xi-\bar W)$ in~\eqref{en:fullnormal}, where $G(r)$ includes the two real $\tau$ derivatives as well as derivatives in the fiber variables.

\medskip
To bound $\mu_{\mathcal A}$ and $K_{\mathcal A}$, we estimate the spatial multipliers $S_j$ and $T_j$ on the cap. For $z=(T,\tau)$ with $2a\le|\tau|\le3$ and fixed $(W,Q)\in\overline D\times D_Q$, the multipliers $S_j,T_j$, $j=0,1$, act in $X$ on trigonometric polynomials with coefficients in $L^4(\mathcal A)$ and frequencies in the chosen cap. Each preserves this space. For $\eta\in\R^2$ with $\eta_{b_0}>0$, define $F_W(\eta)\in[0,\infty)$ by
\begin{equation}\label{en:gaussian}
 F_W(\eta)=\frac{\kappa_\eta}{2}\frac{|Y|}{|\tau|}|W-v_\eta|^2,
 \qquad
 m_\xi=\left(\frac{\kappa_\xi}{|Y|}\right)^{1/2}
          e^{-N|\tau|F_W(\xi/N)}\,.
\end{equation}
The second equality rewrites~\eqref{cur:scalar}, since $\kappa_\xi=N\kappa_{\xi/N}$. By~\eqref{en:measure}, $|Y|/|\tau|$ and its inverse are uniformly bounded. Its first derivatives in the two real $W$ coordinates are bounded at fixed $(Q,z)$, since
\begin{equation*}
 \frac{|Y|}{|\tau|}
 =\frac{|\sqrt2\tau+aQ|}{|\tau|\,|1+W^2/2|}\,.
\end{equation*}
All $\eta$ derivatives of $F_W$ through order four on the support of $\vartheta_{\rm wide}$ are uniformly bounded, as are their first $W$ derivatives. Each differentiation of $e^{-N|\tau|F_W}$ adds at most one factor $N|\tau|$. The prefactor in~\eqref{en:gaussian}, written as $(N/|\tau|)^{1/2}(\kappa_\eta/(|Y|/|\tau|))^{1/2}$ at $\xi=N\eta$, has scaled derivatives through order four bounded by $C(N/|\tau|)^{1/2}$. Since $N|\tau|\ge2$, the powers $(N|\tau|)^k$ with $k\in\{0,1,2,3,4\}$ are all bounded by $(N|\tau|)^4$. Multiplying each symbol of $S_j$ by $\vartheta_{\rm wide}(\xi/N)$ and applying Leibniz' rule therefore bounds its scaled derivatives through order four by $C(N/|\tau|)^{1/2}(N|\tau|)^4$. For $N^{-1}\partial_{X_i}S_j$, the extra factor is $i\eta_i$. For a first real $W$ derivative of $S_j$, the bound becomes $C(N/|\tau|)^{1/2}(N|\tau|)^5$.

On the wall $|W|=R$, the angular support gives $|W-v_\eta|\ge2R/3$. Choose $c_w>0$, using only $R$ and the fixed lower bounds for $\kappa_\eta$ and $|Y|/|\tau|$, so that $F_W\ge c_w/2$ there. In each differentiated symbol of $T_j$, the Gaussian still occurs as $e^{-N|\tau|F_W}$. The additional denominator $|Y|$ is comparable to $|\tau|$, giving the symbol bound $C|\tau|^{-1}(N/|\tau|)^{1/2}(N|\tau|)^4e^{-c_wN|\tau|/2}$. Squaring in the estimate for the quadratic density $K_{\mathcal A}$ gives the exponential rate $c_w$.

By Lemma~\ref{pre:scalar-kernels}\ref{pre:kernel-bounds}--\ref{pre:convolution}, a bound $A$ on the scaled derivatives of each localized symbol through order four gives an $L^4(\T^2_X;L^4(\mathcal A))$ multiplier norm at most $CA$. The $L^p$ spaces on $\T^2_X$ use probability Haar measure $dX$. Using the pointwise equality $\|u_{\M}\|_{L^4(\M)}=\|u_{\M^{\mathrm{op}}}\|_{L^4(\M^{\mathrm{op}})}$, for $j=0,1$ we obtain
\begin{equation}\label{en:Sbounds}
\begin{aligned}
 &\|S_j u_{\mathcal A}\|_{L^4(\T^2_X;L^4(\mathcal A))}
   +N^{-1}\|\nabla_X(S_j u_{\mathcal A})\|_{L^4(\T^2_X;L^4(\mathcal A))}\\
 &\quad +(N|\tau|)^{-1}\|\nabla_W(S_j u_{\mathcal A})\|_{L^4(\T^2_X;L^4(\mathcal A))}\\
 &\qquad\le C(N/|\tau|)^{1/2}(N|\tau|)^4 q_u(T)^{1/4}\,.
\end{aligned}
\end{equation}
Each gradient norm denotes the sum of the norms of its two real coordinate derivatives, with $W$ derivatives taken at fixed $(Q,z)$. For $D_0\in L^4(\T^2_X;L^4(\mathcal A))$, we have $\|D_0^*D_0\|_{L^2(\T^2_X;L^2(\mathcal A))}=\|D_0\|_{L^4(\T^2_X;L^4(\mathcal A))}^2$. Using the multiplier bounds for $D_0=S_j u_{\mathcal A}$ and $D_0=T_j u_{\mathcal A}$ in~\eqref{cur:density} and~\eqref{en:wallK} and summing over $j=0,1$, we obtain
\begin{align}
 \|\mu_{\mathcal A}\|_{L^2(\T^2_X;L^2(\mathcal A))}
 &\le C (N/|\tau|)(N|\tau|)^8 q_u(T)^{1/2}\,,\label{en:densitybound}\\
 \|K_{\mathcal A}\|_{L^2(\T^2_X;L^2(\mathcal A))}
 &\le C |\tau|^{-2}(N/|\tau|)(N|\tau|)^8 e^{-c_wN|\tau|}q_u(T)^{1/2}
       \qquad (W\in\partial D)\,.\label{en:wallbound}
\end{align}

\medskip
We use~\eqref{en:densitybound} to estimate the Dirichlet source. In~\eqref{en:fourier}, $f_\ell=\widehat\mu_{\mathcal A}(\ell)/(4|\mathscr J|^2)$. Parseval's identity in $X$, $|\mathscr J|^{-2}\lesssim |\tau|^2$, and integration over the fixed disk $D$ give, uniformly in $Q\in D_Q$,
\begin{equation}
 \left(\sum_\ell\|f_\ell\|_{L^2(D;Z_{\mathcal A})}^2\right)^{1/2}
 \le C |\tau|^2(N/|\tau|)(N|\tau|)^8 q_u(T)^{1/2}
 = CN|\tau| (N|\tau|)^8 q_u(T)^{1/2}\,.
 \label{en:fbound}
\end{equation}
By Lemma~\ref{en:fiber-solve}, the norm of $P_\ell^{-1}$ on $L^2(D;Z_{\mathcal A})$ is at most $\lambda_1(D)^{-1}$ for every $\ell\in\Gamma^*$. Using this estimate in~\eqref{en:positive} and applying~\eqref{en:fbound} at each $Q$ gives
\begin{equation*}
 \mathcal E(T,\tau)\le 2a^2\lambda_1(D)^{-1}\sum_{\mathcal A}\int_{D_Q}
 \sum_\ell\|f_\ell\|_{L^2(D;Z_{\mathcal A})}^2\,d^2Q
 \le Ca^2(N|\tau|)^{18}q_u(T)\,.
\end{equation*}
This proves~\eqref{en:energyupper}.

\medskip
It remains to combine~\eqref{en:wallbound} with a normal derivative bound uniform in $\ell$. For a smooth potential $V:\overline D\to[0,\infty)$ and a source $f\in L^2(D;\mathbb C)$, let $w=(-\Delta_D+V)^{-1}f\in H^2(D)\cap H^1_0(D)$. The map $f\mapsto\partial_n w$ from $L^2(D)$ to $L^2(\partial D)$ satisfies
\begin{equation}
 \|\partial_n w\|_{L^2(\partial D)}\le (R/2)^{1/2}\|f\|_{L^2(D)}\,.
 \label{en:normal}
\end{equation}
The normal derivative $\partial_n w$ and the arc length measure on $\partial D$ are those fixed in Section~\ref{en:green-section}. We first prove~\eqref{en:normal} for smooth $f$ by duality. For $\phi\in C^\infty(\partial D;\mathbb C)$, take a smooth extension $\widetilde\phi$ to $\overline D$. Then $b=\widetilde\phi-(-\Delta_D+V)^{-1}(( -\Delta+V)\widetilde\phi)$ solves $(-\Delta+V)b=0$ with $b|_{\partial D}=\phi$. Since $\Delta|b|^2=2|\nabla b|^2+2V|b|^2\ge0$, the circle mean $M_b(t)=(2\pi)^{-1}\int_0^{2\pi}|b(te^{k\theta})|^2\,d\theta$ satisfies $M_b'(t)=(2\pi t)^{-1}\int_{|W|<t}\Delta|b|^2\,d^2W\ge0$ for $0<t<R$, by the divergence theorem. Consequently
\begin{equation*}
 \int_D|b|^2\,d^2W=2\pi\int_0^R t M_b(t)\,dt
 \le \pi R^2 M_b(R)=\frac R2\int_{\partial D}|\phi|^2\,dS_W\,.
\end{equation*}

We use this bound for $b$ to estimate the pairing of $\partial_n w$ with $\phi$. Subtracting $\int_D w\overline{(-\Delta+V)b}=0$ from $\int_D f\bar b$ cancels the potential. Green's formula and $w|_{\partial D}=0$ then yield $\int_D f\bar b=-\int_{\partial D}\partial_n w\,\bar\phi$. Taking the supremum over smooth \(\phi\) with
\(\|\phi\|_{L^2(\partial D)}\le1\) proves~\eqref{en:normal} for
smooth \(f\). Choose smooth $f_j\to f$ in $L^2(D)$, indexed by $j\in\mathbb N$, and put $w_j=(-\Delta_D+V)^{-1}f_j$. The bounded inverse gives $w_j\to w$ in $L^2(D)$. The Dirichlet $H^2$ estimate for this fixed potential $V$ gives
\begin{equation*}
 \|w_j-w\|_{H^2(D)}
 \le C_{D,V}\bigl(\|f_j-f\|_{L^2(D)}
                   +\|w_j-w\|_{L^2(D)}\bigr)\longrightarrow0\,.
\end{equation*}
The normal trace is continuous from $H^2(D)$ to $H^{1/2}(\partial D)$, and hence $\partial_n w_j\to\partial_n w$ in $L^2(\partial D)$. Passing to the limit in~\eqref{en:normal} retains its constant $(R/2)^{1/2}$, independently of $V$.

On $\partial D$ we have $V_{\mathcal A}/r=(R/2)\partial_n V_{\mathcal A}$ by the boundary quotient formula in Lemma~\ref{en:green}, with the derivative taken at fixed $Q$. Apply~\eqref{en:normal} to each coordinate of $-\widehat V_{\mathcal A}(\ell)=P_\ell^{-1}f_\ell$ in an orthonormal basis of $Z_{\mathcal A}$. Summing the squared estimates and using Plancherel in $X$ gives, at each $Q\in D_Q$,
\begin{equation*}
 \|V_{\mathcal A}/r\|_{L^2(\T^2_X\times\partial D;L^2(\mathcal A))}^2
 =(R/2)^2\sum_\ell\|\partial_n\widehat V_{\mathcal A}(\ell)\|_{L^2(\partial D;Z_{\mathcal A})}^2
 \le (R/2)^3\sum_\ell\|f_\ell\|_{L^2(D;Z_{\mathcal A})}^2\,.
\end{equation*}
By~\eqref{en:fbound}, the norm of $V_{\mathcal A}/r$ on $\T^2_X\times\partial D$ is at most $C(N|\tau|)^9 q_u(T)^{1/2}$. In the formula for $\mathcal B$ in Lemma~\ref{en:green}, use $G(r)^*A_HG(r)=\mathrm{tr}_{\mathcal A}(HK_{\mathcal A})$ from~\eqref{en:wallK} and sum over $H$. With $\delta_0(r)d^2W=(R/2)dS_W$, this gives
\begin{equation*}
 \mathcal B(z)=\frac{Ra^2}{2}\sum_{\mathcal A}\int_{D_Q}
 \int_{\T^2_X\times\partial D}|\mathscr J|^{-2}
 \mathrm{tr}_{\mathcal A}\bigl((V_{\mathcal A}/r)K_{\mathcal A}\bigr)\,dX\,dS_W\,d^2Q\,.
\end{equation*}
By~\eqref{en:measure}, $|\mathscr J|^{-2}\le C|\tau|^2$. At each $Q$, apply Cauchy--Schwarz in the trace and in $(X,W)\in\T^2_X\times\partial D$. Using~\eqref{en:wallbound} and the bound for $V_{\mathcal A}/r$, we get
\begin{align*}
 |\mathcal B(T,\tau)|
 &\le Ca^2|\tau|^2\sum_{\mathcal A}\int_{D_Q}
 \|V_{\mathcal A}/r\|_{L^2(\T^2_X\times\partial D;L^2(\mathcal A))}
 \|K_{\mathcal A}\|_{L^2(\T^2_X\times\partial D;L^2(\mathcal A))}\,d^2Q\\
 &\le Ca^2(N|\tau|)^{17}(N/|\tau|)e^{-c_wN|\tau|}q_u(T)\\
 &=C(a/|\tau|)^2(N|\tau|)^{18}e^{-c_wN|\tau|}q_u(T)\,.
\end{align*}
This proves~\eqref{en:upperbounds}.
\end{proof}

\subsection{Recovery of the fourth moment}

The first inequality in~\eqref{en:closure} requires a lower bound for $\mathcal E(T,\tau)$ in terms of $q_u(T)$ at $|\tau|=2a$. We prove the lower bound throughout $2a\le|\tau|\le3$ so that it can also be combined with~\eqref{en:upperbounds} to control $|\mathcal B(T,\tau)|$ by $V_N(|\tau|)\mathcal E(T,\tau)$. We test the Dirichlet equation for $V_{\M}$ against $\mu_{\M}$ multiplied by a smooth function supported near $W=0$. There, we can recover $u$ from $S_0u$ by a spatial multiplier. This gives a lower bound for the contribution of $(S_0u)^*S_0u$ to the test pairing.

\begin{lemma}\label{en:recovery}
Fix the annulus, disks and cutoffs of Lemma~\ref{en:upper}. There is $C_*>0$ such that, for every real $\gamma\in(0,R/8)$, there is $c_\gamma>0$ with the following property. For all real $K\ge1$ and $N\ge2$, every finite Fourier sum $u$ of the form~\eqref{pre:wave} with frequencies in the forward annular cap of Lemma~\ref{en:upper}, and every $(T,\tau)\in\R\times\mathbb C_k$ with $2a\le |\tau|\le3$, we have
\begin{equation}
 \mathcal E(T,\tau)\ge c_\gamma a^2(N|\tau|)^{-52}e^{-C_*\gamma^2N|\tau|}q_u(T)\,.
 \label{en:lower}
\end{equation}
Here $a=K/N$. The constant $C_*$ depends only on the fixed annulus, disks and cutoffs; $c_\gamma$ may also depend on $\gamma$. Both are independent of $K,N$, the cap axis, period lattice, finite frequency set, coefficients, and the algebra and trace.
\end{lemma}

\begin{proof}
We prove~\eqref{en:lower} by controlling the pairing with $\mu=\mu_{\M}$ through the contribution of $V=V_{\M}$ to~\eqref{en:positive}. At fixed $z$ consider smooth self-adjoint test fields $\Xi:\T^2_X\times\overline D\times D_Q\to Z_{\M}$ vanishing at $|W|=R$. Write $W=W_1+kW_2$ with $W_1,W_2\in\R$. We use fields with finite value of the Dirichlet form
\begin{equation}
 \mathfrak d(\Xi)=4a^2\sum_\ell\int_{D_Q}\int_D
 \left(\sum_{j=1}^2\|\partial_{W_j}\widehat\Xi(\ell)\|_{L^2(\M)}^2+
   |\mathscr J|^{-2}|\ell\cdot d|^2\|\widehat\Xi(\ell)\|_{L^2(\M)}^2\right)d^2W\,d^2Q\,.
 \label{en:form}
\end{equation}
Write $\mathfrak d(\cdot,\cdot)$ for the associated real bilinear form. At each $(\ell,Q)$ integrate the gradient term by parts with zero boundary values. Then~\eqref{en:fourier} gives $P_\ell(-\widehat V(\ell))=f_\ell$. Consequently
\[
 \mathfrak d(-V,\Xi)
 =4a^2\sum_\ell\int_{D_Q}\operatorname{Re}
    \langle f_\ell,\widehat\Xi(\ell)\rangle_{L^2(D;Z_\M)}\,d^2Q
 =\int_{\Omega_z}\operatorname{tr}_{\M}(\mu\Xi)\,d\zeta
 =:\langle\mu,\Xi\rangle_{d\zeta}.
\]
The fiber trace pairing follows from Parseval's identity in $X$ and $d\zeta=a^2|\mathscr J|^{-2}\,dX\,d^2W\,d^2Q$. Taking $\Xi=-V$ shows $\mathfrak d(-V)\le2\mathcal E$ since the contribution from $\M^{\mathrm{op}}$ in~\eqref{en:positive} is nonnegative. The Cauchy--Schwarz inequality for the form~\eqref{en:form} now yields
\begin{equation}
 |\langle\mu,\Xi\rangle_{d\zeta}|^2\le 2\mathcal E\,\mathfrak d(\Xi).
 \label{en:dualtest}
\end{equation}

Choose a smooth cutoff $\chi:D\times D_Q\to[0,1]$, supported where $|W|<\gamma$, $|Q|<1/400$, and equal to one on $|W|\le\gamma/2$, $|Q|\le1/800$. With $\Xi=\chi\mu$ in~\eqref{en:dualtest}, we need an upper bound for $\mathfrak d(\chi\mu)$ and a lower bound for $\langle\mu,\chi\mu\rangle_{d\zeta}$. We estimate $\mathfrak d(\chi\mu)$ using the derivatives of $S_j u$ in~\eqref{en:Sbounds}. For the pairing we will use a left inverse for $S_0$ on the spatial Fourier support of $u$.

\medskip
To bound $\mathfrak d(\chi\mu)$, observe that the fiber cutoff $\chi_{\mathrm{fib}}(Q)$ equals one on the support of $\chi$. Thus~\eqref{cur:density} gives $\Xi=\chi\sum_{j=0}^1 D_j^*D_j$ with $D_j=S_j u$ for $j=0,1$. If $\partial$ is a first real $W$ or $X$ derivative, then
\[
 \partial(\chi D_j^*D_j)
   =(\partial\chi)D_j^*D_j+
          \chi(\partial D_j)^*D_j+\chi D_j^*(\partial D_j).
\]
The $W$ derivatives are taken at fixed $(z,Q)$, and $\nabla_X\chi=0$. H\"older's inequality for the trace at each $X$, followed by H\"older's inequality in $X$, bounds each product $(\partial D_j)^*D_j$ or $D_j^*(\partial D_j)$ in $L^2(\T^2_X;L^2(\M))$ by $\|\partial D_j\|_{L^4(\T^2_X;L^4(\M))}\|D_j\|_{L^4(\T^2_X;L^4(\M))}$. Applying~\eqref{en:Sbounds} to the two factors gives the following estimate at each $(W,Q)$, with $C_\gamma$ accounting for the $W$ derivatives of $\chi$,
\[
 \|\nabla_W\Xi\|_{L^2(\T^2_X;L^2(\M))}
 +N^{-1}\|\nabla_X\Xi\|_{L^2(\T^2_X;L^2(\M))}
 \le C_\gamma(N/|\tau|)(N|\tau|)^9 q_u(T)^{1/2}.
\]
Since $|\mathscr J|^{-2}|\ell\cdot d|^2\le C|\tau|^2|\ell|^2$, Parseval's identity in $X$ applied to~\eqref{en:form} bounds the potential term by an integral of $C|\tau|^2\|\nabla_X\Xi\|_{L^2(\T^2_X;L^2(\M))}^2$. The two disk areas are bounded, $|\tau|\ge2K/N\ge N^{-1}$, and $1+(N|\tau|)^2\le2(N|\tau|)^2$. Hence
\begin{equation}
 \mathfrak d(\Xi)\le C_\gamma a^2(N/|\tau|)^2(1+(N|\tau|)^2)(N|\tau|)^{18}q_u(T)
 \le C_\gamma a^2N^4(N|\tau|)^{20}q_u(T).
 \label{en:testcost}
\end{equation}

\medskip
The summand $(S_0u)^*S_0u$ in $\mu$ will give the lower bound for $\langle\mu,\chi\mu\rangle_{d\zeta}$. For $|W|\le\gamma$, we construct a bounded periodic multiplier that maps $S_0u$ to $u$. Fix $\psi\in C_c^\infty((-2,2))$ equal to one on $[-1,1]$. For $\eta\in\R^2$ with $\eta_{b_0}>0$, let
\[
 \vartheta_\gamma(\eta)=\vartheta_{\rm wide}(\eta)\psi(v_\eta/\gamma)\,,
\]
and set $\vartheta_\gamma=0$ for $\eta_{b_0}\le0$. Then $\vartheta_\gamma(\xi/N)=1$ for every $\xi\in\Lambda$, and $|v_\eta|\le2\gamma$ on $\operatorname{supp}\vartheta_\gamma$. For $\eta\in\R^2$ with $\eta_{b_0}>0$, the real factor in the symbol of $S_0$ is
\[
 r_W(\eta)=\operatorname{Re}_k V_T^{v_\eta}
          =(1+v_\eta\operatorname{Re}_k W/2)/\sqrt2 .
\]
For $\eta\in\operatorname{supp}\vartheta_\gamma$, we have $r_W(\eta)\ge(1-\gamma^2)/\sqrt2>0$ and $|W-v_\eta|\le3\gamma$. Let $C_Y$ be the bound for $|Y|/|\tau|$ in~\eqref{en:measure}. By~\eqref{en:gaussian},
\begin{equation*}
 0\le F_W(\eta)\le C_{\rm G}\gamma^2,
 \qquad
 C_{\rm G}=\frac92 C_Y\sup_{\eta\in\operatorname{supp}\vartheta_{\rm wide}}\kappa_\eta\,.
\end{equation*}
In particular $e^{-C_{\rm G}\gamma^2N|\tau|}\le e^{-N|\tau|F_W}\le1$ on $\operatorname{supp}\vartheta_\gamma$.

Extend $m_\xi$ to all $\xi\in\R^2$ with $\xi_{b_0}>0$ by~\eqref{en:gaussian}. For fixed $(W,Q,z)$ with $|W|\le\gamma$, define the real symbol $b:\R^2\to\R$ in scaled frequency by
\begin{equation*}
 b(\eta)=\frac{\vartheta_\gamma(\eta)}{m_{N\eta}r_W(\eta)}
 =\left(\frac{|\tau|}{N}\right)^{1/2}
   \vartheta_\gamma(\eta)
   \left(\frac{|Y|/|\tau|}{\kappa_\eta}\right)^{1/2}
   r_W(\eta)^{-1}e^{N|\tau|F_W(\eta)}\,.
\end{equation*}
Use this formula on a neighborhood of $\operatorname{supp}\vartheta_\gamma$ where $r_W>0$, and extend it by zero outside that neighborhood. Such a neighborhood exists because the support is compactly contained in $|v_\eta|<2\gamma$, where $r_W\ge(1-\gamma^2)/\sqrt2$. Thus $b$ is smooth and compactly supported. Let $I_b$ be the periodic multiplier on $L^4(\T^2_X;L^4(\M))$ with symbol $\xi\mapsto b(\xi/N)$, $\xi\in\Gamma^*$, as in Lemma~\ref{pre:scalar-kernels}\ref{pre:periodic-kernel}--\ref{pre:convolution}. For $\xi\in\Lambda$, $b(\xi/N)m_\xi r_W(\xi/N)=1$, so $I_b(S_0u)=u$.

For a multi-index $\nu\in\{0,1,2,\ldots\}^2$, let $\partial_\eta^\nu$ denote differentiation in the scaled frequency coordinates. The derivatives of $v_\eta$ and $\kappa_\eta$ through order four are uniformly bounded on the fixed wide patch, and the chain rule gives $|\partial_\eta^\nu\vartheta_\gamma|\le C_\nu\gamma^{-|\nu|}$ for $|\nu|\le4$. Differentiating $r_Wr_W^{-1}=1$ gives, for $1\le|\nu|\le4$,
\begin{equation*}
 \partial_\eta^\nu(r_W^{-1})
 =-r_W^{-1}\sum_{0<\alpha\le\nu}\binom{\nu}{\alpha}
       (\partial_\eta^\alpha r_W)
       \partial_\eta^{\nu-\alpha}(r_W^{-1})\,.
\end{equation*}
Here $\alpha\le\nu$ is understood componentwise. The positive lower bound for $r_W$ and induction on $|\nu|$ bound every derivative of $r_W^{-1}$ through order four. On the wide patch, $\kappa_\eta\ge c>0$. The derivatives of $s\mapsto s^{-1/2}$ through order four are bounded for $s\ge c$, so the chain rule bounds $\partial_\eta^\nu(\kappa_\eta^{-1/2})$ for $|\nu|\le4$. Since $|Y|/|\tau|$ is independent of $\eta$, the nonexponential factors in $b$, apart from $(|\tau|/N)^{1/2}$, have derivatives bounded by $C_\gamma$.

Each derivative of $e^{N|\tau|F_W}$ either differentiates a derivative of $F_W$ already present or adds one factor $N|\tau|$ and a first derivative of $F_W$. Thus, for $|\nu|\le4$, the bounded derivatives of $F_W$ and $N|\tau|\ge2$ give
\begin{equation*}
 |\partial_\eta^\nu e^{N|\tau|F_W}|
 \le C_\nu(N|\tau|)^{|\nu|}e^{N|\tau|F_W}\,.
\end{equation*}
Leibniz' rule and $F_W\le C_{\rm G}\gamma^2$ on the support now yield
\begin{equation*}
 \max_{|\nu|\le4}\|\partial_\eta^\nu b\|_{L^\infty(\R^2)}
 \le C_\gamma(|\tau|/N)^{1/2}(N|\tau|)^4
                 e^{C_{\rm G}\gamma^2N|\tau|}\,.
\end{equation*}
The prefactor $C_\gamma$ includes the derivatives of $\vartheta_\gamma$. The constant $C_{\rm G}$ depends only on the fixed wide patch and the uniform bound for $|Y|/|\tau|$ in~\eqref{en:measure}, and hence is independent of $\gamma,K,N$.

Let $K_b$ be the periodic convolution kernel of $I_b$. Since $b$ is supported in a fixed annulus, Lemma~\ref{pre:scalar-kernels}\ref{pre:kernel-bounds}--\ref{pre:convolution} gives
\begin{equation*}
 \|K_b\|_{L^1(\T^2_X)}\le \|\check b\|_{L^1(\R^2)}
 \le C\sum_{|\nu|\le4}\|\partial_\eta^\nu b\|_{L^\infty(\R^2)}\,.
\end{equation*}
Using $I_b(S_0u)=u$ gives $q_u(T)^{1/4}\le\|K_b\|_{L^1(\T^2_X)}\|S_0u\|_{L^4(\T^2_X;L^4(\M))}$. The derivative estimate for $b$ therefore gives
\begin{equation}\label{en:S0lower}
 \|S_0u\|_{L^4(\T^2_X;L^4(\M))}
 \ge c_\gamma(N/|\tau|)^{1/2}(N|\tau|)^{-4}
              e^{-C_{\rm G}\gamma^2N|\tau|}q_u(T)^{1/4}\,.
\end{equation}

\medskip
We use the lower bound for $\|S_0u\|_{L^4(\T^2_X;L^4(\M))}$ to estimate $\langle\mu,\chi\mu\rangle_{d\zeta}$ from below. On the support of $\chi$, set $B_j=(S_j u)^*S_j u\in L^2(\M)_+$ for $j=0,1$ at each $X$. By~\eqref{cur:density}, $\mu=B_0+B_1$ there. Since $B_i^{1/2}\in L^4(\M)$, cyclicity and H\"older's inequality for the trace give
\[
 \operatorname{tr}_{\M}(B_0 B_1)
  =\operatorname{tr}_{\M}(B_0^{1/2}B_1 B_0^{1/2})
  =\|B_1^{1/2} B_0^{1/2}\|_{L^2(\M)}^2\ge0.
\]
The first equality uses cyclicity for $B_1 B_0^{1/2}\in L^{4/3}(\M)$ and $B_0^{1/2}\in L^4(\M)$. Expanding the trace of the square, with $\operatorname{tr}_{\M}(B_1B_0)=\operatorname{tr}_{\M}(B_0B_1)$, now gives
\[
 \|\mu\|_{L^2(\M)}^2
 =\|B_0\|_{L^2(\M)}^2+2\operatorname{tr}_{\M}(B_0B_1)+\|B_1\|_{L^2(\M)}^2
 \ge \|B_0\|_{L^2(\M)}^2=\|S_0u\|_{L^4(\M)}^4 .
\]
We integrate this inequality over $X$ and restrict $(W,Q)$ to $|W|\le\gamma/2$, $|Q|\le1/800$, where $\chi=1$. Using $a^2|\mathscr J|^{-2}\ge ca^2|\tau|^2$ from~\eqref{en:measure} and then~\eqref{en:S0lower}, we obtain
\begin{align*}
 \langle\mu,\Xi\rangle_{d\zeta}
 &=\int_{\Omega_z}\chi\|\mu\|_{L^2(\M)}^2\,d\zeta\\
 &\ge ca^2|\tau|^2\int_{\{|W|\le\gamma/2,\,|Q|\le1/800\}}
 \|S_0u\|_{L^4(\T^2_X;L^4(\M))}^4\,d^2W\,d^2Q\\
 &\ge c_\gamma a^2N^2(N|\tau|)^{-16}
             e^{-4C_{\rm G}\gamma^2N|\tau|}q_u(T)\,.
\end{align*}
Here $a^2|\tau|^2(N/|\tau|)^2=a^2N^2$, and $c_\gamma$ absorbs the area of the disk $|W|\le\gamma/2$. If $q_u(T)>0$, the pairing is positive, so $\mathfrak d(\Xi)>0$ by~\eqref{en:dualtest}. Combining the lower bound for the pairing with~\eqref{en:testcost} gives
\begin{equation*}
 \mathcal E(T,\tau)\ge\frac{|\langle\mu,\Xi\rangle_{d\zeta}|^2}{2\mathfrak d(\Xi)}
 \ge c_\gamma a^2(N|\tau|)^{-32-20}e^{-8C_{\rm G}\gamma^2N|\tau|}q_u(T)\,.
\end{equation*}
This is~\eqref{en:lower} with $C_*=8C_{\rm G}$. If $q_u(T)=0$,~\eqref{en:lower} follows from $\mathcal E\ge0$.

\end{proof}

\subsection{Radial comparisons and the periodic cap estimate}\label{en:comparison}\label{en:completion}

We now prove the comparison between $\int_1^2 H(T,2a)\,dT$ and $H(0,4a)$ in~\eqref{en:closure}. We first derive a differential inequality for~$\mathcal E$ by bounding $|\mathcal B(T,\tau)|$ in terms of~$\mathcal E(T,\tau)$.

For the fixed annulus, disks and cutoffs, Lemmas~\ref{en:upper} and~\ref{en:recovery} give $c_w,C_*>0$ independently of $\gamma,K,N$. Choose $0<\gamma<R/8$ so that $C_*\gamma^2<c_w/2$, and fix $K=10$. Let $N\ge2$ be real and recall $a=K/N$. Since $\mathcal Q\ge0$, \eqref{en:balance} gives $\Delta_z\mathcal E\le|\mathcal B|$. Substituting the upper bound for $q_u(T)$ from~\eqref{en:lower} into the estimate for $|\mathcal B(T,\tau)|$ in~\eqref{en:upperbounds} gives
\begin{equation*}
 |\mathcal B(T,\tau)|
 \le C_\gamma|\tau|^{-2}(N|\tau|)^{70}
       e^{-(c_w-C_*\gamma^2)N|\tau|}\mathcal E(T,\tau)
 \qquad (2a\le|\tau|\le3)\,.
\end{equation*}
The factors $a^2$ cancel, and the powers of $N|\tau|$ add to $18+52=70$. Since $x^{70}e^{-c_w x/4}$ is bounded for $x\ge0$ and $c_w-C_*\gamma^2>c_w/2$, we obtain
\begin{equation}\label{en:sub}
 \Delta_z\mathcal E(T,\tau)\le V_N(|\tau|)\mathcal E(T,\tau),
 \qquad V_N(s)=C_\gamma s^{-2}e^{-c_wNs/4}
 \quad(2a\le s\le3)\,.
\end{equation}

For the annular comparison, assume $2a<1/4$. Then $[1/2,3/2]\subset[2a,3]$ and $V_N(s)\le M:=4C_\gamma$ for $s\in[1/2,3/2]$. The radial barriers below use the bound
\begin{equation}\label{en:moment}
 \int_{2a}^3s\log\!\left(\frac{s}{2a}\right)V_N(s)\,ds
 \le C_\gamma\int_{2K}^\infty
    e^{-c_wx/4}\log\!\left(\frac{x}{2K}\right)\frac{dx}{x}
 \le C\,,
\end{equation}
where we used $x=Ns$ and $2Na=2K$. With $\gamma$ and $K$ fixed, both $M$ and the bound in~\eqref{en:moment} are independent of~$N$.

\begin{lemma}\label{en:annulus}
For $a=K/N$ with $2a<1/4$, the circle mean $H$ defined in~\eqref{en:circlemean} satisfies
\begin{equation}
 \int_1^2 H(T,2a)\,dT\le C\log^2(3/(2a)) H(0,4a)\,.
 \label{en:ann}
\end{equation}
\end{lemma}

\begin{proof}
To prove~\eqref{en:ann}, define $j:[0,3]\to[0,1]$ by $j(T)=\sin(\pi T/3)$ and $m:[2a,3]\to[0,\infty)$ by $m(s)=\int_0^3 j(T)H(T,s)\,dT$. Since $j\ge\sqrt3/2$ on $[1,2]$, it suffices to bound $m(2a)$ in terms of~$H(0,4a)$. Write $h_*=\inf_{|T|\le1}H(T,1)$ and $L_S=\log(S/(2a))$ for $S\in\{1,3\}$. We first bound $h_*$ from below by $c\,m(2a)/L_3$, then bound $H(0,4a)$ from below by $c\,h_*/L_1$.

In both comparisons we use the differential inequality for the circle mean~\eqref{en:circlemean}. In polar coordinates for~$\tau$, the angular second derivative integrates to zero. Averaging~\eqref{en:sub} therefore gives
\begin{equation*}
 \mathscr L_N H\le0,\qquad
 \mathscr L_N=\partial_T^2+\partial_s^2+s^{-1}\partial_s-V_N(s)
 \qquad(T\in\R,\ 2a\le s\le3)\,.
\end{equation*}

\medskip
For both comparisons we need radial solutions with prescribed endpoint values. For $S\in\{1,3\}$ and real $\lambda>0$, we seek real-valued functions $b_{01}^{S,\lambda},b_{10}^{S,\lambda}\in C^2([2a,S])$ solving
\begin{equation*}
 b''+s^{-1}b'=(V_N(s)+\lambda^2)b\qquad(2a<s<S)
\end{equation*}
with endpoint values $(0,1)$ and $(1,0)$, respectively. Set $A=L_S$ and $x=\log(s/(2a))\in[0,A]$. Define $P:[0,A]\to(0,\infty)$ by $P(x)=s^2(V_N(s)+\lambda^2)$ with $s=2ae^x$. Since $\partial_x=s\partial_s$, the equation for $w(x)=b(2ae^x)$ becomes $w''=Pw$. By~\eqref{en:moment},
\begin{equation*}
 \int_0^A \min(x,A-x)P(x)\,dx
 \le C+\lambda^2S^2\int_0^\infty t e^{-2t}\,dt=:C_{\lambda,S}\,.
\end{equation*}
For the $V_N$ term, use $dx=ds/s$ and $\min(x,A-x)\le\log(s/(2a))$. For the $\lambda^2$ term, use $t=A-x$, $s=Se^{-t}$ and $\min(x,A-x)\le t$.

Let $w$ solve $w''=Pw$ with $w(0)=0$ and $w'(0)=1$. Integrating twice gives
\begin{equation*}
 w(x)=x+\int_0^x(x-t)P(t)w(t)\,dt\,.
\end{equation*}
The solution $w$ is positive near zero. As long as $w\ge0$, we have $w''\ge0$ and $w(x)\ge x$, which precludes a first positive zero. Define $v:[0,A]\to[1,\infty)$ by $v(x)=w(x)/x$ for $x>0$ and $v(0)=1$. For $0<x\le A$ and $0\le t\le x$, we have $0\le (x-t)t/x\le\min(t,A-t)$. The integral equation and Gronwall's inequality give
\begin{equation*}
 1\le v(x)\le1+\int_0^x\min(t,A-t)P(t)v(t)\,dt
       \le\exp\!\left(\int_0^x\min(t,A-t)P(t)\,dt\right)
       \le e^{C_{\lambda,S}}\,.
\end{equation*}
Thus $b_{01}^{S,\lambda}(s)=w(x)/w(A)$ has endpoint values $(0,1)$, and $w(x)\ge x$ and $w(A)\le Ae^{C_{\lambda,S}}$ give its lower bound. For endpoint values $(1,0)$, let $\widetilde w$ solve $\widetilde w''(x)=P(A-x)\widetilde w(x)$ with $\widetilde w(0)=0$ and $\widetilde w'(0)=1$. The substitution $t=A-x$ gives $\int_0^A\min(x,A-x)P(A-x)\,dx\le C_{\lambda,S}$. Applying the estimate for $w$ with this reflected potential gives $\widetilde w(x)\ge x$ and $\widetilde w(A)\le Ae^{C_{\lambda,S}}$. For the profile $b_{10}^{S,\lambda}(s)=\widetilde w(A-x)/\widetilde w(A)$, the numerator is at least $A-x=\log(S/s)$. Thus for every $s\in[2a,S]$,
\begin{equation}\label{en:radialbounds}
 b_{01}^{S,\lambda}(s)\ge e^{-C_{\lambda,S}}\frac{\log(s/(2a))}{L_S},
 \qquad
 b_{10}^{S,\lambda}(s)\ge e^{-C_{\lambda,S}}\frac{\log(S/s)}{L_S}\,.
\end{equation}

\medskip
To get a lower bound for $H(T,1)$ when $|T|\le1$, we first estimate $m(s)$ on $3/4\le s\le5/4$. Twice integrating by parts, using $j(0)=j(3)=0$ and $j''=-(\pi/3)^2j$, gives
\begin{equation*}
 \int_0^3 j(T)H_{TT}(T,s)\,dT
  =(\pi/3)\big(H(0,s)+H(3,s)\big)-(\pi/3)^2m(s)\,.
\end{equation*}
Since $H\ge0$ and $j\ge0$, integrating $\mathscr L_N H\le0$ against $j$ gives $m''+s^{-1}m'-(V_N+\pi^2/9)m\le0$. The function $m(s)-m(2a)b_{10}^{3,\pi/3}(s)$ satisfies the same differential inequality and has endpoint values $0$ and $m(3)\ge0$. At a negative interior minimum its first derivative would be zero and its second derivative nonnegative. Hence $\partial_s^2+s^{-1}\partial_s-(V_N+\pi^2/9)$ applied to that function would be strictly positive there. This proves $m(s)\ge m(2a)b_{10}^{3,\pi/3}(s)$ on~$[2a,3]$. In particular, \eqref{en:radialbounds} gives
\begin{equation*}
 m(s)\ge c\,m(2a)\frac{\log(3/s)}{L_3}
       \ge c'\frac{m(2a)}{L_3}
 \qquad(3/4\le s\le5/4)\,.
\end{equation*}

We use this lower bound for $m(s)$ on $3/4\le s\le5/4$ to bound $H(T,1)$ from below for every $T\in[-1,1]$. On $U=\R\times\{\tau\in\mathbb C_k:1/2<|\tau|<3/2\}$, put $H_0(T,\tau)=H(T,|\tau|)$. Since $\mathscr L_N H\le0$ and $V_N\le M$ on $[1/2,3/2]$, we have $\Delta_{T,\tau}H_0\le M H_0$. For $\omega=\sqrt{M+1}$, define $F:U\times(-\pi/(3\omega),\pi/(3\omega))\to[0,\infty)$ by $F(T,\tau,y')=H_0(T,\tau)\cos(\omega y')$. In this domain
\begin{equation*}
 \Delta_{T,\tau,y'}F
 =\cos(\omega y')(\Delta_{T,\tau}H_0-\omega^2H_0)\le0\,.
\end{equation*}
Write $B^d(p,r)$ for the open Euclidean ball of center $p$ and radius $r>0$ in $\R^d$. For any $z=(T,\tau)\in U$ and $\rho>0$ such that $B^4((z,0),2\rho)$ lies in the domain of~$F$, the four-dimensional mean inequality gives
\begin{equation}\label{en:local-mean}
 H_0(z)\ge\frac{1}{|B^4(0,2\rho)|}
                  \int_{B^4((z,0),2\rho)}F
 \ge c\,\rho^{-3}\int_{B^3(z,\rho)}H_0(z')\,dz'\,.
\end{equation}
The second inequality follows by integrating over $B^3(z,\rho)\times(-\rho,\rho)$ inside $B^4((z,0),2\rho)$, where $\cos(\omega y')\ge1/2$.

Consider the compact connected set $\mathcal C=[-1,3]\times\{\tau\in\mathbb C_k:3/4\le|\tau|\le5/4\}$. Choose $\rho>0$ depending only on $M$ and~$\mathcal C$, small enough that $B^4((z,0),2\rho)$ lies in the domain of $F$ whenever $z$ is within distance $\rho/8$ of~$\mathcal C$. Since $j\le1$, polar integration gives
\begin{equation*}
 \int_{\mathcal C}H_0(z)\,dz
 \ge 2\pi\int_{3/4}^{5/4}s\,m(s)\,ds
 \ge c\,\frac{m(2a)}{L_3}\,.
\end{equation*}
Cover $\mathcal C$ by $L$ open balls $B_1,\ldots,B_L$ of radius $\rho/8$ with centers in~$\mathcal C$. Their intersection graph is connected because every ball meets $\mathcal C$ and $\mathcal C$ is connected. Relabel the balls so that $\int_{B_1}H_0\ge c\,m(2a)/(L L_3)$. If $B_i\cap B_j\ne\varnothing$, then $B_i\subset B^3(z,\rho)$ for every $z\in B_j$ by the triangle inequality. Applying~\eqref{en:local-mean} at such $z$ and integrating over $B_j$ gives
\begin{equation*}
 \int_{B_j}H_0\ge c\rho^{-3}|B_j|\int_{B_i}H_0\,.
\end{equation*}
For a target point $z=(T,\tau)$ with $|T|\le1$ and $|\tau|=1$, choose a ball $B_j$ containing~$z$. A chain of at most $L-1$ overlaps yields $\int_{B_j}H_0\ge c\,m(2a)/L_3$. Since $B_j\subset B^3(z,\rho)$, \eqref{en:local-mean} now gives
\begin{equation*}
 h_*=\inf_{|T|\le1}H(T,1)\ge c\,\frac{m(2a)}{L_3}\,.
\end{equation*}
The number of balls $L$, their common volume and $\rho$ depend only on $M$ and~$\mathcal C$, so the constant is independent of~$N$.

\medskip
The lower bound at $s=1$ gives boundary data for a comparison that reaches $s=4a$. On the rectangle $\mathscr R=[-1,1]\times[2a,1]$, define $G:\mathscr R\to[0,\infty)$ by
\begin{equation*}
 G(T,s)=h_*\cos(\pi T/2)b_{01}^{1,\pi/2}(s)\,.
\end{equation*}
The equation for $b_{01}^{1,\pi/2}$ gives $\mathscr L_NG=0$. On $T=\pm1$ and on $s=2a$ we have $G=0\le H$. On $s=1$ we have $G(T,1)\le h_*\le H(T,1)$. Thus $H-G\ge0$ on the boundary of $\mathscr R$ and $\mathscr L_N(H-G)\le0$ in its interior. At a negative interior minimum of $H-G$, the second derivatives in $T$ and $s$ would be nonnegative, the first derivatives would vanish, and $-V_N(H-G)$ would be strictly positive. Hence $H\ge G$ on~$\mathscr R$. At $(T,s)=(0,4a)$, \eqref{en:radialbounds} gives
\begin{equation*}
 H(0,4a)\ge h_*b_{01}^{1,\pi/2}(4a)
 \ge c\,h_*\frac{\log2}{L_1}\,.
\end{equation*}
Using $j\ge\sqrt3/2$ on $[1,2]$ and the two lower bounds for $h_*$ and $H(0,4a)$, we conclude that
\begin{equation*}
 \int_1^2 H(T,2a)\,dT
 \le \frac{2}{\sqrt3}m(2a)
 \le C L_3 h_*
 \le C' L_3 L_1 H(0,4a)
 \le C' L_3^2 H(0,4a)\,,
\end{equation*}
which is~\eqref{en:ann}.
\end{proof}

To complete the proof of Theorem~\ref{thm:four} via the reduction in Section~\ref{pre:periodic-cap}, we combine Lemma~\ref{en:annulus} with the energy bounds~\eqref{en:energyupper} and~\eqref{en:lower} for $2a<1/4$ and use a direct multiplier bound for $2a\ge1/4$.

\begin{proof}[Proof of Proposition~\ref{en:cap}]
Use the aperture chosen before~\eqref{en:sub} and $a=K/N$ with $K=10$. We first prove~\eqref{en:capestimate} when $2a\ge1/4$, or equivalently when $2\le N\le8K$. Choose $\chi\in C_c^\infty(\R^2\setminus\{0\})$ equal to one on $\{\eta\in\R^2:c_0\le|\eta|\le C_0\}$ and define $p_{T,N}:\R^2\to\mathbb C$ by $p_{T,N}(\eta)=\chi(\eta)e^{-iTN|\eta|}$. By~\eqref{pre:wave}, the periodic multiplier with symbol $p_{T,N}(\xi/N)$ for $\xi\in\Gamma^*$ sends $u(0,\cdot)$ to $u(T,\cdot)$. For $T\in[1,2]$ and $2\le N\le8K$, the derivatives of $p_{T,N}$ through order four are uniformly bounded. Lemma~\ref{pre:scalar-kernels} gives $q_u(T)\le Cq_u(0)$ with $C$ independent of the period lattice and the coefficient algebra. Integrating over $[1,2]$ proves~\eqref{en:capestimate} in this range.

If $2a<1/4$, apply~\eqref{en:energyupper} at $|\tau|=4a$ and~\eqref{en:lower} at $|\tau|=2a$. Since $N|\tau|=4K$ and $N|\tau|=2K$, respectively, taking circle means gives
\begin{equation*}
 H(0,4a)\le C a^2q_u(0),\qquad c a^2q_u(T)\le H(T,2a)
 \qquad (T\in\R)\,,
\end{equation*}
with constants independent of~$N$. Inserting the upper bound for $H(0,4a)$ and the lower bound for $H(T,2a)$ into~\eqref{en:ann} gives the chain of estimates in~\eqref{en:closure}, and hence
\begin{equation*}
 a^2\int_1^2 q_u(T)\,dT\le C a^2\log^2(3/(2a))q_u(0)\,.
\end{equation*}
Cancelling $a^2$ and using $\log(3/(2a))\lesssim\log(2+N)$ proves~\eqref{en:capestimate}. The scaled symbol bounds use the same cutoff profiles in each orthonormal frame $(e_x,e_{b_0})$. The periodic kernel estimates in Lemma~\ref{pre:scalar-kernels} and the spatial Parseval identities with probability Haar measure are uniform in the period lattice~$\Gamma$. These facts give the asserted uniformity in the frame and lattice.
\end{proof}

By the reduction in Section~\ref{pre:periodic-cap}, this completes the proof of Theorem~\ref{thm:four}.

\section{Local smoothing and quantum wave estimates}\label{sec:lp}

Section~\ref{lp:annular} interpolates the logarithmic \(L^4\) bound to obtain the single-annulus estimate~\eqref{intro:annular}, then rescales it to long time intervals for Sections~\ref{sec:cone} and~\ref{sec:stein}. Section~\ref{lp:sum-annuli} sums over annuli to obtain the Sobolev estimates, including the logarithmically weighted estimate at the critical order. At the critical order, the proof uses the annular Fourier supports of the spatial pieces of \(F\) and their preservation under \(U(t)\). After examining scalar obstructions to the unweighted Sobolev estimates, we introduce the quantum Fourier models and transfer the wave estimates to quantum tori and Moyal planes.

\subsection{Annular local smoothing estimates}\label{lp:annular}
Fix a real \(N\ge2\), a cutoff \(\psi\in C_c^\infty((1/2,2))\), and \(\nu\in\{-1,1\}\).
On test tensors define \((T_NF)(t)=U(\nu t)\psi(|D|/N)F\) for \(t\in I\).
To prove \eqref{intro:annular} we interpolate the operator
\[
 T_N:L^p(\N_{\R^2})\longrightarrow L^p(\N_{I\times\R^2})
\]
on each side of \(p=4\). The \(L^4\) bound for \(T_N\) follows from~\eqref{intro:l4}
and the adjoint identity for \(U(-t)\).
Since \(L^2(\M)\) is Hilbert, Plancherel's identity and \(|I|=1\)
give
\[
 \|T_N\|_{L^2(\N_{\R^2})\to L^2(\N_{I\times\R^2})}
 \le \|\psi\|_{L^\infty((1/2,2))}.
\]
For \(p=\infty\), the spatial kernel satisfies the uniform bound
\begin{equation}\label{lp:infty}
 \sup_{1\le t\le2}\| (e^{i\nu t|\xi|}\psi(|\xi|/N))^\vee\|_{L^1(\R^2)}
   \le C_\psi N^{1/2}.
\end{equation}
On the annulus \(\xi=N\eta\), use an angular partition into
\(O(N^{1/2})\) caps of aperture \(N^{-1/2}\). Such a partition is obtained
from translates of one smooth bump in the angle divided by their sum,
with uniformly bounded overlap and derivatives of each integer order \(k\ge0\) bounded by
\(C_kN^{k/2}\). In coordinates along a cap write
\(\eta=(N^{-1/2}z,r)\), where \((z,r)\in\R^2\), \(r\ge c>0\), and the pair ranges in a fixed
compact set. The cutoffs have bounded derivatives through order four
in these coordinates, since the scaled angular variable is
\(N^{1/2}\arctan(N^{-1/2}z/r)\). Removing the linear phase \(\nu tNr\)
translates the kernel. The remaining phase is
\[
 \nu tN\big(\sqrt{r^2+N^{-1}z^2}-r\big)
       =\nu t z^2/(\sqrt{r^2+N^{-1}z^2}+r).
\]
Its first four derivatives are bounded independently of the cap and \(N\).
By \eqref{pre:kern}, the kernel for each cap in the \((z,r)\)
coordinates has bounded \(L^1(\R^2)\) norm. The inverse linear changes of
frequency coordinates preserve this norm, as does the translation
from the removed phase. Summing over the \(O(N^{1/2})\) caps proves
\eqref{lp:infty}.
Let \(K_t=(e^{i\nu t|\xi|}\psi(|\xi|/N))^\vee\) for \(t\in I\). For fixed \(N\) the symbol depends continuously on \(t\) in the derivatives through order four on its compact support, so \eqref{pre:kern} gives continuity of \(t\mapsto K_t\) in \(L^1(\R^2)\). For \(H\in L^1(\N_{I\times\R^2})\), the Bochner integral \(P_KH=\int_I\widetilde K_t*H(t)\,dt\) defines an element of \(L^1(\N_{\R^2})\). Indeed, Lemma~\ref{pre:scalar-kernels}\ref{pre:convolution} gives
\begin{equation*}
 \|P_KH\|_{L^1(\N_{\R^2})}
 \le \sup_{t\in I}\|K_t\|_{L^1(\R^2)}\,\|H\|_{L^1(\N_{I\times\R^2})}\,.
\end{equation*}
The adjoint \(P_K^*:\N_{\R^2}\to\N_{I\times\R^2}\) defines \(T_N\) at \(p=\infty\). For distinct \(r_0,r_1\in\{2,4,\infty\}\), \(F\in L^{r_0}(\N_{\R^2})\cap L^{r_1}(\N_{\R^2})\) and a spacetime test tensor \(H\), spatial convolution and Fubini give
\(\langle H,T_N F\rangle_{I\times\R^2}=\langle P_KH,F\rangle_{\R^2}\).
For the finite endpoint exponents the kernel action extends from test tensors by the bound for spatial convolution, uniformly in \(t\). The pairing identity therefore identifies the endpoint actions on their intersections.

For \(\mathcal B=\N_{\R^2}\) and \(\mathcal B=\N_{I\times\R^2}\), we use the tracial interpolation identity
\([L^{p_0}(\mathcal B),L^{p_1}(\mathcal B)]_\theta=L^p(\mathcal B)\) isometrically, with
\(1\le p_0<p_1\le\infty\), \(0<\theta<1\), and \(1/p=(1-\theta)/p_0+\theta/p_1\)
\cite[Section 2, (2.1)]{px}. Between \(2\) and \(4\) the weight of
Theorem~\ref{thm:four} is \(2-4/p\), so the logarithmic power is
\(1-2/p\). Between \(4\) and \(\infty\) the weight of \eqref{lp:infty}
is \(1-4/p\). This gives the power \(1/2-2/p\) of \(N\) and the power
\(2/p\) of the logarithm. Thus \(\sigma(p)\) and \(\beta(p)\) are the powers in \eqref{intro:annular} for \(U(\nu t)\).
The endpoint constants are independent of \(\nu\) and of the coefficient algebra.

We extend~\eqref{intro:annular} and~\eqref{lp:infty} to \(U(\nu t)m(D/N)\) for a fixed smooth radial symbol \(m:\R^2\to\mathbb C\) supported in a compact annulus away from zero. Choose finitely many \(d_l>0\), \(1\le l\le M\) with \(M\in\mathbb N\), such that \((d_l/2,2d_l)\) cover the radial support of \(m\). A subordinate smooth partition splits \(m\) into symbols \(m_l\) supported in \(d_l/2<|\eta|<2d_l\), where \(\eta=\xi/N\). Choose a fixed \(\psi_l\in C_c^\infty((1/2,2))\) equal to one on the rescaled radial support of \(m_l\). Then
\begin{equation*}
 e^{i\nu t|\xi|}m_l(\xi/N)
 =e^{i\nu t|\xi|}\psi_l(|\xi|/(d_lN))\,m_l(\xi/N)\,.
\end{equation*}
If \(d_lN\ge2\), apply the annular bound for \(U(\nu t)\) to \(m_l(D/N)F\) with frequency parameter \(d_lN\). The kernel of the product is the convolution of the kernel in~\eqref{lp:infty} with \((m_l(\cdot/N))^\vee\). The latter has \(L^1\) norm \(\|\check m_l\|_{L^1(\R^2)}\) by dilation. Lemma~\ref{pre:scalar-kernels} thus bounds the localizing multipliers uniformly in \(N\), and the product kernel has \(L^1\) norm at most \(C_{m_l,\psi_l}(d_lN)^{1/2}\). If \(d_lN<2\), then \(N<2/d_l\) and the kernel bound~\eqref{pre:kern} applies directly to \(e^{i\nu tN|\eta|}m_l(\eta)\) in the scaled variable. Summing the finitely many pieces gives~\eqref{intro:annular} for \(U(\nu t)m(D/N)\) with the same powers of \(N\) and \(\log(2+N)\), along with the uniform kernel bound
\begin{equation*}
 \sup_{t\in I}\bigl\|(e^{i\nu t|\xi|}m(\xi/N))^\vee\bigr\|_{L^1(\R^2)}
 \le C_m N^{1/2}\,.
\end{equation*}

Rescaling~\eqref{intro:annular} gives a bound on \([S,2S]\) for waves whose spatial frequencies lie in a fixed annulus.

\begin{lemma}\label{cone:long-wave}
Fix \(0<a<b<\infty\) and \(2\le p<\infty\). If \(h\in L^p(\N_{\R^2})\) has spatial Fourier support in \(a\le |\xi|\le b\), then, for each \(\nu\in\{-1,1\}\) and every real \(S\ge1\),
\begin{equation}\label{cone:long-time}
 \left(\int_S^{2S}\|U(\nu s)h\|_{L^p(\N_{\R^2})}^p\,ds\right)^{1/p}
 \le C_{a,b,p}S^{1/p+\sigma(p)}[\log(2+S)]^{\beta(p)}\|h\|_{L^p(\N_{\R^2})} .
\end{equation}
For \(|s|\le2\) one has \(\|U(s)h\|_{L^p(\N_{\R^2})}\le C_{a,b}\|h\|_{L^p(\N_{\R^2})}\).
\end{lemma}

\begin{proof}
We first bound \(U(s)h\) for bounded \(s\), then rescale \([S,2S]\) to \([1,2]\) to prove \eqref{cone:long-time}. Choose \(\chi\in C_c^\infty((0,\infty))\) equal to one on a
neighborhood of \([a,b]\). Then \(\chi(|D|)h=h\) by Lemma~\ref{pre:localized-approximation}. For bounded \(s\), the
derivatives through order four of
\(\xi\mapsto e^{is|\xi|}\chi(|\xi|)\) are bounded on its
support. Lemma~\ref{pre:scalar-kernels} proves the assertion
for \(|s|\le2\) and gives a uniform bound on any fixed bounded
interval, with a constant depending on that interval.

For real \(S\ge1\), define \(h_S\in L^p(\N_{\R^2})\) by \(h_S(y)=h(Sy)\) for \(y\in\R^2\).
To apply \eqref{intro:annular} on \([1,2]\) we need the frequency
support and the norm of \(h_S\). Its distributional Fourier transform
is \(S^{-2}\widehat h(\eta/S)\), supported in \(aS\le|\eta|\le bS\).
For \(\tau\in\R\), substitution of \(\eta=S\xi\) in the half-wave formula and of \(x=Sy\) in the spatial norm gives
\[
 (U(\nu S\tau)h)(Sy)=(U(\nu\tau)h_S)(y),\qquad
 \|h_S\|_{L^p(\N_{\R^2})}^p=S^{-2}\|h\|_{L^p(\N_{\R^2})}^p .
\]
The wave identity follows first for Schwartz tensors and then for \(h\)
by the localized approximation \eqref{pre:localized-core}, using an
enlarged compact annulus for the approximants.

For \(S\ge2\), take \(N=S\), \(m(\eta)=\chi(|\eta|)\), and \(F=h_S\) in the fixed-annulus form of~\eqref{intro:annular} for \(U(\nu\tau)\). Since \(\chi(|D|/S)h_S=h_S\) and the radial partition depends only on \(\chi\), this gives
\[
 \|U(\nu\tau)h_S\|_{L^p([1,2];L^p(\N_{\R^2}))}
 \le C S^{\sigma(p)}[\log(2+S)]^{\beta(p)}\|h_S\|_{L^p(\N_{\R^2})} .
\]
In the norm integral, substitute \(s=S\tau\) with \(\tau\in[1,2]\), and \(x=Sy\) with \(y\in\R^2\). The time Jacobian \(S\) and spatial Jacobian \(S^2\) give
\begin{align*}
 \int_S^{2S}\|U(\nu s)h\|_{L^p(\N_{\R^2})}^p\,ds
 &= S^3\int_1^2\|U(\nu\tau)h_S\|_{L^p(\N_{\R^2})}^p\,d\tau\\
 &\le C S^{3+p\sigma(p)}[\log(2+S)]^{p\beta(p)}
                 S^{-2}\|h\|_{L^p(\N_{\R^2})}^p .
\end{align*}
Taking the \(p\)-th root gives~\eqref{cone:long-time}. For \(1\le S<2\), the localized kernel bound on \(|s|\le4\) gives~\eqref{cone:long-time} directly. With \(\nu=-1\), changing variables in \eqref{cone:long-time} gives the same bound for \(U(s)h\) on \([-2S,-S]\).
\end{proof}

\subsection{Sobolev estimates and logarithmic weights}\label{lp:sum-annuli}
We prove parts \ref{thm:lp-sob} and \ref{thm:lp-weight} of Theorem~\ref{thm:lp}. For \(s>\sigma(p)\), the \(j\)-th annular term has \(L^p(\N_{I\times\R^2})\) norm at most \(C2^{-j(s-\sigma(p))}(1+j)^{\beta(p)}\|J^sF\|_{L^p(\N_{\R^2})}\), and these coefficients are summable in \(j\).
At \(s=\sigma(p)\) we will instead combine an \(\ell^p\) estimate for
the weighted frequency pieces of \(F\) with an \(\ell^2\) estimate for the sum of the half-waves of the unweighted pieces, as in~\eqref{lp:analysis} and~\eqref{lp:synthesis} below.

We first take \(F\) in the Schwartz tensor class of Section~\ref{pre:tracial}. Choose an inhomogeneous radial partition
\(1=\chi_0+\sum_{j\ge1}\varphi(2^{-j}\cdot)\) of \(\R^2\) with real-valued smooth cutoffs \(\chi_0,\varphi\), \(\chi_0\) supported in a ball and \(\varphi\) in a compact annulus in \(\R^2\setminus\{0\}\). For each integer \(j\ge1\) write \(\Pi_j=\varphi(2^{-j}D)\)
in this subsection, and \(\widetilde\Pi_j=\widetilde\varphi(2^{-j}D)\)
where \(\widetilde\varphi\) is smooth, radial, supported in a compact annulus and equal to one on a neighborhood of \(\operatorname{supp}\varphi\). Both spatial multipliers act on \(L^p(\N_{\R^2})\), and on spacetime functions at each time.
For a finite set \(E\subset\{1,2,\ldots\}\), the triangle inequality gives
\[
 \begin{aligned}
 \left\|U(t)\left(\chi_0(D)F+\sum_{j\in E}\Pi_jF\right)
                \right\|_{L^p(\N_{I\times\R^2})}
 &\le \|U(t)\chi_0(D)F\|_{L^p(\N_{I\times\R^2})}\\
 &\quad+\sum_{j\in E}\|U(t)\Pi_jF\|_{L^p(\N_{I\times\R^2})}.
 \end{aligned}
\]
We bound the annular sum uniformly in \(E\) and then treat the term
with multiplier \(\chi_0(D)\).

Let \(G=J^sF\in L^p(\N_{\R^2})\). In the variable \(\eta=2^{-j}\xi\), the symbol of \(2^{js}\Pi_jJ^{-s}\) is
\begin{equation*}
 \varphi(\eta)(2^{-2j}+|\eta|^2)^{-s/2}\,.
\end{equation*}
Its derivatives through order four are uniformly bounded on the fixed support of \(\varphi\), which lies away from the origin. Applying~\eqref{intro:annular} to \(U(t)\widetilde\Pi_j\) in the fixed-annulus form of Section~\ref{lp:annular}
and Lemma~\ref{pre:scalar-kernels} to \(2^{js}\Pi_jJ^{-s}\) gives
\begin{equation}\label{lp:sob-annular-piece}
 \|U(t)\Pi_jF\|_{L^p(\N_{I\times\R^2})}
 \le C_{p,s}2^{-j(s-\sigma(p))}(1+j)^{\beta(p)}
          \|G\|_{L^p(\N_{\R^2})}.
\end{equation}
The coefficients are summable for \(s>\sigma(p)\). For the low-frequency term put \(b(\xi)=\chi_0(\xi)\langle\xi\rangle^{-s}\) for \(\xi\in\R^2\). The symbol of \(U(t)\chi_0(D)J^{-s}\) is \(b(\xi)+b(\xi)(e^{it|\xi|}-1)\). The kernel \(\check b\) is integrable. To treat the residual, choose a smooth \(\rho:\R^2\to\R\) supported in a compact annulus such that \(\sum_{k\in\mathbb Z}\rho(2^{-k}\xi)=1\) for \(\xi\ne0\). For \(k\in\mathbb Z\) let
\begin{equation*}
 b_{k,t}(\xi)=\rho(2^{-k}\xi)b(\xi)(e^{it|\xi|}-1)\,.
\end{equation*}
In the variable \(\eta=2^{-k}\xi\), the factor \(e^{it2^k|\eta|}-1\) and its derivatives through order four are \(O(2^k)\) on the support of \(\rho\) for \(k\le0\), uniformly in \(t\in I\). The integration-by-parts estimate used in~\eqref{pre:kern} therefore gives
\begin{equation*}
 |b_{k,t}^{\vee}(x)|\le C 2^{3k}(1+2^k|x|)^{-4}\qquad(k\le0,\ x\in\R^2)\,.
\end{equation*}
The \(L^1\) norms of these pointwise bounds sum to a finite constant. For \(k>0\), only finitely many \(b_{k,t}\) can be nonzero because \(b\) is compactly supported. For each such \(k\), the derivatives of \(b_{k,t}\) through order four are bounded uniformly for \(t\in I\) on a fixed compact annulus. Applying \((1-\Delta_\xi)^2\) in Fourier inversion gives \(|b_{k,t}^\vee(x)|\le C_k(1+|x|^2)^{-2}\). Hence the sum of the kernels \(b_{k,t}^{\vee}\) converges in \(L^1\), uniformly in \(t\in I\). Its symbol equals \(b(\xi)(e^{it|\xi|}-1)\) for \(\xi\ne0\) by the partition, and both sides vanish at \(\xi=0\).

For finite \(p\), the integrable kernel of \(U(t)\chi_0(D)J^{-s}\) controls \(U(t)\chi_0(D)F\) uniformly for \(t\in I\). Summing~\eqref{lp:sob-annular-piece} in \(j\) gives convergence of \(U(t)\chi_0(D)F+\sum_{j=1}^m U(t)\Pi_jF\) in spacetime \(L^p\) as \(m\to\infty\). For each scalar Schwartz factor \(f\) of a test tensor, the Fourier symbol \(\chi_0(\xi)+\sum_{j=1}^m\varphi(2^{-j}\xi)\) tends to one and is uniformly bounded by finite overlap. Since \(\widehat f\in L^1(\R^2)\), dominated convergence in Fourier inversion identifies the limit with \(U(t)F\) on test tensors. Completion in the norm \(\|J^sF\|_{L^p(\N_{\R^2})}\) now gives~\eqref{intro:sob}. At \(p=2\) and \(s=0\), spatial Plancherel and \(|I|=1\) give equality directly.

For \(p=\infty\), we construct the kernel \(K_t^{(s)}\) of \(U(t)J^{-s}\) by summing its localized kernels. On the \(j\)-th high-frequency annulus, the fixed-annulus extension of~\eqref{lp:infty} bounds the kernel norm for \(U(t)\widetilde\Pi_j\) by \(C2^{j/2}\). Applying~\eqref{pre:kern} to \(\varphi(\eta)(2^{-2j}+|\eta|^2)^{-s/2}\) bounds the \(L^1\) norm of the kernel of \(\Pi_jJ^{-s}\) by \(C_s2^{-js}\). Their convolution has norm at most \(C_s2^{-j(s-1/2)}\). For \(s>1/2\), summing \(C_s2^{-j(s-1/2)}\) over \(j\ge1\) and adding the kernel of \(U(t)\chi_0(D)J^{-s}\) gives uniform convergence of the kernel series for \(K_t^{(s)}\) in \(L^1(\R^2)\) for \(t\in I\). Each localized kernel is \(L^1\)-continuous in \(t\) by~\eqref{pre:kern}, so \(t\mapsto K_t^{(s)}\) is \(L^1\)-continuous. The adjoint of \(H\mapsto\int_I\widetilde K_t^{(s)}*H(t)\,dt\) defines the spacetime operator on \(G\in\N_{\R^2}\), with norm at most \(\sup_{t\in I}\|K_t^{(s)}\|_{L^1(\R^2)}\) by Lemma~\ref{pre:scalar-kernels}\ref{pre:convolution}. The partition of unity identifies the symbol of \(K_t^{(s)}\) as \(e^{it|\xi|}\langle\xi\rangle^{-s}\). Thus the adjoint acts as \(U(t)F\) on test tensors for \(F=J^{-s}G\), and gives~\eqref{intro:sob} for \(p=\infty\).

For the critical estimate \eqref{intro:weighted}, fix \(2<p<\infty\) and \(B,C\in\R\) as in \eqref{intro:logrange}. The multiplier \(D_w=J^{\sigma(p)}L_{B,C}(J)\) maps the logarithmically weighted Sobolev space in \eqref{intro:weighted} into \(L^p(\N_{\R^2})\). Write \(G=D_w F\). For each integer \(j\ge1\), define \(v_j\in L^p(\N_{I\times\R^2})\) by
\(v_j(t)=U(t)\Pi_jF\). To estimate the sum of the \(v_j\) in terms of \(G\), define the spatial multipliers \(A_j:L^p(\N_{\R^2})\to L^p(\N_{\R^2})\) and the real weights \(w_j>0\) by
\[
 A_j=2^{j\sigma(p)}(1+j)^B[\log(2+j)]^C \Pi_jD_w^{-1},\qquad
 w_j=(1+j)^{\beta(p)-B}[\log(2+j)]^{-C}\qquad (j\ge1).
\]
Since \(\widetilde\Pi_jA_j=A_j\), the annular estimate \eqref{intro:annular} for \(U(t)\widetilde\Pi_j\), applied to \(A_jG\), gives
\[
 \|v_j\|_{L^p(\N_{I\times\R^2})}\le C_p w_j\|A_jG\|_{L^p(\N_{\R^2})}.
\]
For \eqref{intro:weighted} we need to sum \(v_j\) at the critical
Sobolev order. We shall bound the \(\ell^p\) norm of
\((\|A_jG\|_{L^p(\N_{\R^2})})_{j\ge1}\) and estimate the spacetime norm of
\(\sum_jv_j\) by the \(\ell^2\) norm of \((\|v_j\|_{L^p(\N_{I\times\R^2})})_{j\ge1}\). Thus H\"older's inequality will require
\((w_j)_{j\ge1}\in\ell^q\) with \(1/q=1/2-1/p=s_p\).
For the spatial multipliers \(A_j\) we shall prove
\begin{equation}\label{lp:analysis}
 \Big(\sum_{j\ge1}\|A_jG\|_{L^p(\N_{\R^2})}^p\Big)^{1/p}
        \le C_{p,B,C}\|G\|_{L^p(\N_{\R^2})}
\end{equation}
and for the spacetime sum we shall prove
\begin{equation}\label{lp:synthesis}
 \Big\|\sum_{j\in E}v_j\Big\|_{L^p(\N_{I\times\R^2})}
 \le C_p\Big(\sum_{j\in E}\|v_j\|_{L^p(\N_{I\times\R^2})}^2\Big)^{1/2}.
\end{equation}
In \eqref{lp:synthesis}, \(E\) may be any finite subset of \(\{1,2,\ldots\}\) and each \(v_j\) may be any spacetime \(L^p\) function with spatial Fourier support in \(2^j\operatorname{supp}\varphi\). We need the constants in~\eqref{lp:analysis} and~\eqref{lp:synthesis} to be uniform in \(E\) and the coefficient algebra.

\medskip
To prove~\eqref{lp:analysis}, we first bound the $L^1$ norms of the kernels of the spatial multipliers $A_j$ uniformly in $j$. For \(\eta=2^{-j}\xi\in\R^2\) write \(r_j=\langle 2^j\eta\rangle\). The symbol of \(A_j\) there is \(\varphi(\eta)\) times
\[
 (2^j/r_j)^{\sigma(p)}
 \big((1+j)/\ell_1(r_j)\big)^B
 \big(\log(2+j)/\ell_2(r_j)\big)^C .
\]
On the fixed annulus \(r_j/2^j=(2^{-2j}+|\eta|^2)^{1/2}\) ranges in a compact subinterval of \((0,\infty)\) and has bounded derivatives through order four. Writing \(\ell_1(r_j)=j\log 2+\log(2^{1-j}+r_j/2^j)\) shows that \(\ell_1(r_j)\asymp 1+j\) there and its positive-order derivatives through order four are bounded. Differentiation of \(\ell_2(r_j)=\log(2+\ell_1(r_j))\) gives the same derivative bounds and \(\ell_2(r_j)\asymp\log(2+j)\). The normalized factors and their reciprocals consequently have uniform derivatives through order four. By \eqref{pre:kern}, each \(A_j\) has kernel of \(L^1\) norm at most \(C_{p,B,C}\).

Keep the target exponent \(p\) fixed in the symbols of \(A_j\). For a finite set \(E\subset\{1,2,\ldots\}\) and \(2\le r\le\infty\), consider the same map
\begin{equation*}
 A_E:L^r(\N_{\R^2})\to L^r(\ell^\infty(E)\bar\otimes\N_{\R^2}),
 \qquad A_EG=(A_jG)_{j\in E}\,,
\end{equation*}
where the trace on \(\ell^\infty(E)\) is counting measure. For finite \(r\) we have
\[
 \|A_EG\|_{L^r(\ell^\infty(E)\bar\otimes\N_{\R^2})}=\Big(\sum_{j\in E}\|A_jG\|_{L^r(\N_{\R^2})}^r\Big)^{1/r}.
\]
For \(r=\infty\) the norm of \(A_EG\) equals \(\sup_{j\in E}\|A_jG\|_{L^\infty(\N_{\R^2})}\), bounded by \(C_{p,B,C}\|G\|_{L^\infty(\N_{\R^2})}\) using the individual kernel norms. At \(r=2\), let \(a_j:\R^2\to\mathbb C\) denote the symbol of \(A_j\). Plancherel's identity bounds \(\sum_{j\in E}\|A_jG\|_{L^2(\N_{\R^2})}^2\) by
\[
 \sup_\xi\sum_{j\in E}|a_j(\xi)|^2\
               \|G\|_{L^2(\N_{\R^2})}^2.
\]
Only a bounded number of annuli contain a given \(\xi\). Interpolating \(A_E\) between \(r=2\) and \(r=\infty\) gives the bound at \(r=p\), which is \eqref{lp:analysis} with the sum restricted to \(E\). Taking \(E=\{1,\ldots,m\}\) and letting \(m\to\infty\) gives \eqref{lp:analysis} by monotone convergence of the scalar sums.

\medskip
For~\eqref{lp:synthesis}, we use random signs to estimate the spacetime sum. Fix a finite \(E\subset\{1,2,\ldots\}\). For \(\epsilon=(\epsilon_j)_{j\in E}\in\{-1,1\}^E\) let \(S_\epsilon=\sum_{j\in E}\epsilon_j\widetilde\Pi_j\). We will recover \(\sum_{j\in E}v_j\) by applying \(S_\epsilon\)
to \(\sum_{j\in E}\epsilon_jv_j\) and averaging in \(\epsilon\).
For this reduction we first prove an \(L^p\) multiplier bound for
\(S_\epsilon\) uniform in \(E\) and \(\epsilon\). Its kernel \(K_\epsilon:\R^2\to\mathbb C\) is given by
\(K_\epsilon(x)=\sum_{j\in E}\epsilon_j2^{2j}\check{\widetilde\varphi}(2^j x)\).
For \(m=0,1\), its derivatives of order \(m\) obey
\[
 |\nabla^m K_\epsilon(x)|
 \le C_m\sum_{j\in E}2^{(2+m)j}(1+2^j|x|)^{-4}
 \le C'_m |x|^{-(2+m)}\qquad(x\ne0).
\]
Indeed the terms with \(2^j|x|\le1\) form an increasing geometric sum, while the terms with \(2^j|x|>1\) are bounded by \(C|x|^{-4}2^{(m-2)j}\). The gradient estimate gives the Lipschitz kernel difference bound in either variable of \(K_\epsilon(x-y)\) at displacements at most \(|x-y|/2\). By Hilbert-valued Plancherel on \(L^2(\R^2;L^2(\M))\), the \(L^2\) norm of \(S_\epsilon\) is bounded by \(\sup_\xi\sum_{j\in E}|\widetilde\varphi(2^{-j}\xi)|\), uniformly by finite overlap. Parcet's scalar-kernel Calder\'on--Zygmund theorem \cite[Theorem A]{parcet} gives a bound on \(L^p(L^\infty(\R^2)\bar\otimes\M)\) depending only on \(p\), the dimension, the kernel bounds and the \(L^2\) bound. All these bounds are uniform in \(E,\epsilon,\M\). Integration in \(t\) gives the same operator norm on \(L^p(\N_{I\times\R^2})\).

Choose the signs \(\epsilon_j\), \(j\in E\), independently, each taking the values \(-1\) and \(1\) with probability \(1/2\). Lemma~\ref{pre:localized-approximation}, with coefficient algebra \(\N_I\), gives \(\widetilde\Pi_j v_j=v_j\) by the spatial Fourier support of \(v_j\). Since \(\mathbb E(\epsilon_j\epsilon_k)=1_{j=k}\), we obtain
\[
 \mathbb E_\epsilon S_\epsilon\Big(\sum_{k\in E}\epsilon_kv_k\Big)
 =\sum_{j,k\in E}\mathbb E_\epsilon(\epsilon_j\epsilon_k)\widetilde\Pi_jv_k
 =\sum_{j\in E}v_j .
\]
Using the triangle inequality for the expectation, the multiplier
bound just proved and H\"older's inequality in the sign variables, we obtain
\[
 \Big\|\sum_{j\in E}v_j\Big\|_{L^p(\N_{I\times\R^2})}
  \le C_p\mathbb E_\epsilon\Big\|\sum_{j\in E}\epsilon_jv_j\Big\|_{L^p(\N_{I\times\R^2})}
 \le C_p\Big(\mathbb E_\epsilon
                 \Big\|\sum_{j\in E}\epsilon_jv_j\Big\|_{L^p(\N_{I\times\R^2})}^p\Big)^{1/p}.
\]
The noncommutative Khintchine inequality and the comparison of
Rademacher moments \cite[Theorem 6.1 and (6.2)]{px} bound
\((\mathbb E_\epsilon\|\sum_{j\in E}\epsilon_jv_j\|_{L^p(\N_{I\times\R^2})}^p)^{1/p}\) by
\[
 C_p\max\left\{
 \Big\|\big(\sum_{j\in E}v_j^*v_j\big)^{1/2}\Big\|_{L^p(\N_{I\times\R^2})},\
 \Big\|\big(\sum_{j\in E}v_jv_j^*\big)^{1/2}\Big\|_{L^p(\N_{I\times\R^2})}\right\}
 .
\]
For the column square function its squared norm is
\[
 \Big\|\sum_{j\in E}v_j^*v_j\Big\|_{L^{p/2}(\N_{I\times\R^2})}
 \le\sum_{j\in E}\|v_j^*v_j\|_{L^{p/2}(\N_{I\times\R^2})}
 =\sum_{j\in E}\|v_j\|_{L^p(\N_{I\times\R^2})}^2 .
\]
For the row square function apply the triangle inequality in \(L^{p/2}(\N_{I\times\R^2})\) to \(v_jv_j^*\), using
\[
 \|v_jv_j^*\|_{L^{p/2}(\N_{I\times\R^2})}
 =\|v_j^*\|_{L^p(\N_{I\times\R^2})}^2
 =\|v_j\|_{L^p(\N_{I\times\R^2})}^2.
\]
Taking square roots in the row and column estimates proves \eqref{lp:synthesis}.

\medskip
It remains to combine~\eqref{lp:analysis} and~\eqref{lp:synthesis} at $s=\sigma(p)$. Let \(q=1/s_p=2p/(p-2)\in(2,\infty)\), so \(1/2=1/p+1/q\). Combining \eqref{lp:synthesis} with the bound for each \(v_j\), then applying H\"older's inequality to the scalar sequences \(w_j^2\) and \(\|A_jG\|_{L^p(\N_{\R^2})}^2\) with exponents \(q/2,p/2\), gives
\[
 \Big\|\sum_{j\in E}U(t)\Pi_jF\Big\|_{L^p(\N_{I\times\R^2})}
 \lesssim \Big(\sum_{j\in E}w_j^q\Big)^{1/q}
               \Big(\sum_{j\in E}\|A_jG\|_{L^p(\N_{\R^2})}^p\Big)^{1/p}.
\]
By~\eqref{lp:analysis} the \(\ell^p\) norm of \((\|A_jG\|_{L^p(\N_{\R^2})})_{j\in E}\) is bounded uniformly in \(E\). The \(\ell^q\) norm of \((w_j)_{j\in E}\) is uniformly bounded if
\(\sum_{j\ge1}(1+j)^{-q(B-\beta(p))}[\log(2+j)]^{-qC}\) converges.
The integral test gives convergence when \(q(B-\beta(p))>1\), or when \(q(B-\beta(p))=1\) and \(qC>1\). Since \(1/q=s_p\), these are precisely the two alternatives in \eqref{intro:logrange}. If \(E\subset\{m,m+1,\ldots\}\), the \(\ell^q\) factor is at most
\((\sum_{j\ge m}w_j^q)^{1/q}\), which tends to zero as \(m\to\infty\).
By~\eqref{lp:analysis}, \(\big(\sum_{j\in E}\|A_jG\|_{L^p(\N_{\R^2})}^p\big)^{1/p}\) is at most
\(C_{p,B,C}\|G\|_{L^p(\N_{\R^2})}\), independently of \(E\). Hence the partial
sums \(\sum_{j=1}^mU(t)\Pi_jF\) are Cauchy in spacetime \(L^p\). For the low-frequency term the proof of~\eqref{intro:sob} applies with \(b(\xi)=\chi_0(\xi)\langle\xi\rangle^{-\sigma(p)}L_{B,C}(\langle\xi\rangle)^{-1}\), which is smooth with compact support. On Schwartz test tensors the partition of unity and dominated convergence in Fourier inversion identify \(U(t)\chi_0(D)F+\sum_{j\ge1}U(t)\Pi_jF\) with \(U(t)F\). Completion in \(\|J^{\sigma(p)}L_{B,C}(J)F\|_{L^p(\N_{\R^2})}\) now proves \eqref{intro:weighted}. In particular \(p=4\) allows \(B=3/4,\ C>1/4\).

\subsection{Scalar obstructions to unweighted Sobolev estimates}
Section~\ref{lp:sum-annuli} proves the unweighted estimate~\eqref{intro:sob} by summing over frequency annuli. To see how far the Sobolev order \(s\) can be lowered, we test the estimate with scalar functions. Any bound valid for every coefficient algebra must in particular hold for \(\M=\mathbb C\).

\begin{proposition}\label{lp:sob-necessary}
Fix \(2\le p\le\infty\) and \(s\in\R\). Suppose the unweighted estimate~\eqref{intro:sob} of Theorem~\ref{thm:lp}\ref{thm:lp-sob} holds, with a constant independent of \(F\), for every scalar \(F\in\mathscr S(\R^2)\). Then
\begin{equation*}
 \begin{array}{ll}
 s\ge0             & \text{if }p=2,\\
 s>0               & \text{if }2<p\le4,\\
 s\ge 1/2-2/p      & \text{if }4<p<\infty,\\
 s>1/2             & \text{if }p=\infty.
 \end{array}
\end{equation*}
\end{proposition}

Together with Theorem~\ref{thm:lp}\ref{thm:lp-sob}, these conditions give the exact range of \(s\) for \(2\le p\le4\) and for \(p=\infty\). For \(4<p<\infty\), the unweighted estimate at \(s=1/2-2/p\) remains open for general coefficient algebras.

\begin{proof}
We use \(U(-t)\) because complex conjugation interchanges \(U(t)\) and \(U(-t)\) and preserves the norms and radial frequency supports. Fix a smooth annular cutoff \(\chi\) equal to one for \(5/4\le|\eta|\le7/4\). In the variable \(\eta=\xi/N\), the derivatives of \(\chi(\eta)(N^{-2}+|\eta|^2)^{s/2}\) through order four are bounded uniformly for \(N\ge2\). Lemma~\ref{pre:scalar-kernels} therefore gives \(\|J^sF_N\|_p\lesssim_s N^s\|F_N\|_p\) whenever \(F_N\in\mathscr S(\R^2)\) has Fourier support in \(5N/4<|\xi|<7N/4\).

At \(p=2\), Plancherel's identity and \(|I|=1\) give \(\|U(-t)F_N\|_{L^2(I\times\R^2)}=\|F_N\|_2\). Taking a nonzero \(F_N\) on each annulus in~\eqref{intro:sob} yields \(1\lesssim_s N^s\), so \(s\ge0\). For \(2<p\le4\), \cite[Proposition 8.2]{openai} gives, for arbitrarily large \(N\ge2\), nonzero Schwartz functions \(F_N\in\mathscr S(\R^2)\) with Fourier support in \(5N/4<|\xi|<7N/4\) and
\begin{equation*}
 \|U(-t)F_N\|_{L^{p,\infty}(I\times\R^2)}
 \ge c_p(\log N)^{1/2-1/p}\|F_N\|_{L^p(\R^2)}\,.
\end{equation*}
The strong spacetime \(L^p\) norm dominates the weak norm, so~\eqref{intro:sob} gives \((\log N)^{1/2-1/p}\lesssim_{p,s}N^s\). Since \(1/2-1/p>0\), this requires \(s>0\).

For \(4<p\le\infty\), we use the focusing example of~\cite[Proposition 8.1]{openai} to show that \eqref{intro:sob} requires \(s\ge\sigma(p)\). Choose \(t_0=3/2\) and a nonzero radial function \(\phi\in C_c^\infty(\R^2)\) with \(\phi\ge0\) supported in \(5/4<|\xi|<7/4\). For real \(N\ge2\) define \(F_N\in\mathscr S(\R^2)\) by \(\widehat F_N(\xi)=e^{it_0|\xi|}\phi(\xi/N)\). Fourier inversion with \(\xi=N\eta\) gives
\begin{equation*}
 U(-t)F_N(x)=\frac{N^2}{(2\pi)^2}\int_{\R^2}e^{iN[x\cdot\eta+(t_0-t)|\eta|]}\phi(\eta)\,d\eta\,.
\end{equation*}
On \(E_N=\{(t,x)\in I\times\R^2:|x|+|t-t_0|<c/N\}\), the real part of the exponential is at least \(1/2\) for all \(\eta\in\operatorname{supp}\phi\) when \(c>0\) is sufficiently small. Since \(\phi\ge0\) and \(\int\phi>0\), we have \(|U(-t)F_N(x)|\gtrsim N^2\) on \(E_N\), which has measure comparable to \(N^{-3}\).

To estimate \(\|F_N\|_{L^p(\R^2)}\), we first use Plancherel's identity to obtain \(\|F_N\|_{L^2(\R^2)}=O(N)\). In polar coordinates \(\xi=N\rho\omega\), with \(\rho>0\) and \(\omega\in\mathbb S^1\), the radial phase in the inverse Fourier formula for \(F_N\) is \(N\rho(t_0+x\cdot\omega)\). For \(|x|\le t_0/2\), its \(\rho\)-derivative is at least \(Nt_0/2\). Integrating by parts \(M\) times in \(\rho\) against the compactly supported amplitude \(\rho\phi(\rho\omega)\) gives \(F_N(x)=O_M(N^{2-M})\) for every integer \(M\ge1\). For \(|x|>t_0/2\), the angular integral is \(2\pi J_0(N\rho|x|)\), where \(J_0\) is the Bessel function of the first kind of order zero. It is \(O(N^{-1/2})\) on the support of \(\phi\) by~\cite[(10.17.3)]{dlmf}. The prefactor \(N^2\) then gives \(\|F_N\|_{L^\infty(\R^2)}=O(N^{3/2})\). Interpolating the \(L^2\) and \(L^\infty\) norms yields \(\|F_N\|_{L^p(\R^2)}\lesssim N^{3/2-1/p}\) for \(2\le p<\infty\).

For \(4<p<\infty\), the lower bound on \(E_N\) and the annular bound \(\|J^sF_N\|_p\lesssim_s N^s\|F_N\|_p\) in~\eqref{intro:sob} give \(N^{2-3/p}\lesssim N^{s+3/2-1/p}\), hence \(s\ge1/2-2/p\). For \(p=\infty\), the essential supremum gives \(N^2\lesssim N^{s+3/2}\), hence \(s\ge1/2\).

To rule out equality at \(p=\infty\), define \(K\in\mathscr S'(\R^2)\) by \(\widehat K(\xi)=e^{-it_0|\xi|}\langle\xi\rangle^{-1/2}\). If \eqref{intro:sob} held for \(s=1/2\), then taking \(F=J^{-1/2}G\) with \(G\in\mathscr S(\R^2)\) would give
\begin{equation*}
 |\langle K,G(-\cdot)\rangle|
 =|(U(-t_0)J^{-1/2}G)(0)|\le C\|G\|_{L^\infty(\R^2)}\,.
\end{equation*}
Fourier inversion gives continuity of \((t,x)\mapsto U(-t)J^{-1/2}G(x)\), so the essential supremum controls evaluation at \((t_0,0)\). By density of \(\mathscr S(\R^2)\) in \(C_0(\R^2)\) and the representation theorem for its dual, \(K\) would be a finite measure. We show that \(K\) has a nonzero principal value term across the circle \(|x|=t_0\).

For \(\epsilon>0\) let \(K_\epsilon\) be the inverse Fourier transform of \(e^{-\epsilon|\xi|}\widehat K(\xi)\). Then \(K_\epsilon\to K\) in \(\mathscr S'(\R^2)\) by dominated convergence. Writing \(r=|x|\), we have
\begin{equation*}
 K_\epsilon(r)=\frac{1}{2\pi}\int_0^\infty e^{-(\epsilon+it_0)\rho}\langle\rho\rangle^{-1/2}J_0(r\rho)\,\rho\,d\rho\,.
\end{equation*}
The integral over \(0\le\rho\le1\) converges locally uniformly in \(r\) as \(\epsilon\downarrow0\). For \(\rho\ge1\), the large-argument expansion \cite[(10.17.3) and Section 10.17(iii)]{dlmf} gives
\begin{equation*}
 \langle\rho\rangle^{-1/2}\rho J_0(r\rho)
 =\sum_{\pm}e^{\pm ir\rho}\left(c_\pm(r)+\frac{d_\pm(r)}{\rho}\right)+O(\rho^{-2})\,,
\end{equation*}
with smooth \(c_\pm,d_\pm\), and \(c_+(r)=e^{-i\pi/4}/\sqrt{2\pi r}\). The remainder is uniform for \(r\) near \(t_0\) and contributes a locally bounded function in the limit.

The terms with phase \(e^{-i(r+t_0)\rho}\) have smooth Abel limits near \(r=t_0\). For the constant amplitude this follows from the formula \(e^{-i(r+t_0)}/[i(r+t_0)]\) for the limit of the integral over \([1,\infty)\). For the amplitude \(1/\rho\), changing variables to \(v=(r+t_0)\rho\) gives the limit \(\int_{r+t_0}^\infty e^{-iv}\,dv/v\), whose lower endpoint stays positive. For the term with phase \(e^{i(r-t_0)\rho}\) and amplitude \(1/\rho\), set \(u=r-t_0\). If \(0<|u|<1\), split the integral at \(A=|u|^{-1}\). The integral up to \(A\) has modulus at most \(\log A\). One integration by parts bounds the tail, uniformly in \(\epsilon>0\), by
\begin{equation*}
 \left|\int_A^\infty e^{iu\rho}\frac{e^{-\epsilon\rho}}{\rho}\,d\rho\right|
 \le \frac{2e^{-\epsilon A}}{|u|A}\le 2\,.
\end{equation*}
Here \(e^{-\epsilon\rho}/\rho\) decreases to zero. Dominated convergence with the locally integrable bound \(2+|\log|u||\) shows that this Abel limit is locally integrable in \(r\).

For the constant amplitude \(c_+(r)\) with phase \(e^{iu\rho}\), we use the Abel limit
\begin{equation*}
 \int_0^\infty e^{iu\rho-\epsilon\rho}\,d\rho=(\epsilon-iu)^{-1}
 \longrightarrow \pi\delta_0(u)+i\,\operatorname{pv}(1/u)
 \qquad\text{in }\mathscr S'(\R)\,.
\end{equation*}
Removing \(0<\rho<1\) changes the limit by a smooth function. Combining the Abel limits thus gives, in the radial variable on an annulus about \(r=t_0\),
\begin{equation*}
 K(r)=\frac{i c_+(r)}{2\pi}\operatorname{pv}\frac{1}{r-t_0}
     +\frac{c_+(t_0)}{2}\delta_{t_0}(r)+h(r)\,,
 \qquad h\in L^1_{\mathrm{loc}}\,.
\end{equation*}
The polar Jacobian \(r\) is smooth and bounded away from zero on this annulus. Away from the circle the principal value term is the function \(i c_+(r)/(2\pi(r-t_0))\), whose absolute integral diverges near \(t_0\) since \(c_+(t_0)\ne0\). The delta term and \(h\) define finite measures on any smaller annulus. Thus \(K\) cannot be a finite measure, which rules out \(s=1/2\) at \(p=\infty\).
\end{proof}

\subsection{Quantum Fourier models}\label{q:models}

We describe translations in Fourier coordinates on the quantum torus and the generalized Moyal plane and verify that their orbits are norm continuous in \(L^p\) for \(1\le p<\infty\). Section~\ref{q:secspacetime} uses these actions to transfer the Euclidean spacetime estimates to both models in the proof of Theorem~\ref{thm:quantum}.

The frequency \(\xi\) ranges over \(\mathbb Z^2\) on
the torus and \(\R^2\) on the plane. In either case translation by
\(z\in\R^2\) multiplies the coefficient at \(\xi\) by \(e^{iz\cdot\xi}\). Thus integrating the orbit at \(-y\) against a scalar kernel in \(y\) with Fourier transform \(m\) realizes the quantum multiplier with symbol \(m\).

Fix a real skew-symmetric \(2\times2\) matrix \(\Theta\). We use two
von Neumann algebras \(\mathcal Q_\Theta\), each equipped with a trace and an
action \(\gamma\) of \(\R^2\) by trace-preserving automorphisms \(\gamma_z:\mathcal Q_\Theta\to\mathcal Q_\Theta\) for \(z\in\R^2\).

On the quantum torus the two unitary generators satisfy
\(U_jU_k=e^{2\pi i\Theta_{jk}}U_kU_j\) for \(j,k\in\{1,2\}\).
For \(n=(n_1,n_2)\in\mathbb Z^2\) the ordered monomial is
\(U^n=U_1^{n_1}U_2^{n_2}\), with the first generator always on the left,
\(U_j^0=1\) and \(U_j^{-k}=(U_j^*)^k\) for positive integers \(k\). For a symbol \(m:\mathbb Z^2\to\mathbb C\), the multiplier \(m(D_\Theta)\) acts on the span of the monomials with values in the same span. On the monomials these conventions give
\[
 \operatorname{tr}_\Theta(U^n)=1_{n=0},\qquad
 \gamma_z(U^n)=e^{iz\cdot n}U^n,\qquad m(D_\Theta)U^n=m(n)U^n .
\]
The weak closure of the polynomial algebra in the trace representation is the quantum torus. The gauge action factors through \(\R^2/(2\pi\mathbb Z)^2\) and preserves the
trace and hence each \(L^p\) norm. For \(1\le p<\infty\) the map
\(z\mapsto\gamma_z x\) is norm continuous for each \(x\in L^p\), by the
character formula on polynomials and density. The gauge action and the identification of its orbit in a semicommutative
\(L^p\) space are described by Chen, Xu and Yin
\cite[Proposition 2.1 and Corollary 2.2]{cxy}.

\medskip
On the generalized Moyal plane, \(\mathcal Q_\Theta\) is the twisted group von Neumann algebra of \(\R^2\) in its regular representation, with generating unitaries \(\lambda_\Theta(\xi)\in\mathcal Q_\Theta\). For \(\xi,\eta\in\R^2\) they satisfy
\[
 \lambda_\Theta(\xi)\lambda_\Theta(\eta)
 =e^{i\langle\eta,\Theta\xi\rangle/2}\lambda_\Theta(\xi+\eta).
\]
For \(f\in L^1(\R^2)\), the weak operator integral \(\lambda_\Theta(f)=\int_{\R^2} f(\xi)\lambda_\Theta(\xi)\,d\xi\) belongs to \(\mathcal Q_\Theta\). For \(f\in\mathscr S(\R^2)\), the Plancherel trace satisfies
\(\operatorname{tr}_\Theta(\lambda_\Theta(f))=f(0)\) and
\(\|\lambda_\Theta(f)\|_{L^2(\mathcal Q_\Theta)}=\|f\|_{L^2(\R^2)}\). We call \(\mathscr S(\mathcal Q_\Theta)=\lambda_\Theta(\mathscr S(\R^2))\) the quantum Schwartz class. For a measurable symbol \(m:\R^2\to\mathbb C\), the multiplier \(m(D_\Theta)\) on \(L^2(\mathcal Q_\Theta)\) has domain \(\{\lambda_\Theta(f):f,mf\in L^2(\R^2)\}\), using the \(L^2\) extension of \(\lambda_\Theta\) given by Plancherel's identity. In Fourier coordinates the dual action and the multiplier are
\[
 \gamma_z(\lambda_\Theta(f))=\lambda_\Theta(e^{iz\cdot\xi}f(\xi)),
 \qquad m(D_\Theta)\lambda_\Theta(f)=\lambda_\Theta(mf).
\]
The formula for \(\gamma_z\) preserves the twisted product and the trace,
and extends to a normal trace-preserving automorphism. For these models
see \cite[Section 2]{restriction} and \cite[Proposition 6.1]{hlw}.
For Schwartz \(f\), the coefficient \((e^{iz\cdot\xi}-1)f(\xi)\) tends to zero in both \(L^2_\xi\) and \(L^1_\xi\) as \(z\to0\). Plancherel's identity and the unitary integral bound \(\|\lambda_\Theta(f)\|_{L^\infty(\mathcal Q_\Theta)}\le\|f\|_{L^1(\R^2)}\) thus give orbit continuity in \(L^2\) and in operator norm on these elements. Interpolation gives continuity in \(L^p\) for \(2<p<\infty\).

\medskip
For either model put \(\mathcal D=\mathcal C\bar\otimes\mathcal Q_\Theta\)
with product trace, where \(\mathcal C\) is any semifinite tracial
coefficient algebra. For \(z\in\R^2\) write \(\beta_z=\mathrm{id}_{\mathcal C}\otimes\gamma_z:\mathcal D\to\mathcal D\).
Since \(\beta_z\) preserves the product trace, \(\|\beta_zx\|_{L^p(\mathcal D)}=\|x\|_{L^p(\mathcal D)}\) for every \(z\in\R^2\), \(1\le p\le\infty\), and \(x\in L^p(\mathcal D)\). For \(2\le p<\infty\), the quantum polynomials are dense in \(L^p(\mathcal Q_\Theta)\) on the torus, and the quantum Schwartz class is dense in \(L^p(\mathcal Q_\Theta)\) on the plane by~\cite[Section 2.2]{restriction}. Together with the density of finite sums of elementary tensors in Section~\ref{pre:tracial}, this extends orbit continuity to every \(x\in L^p(\mathcal D)\).
Continuity on the trace predual \(L^1(\mathcal D)\) will allow us
to define \eqref{q:orbit} at \(p=\infty\) by duality. Given \(d=u|d|\in L^1(\mathcal D)\), put \(a=u|d|^{1/2}\) and \(b=|d|^{1/2}\). Both lie in \(L^2(\mathcal D)\) and
\[
 \|\beta_zd-d\|_{L^1(\mathcal D)}\le
 \|\beta_z a-a\|_{L^2(\mathcal D)}\|b\|_{L^2(\mathcal D)}+\|a\|_{L^2(\mathcal D)}\|\beta_z b-b\|_{L^2(\mathcal D)}\longrightarrow0 .
\]
For bounded \(d\in\mathcal D\) whose support projections have finite trace, interpolation between the
\(L^1\) and \(L^2\) norms gives \(\|\beta_zd-d\|_{L^p(\mathcal D)}\to0\)
when \(1<p<2\). Since \(\beta_z\) preserves the \(L^p(\mathcal D)\) norm,
density of these elements extends continuity to all of \(L^p(\mathcal D)\). In particular, \(z\mapsto\beta_z x\) is norm continuous from \(\R^2\) to \(L^p(\mathcal D)\) for every \(x\in L^p(\mathcal D)\) and \(1\le p<\infty\).

Let \((\Lambda_\Theta,\mu_\Theta)\) be \(\mathbb Z^2\) with counting measure on the torus or \(\R^2\) with Lebesgue measure on the plane. The torus monomials form an orthonormal basis of \(L^2(\mathcal Q_\Theta)\), since \(\operatorname{tr}_\Theta((U^m)^*U^n)=1_{m=n}\) for \(m,n\in\mathbb Z^2\) and the polynomials are dense. On the plane, Plancherel's identity identifies \(\lambda_\Theta(f)\) with \(f\) in \(L^2(\R^2)\). Tensoring either identification with the identity on \(L^2(\mathcal C)\) gives a unitary map \(x\mapsto\widehat x_\Theta\) from \(L^2(\mathcal D)\) onto \(L^2(\Lambda_\Theta,\mu_\Theta;L^2(\mathcal C))\).

We write \(|D_\Theta|=\sqrt{-\Delta_\Theta}\) and \(J_\Theta=(1-\Delta_\Theta)^{1/2}\) for the multipliers of symbols \(|\xi|\) and \(\langle\xi\rangle\). With \(\operatorname{Dom}_2\) denoting the domain in \(L^2(\mathcal D)\),
\begin{align*}
 \operatorname{Dom}_2(|D_\Theta|)
 &=\left\{x\in L^2(\mathcal D):
   \int_{\Lambda_\Theta}|\xi|^2\|\widehat x_\Theta(\xi)\|_{L^2(\mathcal C)}^2\,d\mu_\Theta(\xi)<\infty\right\}\,,\\
 \operatorname{Dom}_2(J_\Theta)
 &=\left\{x\in L^2(\mathcal D):
   \int_{\Lambda_\Theta}(1+|\xi|^2)\|\widehat x_\Theta(\xi)\|_{L^2(\mathcal C)}^2\,d\mu_\Theta(\xi)<\infty\right\}\,.
\end{align*}

\subsection{Spacetime transference and quantum waves}\label{q:secspacetime}

We transfer the Euclidean spacetime estimate by restricting the orbit \(z\mapsto\beta_zx\) to a square. In comparing the Euclidean convolution with the operator on the orbit, the error from a translation by \(y\in\R^2\) is supported on the symmetric difference of the square and its translate by \(y\). For a semifinite tracial von Neumann algebra \(\mathcal D\), write \(\mathcal S_{\mathcal D}=L^\infty(\R^2)\bar\otimes\mathcal D\) and \(\mathcal D_I=L^\infty(I)\bar\otimes\mathcal D\), both with the product trace. The identification~\eqref{pre:fubini} gives \(L^p(\mathcal S_{\mathcal D})=L^p(\R^2;L^p(\mathcal D))\) and \(L^p(\mathcal D_I)=L^p(I;L^p(\mathcal D))\) isometrically for \(1\le p<\infty\).

\begin{proposition}\label{q:spacetime}
Let \(1\le p<\infty\) and let \(\beta\) be an action of \(\R^2\) by
normal trace-preserving \(*\)-automorphisms of \(\mathcal D\), with
norm-continuous orbits on \(L^p(\mathcal D)\). Suppose the jointly measurable kernel \((t,y)\mapsto K_t(y)\) from \(I\times\R^2\) to \(\mathbb C\), with \(I=[1,2]\), satisfies
\(|K_t(y)|\le k(y)\) for a common \(k\in L^1(\R^2)\), and that
\[
 \|K_t*F\|_{L^p(I;L^p(\mathcal S_{\mathcal D}))}
 \le C\|F\|_{L^p(\mathcal S_{\mathcal D})}\qquad
 (F\in L^p(\mathcal S_{\mathcal D})).
\]
For \(t\in I\), define \(Q_t:L^p(\mathcal D)\to L^p(\mathcal D)\) by
\begin{equation}\label{q:orbit}
 Q_tx=\int_{\R^2}K_t(y)\beta_{-y}x\,dy
\end{equation}
for \(x\in L^p(\mathcal D)\). The spacetime map \(x\mapsto(t\mapsto Q_tx)\) satisfies
\begin{equation}\label{q:transfer}
 \|Q_tx\|_{L^p(I;L^p(\mathcal D))}\le C\|x\|_{L^p(\mathcal D)} .
\end{equation}
\end{proposition}

\begin{proof}
We prove~\eqref{q:transfer} by applying the Euclidean estimate to a localized orbit of \(x\in L^p(\mathcal D)\). The continuous orbit \(y\mapsto\beta_{-y}x\) has separable closed linear span in \(L^p(\mathcal D)\). Since \(K_t(y)\) is jointly measurable, the map \((t,y)\mapsto K_t(y)\beta_{-y}x\) is strongly measurable. The bound \(\|K_t(y)\beta_{-y}x\|_{L^p(\mathcal D)}\le k(y)\|x\|_{L^p(\mathcal D)}\) justifies~\eqref{q:orbit} as a Bochner integral. Bochner Fubini gives strong measurability of \(t\mapsto Q_tx\), and \(\|Q_tx\|_{L^p(\mathcal D)}\le\|k\|_{L^1(\R^2)}\|x\|_{L^p(\mathcal D)}\). For \(R>0\), put \(B_R=[-R,R]^2\), \(w_R=|B_R|^{-1/p}1_{B_R}\), and define \(F_R:\R^2\to L^p(\mathcal D)\) by \(F_R(z)=w_R(z)\beta_zx\). Since \(\|w_R\|_{L^p(\R^2)}=1\) and \(\beta_z\) preserves the \(L^p(\mathcal D)\) norm,
\begin{equation*}
 \|F_R\|_{L^p(\mathcal S_{\mathcal D})}=\|x\|_{L^p(\mathcal D)},\qquad
 \|w_R(z)\beta_z(Q_tx)\|_{L^p(I\times\R^2;L^p(\mathcal D))}
 =\|Q_tx\|_{L^p(\mathcal D_I)}\,.
\end{equation*}
Indeed, \(\beta_z(Q_tx)=\int_{\R^2}K_t(y)\beta_{z-y}x\,dy\), which also gives joint strong measurability in \((t,z)\). Subtracting the two integrals gives
\begin{equation*}
 (K_t*F_R)(z)-w_R(z)\beta_z(Q_tx)
 =\int_{\R^2} K_t(y)[w_R(z-y)-w_R(z)]\beta_{z-y}x\,dy\,.
\end{equation*}
For each \(y\in\R^2\), trace preservation and \(|I|=1\) give
\begin{equation*}
 \bigl\|K_t(y)[w_R(z-y)-w_R(z)]\beta_{z-y}x\bigr\|_{L^p(I\times\R^2;L^p(\mathcal D))}
 \le k(y)\|x\|_{L^p(\mathcal D)}\|w_R(\cdot-y)-w_R\|_{L^p(\R^2)}\,.
\end{equation*}
Applying Minkowski's integral inequality and the hypothesis on \(K_t*F_R\) gives
\begin{equation*}
 \|Q_tx\|_{L^p(\mathcal D_I)}
 \le C\|x\|_{L^p(\mathcal D)}+\|x\|_{L^p(\mathcal D)}\int_{\R^2} k(y)
       \left(\frac{|(B_R+y)\triangle B_R|}{|B_R|}\right)^{1/p}dy\,,
\end{equation*}
where the area ratio is \(\|w_R(\cdot-y)-w_R\|_{L^p(\R^2)}^p\). For each fixed \(y\) this ratio tends to zero as \(R\to\infty\), and it is at most two. Dominated convergence proves~\eqref{q:transfer}.
\end{proof}

We now prove Theorem~\ref{thm:quantum}, beginning with the annular estimate in part~\ref{thm:q-annular}, using it to prove part~\ref{thm:q-wave}, and then establishing the critical estimates in part~\ref{thm:q-weight}.

For the annular estimate we apply Proposition~\ref{q:spacetime} to \(\beta_z=\mathrm{id}_{\mathcal C}\otimes\gamma_z\) on either quantum model. For \(1\le p<\infty\) and a scalar kernel \(K\in L^1(\R^2)\), denote the operator in~\eqref{q:orbit} with kernel \(K\) by \(A_K:L^p(\mathcal D)\to L^p(\mathcal D)\). For \(x\in L^2(\mathcal D)\), taking the quantum Fourier transform in the Bochner integral for \(A_Kx\) gives
\begin{equation*}
 \widehat{A_Kx}_\Theta(\xi)
 =\left(\int_{\R^2}K(y)e^{-iy\cdot\xi}\,dy\right)\widehat x_\Theta(\xi)
 =\widehat K(\xi)\widehat x_\Theta(\xi)\,.
\end{equation*}
The first equality follows from the action formula for \(\beta_{-y}\) and Fubini after pairing with an element of \(L^2(\Lambda_\Theta,\mu_\Theta;L^2(\mathcal C))\). If \(x\in L^2(\mathcal D)\cap L^p(\mathcal D)\), the orbit is continuous in the intersection with norm \(\|\cdot\|_{L^2(\mathcal D)}+\|\cdot\|_{L^p(\mathcal D)}\). The integral in that norm agrees with the integrals in \(L^2(\mathcal D)\) and \(L^p(\mathcal D)\).

Fix a real \(N\ge2\) and a cutoff \(\psi\in C_c^\infty((0,\infty))\). For \(t\in I\) let \(m_t:\R^2\to\mathbb C\) be \(m_t(\xi)=e^{it|\xi|}\psi(|\xi|/N)\). For fixed \(N\), the symbols \(m_t\) have common compact support away from zero, with all derivatives bounded uniformly in \(t\). Fourier inversion gives jointly continuous kernels \(K_t=\check m_t\), and integration by parts as in~\eqref{pre:kern} gives \(|K_t(y)|\le C_{N,\psi}(1+|y|^2)^{-2}\). The Fourier formula for \(A_{K_t}\) identifies \(Q_t\) with \(U_\Theta(t)\psi(|D_\Theta|/N)\) on \(L^2(\mathcal D)\cap L^p(\mathcal D)\). Applying~\eqref{q:transfer} with the Euclidean annular bound for \(U(t)\) proves Theorem~\ref{thm:quantum}\ref{thm:q-annular} for \(2\le p<\infty\) and \(\psi\) supported in \((1/2,2)\). For \(\psi\) supported in any fixed compact subinterval of \((0,\infty)\), use the fixed-annulus form of~\eqref{intro:annular} from Section~\ref{lp:annular} with \(m(\eta)=\psi(|\eta|)\). Substituting its bound into~\eqref{q:transfer} for \(K_t=\check m_t\) gives the factor \(N^{\sigma(p)}[\log(2+N)]^{\beta(p)}\) for \(Q_t\).

For \(p=\infty\) and \(\psi\in C_c^\infty((1/2,2))\), define \(Q_t\) by~\eqref{q:orbit} in the weak-star sense. For the pairing \(\langle d,x\rangle=\operatorname{tr}_{\mathcal D}(d^*x)\) with \(d\in L^1(\mathcal D)\), its preadjoint is \(d\mapsto\int_{\R^2}\overline{K_t(y)}\beta_y d\,dy\). This Bochner integral exists by \(L^1\) orbit continuity and has norm at most \(\|K_t\|_{L^1(\R^2)}\|d\|_{L^1(\mathcal D)}\). Testing against bounded \(d\) whose support projections have finite trace shows that the weak-star integral agrees with the Fourier multiplier on \(\mathcal D\cap L^2(\mathcal D)\). The bound~\eqref{lp:infty} therefore gives the annular estimate in operator norm, uniformly in \(t\). For fixed \(N\), applying~\eqref{pre:kern} to \(m_u-m_t\) as \(u\to t\) in \(I\) shows that \(t\mapsto K_t\) is continuous in \(L^1(\R^2)\). Hence \(t\mapsto Q_tx\) is continuous in operator norm and defines an element of \(\mathcal D_I\) by uniform approximation by step functions with values in \(\mathcal D\).

Since each \(\beta_y\), \(y\in\R^2\), preserves adjoints, the orbit formula~\eqref{q:orbit} gives \((A_Kx)^*=A_{\overline K}(x^*)\) for \(K\in L^1(\R^2)\). The equality \(\widehat{\overline K}(\xi)=\overline{\widehat K(-\xi)}\) then gives the estimates for \(U_\Theta(-t)\) on \(I\) by applying the bound for \(U_\Theta(t)\) to \(x^*\) with cutoff \(\overline\psi\). For the weak-star integral at \(p=\infty\), the same identity follows by testing against the trace predual and conjugating the kernel.

\medskip
We use the annular estimate to bound the solution of the quantum Cauchy problem. Fix \(2\le p<\infty\). Let \(x,y\) lie in the span of \(c\otimes U^n\), \(n\in\mathbb Z^2\), on the torus, or of \(c\otimes\lambda_\Theta(f)\), \(f\in\mathscr S(\R^2)\), on the plane. In both cases \(c\in\mathcal C\) is bounded and its support projections have finite trace. The Fourier multiplier formula for \(u_\Theta(t;x,y)\) preceding Theorem~\ref{thm:quantum} gives a solution of
\[
 u_{tt}-\Delta_\Theta u=0,\qquad u(0)=x,\qquad u_t(0)=y .
\]
Indeed the scalar Fourier equation is \(u_{tt}+|\xi|^2u=0\), with solution
multipliers \(\cos(t|\xi|)\) and \(\sin(t|\xi|)/|\xi|\). To prove~\eqref{intro:quantum} for \(s>\sigma(p)\), take the partition from Section~\ref{lp:sum-annuli}. We use \(\chi_0\) for the low-frequency terms and \(\varphi(2^{-j}\cdot)\), \(j\ge1\), for the annular terms. The low-frequency symbols acting on \(J_\Theta^s x\) and \(J_\Theta^{s-1}y\) are
\[
 \chi_0(\xi)\cos(t|\xi|)\langle\xi\rangle^{-s},\qquad
 \chi_0(\xi)\frac{\sin(t|\xi|)}{|\xi|}\langle\xi\rangle^{1-s}.
\]
The cosine and sine quotient are smooth functions of \(|\xi|^2\) near zero, and the sine quotient has value \(t\) there. Compact support and the derivative bound \eqref{pre:kern} give uniform kernel norms for \(t\in I\). Using \(\|\beta_z a\|_{L^p(\mathcal D)}=\|a\|_{L^p(\mathcal D)}\) in \eqref{q:orbit}, their
contributions to \eqref{intro:quantum} are bounded by
\(C_{p,s}\|J_\Theta^s x\|_{L^p(\mathcal D)}\) and \(C_{p,s}\|J_\Theta^{s-1}y\|_{L^p(\mathcal D)}\),
respectively.

On each higher annulus apply to \(J_\Theta^s x\) and \(J_\Theta^{s-1}y\), respectively, the multipliers with symbols
\begin{equation*}
 \varphi(2^{-j}\xi)\langle\xi\rangle^{-s},\qquad
 \varphi(2^{-j}\xi)|\xi|^{-1}\langle\xi\rangle^{1-s},
 \qquad j\ge1\,.
\end{equation*}
After multiplication by \(2^{js}\), their derivatives through order four in the scaled variable \(\eta=2^{-j}\xi\) are bounded uniformly in \(j\). By~\eqref{pre:kern} and~\eqref{q:orbit}, both multipliers have operator norm at most \(C_s2^{-js}\). Apply the quantum annular estimate with \(\widetilde\varphi(2^{-j}\cdot)\), where the fixed radial cutoff \(\widetilde\varphi\) is supported on a compact annulus away from zero and equals one on \(\operatorname{supp}\varphi\). Applying this estimate to \(U_\Theta(t)\) and \(U_\Theta(-t)\), using the triangle inequality for each \(j\), and summing gives
\begin{equation*}
 \begin{aligned}
 &\sum_{j\ge1}\|\varphi(2^{-j}D_\Theta)u_\Theta(t;x,y)\|_{L^p(\mathcal D_I)}\\
 &\quad\le C_{p,s}\big(\|J_\Theta^s x\|_{L^p(\mathcal D)}+\|J_\Theta^{s-1}y\|_{L^p(\mathcal D)}\big)
        \sum_{j\ge1}2^{-j(s-\sigma(p))}(1+j)^{\beta(p)}\,.
 \end{aligned}
\end{equation*}
Thus \(\sum_{j\ge1}\varphi(2^{-j}D_\Theta)u_\Theta(t;x,y)\) converges absolutely in \(L^p(\mathcal D_I)\) for \(s>\sigma(p)\). To identify its sum after adding \(\chi_0(D_\Theta)u_\Theta(t;x,y)\), recall that \(x,y\in L^2(\mathcal D)\) on the tensor span. Keeping the annuli \(1\le j\le M\) and the low-frequency term multiplies the Fourier coefficients of \(u_\Theta(t;x,y)\) by \(\chi_0+\sum_{j=1}^M\varphi(2^{-j}\cdot)\). The multipliers are uniformly bounded by finite overlap and converge pointwise to one. The cosine and sine quotient are bounded for \(t\in I\), so Plancherel and dominated convergence give convergence to \(u_\Theta(t;x,y)\) in \(L^2(\mathcal D_I)\). The \(L^2\) and \(L^p\) limits coincide, since both imply convergence in measure for the product trace. The low- and high-frequency bounds now prove~\eqref{intro:quantum} on the tensor span. The estimate extends the solution map to the completions in \(\|J_\Theta^s x\|_{L^p(\mathcal D)}\) and \(\|J_\Theta^{s-1}y\|_{L^p(\mathcal D)}\), respectively.

At the endpoint \(p=2\), \(s=0\), use \(|\cos(tr)|\le1\), \(\langle r\rangle|\sin(tr)|/r\le C_I\)
and Plancherel's identity for the quantum Fourier transform tensored with the identity on
\(L^2(\mathcal C)\). Since \(|I|=1\), the cosine contribution has spacetime norm at most \(\|x\|_{L^2(\mathcal D)}\), and the sine-quotient contribution has spacetime norm at most \(C_I\|J_\Theta^{-1}y\|_{L^2(\mathcal D)}\).

\medskip
It remains to prove the estimates at Sobolev order $\sigma(p)$. For \(2<p<\infty\) and \(B,C\) as in \eqref{intro:logrange}, put \(D_w=J^{\sigma(p)}L_{B,C}(J)\). The half-wave assertion of Theorem~\ref{thm:quantum}\ref{thm:q-weight}
follows once we show that \(G\mapsto (U_\Theta(t)D_{w,\Theta}^{-1}G)_{t\in I}\) is bounded
from \(L^p(\mathcal D)\) to \(L^p(\mathcal D_I)\). We first transfer operators with compact frequency support.
Write
\[
 d_w(\xi)=\langle\xi\rangle^{\sigma(p)}L_{B,C}(\langle\xi\rangle),
 \qquad \xi\in\R^2,
\]
for the symbol of \(D_w\). The functions \(d_w\) and \(d_w^{-1}\) are smooth, and each derivative has at most polynomial growth. Hence multiplication by either function preserves the Schwartz class on the plane; on the torus it preserves finite Fourier sums. Thus \(D_{w,\Theta}\) maps the span of these elementary tensors bijectively onto itself. Its image is dense in \(L^p(\mathcal D)\), so \(D_{w,\Theta}\) extends to an isometric isomorphism from the completion in \(\|D_{w,\Theta}x\|_{L^p(\mathcal D)}\) onto \(L^p(\mathcal D)\).

Choose \(\chi\in C_c^\infty(\R^2)\)
equal to one near zero. For each integer \(m\ge1\) the localized
Euclidean operator has symbol
\(e^{it|\xi|}\chi(2^{-m}\xi)/d_w(\xi)\). The Euclidean spacetime norm stays uniformly bounded
in \(m\) by the Euclidean weighted estimate for \(U(t)\) and the convolution bound \(\|\check\chi\|_{L^1(\R^2)}\).
For fixed \(m\) the finitely many shells with \(1\lesssim|\xi|\lesssim 2^m\)
have a common integrable kernel bound uniform in \(t\). At low frequencies we write
\(e^{it|\xi|}=1+(e^{it|\xi|}-1)\) and separate the smooth
term \(\chi_0(\xi)\chi(2^{-m}\xi)/d_w(\xi)\). On a shell \(|\xi|\simeq 2^k\),
with \(k\in\mathbb Z\), \(k\le0\), the rescaled factor \(e^{it2^k|\eta|}-1\), \(\eta=2^{-k}\xi\), and its first four derivatives are \(O(2^k)\) uniformly for \(t\in I\). With the annular partition \(\rho\) from Section~\ref{lp:sum-annuli}, the symbol on this shell is
\begin{equation*}
 \rho(2^{-k}\xi)\chi_0(\xi)\chi(2^{-m}\xi)(e^{it|\xi|}-1)/d_w(\xi)\,.
\end{equation*}
Its inverse Fourier transform at \(z\in\R^2\) is bounded by
\(C2^{3k}(1+2^k|z|)^{-4}\) by the low-frequency calculation in
Section~\ref{lp:sum-annuli}. The \(L^1\) norms of these majorants sum over \(k\le0\).
The sum of the localized kernels gives a jointly measurable kernel \(K_{m,t}\) with Fourier transform \(e^{it|\xi|}\chi(2^{-m}\xi)/d_w(\xi)\) and a common integrable majorant for \(t\in I\) at each fixed \(m\). Let \(T_m:L^p(\mathcal D)\to L^p(\mathcal D_I)\) be defined by \((T_mG)(t)=A_{K_{m,t}}G\). By~\eqref{q:transfer}, the norms of \(T_m\) are bounded independently of \(m\).

To pass to the limit in \(T_m\) as \(m\to\infty\), for each integer \(l\ge1\) write \(\chi_l(\xi)=\chi(2^{-l}\xi)\) and \(P_l=\chi_l(D_\Theta)=A_{\check\chi_l}\). The rescaling \(\check\chi_l(y)=2^{2l}\check\chi(2^ly)\) and \(\int_{\R^2}\check\chi(y)\,dy=\chi(0)=1\) give, for every \(G\in L^p(\mathcal D)\),
\begin{equation*}
 \|P_lG-G\|_{L^p(\mathcal D)}
 \le\int_{\R^2}|\check\chi(y)|\,\|\beta_{-2^{-l}y}G-G\|_{L^p(\mathcal D)}\,dy
 \longrightarrow0
 \qquad(l\to\infty)
\end{equation*}
by orbit continuity and dominated convergence.

For \(K,L\in L^1(\R^2)\) and \(G\in L^p(\mathcal D)\), Bochner Fubini and the action law give
\begin{equation*}
 A_K A_L G
 =\int_{\R^2}\int_{\R^2}
      K(y)L(v)\beta_{-(y+v)}G\,dv\,dy
 =A_{K*L}G\,.
\end{equation*}
The integrand has norm at most \(|K(y)L(v)|\|G\|_{L^p(\mathcal D)}\), so the changes in the order of integration are valid. For fixed \(l\), we have \(\chi_m\chi_l=\chi_l\) for every sufficiently large integer \(m\). The Fourier transform of \(K_{m,t}*\check\chi_l\) is then \(e^{it|\xi|}\chi_l(\xi)/d_w(\xi)\), independently of \(m\). By uniqueness of the Fourier transform on \(L^1\) and the composition rule, \(T_mP_lG=T_nP_lG\) for all sufficiently large \(m,n\). The common bound \(\|T_m\|\le C\) now gives
\begin{equation*}
 \|T_mG-T_nG\|_{L^p(\mathcal D_I)}
 \le 2C\|G-P_lG\|_{L^p(\mathcal D)}\,.
\end{equation*}
Since \(\|G-P_lG\|_{L^p(\mathcal D)}\to0\), the sequence \(T_mG\) is Cauchy in \(L^p(\mathcal D_I)\). Denote its limit by \(T_1 G\). Then \(\|T_1 G\|_{L^p(\mathcal D_I)}\le C\|G\|_{L^p(\mathcal D)}\).

For \(x\) in the initial tensor span, \(G=D_{w,\Theta}x\) belongs to \(L^2(\mathcal D)\cap L^p(\mathcal D)\). The \(L^2\) multiplier formula for \(A_{K_{m,t}}\) shows that \(T_mG=U_\Theta(t)\chi_m(D_\Theta)x\). By Plancherel and dominated convergence, \(\chi_m(D_\Theta)x\to x\) in \(L^2(\mathcal D)\). Since \(U_\Theta(t)\) is unitary on \(L^2(\mathcal D)\) and \(|I|=1\), we have \(T_mG\to U_\Theta(t)x\) in \(L^2(\mathcal D_I)\). The \(L^2\) and \(L^p\) limits agree by uniqueness of the limit in measure. Thus \((T_1(D_{w,\Theta}x))(t)=U_\Theta(t)x\) on the tensor span, and completion in \(\|D_{w,\Theta}x\|_{L^p(\mathcal D)}\) proves the weighted half-wave bound.

For \(G\in L^p(\mathcal D)\), define \(T_{-1}G=(T_1(G^*))^*\) where the outer adjoint acts on \(T_1(G^*)\) in \(L^p(\mathcal D_I)\). Then \(\|T_{-1}G\|_{L^p(\mathcal D_I)}\le C\|G\|_{L^p(\mathcal D)}\). Since \(d_w\) is real and even, \((D_{w,\Theta}x)^*=D_{w,\Theta}x^*\) on the initial tensor span. For \(\xi\in\R^2\), set \(\chi_m^\sharp(\xi)=\overline{\chi_m(-\xi)}\) as in~\eqref{intro:half-wave-adjoint}. The orbit adjoint identity applied to \(K_{m,t}\) gives
\begin{equation*}
 \bigl[(T_m(D_{w,\Theta}x^*))(t)\bigr]^*
 =U_\Theta(-t)\chi_m^\sharp(D_\Theta)x\,.
\end{equation*}
As \(m\to\infty\), the left side converges to \(T_{-1}(D_{w,\Theta}x)\) in \(L^p(\mathcal D_I)\). The right side converges to \(U_\Theta(-t)x\) in \(L^2(\mathcal D_I)\) because the bounded symbols \(\chi_m^\sharp\) converge pointwise to one. The two limits agree in measure, so \(T_{-1}\) gives the weighted estimate for \(U_\Theta(-t)\) on \(I\).

To prove the weighted estimate for \(u_\Theta(t;x,y)\), put \(G_x=D_{w,\Theta}x\) and \(G_y=J_\Theta^{-1}D_{w,\Theta}y\) for \(x,y\) in the initial tensor span. The displacement term is \(\tfrac12(T_1+T_{-1})G_x\). Split the velocity term using \(\chi_0\) and \(1-\chi_0\). The low-frequency operator \(S_{\mathrm{low}}:L^p(\mathcal D)\to L^p(\mathcal D_I)\) acting on \(G_y\) has symbol
\begin{equation*}
 \chi_0(\xi)\frac{\sin(t|\xi|)}{|\xi|}\frac{\langle\xi\rangle}{d_w(\xi)}
 \qquad(t\in I,\ \xi\in\R^2)\,.
\end{equation*}
The symbol is smooth and compactly supported, with uniformly bounded derivatives for \(t\in I\). By~\eqref{pre:kern} and~\eqref{q:orbit}, \(\|S_{\mathrm{low}}G_y\|_{L^p(\mathcal D_I)}\le C\|G_y\|_{L^p(\mathcal D)}\).

For the high-frequency velocity term define \(a:\R^2\to\R\) by \(a(\xi)=(1-\chi_0(\xi))\langle\xi\rangle/|\xi|\) for \(\xi\ne0\) and \(a(0)=0\). The partition of unity gives
\begin{equation*}
 a(\xi)=1-\chi_0(\xi)+\sum_{j\ge1}\varphi(2^{-j}\xi)
                        \left(\frac{\langle\xi\rangle}{|\xi|}-1\right)\,,
\end{equation*}
where each annular symbol is zero near the origin. In the scaled variable \(\eta=2^{-j}\xi\) on \(\operatorname{supp}\varphi\),
\begin{equation*}
 \frac{\langle2^j\eta\rangle}{2^j|\eta|}-1
 =\frac{2^{-2j}|\eta|^{-2}}
        {\sqrt{1+2^{-2j}|\eta|^{-2}}+1},\qquad j\ge1\,.
\end{equation*}
This expression and its derivatives in \(\eta\) through order four are \(O(2^{-2j})\) on the fixed annulus. By~\eqref{pre:kern} the annular kernels have summable \(L^1\) norms. Using~\eqref{q:orbit} for the annular kernels and for \(\check\chi_0\) therefore defines \(a(D_\Theta)\) as a bounded operator on \(L^p(\mathcal D)\). On the tensor span it agrees with the multiplier of symbol \(a\) by the \(L^2\) formula for \(A_K\). The symbol \(a\) is smooth, with derivatives of at most polynomial growth, so this multiplier preserves the tensor span.

On the tensor span, the identity \(\sin(t|\xi|)=(e^{it|\xi|}-e^{-it|\xi|})/(2i)\) gives the full solution
\begin{equation*}
 u_\Theta(\cdot;x,y)=\frac{T_1+T_{-1}}{2}G_x+S_{\mathrm{low}}G_y
                 +\frac{T_1-T_{-1}}{2i}a(D_\Theta)G_y\,.
\end{equation*}
The bounds for \(T_1,T_{-1},S_{\mathrm{low}}\) and \(a(D_\Theta)\) give \(\|u_\Theta(\cdot;x,y)\|_{L^p(\mathcal D_I)}\le C(\|G_x\|_{L^p(\mathcal D)}+\|G_y\|_{L^p(\mathcal D)})\). This proves the wave estimate of Theorem~\ref{thm:quantum}\ref{thm:q-weight} after completion in \(\|D_{w,\Theta}x\|_{L^p(\mathcal D)}\) and \(\|J_\Theta^{-1}D_{w,\Theta}y\|_{L^p(\mathcal D)}\), respectively.

\section[Maximal inequalities and convergence for generalized circular means]{Maximal inequalities and convergence\\for generalized circular means}\label{max:sec}

We prove the all-radius estimate \eqref{intro:maximal} on
$\mathcal A=L^\infty(\R^2)\,\bar\otimes\,\M$, equipped with the product trace.
For $2\le p<\infty$ and $t|\xi|\asymp 2^j$, local smoothing and the decay of
the Bessel symbol give a bound for the mean and its radius derivative
on $1\le t\le2$. From this pair of estimates we construct positive
majorants whose norms are summable in \(j\ge1\) over all radii when \(\Re\alpha>-g(p)\).
Section~\ref{q:secmaximal} transfers the maximal estimates to trace-preserving actions. Section~\ref{q:secconvergence} uses this control to prove bilateral almost uniform convergence as the radius tends to zero.

For $t>0$ and $\alpha\in\mathbb C$, the circular mean $A_t^\alpha$ preserves the linear span of elementary tensors $f\otimes a$, where $f\in\mathscr S(\R^2)$ and $a\in\M$ is bounded with support projections of finite trace. Its spatial Fourier symbol is $\xi\mapsto m_\alpha(t|\xi|)$ for $\xi\in\R^2$, where $m_\alpha:[0,\infty)\to\mathbb C$ is determined for $r>0$ by the Bessel function $J_\alpha$ of the first kind and its series expansion
\begin{equation}\label{max:symbol}
 m_\alpha(r)=\pi 2^\alpha r^{-\alpha}J_\alpha(r)
 =\pi\sum_{n=0}^\infty\frac{(-r^2/4)^n}{n!\,\Gamma(n+1+\alpha)}\,.
\end{equation}
We use the real logarithm in $r^{-\alpha}$. The entire reciprocal gamma function in the series gives $m_\alpha(0)=\pi/\Gamma(1+\alpha)$, including at poles of $\Gamma$. When $\Re\alpha>0$, the kernel $K_t^\alpha\in L^1(\R^2)$ is given for almost every $y\in\R^2$ by
\begin{equation}\label{max:kernel}
 K_t^\alpha(y)=\frac{1}{t^2\Gamma(\alpha)}
      (1-|y|^2/t^2)_+^{\alpha-1}\,.
\end{equation}
Indeed, at $r=t|\xi|$ polar integration gives
\[
 \widehat K_t^\alpha(\xi)
 =\frac{2\pi}{\Gamma(\alpha)}\int_0^1(1-s^2)^{\alpha-1}J_0(sr)s\,ds .
\]
Insert the power series for $J_0(sr)$, which is absolutely and uniformly convergent on $0\le s\le1$ for each fixed $r$. The substitution $v=s^2$ gives
\[
 \int_0^1(1-s^2)^{\alpha-1}s^{2n+1}\,ds
 =\frac{\Gamma(n+1)\Gamma(\alpha)}{2\Gamma(n+1+\alpha)} .
\]
This proves~\eqref{max:symbol} in the integrable range. Near $|y|=t$ write $y=t\sqrt{1-v}\,\omega$, where $v=1-|y/t|^2$ and $\omega\in S^1$. For a smooth test function $\phi$ supported in a thin annulus around the circle, put $\Phi(v)=\int_{S^1}\phi(t\sqrt{1-v}\,\omega)\,d\omega$, where $d\omega$ is arc length on the unit circle. The function $\Phi$ is a smooth test function of $v$, extended by zero outside the coordinate interval. Writing $r=|y|$ gives $r\,dr=-\tfrac12 t^2\,dv$, so for $\Re\alpha>0$ we obtain
\begin{equation*}
 \langle K_t^\alpha,\phi\rangle
 =\frac12\left\langle \frac{v_+^{\alpha-1}}{\Gamma(\alpha)},\Phi\right\rangle\,.
\end{equation*}
The identity
\[
 \frac{v_+^{\alpha-1}}{\Gamma(\alpha)}
 =\partial_v^k\frac{v_+^{\alpha+k-1}}{\Gamma(\alpha+k)}
 \qquad (k\in\mathbb N,\ \Re(\alpha+k)>0)
\]
continues the one-variable pairing to $\Re(\alpha+k)>0$. The continuations for different $k$ agree by the identity theorem. The polar pairing therefore defines an entire distribution family near the circle. It agrees with the analytic kernels on the interior and exterior overlaps by the identity theorem, since both constructions agree with~\eqref{max:kernel} when $\Re\alpha>0$. Testing outside the closed disk gives zero throughout the continuation, so $\operatorname{supp}K_t^\alpha\subset\overline B(0,t)$ for all $\alpha$. For fixed \(t>0\), take \(\zeta_t\in C_c^\infty(\R^2)\) equal to
one on a neighborhood of the closed disk. Testing against
\(\zeta_t(y)e^{-iy\cdot\xi}\) extends the Fourier identity to every \(\alpha\in\mathbb C\). At $\alpha=0$, the inverse transform of $\pi J_0(t|\xi|)$ is $\pi$ times the probability measure on $|y|=t$.

To measure a family of radii, fix $1\le p\le\infty$ and an index set $E$.
For a family $Y_E=(Y_t)_{t\in E}$ in $L^p(\mathcal A)$, the maximal norm is
\begin{equation}\label{max:norm}
 \|Y_E\|_{L^p(\mathcal A;\ell^\infty(E))}
   =\inf_{Y_t=a z_t b}\|a\|_{L^{2p}(\mathcal A)}\sup_{t\in E}\|z_t\|_{L^\infty(\mathcal A)}\|b\|_{L^{2p}(\mathcal A)}\,,
\end{equation}
with common factors $a,b\in L^{2p}(\mathcal A)$ and $z_t\in\mathcal A$. An infimum over an empty set is $+\infty$. We suppress $E$ in $\ell^\infty(E)$ when the indices are displayed on the family. The norm can be computed on finite sets of indices
\cite[\S2.2, (2.1) and the arbitrary-index formulation on p.~8]{hlw}.
For arbitrary complex families,
\begin{equation}\label{max:finite}
 \|Y_E\|_{L^p(\mathcal A;\ell^\infty(E))}
 =\sup_{\substack{E_0\subset E\\ E_0\ \mathrm{finite}}}
       \|Y_{E_0}\|_{L^p(\mathcal A;\ell^\infty(E_0))}.
\end{equation}
For self-adjoint $Y_t$, the order characterization in the same reference gives
\[
 \|Y_E\|_{L^p(\mathcal A;\ell^\infty(E))}=\inf\{\|B\|_{L^p(\mathcal A)}:B\in L^p(\mathcal A)_+,\ -B\le Y_t\le B\ \text{for all }t\in E\}.
\]
To see the order description for uncountably many indices in the range $1<p<\infty$ used below, suppose the supremum in~\eqref{max:finite} is at most $C$ for a self-adjoint family. In the positive ball of radius $C+\delta$, $\delta>0$, each condition $B\ge\pm Y_t$ is weakly closed. The ball is weakly compact by reflexivity, and the finite-family order formula gives the finite intersection property. Thus a single $B$ satisfies all the inequalities. Letting $\delta\downarrow0$ gives the same bound for the infimum.

For general $Y_t$ we can control $\Re Y_t$ and $\Im Y_t$ separately. Taking adjoints in a factorization gives $Y_t^*=b^*z_t^*a^*$ and preserves the norm. Consequently each of the two self-adjoint families has norm at most $\|Y_E\|_{L^p(\mathcal A;\ell^\infty(E))}$. Conversely, majorants of norms $C_1$ and $C_2$ for both signs of the real and imaginary parts give the upper bound $C_1+C_2$ in~\eqref{max:norm}, by the triangle inequality. These formulations also apply with $\mathcal A$ replaced by any semifinite tracial algebra.

\subsection{A frequency band on a bounded interval of radii}

Fix $2\le p<\infty$. Choose $\varphi\in C_c^\infty((1/2,2))$ with $\sum_{j\in\mathbb Z}\varphi(2^{-j}r)=1$ for $r>0$, and define the low-frequency cutoff $\eta$ on $[0,\infty)$ by $\eta(r)=1-\sum_{j\ge1}\varphi(2^{-j}r)$. The decomposition in the variable $t|\xi|$, used in~\cite[\S3]{hlw}, reads
\begin{equation}\label{max:decomp}
 A_t^\alpha=A_{t,0}+\sum_{j\ge1} A_{t,j},\qquad
 A_{t,0}=(\eta m_\alpha)(t|D|),\quad
 A_{t,j}=m_\alpha(t|D|)\varphi(2^{-j}t|D|)\,.
\end{equation}
For fixed $t>0$ and integer $j\ge0$, the smooth compactly supported
symbol of $A_{t,j}$ defines a bounded operator on $L^p(\mathcal A)$.
The triangle inequality in the maximal space reduces
\eqref{intro:maximal} for \(2\le p<\infty\) to bounding
\[
 \|(A_{t,0}F)_{t>0}\|_{L^p(\mathcal A;\ell^\infty)}
 +\sum_{j\ge1}\|(A_{t,j}F)_{t>0}\|_{L^p(\mathcal A;\ell^\infty)} .
\]
We treat the low-frequency family by positive ball averages. For
$j\ge1$ we first restrict $t$ to $I=[1,2]$ to apply local smoothing,
then combine the dilates of $I$. The range $1<p<2$ will follow by
interpolation in $\alpha$ after the summation at $p=2$.

The inverse transform $k_0$ of $\eta(|\xi|)m_\alpha(|\xi|)$ is Schwartz by~\eqref{max:symbol}, with seminorms bounded for $\alpha$ in compact subsets of~$\mathbb C$. In particular
\begin{equation}\label{max:low-kernel}
 t^{-2}|k_0(y/t)|\le C_\alpha
 \sum_{l\ge0}2^{-2l}|B(0,2^lt)|^{-1}\mathbf1_{B(0,2^lt)}(y)\,.
\end{equation}
For $F\in L^p(\mathcal A)_+$, the kernel majorant gives an order bound for both signs of $\Re A_{t,0}F$ and $\Im A_{t,0}F$: each is bounded above by
\[
 C_\alpha\sum_{l\ge0}2^{-2l}|B(0,2^lt)|^{-1}
       \int_{|y|<2^lt}F(\cdot-y)\,dy.
\]
By the noncommutative Hardy--Littlewood maximal inequality for balls, in the form used in~\cite[proof of Lemma~3.2]{hlw}, these averages, over all ball radii, have a common positive majorant of $L^p(\mathcal A)$ norm at most $C_p\|F\|_{L^p(\mathcal A)}$. Splitting the real and imaginary parts of general $F$ into positive and negative parts thus gives
\begin{equation}\label{max:low}
 \|(A_{t,0}F)_{t>0}\|_{L^p(\mathcal A;\ell^\infty)}
       \le C_{p,\alpha}\|F\|_{L^p(\mathcal A)}\,.
\end{equation}

The annular estimate~\eqref{intro:annular} controls the $L^p(I;L^p(\mathcal A))$ norm. We use the following form of the time Sobolev estimate of Hong, the first author and Wang~\cite[Lemma~3.6]{hlw} to pass to a maximal estimate. For self-adjoint families, the lemma also gives an explicit positive majorant.

\begin{lemma}\label{max:timelemma}
Let $(\mathcal B,\tau)$ be a semifinite tracial von Neumann algebra, $2\le p<\infty$, and $I=[1,2]$. For $H\in C^1(I;L^p(\mathcal B))$ and each integer $N\ge1$,
\begin{equation}\label{max:order}
 \|(H(u))_{u\in I}\|_{L^p(\mathcal B;\ell^\infty)}
 \le 2 N^{1/p}\Big(\|H\|_{L^p(I;L^p(\mathcal B))}
                 + N^{-1}\|H'\|_{L^p(I;L^p(\mathcal B))}\Big)\,.
\end{equation}
If $H$ is self-adjoint, there is a single $B\in L^p(\mathcal B)_+$ such that $-B\le H(u)\le B$ for all $u\in I$, and $\|B\|_{L^p(\mathcal B)}$ satisfies the right-hand bound in~\eqref{max:order} with factor $1$ in place of $2$.
\end{lemma}
\begin{proof}
For self-adjoint $H$, we construct a positive majorant on each subinterval
of length $N^{-1}$ and combine these majorants by summing their $p$th powers.
Partition $I$ into $N$ intervals $J$ of length $N^{-1}$ with disjoint interiors. For $u,v\in J$, the identity $H(u)=H(v)+\int_v^u H'(s)\,ds$ and the order bound $\pm H'(s)\le |H'(s)|$ give
\begin{equation*}
 \pm H(u)\le |H(v)|+\int_J |H'(s)|\,ds\,.
\end{equation*}
Averaging over $v\in J$ gives
\begin{equation*}
 -B_J\le H(u)\le B_J\quad(u\in J),\qquad
 B_J=N\int_J |H(v)|\,dv+\int_J |H'(v)|\,dv\,.
\end{equation*}
To justify the Bochner integrals defining $B_J$, we check strong measurability of $|H|$ and $|H'|$ in $L^p(\mathcal B)$. For a continuous curve $X:J\to L^p(\mathcal B)$ with self-adjoint values, we use the representation
\[
 |X(v)|=\lim_{n\to\infty}\frac1\pi
  \int_{1/n}^{n}\lambda^{-1/2}X(v)^2(X(v)^2+\lambda)^{-1}\,d\lambda .
\]
The spectral theorem gives, for $\lambda>0$ and $v\in J$,
\begin{equation*}
 0\le X(v)^2(X(v)^2+\lambda)^{-1}
 \le \frac{|X(v)|}{2\sqrt\lambda}\,.
\end{equation*}
In particular the integrand in the truncated integral belongs to $L^p(\mathcal B)$, and spectral dominated convergence gives continuity in $\lambda$ in $L^p(\mathcal B)$ on compact subintervals of $(0,\infty)$. For fixed $\lambda>0$, the scalar function $s\mapsto s^2/(s^2+\lambda)$ on $\R$ is a constant plus a linear combination of $(s\pm i\sqrt\lambda)^{-1}$. The resolvent identity gives
\[
 \big\|X(v)^2(X(v)^2+\lambda)^{-1}
       -X(w)^2(X(w)^2+\lambda)^{-1}\big\|_{L^p(\mathcal B)}
 \le C_\lambda\|X(v)-X(w)\|_{L^p(\mathcal B)}.
\]
Here $v,w\in J$, and $C_\lambda$ is bounded when $\lambda$ ranges over a compact subinterval of $(0,\infty)$. Thus each truncated integral is continuous in $v$. The spectral theorem and dominated convergence against the trace measure of $|X(v)|^p$ give convergence in $L^p(\mathcal B)$ for each $v$. Thus $v\mapsto|X(v)|$ is strongly measurable. Applying this conclusion to $H$ and $H'$, and using $\||X(v)|\|_{L^p(\mathcal B)}=\|X(v)\|_{L^p(\mathcal B)}$, proves that both integrands defining $B_J$ are Bochner integrable.

It remains to dominate $H(u)$ on all of $I$ by one positive operator.
Put $S=\sum_J B_J^p\in L^1(\mathcal B)_+$ and $B=S^{1/p}$. From $B_J^p\le S$ and operator monotonicity of the $1/p$ power we obtain $B_J\le B$ for each $J$. For positive measurable operators this monotonicity follows from the integral representation of $s^{1/p}$ by positive multiples of $s/(s+\lambda)$ and reversal of order by the resolvent $(s+\lambda)^{-1}$ for $\lambda>0$. Taking the trace gives $\|B\|_{L^p(\mathcal B)}^p=\sum_J\|B_J\|_{L^p(\mathcal B)}^p$. On each $J$, H\"older's inequality in time gives
\begin{equation*}
 \|B_J\|_{L^p(\mathcal B)}\le N^{1/p}\|H\|_{L^p(J;L^p(\mathcal B))}
               +N^{-1+1/p}\|H'\|_{L^p(J;L^p(\mathcal B))}\,.
\end{equation*}
Taking the $\ell^p$ norm over the partition and using
$\|B\|_{L^p(\mathcal B)}^p=\sum_J\|B_J\|_{L^p(\mathcal B)}^p$ gives \eqref{max:order} with factor one. For general $H$, apply the self-adjoint conclusion separately to $\Re H$ and $\Im H$. Taking real or imaginary parts commutes with differentiation and decreases the $L^p(\mathcal B)$ norm. The triangle inequality in~\eqref{max:norm} now gives~\eqref{max:order}.
\end{proof}

We apply Lemma~\ref{max:timelemma} at time resolution $2^{-j}$, for integers $j\ge1$, to obtain a maximal estimate for the frequency band $A_{u,j}$.

\begin{lemma}\label{max:local}
For $2\le p<\infty$, $\alpha\in\mathbb C$, integers $j\ge1$, and $F\in L^p(\mathcal A)$,
\begin{equation}\label{max:band}
 \|(A_{u,j}F)_{u\in I}\|_{L^p(\mathcal A;\ell^\infty)}
 \le C_{p,\alpha}\, 2^{-j(\Re\alpha+g(p))}(1+j)^{2g(p)}\|F\|_{L^p(\mathcal A)}\,.
\end{equation}
For fixed $p$, the constants may be taken uniformly on compact sets of orders~$\alpha$.
\end{lemma}
\begin{proof}
Put $N=2^j$. By Lemma~\ref{max:timelemma}, it suffices to bound $A_{u,j}F$ and $N^{-1}\partial_u A_{u,j}F$ in $L^p(I;L^p(\mathcal A))$. We estimate the terms with $U(u)$ and $U(-u)$ by the annular estimate of Theorem~\ref{thm:lp}.

For $r\ge1$, the identity $2J_\alpha=H_\alpha^{(1)}+H_\alpha^{(2)}$ gives the exact decomposition
\begin{equation}\label{max:hankel}
 m_\alpha(r)=\sum_{\nu=\pm1}e^{i\nu r}r^{-1/2-\alpha}b_{\alpha,\nu}(r),\qquad
 |\partial_r^l b_{\alpha,\nu}(r)|\le C_{l,\alpha}\,r^{-l}
 \quad(l=0,1,2,\ldots)\,.
\end{equation}
Indeed, take $b_{\alpha,+}(z)=\pi2^{\alpha-1} z^{1/2}e^{-iz}H_\alpha^{(1)}(z)$ and $b_{\alpha,-}(z)=\pi2^{\alpha-1} z^{1/2}e^{iz}H_\alpha^{(2)}(z)$ for $z$ in a complex sector about the positive real axis, using the branch of $z^{1/2}$ positive on that axis. The large-argument expansions with their complex-order remainder bounds \cite[\S10.17(i),(iv)]{dlmf} bound both functions there, uniformly for $\alpha$ in a compact set. Cauchy's formula on disks of radius $cr$ about $r$ gives~\eqref{max:hankel}, with the bounded range of $r$ covered by continuity.

To apply \eqref{intro:annular}, we separate the oscillatory phase from a
spatial multiplier whose kernel has bounded $L^1(\R^2)$ norm. Take
$\chi\in C_c^\infty((0,\infty))$ equal to one on $[1/4,2]$. For $u\in I$ and $\xi\in\R^2$, put $w=u|\xi|/N$. For $\nu\in\{-1,1\}$ define $a_\nu:\R\to\mathbb C$ by $a_\nu(w)=\varphi(w)w^{-1/2-\alpha} b_{\alpha,\nu}(Nw)$ on the support of $\varphi$, extended by zero. The term with phase $e^{i\nu u|\xi|}$ has symbol $N^{-1/2-\alpha}e^{i\nu u|\xi|}a_\nu(w)$. Since $\partial_u w=|\xi|/N$, differentiation acts both on the phase and on $a_\nu$. After extracting $N^{-1/2-\alpha}e^{i\nu u|\xi|}\chi(|\xi|/N)$, the residual symbols for $A_{u,j}$ and $N^{-1}\partial_u A_{u,j}$ are, respectively,
\begin{equation*}
 a_\nu(w),\qquad \frac{|\xi|}{N}\big(i\nu a_\nu(w)+N^{-1}a_\nu'(w)\big)\,.
\end{equation*}
For $u\in I$, the residual symbols in this display have fixed compact annular support in $\zeta=\xi/N$. Their derivatives in $\zeta$ through order four are uniformly bounded by~\eqref{max:hankel}. Applying $(1-\Delta_\zeta)^2$ in Fourier inversion gives an integrable majorant $C(1+|x|^2)^{-2}$ in the scaled spatial variable, hence a uniform $L^1(\R^2)$ kernel bound in the original variable.

We can therefore estimate both $A_{u,j}F$ and $N^{-1}\partial_u A_{u,j}F$ by applying these kernels with uniformly bounded $L^1(\R^2)$ norms to $U(\nu u)\chi(|D|/N)F$. Taking the $L^p(I)$ norm in $u$, we use~\eqref{intro:annular} for both time signs and for the fixed cutoff $\chi$, as justified by the finite annular covering in Section~\ref{lp:annular}. Restoring the factor $N^{-1/2-\alpha}$ and using $\beta(p)=2g(p)$ gives
\begin{equation}\label{max:time}
 \sum_{s=0}^1 N^{-s}\|\partial_u^s A_{u,j}F\|_{L^p(I;L^p(\mathcal A))}
 \le C_{p,\alpha}N^{-1/2-\Re\alpha+\sigma(p)}(1+j)^{2g(p)}\|F\|_{L^p(\mathcal A)}\,.
\end{equation}

For fixed $j$ the map $u\mapsto A_{u,j}F$ is $C^1$ in $L^p(\mathcal A)$ because the compact annular symbol and its $u$-derivative have Schwartz kernels continuous in $L^1(\R^2)$. Apply Lemma~\ref{max:timelemma} to this map with $N=2^j$. Multiplying the right side of~\eqref{max:time} by $2N^{1/p}$ gives~\eqref{max:band}, since $1/p-1/2+\sigma(p)=-g(p)$.
\end{proof}

\subsection{All radii and the critical order}

We combine the estimate on \(I=[1,2]\) with its dyadic dilates, using the bounded overlap of the spatial frequency annuli selected by \(A_{2^ku,j}\) as \(k\in\mathbb Z\) varies with \(u\in I\) and \(j\ge1\) fixed. This gives \eqref{intro:maximal} for \(2\le p<\infty\); interpolation in \(\alpha\) will then give the range \(1<p<2\).

\begin{proof}[Proof of the Euclidean assertion of Theorem~\ref{thm:maximal}]
For $2\le p<\infty$, it remains to pass from~\eqref{max:band} to a summable estimate on $(A_{t,j}F)_{t>0}$. For $j\ge1$ put $D_j=2^{-j(\Re\alpha+g(p))}(1+j)^{2g(p)}$. If $t=2^ku$ with $k\in\mathbb Z$ and $u\in[1,2]$, the support of $\varphi(2^{-j}t|\xi|)$ forces $2^{k-j}|\xi|\in[1/4,2]$. At fixed $j$ these annuli have bounded overlap as $k$ varies. For $l\in\mathbb Z$ let $\Delta_l=\chi(2^{-l}|D|)$ on $L^p(\mathcal A)$, where $\chi$ is the cutoff in Lemma~\ref{max:local}. Since $\chi=1$ on $[1/4,2]$, we have
\begin{equation*}
 A_{2^ku,j}F=A_{2^ku,j}\Delta_{j-k}F\,.
\end{equation*}
Apply~\eqref{max:band} to $G_{k,j}(x)=(\Delta_{j-k}F)(2^kx)$ for $x\in\R^2$. For each $\bullet\in\{\Re,\Im\}$, the order description of~\eqref{max:norm} gives a positive majorant $\widetilde B_{k,j}^\bullet\in L^p(\mathcal A)$ for both signs of $\bullet(A_{u,j}G_{k,j})$ on $u\in I$, with $\|\widetilde B_{k,j}^\bullet\|_{L^p(\mathcal A)}\le C_{p,\alpha}D_j\|G_{k,j}\|_{L^p(\mathcal A)}$. By Fourier dilation,
\begin{equation*}
 (A_{u,j}G_{k,j})(x)=(A_{2^ku,j}\Delta_{j-k}F)(2^kx)\,.
\end{equation*}
Set $B_{k,j}^\bullet(x)=\widetilde B_{k,j}^\bullet(2^{-k}x)$. The dilation preserves order. The two changes of spatial variables give
\begin{equation*}
 \|G_{k,j}\|_{L^p(\mathcal A)}=2^{-2k/p}\|\Delta_{j-k}F\|_{L^p(\mathcal A)},\qquad
 \|B_{k,j}^\bullet\|_{L^p(\mathcal A)}=2^{2k/p}\|\widetilde B_{k,j}^\bullet\|_{L^p(\mathcal A)}.
\end{equation*}
The powers of \(2\) cancel, giving
\begin{equation*}
 \pm\bullet(A_{2^ku,j}F)\le B_{k,j}^\bullet\quad(u\in I),\qquad
 \|B_{k,j}^\bullet\|_{L^p(\mathcal A)}\le C_{p,\alpha}D_j\|\Delta_{j-k}F\|_{L^p(\mathcal A)}\,.
\end{equation*}

For fixed $j$ the index $j-k$ runs through all integers as $k$ varies. To sum the norms of the majorants we use
\begin{equation}\label{max:analysis}
 \sum_{l\in\mathbb Z}\|\Delta_lF\|_{L^p(\mathcal A)}^p\le C_p^p\|F\|_{L^p(\mathcal A)}^p\,.
\end{equation}
Indeed, consider $F\mapsto(\Delta_l F)_l$ with target the ordinary tracial $L^p$ space of the direct sum algebra $\ell^\infty(\mathbb Z)\,\bar\otimes\,\mathcal A$ with counting trace in the first factor. Plancherel's identity in the spatial variables and bounded overlap of the annuli give its $L^2$ bound. The $L^\infty$ norm in this direct sum is $\sup_l\|\Delta_l F\|_{L^\infty(\mathcal A)}$,
bounded by $\|\chi(|\cdot|)^\vee\|_{L^1(\R^2)}\|F\|_{L^\infty(\mathcal A)}$ since dilation
preserves the $L^1(\R^2)$ norm of the kernels. Complex interpolation, first on finite index sets and then by monotone convergence of the sums, gives~\eqref{max:analysis}.

For each $\bullet\in\{\Re,\Im\}$, the bound for $B_{k,j}^\bullet$ and~\eqref{max:analysis} imply convergence of $S_j^\bullet=\sum_{k\in\mathbb Z}(B_{k,j}^\bullet)^p$ in $L^1(\mathcal A)_+$. Define $B_j^\bullet=(S_j^\bullet)^{1/p}$. Taking the trace gives
\begin{equation*}
 \|B_j^\bullet\|_{L^p(\mathcal A)}^p=\sum_{k\in\mathbb Z}\|B_{k,j}^\bullet\|_{L^p(\mathcal A)}^p
 \le C_{p,\alpha}^pD_j^p\|F\|_{L^p(\mathcal A)}^p\,.
\end{equation*}
The positive cone is closed in $L^1(\mathcal A)$, so $(B_{k,j}^\bullet)^p\le S_j^\bullet$ for each $k$. As in the proof of Lemma~\ref{max:timelemma}, operator monotonicity of the $1/p$ power gives $B_{k,j}^\bullet\le B_j^\bullet$. Every $t>0$ has the form $2^ku$ with $u\in I$, so $\pm\bullet(A_{t,j}F)\le B_j^\bullet$ for all $t>0$. The triangle inequality in~\eqref{max:norm} now gives
\begin{equation}\label{max:fullband}
 \|(A_{t,j}F)_{t>0}\|_{L^p(\mathcal A;\ell^\infty)}
      \le C_{p,\alpha}D_j\|F\|_{L^p(\mathcal A)}\,.
\end{equation}
For $\Re\alpha+g(p)>0$, the triangle inequality in~\eqref{max:norm} bounds the norm of the family in~\eqref{max:decomp} by the sum of~\eqref{max:low} and~\eqref{max:fullband} over $j\ge1$. Indeed $\sum_{j\ge1}D_j<\infty$, so the family-valued series
$\sum_{j\ge1}(A_{t,j}F)_{t>0}$ converges in
$L^p(\mathcal A;\ell^\infty)$. For a family $Y_E$ and $t\in E$, tracial H\"older's inequality applied to the factorization in~\eqref{max:norm} gives $\|Y_t\|_{L^p(\mathcal A)}\le\|Y_E\|_{L^p(\mathcal A;\ell^\infty(E))}$. Thus for each $t>0$ the family-norm limit gives the limit of its $t$-coordinate in $L^p(\mathcal A)$. To identify this coordinate limit, take a Schwartz tensor with scalar factor $f\in\mathscr S(\R^2)$. The series~\eqref{max:symbol} and the Hankel estimate~\eqref{max:hankel} show that $m_\alpha(t|\xi|)\widehat f(\xi)$ is integrable in $\xi$. At fixed $t$ the partial cutoff sums $\eta(t|\xi|)+\sum_{j=1}^m\varphi(2^{-j}t|\xi|)$ are uniformly bounded by finite overlap and converge to one by~\eqref{max:decomp}. Dominated convergence in Fourier inversion therefore identifies the coordinate limit of $A_{t,0}F+\sum_{j=1}^m A_{t,j}F$ with $A_t^\alpha F$ on Schwartz tensors. Density in $L^p(\mathcal A)$ gives the asserted extension and estimate for all~$F$ in this range. In particular, $g(2)=0$, so the $L^2$ band norms sum for $\Re\alpha>0$.

\medskip
For $1<p<2$ the interpolation argument requires a bound for the
$L^2$ maximal norm along the entire line $\Re\alpha=h$, including its
growth in $\Im\alpha$. Fix $h>0$ and, for $v\in\R$, put $\alpha=h+iv$. For $r\ge1$ we claim
\begin{equation}\label{max:vertical-symbol}
 |m_\alpha(r)|+|m_\alpha'(r)|\le C_h e^{C_h(1+|v|)}r^{-h-1/2}\,.
\end{equation}
Integrating the kernel~\eqref{max:kernel} at $t=1$ in its second coordinate and using the beta integral gives
\begin{equation*}
 m_\alpha(r)=\frac{\sqrt\pi}{\Gamma(\alpha+1/2)}
    \int_{-1}^1 e^{-irs}(1-s^2)^{\alpha-1/2}\,ds\,.
\end{equation*}
Differentiation in $r$ inserts the bounded smooth factor $-is$, so it suffices to estimate the displayed integral with amplitude $(1-s^2)^{\alpha-1/2}$ or $-is(1-s^2)^{\alpha-1/2}$. Near each endpoint $s=1$ and $s=-1$, we put $d=1-|s|$ and use a smooth dyadic partition $d\asymp D$, $D=2^{-l}$ with integers $l\ge1$. On such a piece $1-s^2=d(2-d)$ and the length of the interval is $O(D)$. Each derivative of $d^{\alpha-1/2}$ costs at most $C_{h,M}(1+|v|)D^{-1}$ relative to its size. The dyadic cutoff has derivatives of order $m$ bounded by $C_m D^{-m}$. Thus the derivatives of order $m$ of the localized amplitude are bounded by $C_{h,M}(1+|v|)^M D^{h-1/2-m}$ for $0\le m\le M$. The trivial estimate and $M$ integrations by parts therefore bound its integral by
\[
 C_{h,M}(1+|v|)^M D^{h+1/2}\min\{1,(rD)^{-M}\}.
\]
Choose an integer $M>h+1/2$. For $D\le r^{-1}$ the geometric sum of $D^{h+1/2}$ is $O(r^{-h-1/2})$. For $D>r^{-1}$ we sum $r^{-M}D^{h+1/2-M}$, whose exponent of $D$ is negative, and obtain the same bound. On the compact subinterval away from $s=\pm1$, integrating by parts
$M$ times gives $O_{h,M}((1+|v|)^M r^{-M})$. Stirling's formula on $\Re(\alpha+1/2)=h+1/2$ now proves~\eqref{max:vertical-symbol}.

For $N=2^j$, $u\in I$ and $u|\xi|\asymp N$,~\eqref{max:vertical-symbol} bounds the symbols of $A_{u,j}$ and $N^{-1}\partial_u A_{u,j}$ by $C_h e^{C_h(1+|v|)}N^{-h-1/2}$, since $|\xi|/N$ is bounded on their supports. Plancherel's identity in the spatial variables on the Hilbert space $L^2(\M)$ gives this same bound, times $\|F\|_{L^2(\mathcal A)}$, for both $L^2(I;L^2(\mathcal A))$ norms because $|I|=1$. Lemma~\ref{max:timelemma} costs a factor $2N^{1/2}$ to control all $u\in I$ simultaneously. The resulting bound on the maximal norm over $u\in I$ therefore contains the decaying factor $2N^{1/2}N^{-h-1/2}=2N^{-h}$. Summing $\|\Delta_{j-k}F\|_{L^2(\mathcal A)}^2$ over $k$ by~\eqref{max:analysis} now gives~\eqref{max:fullband} with $D_j=2^{-jh}$ and constant $C_he^{C_h(1+|v|)}$.

For the low-frequency bound~\eqref{max:low}, the spatial derivatives in $\xi$ through order four of $m_\alpha(|\xi|)$ on the support of $\eta$ are bounded by integrating $|K_1^\alpha(y)|(1+|y|)^4$. By~\eqref{max:kernel} this is at most $C_h/|\Gamma(h+iv)|$. Integrating by parts four times in Fourier inversion gives~\eqref{max:low-kernel} with \(C_\alpha\le C_h e^{C_h(1+|v|)}\). Adding the low frequencies and summing $N^{-h}$ yields
\begin{equation}\label{max:vertical}
 \|(A_t^{h+iv} F)_{t>0}\|_{L^2(\mathcal A;\ell^\infty)}
 \le C_h e^{C_h(1+|v|)}\|F\|_{L^2(\mathcal A)}\,.
\end{equation}

\medskip
To interpolate with the \(L^2\) estimate~\eqref{max:vertical}, we use ball averages on the line \(\Re\alpha=1\). On $\Re\alpha=1$ the absolute value of~\eqref{max:kernel} is constant on the disk $|y|<t$. For $F\ge0$ denote its normalized ball average by $B_tF=|B(0,t)|^{-1}\int_{|y|<t}F(\cdot-y)\,dy$. The scalar inequalities $\pm\Re K_t^{1+iv},\ \pm\Im K_t^{1+iv}\le |K_t^{1+iv}|$ give
\begin{equation*}
 \pm\Re A_t^{1+iv}F,\quad\pm\Im A_t^{1+iv}F
 \ \le\ \frac{\pi}{|\Gamma(1+iv)|}B_tF .
\end{equation*}
For $1<q<\infty$ and $F\in L^q(\mathcal A)_+$, the ball maximal inequality used in~\eqref{max:low} supplies a positive element of $L^q(\mathcal A)$ dominating every $B_tF$ with norm at most $C_q\|F\|_{L^q(\mathcal A)}$. Decomposing the real and imaginary parts of general $F$ into their positive and negative parts thus gives
\begin{equation*}
 \|(A_t^{1+iv}F)_{t>0}\|_{L^q(\mathcal A;\ell^\infty)}
 \le C_q |\Gamma(1+iv)|^{-1}\|F\|_{L^q(\mathcal A)}\,.
\end{equation*}
Fix $1<p<2$ and $\alpha=s+ib$ with $s>2/p-1$ and $b\in\R$. We must interpolate to the exponent $p$ and to the order $s+ib$ at the same point of the strip. For $1<q<p$, the relation $1/p=(1-\vartheta)/q+\vartheta/2$ determines
\begin{equation*}
 1-\vartheta=\frac{1/p-1/2}{1/q-1/2}\;\longrightarrow\;2/p-1
 \quad(q\downarrow1).
\end{equation*}
Choose $q$ close enough to one that $1-\vartheta<s$. To realize $s=(1-\vartheta)\cdot1+\vartheta h$ with a positive order $h$ at the $L^2$ boundary, set
\begin{equation*}
 h=\frac{s-1+\vartheta}{\vartheta}>0,\qquad
 \alpha(z)=1+ib+(h-1)z\quad(0\le\Re z\le1).
\end{equation*}
Then $\alpha(\vartheta)=\alpha$. For a finite set $E\subset(0,\infty)$ take $T_zF=(A_t^{\alpha(z)}F)_{t\in E}$. Since $\Re\alpha(z)\ge\min(1,h)>0$, the kernel in~\eqref{max:kernel} depends analytically on $z$ in $L^1(\R^2)$. Differentiation produces powers of $\log(1-|y/t|^2)$ in the integrable weight, together with derivatives of the reciprocal gamma factor. For $\Re\alpha>0$, polar integration gives
\begin{equation*}
 \|K_t^\alpha\|_{L^1(\R^2)}=\frac{2\pi}{|\Gamma(\alpha)|}\int_0^1
        (1-s^2)^{\Re\alpha-1}s\,ds
 =\frac{\pi}{\Re\alpha\,|\Gamma(\alpha)|}\,.
\end{equation*}
Thus $\|K_t^{\alpha(z)}\|_{L^1(\R^2)}$ grows at most exponentially in $|\Im z|$ by Stirling's formula, uniformly in $t$. For fixed finite $E$, the maximal norm of a tuple is bounded by the sum of its coordinate norms. The coordinate kernel estimates therefore give analyticity and the requisite growth of $T_zF$ for each finite $E$ on the dense space of Schwartz tensors from Section~\ref{pre:tracial}. The boundary maximal bounds below are uniform in $E$.

Multiply $T_z$ by $e^{(z-\vartheta)^2}$. For $z=j+iy$ with $j\in\{0,1\}$ and $y\in\R$, the Gaussian factor has modulus $e^{(j-\vartheta)^2-y^2}$, while the imaginary part of the order is $b+(h-1)y$. Thus $e^{(z-\vartheta)^2}T_z$ is uniformly bounded on each boundary line: by the ball estimate from $L^q(\mathcal A)$ to $L^q(\mathcal A;\ell^\infty(E))$ when $j=0$, and by~\eqref{max:vertical} from $L^2(\mathcal A)$ to $L^2(\mathcal A;\ell^\infty(E))$ when $j=1$. These bounds are independent of $y$ and $E$. To interpolate the target spaces we use the isometric identity~\cite[Lemma~2.7(i)]{hlw}
\begin{equation*}
 [L^q(\mathcal A;\ell^\infty(E)),L^2(\mathcal A;\ell^\infty(E))]_\vartheta
   = L^p(\mathcal A;\ell^\infty(E))\,.
\end{equation*}
The Gaussian equals one at $\vartheta$, giving $\|(A_t^\alpha F)_{t\in E}\|_{L^p(\mathcal A;\ell^\infty(E))}\le C_{p,\alpha}\|F\|_{L^p(\mathcal A)}$. Since $C_{p,\alpha}$ is independent of $E$, taking the supremum in
\eqref{max:finite} proves~\eqref{intro:maximal} for every $1<p<2$ in the asserted range.
\end{proof}

\begin{corollary}\label{max:growth}
Fix $2<p<\infty$ and $0\le L<\infty$. For $\alpha\in\mathbb C$ with $|\Im\alpha|\le L$ and $0<\varepsilon=\Re\alpha+g(p)\le1$, the least constant in~\eqref{intro:maximal}, uniform over all semifinite tracial coefficient algebras, is at most $C_{p,L}\varepsilon^{-1-2g(p)}$.
\end{corollary}
\begin{proof}
The proof of \eqref{intro:maximal} bounds its constant by a uniform multiple of $1+\sum_{j\ge1}D_j$. We therefore need to bound this sum as $\varepsilon\downarrow0$. The compact rectangle $-g(p)\le\Re\alpha\le1-g(p)$, $|\Im\alpha|\le L$ contains all orders under consideration, and on this rectangle the kernel bounds in~\eqref{max:low} and Lemma~\ref{max:local} are uniform.

To estimate $\sum_{j\ge1}D_j$, we compare it with an integral. For $j\le v\le j+1$, $2^{-\varepsilon j}(1+j)^{2g(p)}\le 2^{\varepsilon}2^{-\varepsilon v}(1+v)^{2g(p)}$. Integrating on each interval and setting $w=\varepsilon v$ gives
\[
 \sum_{j\ge1}D_j
 \le 2^\varepsilon\varepsilon^{-1-2g(p)}
       \int_{\varepsilon}^{\infty}2^{-w}(\varepsilon+w)^{2g(p)}\,dw
 \le C_p\varepsilon^{-1-2g(p)} .
\]
Together with \eqref{max:low} and \eqref{max:fullband} this gives the claimed bound.
\end{proof}

Liu, Shen, Song and Yan~\cite[Theorem~1.1]{lssy} proved that scalar strong type requires $\Re\alpha\ge\max(1/p-1/2,-1/p)$. On the boundary segment $2<p\le4$, Nowak, Roncal and Szarek~\cite[Corollaries 1.5 and 1.6]{nrs} obtained weak-type failure for the full-radius operator on radial functions in $L^p(\R^2)$ at real order $\alpha=1/p-1/2$. The following one-coordinate test works on the bounded interval $1\le t\le2$, allows every imaginary part, and also applies inside each fixed quantum torus.

\begin{proposition}\label{max:endpoint}
Let $2<p\le4$ and $\alpha=1/p-1/2+ib$, $b\in\R$. The scalar operator $F\mapsto\sup_{1\le t\le2}|A_t^\alpha F|$ fails to map $L^p(\R^2)$ to weak $L^p(\R^2)$. The corresponding strong maximal inequality also fails on each fixed two-dimensional quantum torus of period~$2\pi$.
\end{proposition}
\begin{proof}
We construct a function $F_\epsilon$ with $\|F_\epsilon\|_{L^p(\R^2)}=O((\log(1/\epsilon))^{1/p})$ for which $|A_{x_1}^\alpha F_\epsilon(x_1,x_2)|\gtrsim\log(1/\epsilon)$ on a fixed rectangle. Taking $t=x_1$ there lets us use the singularity in the one-dimensional marginal of $K_t^\alpha$.

Integrating~\eqref{max:kernel} in the second coordinate for $\Re\alpha>0$ yields the following formula for the marginal density $k_t^\alpha\in L^1(\R)$
\begin{equation*}
 k_t^\alpha(y_1)=\frac{\sqrt\pi}{t\,\Gamma(\alpha+1/2)}
                      (1-y_1^2/t^2)_+^{\alpha-1/2}\,.
\end{equation*}
For $|y_1|<t$, substitute $y_2=t(1-y_1^2/t^2)^{1/2}z$ in the kernel. The integral over $-1<z<1$ is $\int_{-1}^1(1-z^2)^{\alpha-1}dz=\sqrt\pi\,\Gamma(\alpha)/\Gamma(\alpha+1/2)$, giving the displayed density. For $|y_1|>t$ the support of $K_t^\alpha$ makes the marginal zero. The marginal formula extends to $\Re\alpha>-1/2$ by analytic continuation against smooth test functions, since its right-hand side is locally integrable there.

To match the phase of the marginal at $y_1=t$, fix $0<c<1/4$, choose $0<\epsilon<c/4$, and take a smooth function $\omega_\epsilon:\R\to[0,1]$ supported in $[\epsilon,c]$ and equal to one on $[2\epsilon,c/2]$. Define $h_\epsilon:\R\to\mathbb C$ by $h_\epsilon(v)=v^{-1/p-ib}\omega_\epsilon(v)$ for $v>0$, zero elsewhere. Then $\|h_\epsilon\|_{L^p(\R)}^p\le \int_\epsilon^c v^{-1}dv=\log(c/\epsilon)$. At $x_1=t\in[1,2]$, set $v=t-y_1$ in the marginal integral. Factoring $1-(1-v/t)^2=(2v/t)(1-v/(2t))$ gives
\begin{equation*}
 (k_t^\alpha*h_\epsilon)(t)
 =\frac{\sqrt\pi\,2^{\alpha-1/2}t^{-\alpha-1/2}}{\Gamma(\alpha+1/2)}
   \int_\epsilon^c \omega_\epsilon(v)(1-v/(2t))^{\alpha-1/2}\,\frac{dv}{v}\,.
\end{equation*}
The integral in this display equals $\log(1/\epsilon)+O_{\alpha,c}(1)$ uniformly for $t\in[1,2]$. Indeed, the factor $(1-v/(2t))^{\alpha-1/2}$ differs from one by $O_\alpha(v)$, and the transition regions $[\epsilon,2\epsilon]$ and $[c/2,c]$ of $\omega_\epsilon$ have bounded logarithmic measure. The prefactor has modulus bounded away from zero, so $|(k_t^\alpha*h_\epsilon)(t)|\gtrsim_{\alpha,c}\log(1/\epsilon)$ for all $t\in[1,2]$ once $\epsilon$ is sufficiently small.

Define $F_\epsilon:\R^2\to\mathbb C$ by $F_\epsilon(x_1,x_2)=h_\epsilon(x_1)\rho(x_2)$ with $\rho\in C_c^\infty(\R)$ identically one on $[-4,4]$. For $1\le x_1\le2$, $|x_2|\le1$, the support of $K_{x_1}^\alpha$ implies $A_{x_1}^\alpha F_\epsilon(x_1,x_2)=(k_{x_1}^\alpha*h_\epsilon)(x_1)$. Thus $\sup_{1\le t\le2}|A_t^\alpha F_\epsilon|\gtrsim\log(1/\epsilon)$ on $[1,2]\times[-1,1]$. Taking a level equal to half this lower bound and using $\|F_\epsilon\|_{L^p(\R^2)}\lesssim(\log(1/\epsilon))^{1/p}$ gives
\begin{equation*}
 \frac{\big\|\sup_{1\le t\le2}|A_t^\alpha F_\epsilon|\big\|_{L^{p,\infty}(\R^2)}}
      {\|F_\epsilon\|_{L^p(\R^2)}}
 \gtrsim(\log(1/\epsilon))^{1-1/p}.
\end{equation*}
The quotient tends to infinity as $\epsilon\downarrow0$.

For the torus assertion, use the $2\pi$-periodic extension $h_\epsilon^{\mathrm{per}}$ of $h_\epsilon$ in the abelian subalgebra generated by the first coordinate unitary $U_1$. Its trace is probability Haar measure $dy/(2\pi)$. The restricted convolution kernel for this measure is
\begin{equation*}
 k_{t,\mathrm{per}}^\alpha(y)=2\pi\sum_{\ell\in\mathbb Z}
                             k_t^\alpha(y+2\pi\ell)\,.
\end{equation*}
Indeed, its Fourier coefficient at $n\in\mathbb Z$ in probability Haar measure is $\widehat k_t^\alpha(n)=m_\alpha(t|n|)$, the multiplier on $U_1^n$. Since $k_t^\alpha$ is supported on $[-t,t]$, for $x_1=t\in[1,2]$ its convolution only uses the values of $h_\epsilon^{\mathrm{per}}$ on $[0,4]$. There $h_\epsilon^{\mathrm{per}}=h_\epsilon$ since $4<2\pi$ and $h_\epsilon$ is supported in $[\epsilon,c]$. Thus $|(k_{t,\mathrm{per}}^\alpha*h_\epsilon^{\mathrm{per}})(t)|\gtrsim_{\alpha,c}\log(1/\epsilon)$ for $1\le t\le2$.

For $1\le t\le2$, write $z_t=m_\alpha(t|D_\Theta|)h_\epsilon^{\mathrm{per}}$. If the strong bound held on the torus, the order characterization would give positive majorants $B_{\Re},B_{\Im}$ for both signs of $\Re z_t$ and $\Im z_t$, respectively, with $L^p(\mathcal Q_\Theta)$ norms bounded by a constant times $\|h_\epsilon^{\mathrm{per}}\|_{L^p(\mathcal Q_\Theta)}$. Let $E_1$ denote the trace-preserving conditional expectation given by averaging rotations of the second generator. Each $z_t$ belongs to this abelian subalgebra, so $|z_t|\le|\Re z_t|+|\Im z_t|$ is bounded by $E_1 B_{\Re}+E_1 B_{\Im}$, whose norm is at most $\|B_{\Re}\|_{L^p(\mathcal Q_\Theta)}+\|B_{\Im}\|_{L^p(\mathcal Q_\Theta)}$. Taking rational $t\in[1,2]$ gives these pointwise bounds on a single set of full measure. For almost every $x_1\in[1,2]$, continuity of $t\mapsto(k_{t,\mathrm{per}}^\alpha*h_\epsilon^{\mathrm{per}})(x_1)$ extends the bound to $t=x_1$, contradicting the marginal lower bound.
\end{proof}

\begin{corollary}\label{max:lowergrowth}
Fix $2<p\le4$ and $0\le L<\infty$. Set $\alpha=1/p-1/2+\varepsilon+ib$ for $b\in[-L,L]$ and $\varepsilon>0$. There are constants $c_{p,L},\varepsilon_0(p,L)>0$ such that, whenever $0<\varepsilon\le\varepsilon_0(p,L)$, the best strong maximal constant for $1\le t\le2$ on the scalar plane is at least $c_{p,L}\varepsilon^{-(1-1/p)}$. The same lower bound holds on each two-dimensional quantum torus of period~$2\pi$, with constants independent of the deformation parameter.
\end{corollary}
\begin{proof}
In the construction of $F_\epsilon$ in Proposition~\ref{max:endpoint}, take $c=1/8$ and $\epsilon=e^{-1/\varepsilon}$. Then $\|F_\epsilon\|_{L^p(\R^2)}\lesssim\varepsilon^{-1/p}$. We shall show that $|A_{x_1}^\alpha F_\epsilon(x_1,x_2)|\gtrsim_{p,L}\varepsilon^{-1}$ on $[1,2]\times[-1,1]$.

In the marginal integral the product of the factor $v^{-1/p-ib}$ from $h_\epsilon$ and the kernel factor $v^{\alpha-1/2}$ is now $v^{\varepsilon-1}$. The function $\omega_\epsilon$ equals one on $[2\epsilon,c/2]$, where
\begin{equation*}
 \int_{2\epsilon}^{c/2} v^{\varepsilon-1}\,dv
   =\frac{(c/2)^\varepsilon-(2\epsilon)^\varepsilon}{\varepsilon}
   \gtrsim \varepsilon^{-1}
\end{equation*}
for small $\varepsilon$, since $(2\epsilon)^\varepsilon\to e^{-1}$. Uniformly for $1\le t\le2$ and $|b|\le L$, we have $(1-v/(2t))^{\alpha-1/2}=1+O_{p,L}(v)$ on $(0,c)$. Replacing $(1-v/(2t))^{\alpha-1/2}$ by one in the marginal integral changes its value by $O_{p,L}(1)$. The prefactor $\sqrt\pi\,2^{\alpha-1/2}t^{-\alpha-1/2}/\Gamma(\alpha+1/2)$ stays bounded away from zero on the indicated compact set of orders since $\Re(\alpha+1/2)\ge1/p$. On the rectangle $1\le x_1\le2$, $|x_2|\le1$ $|A_{x_1}^\alpha F_\epsilon(x_1,x_2)|$ is therefore at least $c_{p,L}\varepsilon^{-1}$. Taking the $L^p$ norm on this rectangle and dividing by $\|F_\epsilon\|_{L^p(\R^2)}$ gives the asserted power $\varepsilon^{-(1-1/p)}$. Convolution with $h_\epsilon^{\mathrm{per}}$ on the circle gives the same marginal integral on the interval $1\le x_1\le2$. Applying the conditional expectation from the proof of Proposition~\ref{max:endpoint} transfers this lower bound to the quantum torus maximal norm.
\end{proof}

For $2<p\le4$, Corollaries~\ref{max:growth} and~\ref{max:lowergrowth} leave a gap between the powers $1-1/p$ and $1+2g(p)$ in the growth of the best uniform constant as $\varepsilon\downarrow0$. On the other segment $p>4$, the full-radius scalar operator is of strong type at the real order $\alpha=-1/p$ on radial functions~\cite[Corollary 1.5]{nrs}. For arbitrary $F\in L^p(\R^2)$ our full-radius summation in~\eqref{max:fullband} requires $\Re\alpha>-1/p$.
\subsection{Transference to trace-preserving actions}\label{q:secmaximal}

We prove the quantum maximal estimate in Theorem~\ref{thm:maximal}\ref{thm:q-max} by first establishing the orbit estimate~\eqref{q:maximal} of Theorem~\ref{q:action} for general trace-preserving actions.

Fix a trace-preserving action \(\beta\) of \(\R^2\) on a semifinite tracial
von Neumann algebra \(\mathcal D\), and write
\(\mathcal S_{\mathcal D}=L^\infty(\R^2)\bar\otimes\mathcal D\) as in
Section~\ref{q:secspacetime}. To transfer \eqref{intro:maximal}, we
seek one positive majorant for both signs of the real parts at all
radii, and another for both signs of the imaginary parts. For a finite set
of radii, compact support of the kernels lets us cut off the orbit
outside one large square without changing any of its convolutions on
a smaller square. Averaging a Euclidean majorant over the smaller
square then gives a majorant in \(L^p(\mathcal D)\). The finite-family
bound will be independent of the radii. This localization was used for
measure averages in \cite[Proposition 8.1]{lailiu}.

\medskip
We define the orbit average for a compactly supported distribution first on elements with smooth $L^p$ orbits. For \(1<p<\infty\), put
\[
 \mathcal D_{\beta,p}^{\infty}
 =\{x\in L^p(\mathcal D):z\mapsto\beta_zx
                     \text{ belongs to } C^\infty(\R^2;L^p(\mathcal D))\}.
\]
For each multi-index \(\nu=(\nu_1,\nu_2)\in\{0,1,\ldots\}^2\), the
derivative \(\partial_z^\nu=\partial_{z_1}^{\nu_1}\partial_{z_2}^{\nu_2}\),
of order \(|\nu|=\nu_1+\nu_2\), must exist and be continuous
in the \(L^p(\mathcal D)\) norm. A compactly supported scalar distribution \(K\)
can act on the orbit of each such \(x\). In detail, choose a closed ball \(D_0\)
whose interior contains \(\operatorname{supp}K\). By continuity of a
distribution on test functions supported in a fixed compact set, there
are \(m\in\{0,1,\ldots\}\) and \(C_K<\infty\) such that
\[
 |\langle K,f\rangle|\le C_K\max_{|\nu|\le m}\sup_{y\in D_0}
                    |\partial_y^\nu f(y)|\qquad(f\in C^\infty(\R^2)).
\]
Here \(\langle K,f\rangle=\langle K,\theta f\rangle\) for a fixed \(\theta\in C_c^\infty(D_0^\circ)\) equal to one near \(\operatorname{supp}K\). The constant \(C_K\) absorbs the derivatives of \(\theta\) through order \(m\). Pairings with scalar test functions are complex linear.

For a complex Banach space \(V\), first take a finite sum \(F=\sum_i f_i v_i\) with \(f_i\in C^\infty(\R^2;\mathbb C)\) and \(v_i\in V\). Define
\(\langle K,F\rangle=\sum_i\langle K,f_i\rangle v_i\). Testing with
linear functionals in the unit ball of \(V^*\) gives
\begin{equation}\label{q:order}
 \|\langle K,F\rangle\|_V
 \le C_K\max_{|\nu|\le m}\sup_{y\in D_0}\|\partial_y^\nu F(y)\|_V .
\end{equation}
For smooth \(V\)-valued \(F\) these finite sums approximate \(F\) in
the seminorm on the right. Indeed, localize \(F\) to a neighborhood of
\(D_0\), mollify in \(y\), and approximate the mollifying integral and
its derivatives through order \(m\) uniformly on \(D_0\) by Riemann sums
of scalar translates of the mollifier times fixed values of \(F\).
Thus \eqref{q:order} defines the pairing \(F\mapsto\langle K,F\rangle\) from \(C^\infty(\R^2;V)\) to \(V\). Its value depends only on \(F\) near \(\operatorname{supp}K\). For any bounded linear map \(T:V\to W\) into a complex Banach space \(W\), \(T\langle K,F\rangle=\langle K,TF\rangle\) by finite-sum approximation.
In particular \(y\mapsto\beta_{-y}x\) is an admissible test function for
\(V=L^p(\mathcal D)\), \(x\in\mathcal D_{\beta,p}^{\infty}\).

\begin{proposition}\label{q:maxtransfer}
Let \(1<p<\infty\), and let \(\beta\) be an action of \(\R^2\) by normal
trace-preserving \(*\)-automorphisms of \(\mathcal D\) with norm-continuous
\(L^p\) orbits. Let \((K_t)_{t>0}\) be compactly supported scalar
distributions on \(\R^2\). Suppose there is a constant \(C\) such that for
every finite set \(E\subset(0,\infty)\),
\[
 \|(K_t*F)_{t\in E}\|_{L^p(\mathcal S_{\mathcal D};\ell^\infty(E))}
 \le C\|F\|_{L^p(\mathcal S_{\mathcal D})}
 \quad
 (F\in C_c^\infty(\R^2)\otimes L^p(\mathcal D)).
\]
For each \(t>0\) define \(Q_t:\mathcal D_{\beta,p}^{\infty}\to L^p(\mathcal D)\) by \(Q_tx=\langle K_t(y),\beta_{-y}x\rangle\).
Then \(x\mapsto(Q_tx)_{t>0}\) extends to all of \(L^p(\mathcal D)\) with
\[
 \|(Q_tx)_{t>0}\|_{L^p(\mathcal D;\ell^\infty)}
 \le 2C\|x\|_{L^p(\mathcal D)} .
\]
\end{proposition}

\begin{proof}
For fixed \(x\in\mathcal D_{\beta,p}^{\infty}\), we construct a positive majorant for both signs of \(\operatorname{Re} Q_t x\), simultaneously for all \(t>0\), with norm at most \(C\|x\|_{L^p(\mathcal D)}\). We first localize to a finite set of radii and average a Euclidean majorant. Weak compactness then gives one majorant for all radii. Applying the construction to the imaginary parts and approximating arbitrary \(L^p\) elements by smooth orbits will complete the proof.

\medskip
To localize, fix a nonempty finite set \(E\subset(0,\infty)\). Write \(B_r=[-r,r]^2\) for \(r>0\) and choose \(S>0\) such that every \(K_t\), \(t\in E\), is supported in \(B_S\). For \(R>0\) choose a smooth \(\rho_R:\R^2\to[0,1]\), identically one on \(B_{R+S+1}\) and
supported in \(B_{R+S+2}\). The function \(F_R(z)=\rho_R(z)\beta_zx\)
is smooth and compactly supported with values in \(L^p(\mathcal D)\).
The hypothesis extends the map \(F\mapsto(K_t*F)_{t\in E}\) from finite tensors to all of \(L^p(\mathcal S_{\mathcal D})\), with values in \(L^p(\mathcal S_{\mathcal D};\ell^\infty(E))\). Indeed, if a sequence of finite tensors \(F_n\) converges to \(F\) in \(L^p(\mathcal S_{\mathcal D})\), the assumed bound makes \((K_t*F_n)_{t\in E}\) Cauchy in \(L^p(\mathcal S_{\mathcal D};\ell^\infty(E))\) as \(n\to\infty\). Completeness gives the extension with norm at most \(C\). Projection onto one coordinate has norm at most one in \eqref{max:norm} by H\"older's inequality for the trace, so each coordinate agrees with the bounded extension of \(F\mapsto K_t*F\).
For smooth \(L^p(\mathcal D)\)-valued functions \(F\), write \((K_t*F)(z)=\langle K_t(y),F(z-y)\rangle\). We verify that this formula gives the coordinates of the bounded extension at \(F_R\). Let \(m\) be the largest of the orders used in~\eqref{q:order} for \(t\in E\). The mollification and finite-sum approximation used there provide finite tensors \(F_j\), \(j\in\mathbb N\), with supports in a common compact set \(B\) containing \(\operatorname{supp}F_R\), and with derivatives through order \(m\) converging uniformly to those of \(F_R\). In particular, \(F_j\to F_R\) in spatial \(L^p\). For each \(t\in E\), applying~\eqref{q:order} to \(y\mapsto(F_j-F_R)(z-y)\) gives
\begin{equation*}
 \sup_{z\in\R^2}\|(K_t*(F_j-F_R))(z)\|_{L^p(\mathcal D)}
 \le C_t\max_{|\nu|\le m}\sup_{w\in\R^2}
       \|\partial^\nu(F_j-F_R)(w)\|_{L^p(\mathcal D)}\longrightarrow0\,.
\end{equation*}
The functions \(K_t*F_j\) and \(K_t*F_R\) are supported in \(B+B_S=\{u+v:u\in B,\ v\in B_S\}\). Consequently \(K_t*F_j\) converges in spatial \(L^p\) to \(K_t*F_R\) as defined by the pairing. The uniform convergence of the continuous functions \(K_t*F_j\) shows that this pairing is the continuous representative of the coordinate of the extension.
Thus for \(z\in B_R\) and \(t\in E\),
\[
 (K_t*F_R)(z)=\langle K_t(y),\beta_{z-y}x\rangle
             =\beta_z(Q_tx),
\]
since \(\rho_R(z-y)=1\) on a neighborhood of \(\operatorname{supp}K_t\).

\medskip
On the square $B_R$ we can now use a positive majorant for the Euclidean convolutions to bound the real parts of $\beta_z(Q_tx)$ for $t\in E$. Fix \(\delta>0\). The adjoint is an isometry on the factorization norm
\eqref{max:norm}, so taking real parts is contractive. By the order
description of this norm there is \(b_R\in L^p(\mathcal S_{\mathcal D})_+\)
with
\[
 \pm\operatorname{Re}(K_t*F_R)\le b_R\quad(t\in E),
 \qquad
 \|b_R\|_{L^p(\mathcal S_{\mathcal D})}\le(C+\delta)|B_{R+S+2}|^{1/p}\|x\|_{L^p(\mathcal D)} .
\]
Here \(\|F_R\|_{L^p(\mathcal S_{\mathcal D})}^p=\|x\|_{L^p(\mathcal D)}^p\int |\rho_R(z)|^p\,dz
\le |B_{R+S+2}|\|x\|_{L^p(\mathcal D)}^p\), since each \(\beta_z\) preserves the
\(L^p(\mathcal D)\) norm. For \(x=0\) take \(b_R=0\). Under \eqref{pre:fubini}, the positive cone of
\(L^p(\mathcal S_{\mathcal D})\) consists of fields positive almost
everywhere (functional calculus for the negative part is pointwise).
Since \(E\) is finite, the exceptional null sets for the displayed
inequalities have a common null union. Thus they hold simultaneously
for almost every \(z\in B_R\).
For these \(z\), applying \(\beta_{-z}\) and using the formula for \((K_t*F_R)(z)\) gives \(\beta_{-z}(b_R(z))\ge\pm\operatorname{Re}Q_tx\) for every \(t\in E\). The map \(z\mapsto\beta_{-z}(b_R(z))\) is strongly measurable. Indeed, approximate \(b_R\) pointwise almost everywhere in \(L^p(\mathcal D)\) by simple functions. Applying \(\beta_{-z}\) to the finitely many values of each simple function gives a strongly measurable function by orbit continuity. The approximations still converge after applying \(\beta_{-z}\) because it preserves the \(L^p\) norm. H\"older's inequality on \(B_R\) gives
\begin{equation*}
 \int_{B_R}\|\beta_{-z}(b_R(z))\|_{L^p(\mathcal D)}\,dz
 \le |B_R|^{1-1/p}\|b_R\|_{L^p(\mathcal S_{\mathcal D})}\,.
\end{equation*}
We may therefore take the Bochner average
\begin{equation*}
 b_{R,\delta}=\frac{1}{|B_R|}\int_{B_R}\beta_{-z}(b_R(z))\,dz\,.
\end{equation*}
For every \(t\in E\), the fields \(\beta_{-z}(b_R(z))\) and \(\beta_{-z}(b_R(z))\pm\operatorname{Re}Q_tx\) are positive for almost every \(z\in B_R\). Their averages lie in the closed convex cone \(L^p(\mathcal D)_+\), so \(b_{R,\delta}\ge0\) and \(b_{R,\delta}\ge\pm\operatorname{Re}Q_tx\). By the integral estimate and the bound for \(b_R\),
\begin{equation*}
 \|b_{R,\delta}\|_{L^p(\mathcal D)}\le (C+\delta)
 \left(\frac{|B_{R+S+2}|}{|B_R|}\right)^{1/p}\|x\|_{L^p(\mathcal D)}\,.
\end{equation*}

\medskip
We pass from the majorants $b_{R,\delta}$ for finite $E$ to a single majorant for all $t>0$. For every \(t>0\), the set
\begin{equation*}
 \mathcal C_t=\{b\in L^p(\mathcal D)_+:
              b\pm\operatorname{Re}Q_tx\in L^p(\mathcal D)_+\}
\end{equation*}
is weakly closed because the positive cone is norm closed and convex. For the fixed finite \(E\), take \(R=n\) and \(\delta=1/n\) with \(n\in\mathbb N\). The volume ratio tends to one. The sequence \(b_{n,1/n}\) is bounded, so reflexivity of \(L^p(\mathcal D)\) gives a weakly convergent subnet. Its limit belongs to \(\bigcap_{t\in E}\mathcal C_t\) and has norm at most \(C\|x\|_{L^p(\mathcal D)}\) by weak lower semicontinuity. Since \(E\) was arbitrary, the intersections of \(\mathcal C_t\) with the weakly compact ball of radius \(C\|x\|_{L^p(\mathcal D)}\) have the finite intersection property. A point in their intersection is a majorant for \(\pm\operatorname{Re}Q_tx\) for all \(t>0\). For the kernels \(-iK_t\), multiplication by \(-i\) preserves the Euclidean family norm in~\eqref{max:norm}, and \(\operatorname{Re}(-iQ_tx)=\operatorname{Im}Q_tx\). The real-part bound just proved therefore gives a positive majorant for \(\pm\operatorname{Im}Q_tx\), for all \(t>0\), of norm at most \(C\|x\|_{L^p(\mathcal D)}\). The triangle inequality in~\eqref{max:norm} gives the asserted factor \(2C\).

\medskip
To pass to all \(x\in L^p(\mathcal D)\), take a nonnegative
\(\rho\in C_c^\infty(\R^2)\) with \(\int\rho=1\), and write
\(\rho_\delta(y)=\delta^{-2}\rho(y/\delta)\). For
\(x_\delta=\int\rho_\delta(y)\beta_y x\,dy\),
\[
 \beta_z x_\delta=\int \rho_\delta(y-z)\beta_yx\,dy,\qquad
 \|x_\delta-x\|_{L^p(\mathcal D)}\le\int \rho_\delta(y)\|\beta_yx-x\|_{L^p(\mathcal D)}\,dy .
\]
Differentiating \(z\mapsto\rho_\delta(\cdot-z)\) in \(L^1\) shows that \(x_\delta\in\mathcal D_{\beta,p}^{\infty}\). The integral bounding \(\|x_\delta-x\|_{L^p(\mathcal D)}\) tends to zero by continuity of \(y\mapsto\beta_yx\) at \(y=0\). Applying the bound for smooth orbits to \(x_\delta-x_{\delta'}\) gives
\begin{equation*}
 \|(Q_tx_\delta-Q_tx_{\delta'})_{t>0}\|_{L^p(\mathcal D;\ell^\infty)}
 \le 2C\|x_\delta-x_{\delta'}\|_{L^p(\mathcal D)}\longrightarrow0
 \qquad(\delta,\delta'\downarrow0)\,.
\end{equation*}
Completeness of the space in~\eqref{max:norm} gives the required extension to \(x\), with the same bound \(2C\|x\|_{L^p(\mathcal D)}\).
\end{proof}

\begin{theorem}\label{q:action}
Fix \(p\in(1,\infty)\) and \(\alpha\in\mathbb C\) with \(\operatorname{Re}\alpha>a(p)\). Let
\(\beta\) be an action of \(\R^2\) by normal trace-preserving
\(*\)-automorphisms on a semifinite tracial von Neumann algebra
\(\mathcal D\) with norm-continuous \(L^p\) orbits. The generalized circular means defined on elements with smooth
\(L^p\) orbit by
\[
 A_{t,\beta}^{\alpha}x=\langle K_t^\alpha(y),\beta_{-y}x\rangle
 \qquad (t>0)
\]
extend to a bounded family map \(L^p(\mathcal D)\to L^p(\mathcal D;\ell^\infty((0,\infty)))\) satisfying
\begin{equation}\label{q:maximal}
 \|(A_{t,\beta}^{\alpha}x)_{t>0}\|_{L^p(\mathcal D;\ell^\infty)}
 \le C_{p,\alpha}\|x\|_{L^p(\mathcal D)}\qquad (x\in L^p(\mathcal D)).
\end{equation}
The constant is uniform in \(\mathcal D\) and \(\beta\). For \(2<p<\infty\) and any \(L\in[0,\infty)\), if
\(|\operatorname{Im}\alpha|\le L\) and
\(0<\varepsilon=\operatorname{Re}\alpha+g(p)\le1\) one may take
\(C_{p,\alpha}\le C_{p,L}\varepsilon^{-1-2g(p)}\).
\end{theorem}

\begin{proof}
For every \(t>0\) the distribution \(K_t^\alpha\) in \eqref{max:kernel}
is supported in the closed disk of radius \(t\). For a finite set of radii all these disks lie in a common square. Theorem~\ref{thm:maximal}, with coefficient algebra \(\mathcal D\), gives a bound independent of the chosen radii for Euclidean convolutions by \(K_t^\alpha\). These are the support and norm hypotheses of Proposition~\ref{q:maxtransfer}. Applying it to the given action yields \eqref{q:maximal} first on smooth orbits and then on all of \(L^p(\mathcal D)\) by orbit mollification. The proposition costs at most a factor of two in the norm, so for \(2<p<\infty\) it also transfers the dependence on \(\varepsilon\) from Corollary~\ref{max:growth}.
\end{proof}

\medskip
To deduce part~\ref{thm:q-max} of Theorem~\ref{thm:maximal} from
\eqref{q:maximal}, it remains to identify the orbit averages with
the Fourier multipliers \(m_\alpha(t|D_\Theta|)\).
For the action \(\beta=\mathrm{id}_{\mathcal C}\otimes\gamma\) in
Section~\ref{q:models}, the Fourier transform of \(K_t^\alpha\) determines
the orbit multiplier explicitly. On the quantum torus, testing \(K_t^\alpha\) against
\(y\mapsto e^{-iy\cdot n}\) gives
\[
 A_{t,\beta}^\alpha(c\otimes U^n)
   =m_\alpha(t|n|)\,c\otimes U^n,\qquad n\in\mathbb Z^2,
\]
for bounded \(c\in\mathcal C\) whose support projections have finite trace.

On the Moyal plane with \(1<p\le2\), we have \(\Re\alpha>a(p)\ge0\), so \(K_t^\alpha\in L^1(\R^2)\) for each \(t>0\). In the notation of Section~\ref{q:secspacetime}, the Bochner operator \(A_K\) associated to any \(K\in L^1(\R^2)\) is given by the integral in~\eqref{q:orbit}. For \(K=K_t^\alpha\) this operator agrees with the distributional pairing on smooth orbits, and hence with its bounded \(L^p\) extension. For \(x\in L^p(\mathcal D)\cap L^2(\mathcal D)\), the integral in the intersection norm agrees with its values in each space. Applying the Fourier formula for \(A_K\) and the identity \(\widehat K_t^\alpha(\xi)=m_\alpha(t|\xi|)\) gives
\begin{equation*}
 \widehat{A_{t,\beta}^\alpha x}_\Theta(\xi)
 =m_\alpha(t|\xi|)\widehat x_\Theta(\xi)
 \qquad\text{for almost every }\xi\in\R^2\,.
\end{equation*}
The intersection contains all bounded elements whose support projections have finite trace and is therefore dense in \(L^p(\mathcal D)\). This identifies the \(L^p\) extension with the claimed multiplier.

For \(2<p<\infty\), the permitted orders can have \(\Re\alpha\le0\), so we use the distributional pairing. Fix \(t>0\) and \(f\in\mathscr S(\R^2)\). Put \(V=L^2(\mathcal Q_\Theta)\cap L^p(\mathcal Q_\Theta)\) with norm \(\|u\|_V=\|u\|_{L^2(\mathcal Q_\Theta)}+\|u\|_{L^p(\mathcal Q_\Theta)}\). Plancherel's identity and \(\|\lambda_\Theta(g)\|_{L^\infty(\mathcal Q_\Theta)}\le\|g\|_{L^1(\R^2)}\), followed by interpolation, make \(\lambda_\Theta:L^1(\R^2)\cap L^2(\R^2)\to V\) bounded. The map \(y\mapsto e^{-iy\cdot\xi}f(\xi)\) is smooth with values in \(L^1\cap L^2\). Its derivative of order \(\nu\) is \((-i\xi)^\nu e^{-iy\cdot\xi}f(\xi)\), by dominated convergence in both norms. Thus the orbit \(F(y)=\lambda_\Theta(e^{-iy\cdot\xi}f(\xi))\) is smooth with values in \(V\).

Take a nonnegative \(\rho\in C_c^\infty(\R^2)\) with \(\int\rho=1\), and for \(\epsilon>0\) set \(\rho_\epsilon(y)=\epsilon^{-2}\rho(y/\epsilon)\), \(\widetilde\rho_\epsilon(y)=\rho_\epsilon(-y)\) and \(K_\epsilon=K_t^\alpha*\rho_\epsilon\). Integrating \(K_\epsilon(y)=\langle K_t^\alpha(s),\rho_\epsilon(y-s)\rangle\) against \(F(y)\) gives
\begin{equation*}
 \int_{\R^2}K_\epsilon(y)F(y)\,dy
 =\langle K_t^\alpha,\widetilde\rho_\epsilon*F\rangle
 \longrightarrow\langle K_t^\alpha,F\rangle\qquad\text{in }V\,.
\end{equation*}
For the equality, the integral can be restricted to a common compact set in \(y\) for all \(s\) in the ball \(D_0\) used in~\eqref{q:order}. Riemann sums for the integral with kernel \(\rho_\epsilon(y-s)\) converge uniformly in \(s\in D_0\) with derivatives through the order \(m\) in~\eqref{q:order}, allowing the interchange. Applying~\eqref{q:order} to \(\widetilde\rho_\epsilon*F-F\) gives convergence of the pairings in \(V\), since all derivatives of this difference through order \(m\) tend uniformly to zero on \(D_0\).

In the trace Hilbert space, the Fourier coefficient of \(\int K_\epsilon(y)F(y)\,dy\) is \(\widehat\rho(\epsilon\xi)m_\alpha(t|\xi|)f(\xi)\). Here \(|\widehat\rho(\epsilon\xi)|\le1\) and \(\widehat\rho(\epsilon\xi)\to1\). The Fourier transform of the compactly supported distribution \(K_t^\alpha\) grows at most polynomially by~\eqref{q:order} applied to characters, so \(m_\alpha(t|\xi|)f(\xi)\) is square integrable. Dominated convergence identifies the Fourier coefficient of the \(V\)-limit as \(m_\alpha(t|\xi|)f(\xi)\). Tensoring these identities
with bounded \(c\in\mathcal C\) whose support projections have finite trace and using density identifies
\(A_{t,\beta}^\alpha\) with \(m_\alpha(t|D_\Theta|)\) in each amplified
model. Theorem~\ref{q:action} now proves Theorem~\ref{thm:maximal}\ref{thm:q-max}.

\subsection{Bilateral almost uniform convergence}\label{q:secconvergence}

To prove the small-radius convergence in Corollary~\ref{q:bau}, we use the following dilation identity to identify the limit on smooth orbits from the mass of \(K_1^\alpha\).
For every \(t>0,\ \alpha\in\mathbb C\), and scalar
\(\varphi\in C_c^\infty(\R^2)\),
\begin{equation}\label{q:dilation}
 \langle K_t^\alpha,\varphi\rangle
 =\langle K_1^\alpha,y\mapsto\varphi(ty)\rangle .
\end{equation}
For \(\Re\alpha>0\), substitute \(y=tv\) in \eqref{max:kernel}.
The Jacobian \(t^2\) cancels the dilation factor \(t^{-2}\).
For fixed \(\varphi,t\) both pairings are entire in \(\alpha\) by the
continuation following \eqref{max:kernel}. The identity theorem proves
\eqref{q:dilation} at every complex order. Cutting off on a neighborhood
of the supports and applying \eqref{q:order} extends the scaling formula
to smooth Banach-valued test functions. In particular,
\(\langle K_1^\alpha,1\rangle=\widehat K_1^\alpha(0)=m_\alpha(0)\).

For normalized geometric averages, the two authors of this paper used a maximal
inequality to obtain bilateral almost uniform convergence along an action
\cite[Proposition 8.1 and Corollary 8.3]{lailiu}. For the complex-order
distributions \(K_t^\alpha\), \eqref{q:dilation} supplies convergence on
smooth orbits, and \eqref{q:maximal} controls the error for every \(x\in L^p(\mathcal D)\).

\begin{corollary}\label{q:bau}
Fix \(p\in(1,\infty)\) and \(\alpha\in\mathbb C\) with \(\operatorname{Re}\alpha>a(p)\). Let \(\beta\)
be an action of \(\R^2\) by normal trace-preserving \(*\)-automorphisms on a semifinite tracial algebra
\(\mathcal D\) with \(L^p\)-continuous orbits. As \(t\downarrow0\),
every \(x\in L^p(\mathcal D)\) satisfies
\[
 A_{t,\beta}^\alpha x\longrightarrow \frac{\pi}{\Gamma(1+\alpha)}x
 \qquad\text{bilaterally almost uniformly}.
\]
\end{corollary}
\begin{proof}
Fix \(x\in L^p(\mathcal D)\) and \(\eta>0\). Put \(c_\alpha=m_\alpha(0)\) and denote the trace on \(\mathcal D\) by \(\operatorname{tr}_{\mathcal D}\). We seek a projection \(e\in\mathcal D\) with \(\operatorname{tr}_{\mathcal D}(1-e)<\eta\) and \(\|e(A_{t,\beta}^\alpha x-c_\alpha x)e\|_{L^\infty(\mathcal D)}\to0\) as \(t\downarrow0\). For approximants \(x_n\) to be chosen, put \(h_n=x-x_n\) for \(n\in\mathbb N\). For any projection \(e\in\mathcal D\) and \(t>0\) we have the identity in \(L^p(\mathcal D)\)
\begin{equation*}
 e(A_{t,\beta}^\alpha x-c_\alpha x)e
 =e(A_{t,\beta}^\alpha x_n-c_\alpha x_n)e
  +e(A_{t,\beta}^\alpha h_n-c_\alpha h_n)e\,.
\end{equation*}
We first prove operator-norm convergence on a dense class from which we can choose the \(x_n\). For the error we use~\eqref{q:maximal} to find one projection \(e\) with \(\operatorname{tr}_{\mathcal D}(1-e)<\eta\) such that, for some \(n_0\in\mathbb N\),
\begin{equation*}
 \sup_{t>0}\|e(A_{t,\beta}^\alpha h_n-c_\alpha h_n)e\|_{L^\infty(\mathcal D)}
 \le 2^{1-n}\qquad(n\ge n_0)\,.
\end{equation*}

Put \(V=\mathcal D\cap L^p(\mathcal D)\) with norm \(\|a\|_V=\max(\|a\|_{L^\infty(\mathcal D)},\|a\|_{L^p(\mathcal D)})\). We construct the dense class by smoothing the orbits of elements of \(V\). For \(\zeta\in C_c^\infty(\R^2)\) and \(a\in V\), put \(S_\zeta a=\int_{\R^2}\zeta(y)\beta_y a\,dy\), initially defined by the \(L^p\) Bochner integral. In addition to
\(\|S_\zeta a\|_{L^p(\mathcal D)}\le\|\zeta\|_{L^1(\R^2)}\|a\|_{L^p(\mathcal D)}\) we have
\(\|S_\zeta a\|_{L^\infty(\mathcal D)}\le\|\zeta\|_{L^1(\R^2)}\|a\|_{L^\infty(\mathcal D)}\). Indeed, Riemann
sums converge in \(L^p\) and have operator norms bounded by
\(\|a\|_{L^\infty(\mathcal D)}\) times the scalar Riemann sums of \(|\zeta|\). The
operator-norm bound for \(S_\zeta a\) follows by testing these sums against bounded \(d\in\mathcal D\) whose support projections have finite trace. We pass to the limit using \(L^p\)-convergence and \(d\in L^{p'}\), where \(1/p+1/p'=1\). The resulting bound is \(|\operatorname{tr}_{\mathcal D}(d^*S_\zeta a)|\le\|\zeta\|_{L^1(\R^2)}\|a\|_{L^\infty(\mathcal D)}\|d\|_{L^1(\mathcal D)}\). To recover the operator norm, let
\(S_\zeta a=u|S_\zeta a|\), put \(M=\|\zeta\|_{L^1(\R^2)}\|a\|_{L^\infty(\mathcal D)}\), and
test with \(d=ue_\varepsilon\), where
\(e_\varepsilon=1_{(M+\varepsilon,\infty)}(|S_\zeta a|)\).
This spectral projection has finite trace by Chebyshev's inequality for the spectral measure of
\(|S_\zeta a|\), so the pairing bound gives
\((M+\varepsilon)\operatorname{tr}_{\mathcal D}(e_\varepsilon)
 \le\operatorname{tr}_{\mathcal D}(e_\varepsilon|S_\zeta a|)\le M\operatorname{tr}_{\mathcal D}(e_\varepsilon)\).
Thus \(e_\varepsilon=0\) for every \(\varepsilon>0\), proving \(\|S_\zeta a\|_{L^\infty(\mathcal D)}\le M\).

For fixed \(a\in V\), taking the maximum of the \(L^p(\mathcal D)\) and \(L^\infty(\mathcal D)\) estimates gives \(\|S_\zeta a\|_V\le\|\zeta\|_{L^1(\R^2)}\|a\|_V\). By density, \(\zeta\mapsto S_\zeta a\) extends to a bounded linear map \(L^1(\R^2)\to V\). Set \(x_0=S_\zeta a\) with \(\zeta\in C_c^\infty(\R^2)\). The identity \(\beta_zx_0=S_{\zeta(\cdot-z)}a\) and differentiability of \(z\mapsto\zeta(\cdot-z)\) of all orders in \(L^1(\R^2)\) show that the orbit of \(x_0\) is smooth in \(V\). The set of such \(x_0\) is dense in \(L^p(\mathcal D)\). Indeed, \(V\) is dense by Section~\ref{pre:tracial}, and we can take \(\rho_\delta(y)=\delta^{-2}\rho(y/\delta)\) for \(\delta>0\) and nonnegative \(\rho\in C_c^\infty(\R^2)\) of integral one. For \(a\in V\), orbit continuity gives
\begin{equation*}
 \|S_{\rho_\delta}a-a\|_{L^p(\mathcal D)}
 \le \int_{\R^2}\rho(y)\|\beta_{\delta y}a-a\|_{L^p(\mathcal D)}\,dy
 \longrightarrow0\qquad(\delta\downarrow0)\,.
\end{equation*}

For such \(x_0\), we estimate \(A_{t,\beta}^\alpha x_0-c_\alpha x_0\) in \(V\). Formula \eqref{q:dilation} for the smooth orbit and the mass of \(K_1^\alpha\) give
\[
 A_{t,\beta}^\alpha x_0-c_\alpha x_0
 =\langle K_1^\alpha(y),\beta_{-ty}x_0-x_0\rangle .
\]
The \(V\)-valued pairing agrees with the defining \(L^p\) pairing by the continuous inclusion \(V\hookrightarrow L^p\). By \eqref{q:order}, with \(K=K_1^\alpha\) and \(D_0=\{|y|\le2\}\), it suffices to show that all derivatives in \(y\) through order \(m\) of \(\beta_{-ty}x_0-x_0\) tend to zero uniformly on \(D_0\) in \(V\). The identity \(\beta_{-ty}x_0=S_{\zeta(\cdot+ty)}a\) gives
\begin{equation*}
 \sup_{y\in D_0}\|\beta_{-ty}x_0-x_0\|_V
 \le \|a\|_V\sup_{y\in D_0}\|\zeta(\cdot+ty)-\zeta\|_{L^1(\R^2)}
 \longrightarrow0 \qquad(t\downarrow0)
\end{equation*}
by continuity of translation in \(L^1(\R^2)\), since \(|ty|\le2t\) for \(y\in D_0\). For \(1\le|\nu|\le m\), differentiating \(\zeta(\cdot+ty)\) gives
\begin{equation*}
 \sup_{y\in D_0}\|\partial_y^\nu\beta_{-ty}x_0\|_V
 \le t^{|\nu|}\|a\|_V\,\|\partial^\nu\zeta\|_{L^1(\R^2)}
 \longrightarrow0 \qquad(t\downarrow0)\,.
\end{equation*}
Hence \(\|A_{t,\beta}^\alpha x_0-c_\alpha x_0\|_{L^\infty(\mathcal D)}\to0\).

We now choose \(x_n=S_{\zeta_n}a_n\) with \(a_n\in V\) and \(\zeta_n\in C_c^\infty(\R^2)\), so that \(\|h_n\|_{L^p(\mathcal D)}=\|x-x_n\|_{L^p(\mathcal D)}\le 2^{-2n}\) for every \(n\in\mathbb N\). For \(t>0\) write \(Y_{n,t}=A_{t,\beta}^\alpha h_n-c_\alpha h_n\). The constant family
\((c_\alpha h_n)_{t>0}\) has maximal norm \(|c_\alpha|\|h_n\|_{L^p(\mathcal D)}\) by factoring \(h_n=u_n|h_n|\) as
\((u_n|h_n|^{1/2})\,1\,|h_n|^{1/2}\) in \eqref{max:norm} and using
the norm of a single coordinate for the reverse inequality. The triangle inequality and \eqref{q:maximal} yield
\[
 \|(Y_{n,t})_{t>0}\|_{L^p(\mathcal D;\ell^\infty)}
 \le (C_{p,\alpha}+|c_\alpha|)\|h_n\|_{L^p(\mathcal D)}.
\]
In~\eqref{max:norm} every \(t>0\) is a coordinate of the family. By the order characterization after~\eqref{max:finite}, choose \(B_n^{\mathrm{Re}},B_n^{\mathrm{Im}}\in L^p(\mathcal D)_+\) such that
\begin{equation*}
 -B_n^{\mathrm{Re}}\le\operatorname{Re}Y_{n,t}\le B_n^{\mathrm{Re}},
 \qquad
 -B_n^{\mathrm{Im}}\le\operatorname{Im}Y_{n,t}\le B_n^{\mathrm{Im}}
 \qquad(t>0)\,.
\end{equation*}
We may take the norm of each majorant to be at most \(\|(Y_{n,t})_{t>0}\|_{L^p(\mathcal D;\ell^\infty)}+2^{-2n}\), since both the real and imaginary parts are contractive in the maximal norm. Thus \(b_n=B_n^{\mathrm{Re}}+B_n^{\mathrm{Im}}\) dominates both signs of both parts for every \(t>0\), and \(\|b_n\|_{L^p(\mathcal D)}\le C'_{p,\alpha}2^{-2n}\).

For \(n\in\mathbb N\), let \(e_n=1_{[0,2^{-n}]}(b_n)\). For \(n_0\in\mathbb N\) put \(e=\bigwedge_{n\ge n_0}e_n\), so \(1-e=\bigvee_{n\ge n_0}(1-e_n)\). By subadditivity of the trace on joins of projections and Chebyshev's inequality for the spectral measure of \(b_n\),
\begin{equation*}
 \operatorname{tr}_{\mathcal D}(1-e)
 \le\sum_{n\ge n_0}\operatorname{tr}_{\mathcal D}(1-e_n)
 \le\sum_{n\ge n_0}2^{np}\|b_n\|_{L^p(\mathcal D)}^p
 \le (C'_{p,\alpha})^p\sum_{n\ge n_0}2^{-np}\,.
\end{equation*}
Choose \(n_0\) so that the right-hand side is smaller than \(\eta\), and fix the resulting \(e\). For each \(n\ge n_0\), we have \(e\le e_n\), whence \(e b_n e=e(e_n b_n e_n)e\le 2^{-n}e\). Compressing the order inequalities gives
\begin{align*}
 -2^{-n}e &\le e(\operatorname{Re}Y_{n,t})e\le 2^{-n}e\,,\\
 -2^{-n}e &\le e(\operatorname{Im}Y_{n,t})e\le 2^{-n}e
 \qquad(t>0)\,.
\end{align*}
Both compressed parts are bounded in operator norm by \(2^{-n}\), and hence \(\sup_{t>0}\|eY_{n,t}e\|_{L^\infty(\mathcal D)}\le 2^{1-n}\) for every \(n\ge n_0\).

Writing \(x=x_n+h_n\) and using \(\sup_{t>0}\|eY_{n,t}e\|_{L^\infty(\mathcal D)}\le2^{1-n}\) in the triangle inequality gives
\begin{equation*}
 \|e(A_{t,\beta}^\alpha x-c_\alpha x)e\|_{L^\infty(\mathcal D)}
 \le \|A_{t,\beta}^\alpha x_n-c_\alpha x_n\|_{L^\infty(\mathcal D)}+2^{1-n}
 \qquad (t>0,\ n\ge n_0)\,.
\end{equation*}
For each \(n\ge n_0\) the first term on the right tends to zero as \(t\downarrow0\), since \(x_n=S_{\zeta_n}a_n\) belongs to the class on which we proved operator-norm convergence. Taking the upper limit as \(t\downarrow0\) and then letting \(n\to\infty\) with \(e\) fixed gives the required convergence.
\end{proof}

\section{Truncated cone multipliers}\label{sec:cone}

We write \((x,t)\in\R^2\times\R\) for the physical variables and
\((\xi,\rho)\in\R^2\times\R\) for their Fourier variables, with phase \(x\cdot\xi+t\rho\).
For a bounded measurable symbol \(m:\R^3\to\mathbb C\), let
\(T_m:L^2(\N_{\R^3})\to L^2(\N_{\R^3})\) be the multiplier with symbol \(m\).
For a real exponent \(1<p<\infty\), put
\[
 \kappa(p)=\max\{0,\,2|1/p-1/2|-1/2\}.
\]
For a real order \(\lambda>0\) and a cutoff \(\chi\in C_c^\infty((0,\infty))\), let \(M_\lambda:\R^3\to\mathbb C\) be zero on \(\rho\le0\) and set
\begin{equation}\label{cone:symbol}
 M_\lambda(\xi,\rho)=\chi(\rho)\left(1-\frac{|\xi|^2}{\rho^2}\right)_+^\lambda
 \qquad(\rho>0).
\end{equation}

In the scalar theory, the sharp square function estimate for the cone in
\(\R^3\) implies both planar local smoothing and the classical open
range of cone multiplier orders \cite[Section~1.1]{gwz}.
Our aim is to bound \(T_{M_\lambda}\) on \(L^p(\N_{\R^3})\) for
\(\lambda>\kappa(p)\), uniformly in the coefficient algebra
(Theorem~\ref{cone:theorem}). Near the cone we decompose \(M_\lambda\)
into shells of thickness \(0<\delta\le1/2\). For a smooth normal profile
\(a((\rho-|\xi|)/\delta)\), \(a\in C_c^\infty(\R)\), Fourier inversion in \(\rho\) gives
the phase \(e^{is|\xi|}\) with weight \(\delta\check a(\delta s)\),
\(s\in\R\).
For a Schwartz test tensor \(F\) and a smooth spatial cutoff \(b\) supported away from the origin, the convolution involves \(U(s)b(D_x)F(\cdot,t-s)\). Changing variables to \(u=t-s\) in the spacetime integral lets us apply Lemma~\ref{cone:long-wave} to \(b(D_x)F(\cdot,u)\) for each fixed \(u\in\R\).
Section~\ref{cone:order} interpolates the shell bounds and sums them with weights \(\delta^\lambda\). A fixed temporal frequency slice also gives scalar necessary conditions on \(\lambda\). Section~\ref{cone:quantum-section} transfers the cone estimates to three-dimensional quantum Euclidean spaces.

\subsection{Thin cone shell estimates}

We bound the shell multiplier on \(L^4(\N_{\R^3})\) by averaging the half-wave. We also estimate the \(L^1(\R^3)\) norm of the kernel by dividing the shell into angular caps. The kernel estimate gives an operator bound on \(L^\infty(\N_{\R^3})\), which we interpolate with the bound on \(L^4(\N_{\R^3})\) in Section~\ref{cone:order}. The profile in the normal variable may depend smoothly on \((\xi,\rho)\).

\begin{proposition}\label{cone:shell}
Fix \(b\in C_c^\infty(\R^2\setminus\{0\})\), \(c\in C_c^\infty(\R)\), and put \(L_\delta=1+\log(\delta^{-1})\) for \(0<\delta\le1/2\). Let \(A_\delta\in C_c^\infty(\R^3\times\R)\) be supported in a fixed compact set and satisfy \(\sup_\delta\|A_\delta\|_{C^8}\le1\). For
\[
 m_\delta(\xi,\rho)=b(\xi)c(\rho)
 A_\delta\big(\xi,\rho,(\rho-|\xi|)/\delta\big)
\]
one has
\begin{equation}\label{cone:shell-estimates}
 \|T_{m_\delta}\|_{L^4(\N_{\R^3})\to L^4(\N_{\R^3})}\le C L_\delta^{1/2},\qquad
 \|\check m_\delta\|_{L^1(\R^3)}\le C\delta^{-1/2}.
\end{equation}
Here and below the operator norms are on \(L^p(\N_{\R^3})\). The constants in \eqref{cone:shell-estimates} depend only on the fixed supports and the \(C^4\) norms of \(b,c\), and are uniform in the semifinite tracial algebra \(\M\).
\end{proposition}

\begin{proof}
To prove \eqref{cone:shell-estimates}, we first bound
\(\|\check m_\delta\|_{L^1(\R^3)}\) by decomposing the circle into angular caps.
For \(\|T_{m_\delta}\|_{L^4(\N_{\R^3})\to L^4(\N_{\R^3})}\), we retain the entire annulus and use
time averaging.

Choose a smooth angular partition on the circle with \(O(\delta^{-1/2})\) caps of aperture comparable to \(\delta^{1/2}\) and less than \(\pi/4\). Write \(q_{\delta,e}:\R^2\setminus\{0\}\to[0,1]\) for the angular cutoff of the cap centered at \(e\in\mathbb S^1\) and set \(m_{\delta,e}(\xi,\rho)=q_{\delta,e}(\xi)m_\delta(\xi,\rho)\). The cutoffs have uniformly bounded derivatives in the scaled angle, since they are obtained by normalizing translates of a fixed bump in that variable. Choose a unit vector \(e^\perp\) perpendicular to \(e\), and define \(B_e:\R^3\to\R^3\) by
\[
 B_e(u,v,w)=(u e+\delta^{1/2}v e^\perp,\,u+\delta w)
 \qquad ((u,v,w)\in\R^3).
\]
Put \(g_{\delta,e}=m_{\delta,e}\circ B_e\). For \(\xi\) in the cap, write \(\theta\) for its signed angle from \(e\). Since \(u=|\xi|\cos\theta\) and \(v=\delta^{-1/2}|\xi|\sin\theta\), the support of \(b\) and the aperture bound give \(0<u_0\le u\le u_1\) and \(|v|\le C\) on the support of \(g_{\delta,e}\), with \(u_0,u_1,C\) independent of \(\delta,e\). There the normal coordinate satisfies
\[
 \frac{\rho-|\xi|}{\delta}
 =w-\frac{v^2}{u+\sqrt{u^2+\delta v^2}} .
\]
The support of \(A_\delta\) bounds \((\rho-|\xi|)/\delta\), so this formula also bounds \(|w|\). Thus \(g_{\delta,e}\) is supported in a fixed box with \(u\) bounded away from zero.

On this box the scaled angle
\[
 \delta^{-1/2}\arctan(\delta^{1/2}v/u)
 =\frac vu\int_0^1\frac{d\tau}{1+\delta\tau^2(v/u)^2}
\]
and the normal coordinate have uniformly bounded derivatives through order four. The latter follows from \(u+\sqrt{u^2+\delta v^2}\ge 2u_0\). The chain rule applied to \(q_{\delta,e}\), \(b\), \(c\), and \(A_\delta\) now bounds the \(C^4\) norm of \(g_{\delta,e}\) uniformly. For \(X=(x,t)\in\R^3\), Fourier inversion gives
\[
 (1+|X|^2)^2|\check g_{\delta,e}(X)|
 \le (2\pi)^{-3}\|(1-\Delta)^2g_{\delta,e}\|_{L^1(\R^3)}\le C .
\]
The change of frequency variables by \(B_e\), whose determinant has absolute value \(\delta^{3/2}\), gives
\[
 \check m_{\delta,e}(X)=\delta^{3/2}\check g_{\delta,e}(B_e^{\mathsf T}X),
 \qquad
 \|\check m_{\delta,e}\|_{L^1(\R^3)}
   =\|\check g_{\delta,e}\|_{L^1(\R^3)}\le C .
\]
Summing \(\|\check m_{\delta,e}\|_{L^1(\R^3)}\) over the caps proves \(\|\check m_\delta\|_{L^1(\R^3)}\lesssim\delta^{-1/2}\).

\medskip
For the \(L^4\) bound, first take
\(m_\delta(\xi,\rho)=b(\xi)c(\rho)a((\rho-|\xi|)/\delta)\), where \(a\in C_c^\infty(\R)\) is supported in a fixed interval. The factor \(c(\rho)\) acts by convolution in \(t\) with the integrable kernel \(\check c\). We may therefore omit this factor at cost \(\|\check c\|_{L^1(\R)}\), writing \(\widetilde m_\delta(\xi,\rho)=b(\xi)a((\rho-|\xi|)/\delta)\). At fixed \(\xi\in\R^2\), the substitution \(r=(\rho-|\xi|)/\delta\), with \(d\rho=\delta\,dr\), gives the inverse transform in \(\rho\)
\[
 \frac1{2\pi}\int_\R e^{is\rho}a\big((\rho-|\xi|)/\delta\big)\,d\rho
 =\delta e^{is|\xi|}\check a(\delta s),\qquad
 \check a(s)=\frac1{2\pi}\int_\R e^{isr}a(r)\,dr .
\]
For \(s\in\R\), differentiating the localized symbol \(\xi\mapsto e^{is|\xi|}b(\xi)\) through order four and using~\eqref{pre:kern} gives
\[
 \|(e^{is|\cdot|}b)^\vee\|_{L^1(\R^2)}
 \le C_b(1+|s|)^4 .
\]
Let \(F\) be a finite sum of Schwartz test tensors. The spatial kernel of \(U(s)b(D_x)\) varies continuously in \(L^1(\R^2)\) on bounded \(s\)-intervals by~\eqref{pre:kern}. Together with norm continuity of the translated slices of \(F\), this gives strong measurability of the integrand below. For fixed \(\delta\), the rapid decay of \(\check a(\delta s)\) and the polynomial kernel bound ensure absolute Bochner integrability in \(L^4(\N_{\R^2})\) at each \(t\) and in \(L^4(\N_{\R^3})\), using~\eqref{pre:fubini} and invariance under time translation. The same integral exists in \(L^2(\N_{\R^3})\) by Plancherel's identity. Fourier inversion in time gives
\[
 T_{\widetilde m_\delta}F(x,t)=
 \int_\R\delta\check a(\delta s)\big[U(s)b(D_x)F(\cdot,t-s)\big](x)\,ds .
\]
Set \(w_\delta(s)=\delta|\check a(\delta s)|\) for \(s\in\R\), \(W=\int_\R w_\delta(s)\,ds=\|\check a\|_{L^1(\R)}\), and \(h_u=b(D_x)F(\cdot,u)\in L^4(\N_{\R^2})\) for \(u\in\R\). To use~\eqref{cone:long-time} with \(h_u\) fixed, apply the triangle inequality in \(L^4(\N_{\R^2})\), then use the scalar H\"older bound \((\int w_\delta g)^4\le W^3\int w_\delta g^4\) with \(g(s)=\|U(s)h_{t-s}\|_{L^4(\N_{\R^2})}\). Integrating in \(t\), changing variables \(u=t-s\) at fixed \(s\), and applying Tonelli's theorem gives
\begin{align*}
 \|T_{\widetilde m_\delta}F\|_{L^4(\N_{\R^3})}^4
 &\le W^3\int_\R w_\delta(s)\int_\R
             \|U(s)h_{t-s}\|_{L^4(\N_{\R^2})}^4\,dt\,ds\\
 &=W^3\int_\R I_\delta(u)\,du,\qquad
 I_\delta(u)=\int_\R w_\delta(s)\|U(s)h_u\|_{L^4(\N_{\R^2})}^4\,ds .
\end{align*}

To bound \(I_\delta(u)\), choose \(0<r_0<r_1<\infty\) such that each \(h_u\) has spatial Fourier support in \(r_0\le|\xi|\le r_1\), using \(\operatorname{supp}\widehat h_u\subset\operatorname{supp}b\). Since \(\sigma(4)=0\) and \(\beta(4)=1/2\), the fourth power of~\eqref{cone:long-time} at \(S=2^j\) gives
\[
 \int_{2^j\le |s|\le2^{j+1}}\|U(s)h_u\|_{L^4(\N_{\R^2})}^4\,ds
 \le C2^j(1+j)^2\|h_u\|_{L^4(\N_{\R^2})}^4
 \qquad (j=0,1,2,\ldots).
\]
For the negative interval, the substitution \(r=-s\) lets us apply~\eqref{cone:long-time} to \(U(-r)h_u\) on \([2^j,2^{j+1}]\). For \(|s|\le1\), Lemma~\ref{cone:long-wave} gives \(\|U(s)h_u\|_{L^4(\N_{\R^2})}\le C\|h_u\|_{L^4(\N_{\R^2})}\). Two integrations by parts, together with the direct bound near \(\delta s=0\), give \(w_\delta(s)\lesssim\|a\|_{C^2}\delta(1+\delta|s|)^{-2}\). Multiplying the integral over each time interval by the supremum of \(w_\delta\) there yields
\[
 I_\delta(u)\le C\|a\|_{C^2}\|h_u\|_{L^4(\N_{\R^2})}^4
 \left(\delta+\sum_{j\ge0}\delta 2^j(1+\delta 2^j)^{-2}(1+j)^2\right).
\]
With \(J=\lfloor\log_2(\delta^{-1})\rfloor\), the part \(0\le j\le J\) of the sum is at most \(2\delta 2^J(1+J)^2\). For \(j=J+n\), \(n\ge1\), the remaining sum is bounded by \(C\sum_{n\ge1}2^{-n}(1+J+n)^2\lesssim(1+J)^2\), since \(1/2<\delta 2^J\le1\). Thus \(I_\delta(u)\le C\|a\|_{C^2}L_\delta^2\|h_u\|_{L^4(\N_{\R^2})}^4\).

Integrating the weight estimate gives \(W\lesssim\|a\|_{C^2}\). Spatial convolution by \(\check b\) gives
\[
 \int_\R\|h_u\|_{L^4(\N_{\R^2})}^4\,du
 \le\|\check b\|_{L^1(\R^2)}^4\|F\|_{L^4(\N_{\R^3})}^4 .
\]
Inserting the bound for \(I_\delta(u)\) into \(W^3\int_\R I_\delta(u)\,du\) and taking fourth roots gives
\[
 \|T_{\widetilde m_\delta}F\|_{L^4(\N_{\R^3})}
 \le C\|\check b\|_{L^1(\R^2)}\|a\|_{C^2}L_\delta^{1/2}\|F\|_{L^4(\N_{\R^3})}.
\]
Convolving in time by \(\check c\) and using density of the test tensors gives the bound for a single normal profile
\[
 \|T_{m_\delta}\|_{L^4(\N_{\R^3})\to L^4(\N_{\R^3})}
 \le C\|\check c\|_{L^1(\R)}\|\check b\|_{L^1(\R^2)}\|a\|_{C^2}L_\delta^{1/2} .
\]

\medskip
To apply the single-profile bound to a general amplitude \(A_\delta\), we separate its dependence on the frequency coordinates from its dependence on the normal variable \(r\). Take a fixed box \(B\) of side \(D>0\) in the three frequency coordinates \(z=(\xi,\rho)\). Its interior contains both the projections of the supports of \(A_\delta\) and \(\operatorname{supp}b\times\operatorname{supp}c\). Expand the periodic extension of \(A_\delta\) in \(z\) as \(\sum_{n\in\mathbb Z^3}e^{2\pi i n\cdot z/D}a_{\delta,n}(r)\), where the coefficients \(a_{\delta,n}:\R\to\mathbb C\) are given by
\[
 a_{\delta,n}(r)=D^{-3}\int_B e^{-2\pi i n\cdot z/D}A_\delta(z,r)\,dz .
\]
They are supported in a common compact interval in \(r\). For \(0\le k\le2\), integrating by parts with \((1-\Delta_z)^2\) gives
\[
 \partial_r^k a_{\delta,n}(r)
 = \frac{D^{-3}}{(1+4\pi^2|n|^2/D^2)^2}
   \int_B e^{-2\pi i n\cdot z/D}(1-\Delta_z)^2\partial_r^k A_\delta(z,r)\,dz .
\]
The derivatives of \(A_\delta\) in the integral have total order at most six. Hence \(\|a_{\delta,n}\|_{C^2}\le C(1+|n|)^{-4}\) uniformly in \(\delta\).

For \(n\in\mathbb Z^3\), set \(m_{\delta,n}(\xi,\rho)=b(\xi)c(\rho)a_{\delta,n}((\rho-|\xi|)/\delta)\). For \(h\in\R^3\) let \(\tau_hG(X)=G(X+h)\) on \(L^4(\N_{\R^3})\). Multiplication of the symbol \(m_{\delta,n}\) by \(e^{2\pi i n\cdot z/D}\) corresponds to \(\tau_{2\pi n/D}T_{m_{\delta,n}}\). Translation is an isometry on \(L^4(\N_{\R^3})\), so the single-profile estimate bounds the norm of this operator by \(CL_\delta^{1/2}(1+|n|)^{-4}\). It follows that
\[
 T_{m_\delta}=\sum_{n\in\mathbb Z^3}\tau_{2\pi n/D}T_{m_{\delta,n}}
\]
in operator norm on \(L^4(\N_{\R^3})\). Indeed, the norm estimates are summable, and the partial symbols converge uniformly to \(m_\delta\) by the decay of \(\|a_{\delta,n}\|_{C^2}\). For \(F\in L^2(\N_{\R^3})\cap L^4(\N_{\R^3})\), Plancherel's identity gives convergence to \(T_{m_\delta}F\) in \(L^2\). Uniqueness of the limit in measure identifies this with the \(L^4\) sum. Density gives the extension to all of \(L^4\), and summing the norms proves the \(L^4\) estimate in~\eqref{cone:shell-estimates}.

\end{proof}

\subsection{Bounds and necessary conditions for the cone order}\label{cone:order}

Interpolating the bounds in~\eqref{cone:shell-estimates} gives the \(L^p\) estimates needed to sum the dyadic pieces of \(M_\lambda\).

\begin{theorem}\label{cone:theorem}
For \(1<p<\infty\) and \(\lambda>\kappa(p)\), \(T_{M_\lambda}\) extends from \(L^p\cap L^2(\N_{\R^3})\) to a bounded operator on \(L^p(\N_{\R^3})\), of norm at most \(C_{p,\lambda,\chi}\), uniformly in \(\M\) and its trace.
\end{theorem}

\begin{proof}
Fix \(1<p<\infty\) and \(\lambda>\kappa(p)\). Choose \(0<\rho_-<\rho_+<\infty\) with \(\operatorname{supp}\chi\subset[\rho_-,\rho_+]\), and an integer \(j_0\ge1\) with \(2^{1-j_0}\le\rho_-/2\). Take a nonnegative \(\varphi\in C_c^\infty((1/2,2))\) with \(\sum_{j\in\mathbb Z}\varphi(2^j u)=1\) for \(u>0\). For \((\xi,\rho)\in\R^3\) put \(u=\rho-|\xi|\) and define
\begin{equation*}
 m_j(\xi,\rho)=M_\lambda(\xi,\rho)\varphi(2^j u)\quad(j\ge j_0),
 \qquad m_{\rm int}=M_\lambda-\sum_{j\ge j_0}m_j\,.
\end{equation*}
For each integer \(J\ge j_0\) and \(F\in L^2\cap L^p(\N_{\R^3})\), the triangle inequality gives
\begin{multline*}
 \Big\|T_{m_{\rm int}}F+\sum_{j=j_0}^J T_{m_j}F\Big\|_{L^p(\N_{\R^3})}
 \le \\
 \Big(\|T_{m_{\rm int}}\|_{L^p(\N_{\R^3})\to L^p(\N_{\R^3})}
          +\sum_{j=j_0}^J\|T_{m_j}\|_{L^p(\N_{\R^3})\to L^p(\N_{\R^3})}\Big)\|F\|_{L^p(\N_{\R^3})}\,.
\end{multline*}
We bound the first operator norm by \(\|\check m_{\rm int}\|_{L^1(\R^3)}\) and show that the norms of \(T_{m_j}\) are summable. We then identify the operator sum with \(T_{M_\lambda}\).

\medskip
To see that $\check m_{\rm int}$ is integrable, we show that $m_{\rm int}$ is smooth and compactly supported. For \(0<u<2^{-j_0}\) all terms with \(j<j_0\) in the partition vanish, so \(\sum_{j\ge j_0}m_j=M_\lambda\) there. When \(u\le0\), both \(M_\lambda\) and the shell symbols vanish. Thus \(m_{\rm int}\) vanishes near \(u=0\). On \(\operatorname{supp}\chi\), all \(m_j\) vanish for \(|\xi|<\rho_-/2\), since their supports require \(u<2^{1-j_0}\). In that neighborhood \(m_{\rm int}=M_\lambda\), which is smooth in \(\xi\) at zero by~\eqref{cone:symbol}. Away from \(\xi=0\) and \(u=0\), the shell sum is locally finite and smooth. Thus \(m_{\rm int}\) is smooth with compact support, so its inverse Fourier transform is integrable and bounds \(\|T_{m_{\rm int}}\|_{L^p(\N_{\R^3})\to L^p(\N_{\R^3})}\) by convolution.

\medskip
We bound the annular terms by applying Proposition~\ref{cone:shell} to the normalized symbols \(2^{j\lambda}m_j\). For each \(j\ge j_0\) let \(\delta=2^{-j}\) and \(v=u/\delta\). On \(\operatorname{supp}m_j\) we have \(\rho\in[\rho_-,\rho_+]\), \(\delta/2<u<2\delta\) and \(\rho_-/2\le|\xi|\le\rho_+\). Choose \(b\in C_c^\infty(\R^2\setminus\{0\})\) equal to one on this fixed annulus and put \(c=\chi\). Take \(\widetilde b\in C_c^\infty(\R^2\setminus\{0\})\) and \(\widetilde c\in C_c^\infty((0,\infty))\) equal to one on \(\operatorname{supp}b\) and \(\operatorname{supp}c\), respectively. Define \(A:\R^3\times\R\to\mathbb C\) for \(\rho>0\), \(\xi\ne0\) and \(v>0\) by
\begin{equation*}
 A(\xi,\rho,v)=
 \widetilde b(\xi)\widetilde c(\rho)\varphi(v)v^\lambda
 \left(\frac{\rho+|\xi|}{\rho^2}\right)^\lambda
\end{equation*}
and extend by zero elsewhere. The factorization \(1-|\xi|^2/\rho^2=u(\rho+|\xi|)/\rho^2\) for \(\rho>0\) gives
\begin{equation*}
 2^{j\lambda}m_j(\xi,\rho)=b(\xi)c(\rho)A(\xi,\rho,2^j u)
 \qquad(j\ge j_0)\,.
\end{equation*}
On the fixed compact support of \(A\), both positive bases \(v\) and \((\rho+|\xi|)/\rho^2\) are smooth and bounded above and away from zero. The cutoffs \(\widetilde b(\xi)\), \(\widetilde c(\rho)\), and \(\varphi(v)\) are supported away from \(\xi=0\), \(\rho=0\), and \(v=0\), respectively. Thus the zero extension is smooth, and differentiating the product gives \(\|A\|_{C^8}\le C_{\chi,\lambda}\). Dividing \(A\) by a bound \(C_{\chi,\lambda}\ge1\) gives the normalization in Proposition~\ref{cone:shell}.

For any normalized symbol \(m_\delta\) in Proposition~\ref{cone:shell}, Plancherel's identity and convolution with \(\check m_\delta\) give, respectively,
\[
 \|T_{m_\delta}\|_{L^2(\N_{\R^3})\to L^2(\N_{\R^3})}\lesssim1,
 \qquad
 \|T_{m_\delta}\|_{L^\infty(\N_{\R^3})\to L^\infty(\N_{\R^3})}
 \lesssim\delta^{-1/2}.
\]
Interpolating \(L^2(\N_{\R^3})\) and \(L^4(\N_{\R^3})\) with weight \(2-4/p\) at \(L^4(\N_{\R^3})\) gives \(C_p L_\delta^{1-2/p}\) for \(2\le p\le4\). Interpolating \(L^4(\N_{\R^3})\) and \(L^\infty(\N_{\R^3})\) with weight \(4/p\) at \(L^4(\N_{\R^3})\) gives \(C_p\delta^{-(1/2-2/p)}L_\delta^{2/p}\) for \(4\le p<\infty\).

If \(1<p<2\), put \(p'=p/(p-1)>2\). Apply the estimate just proved at \(p'\) to \(T_{\overline{m_\delta}}\), whose symbol has the same support and derivative bounds. Under the trace pairing this operator is adjoint to \(T_{m_\delta}\) on \(L^2\), so its \(L^{p'}\) bound gives the \(L^p\) extension of \(T_{m_\delta}\) with the same bound. For \(2\le p<\infty\), the power of \(\delta^{-1}\) in the interpolated estimates is \(\kappa(p)=\max\{0,1/2-2/p\}\). Since \(\kappa(p')=\kappa(p)\) and \(L_\delta^\gamma\le C_{\gamma,\epsilon}\delta^{-\epsilon}\) for fixed \(\gamma\ge0\) and \(\epsilon>0\), this proves
\begin{equation}\label{cone:open-shell}
 \|T_{m_\delta}\|_{L^p(\N_{\R^3})\to L^p(\N_{\R^3})}\le C_{p,\epsilon}\delta^{-\kappa(p)-\epsilon}
 \qquad(1<p<\infty,\ \epsilon>0)\,.
\end{equation}
The constants are uniform in \(\M\) and its trace and have the cutoff dependence of Proposition~\ref{cone:shell}. Apply~\eqref{cone:open-shell} to \(2^{j\lambda}m_j\) with the normalization of \(A\) accounted for. Choosing \(\epsilon=(\lambda-\kappa(p))/2\) gives
\begin{equation*}
 \sum_{j\ge j_0}\|T_{m_j}\|_{L^p(\N_{\R^3})\to L^p(\N_{\R^3})}
 \le C_{p,\lambda,\chi,\epsilon}
       \sum_{j\ge j_0}2^{-j(\lambda-\kappa(p)-\epsilon)}<\infty\,.
\end{equation*}

\medskip
It remains to identify the limit of the bounded operator sums with \(T_{M_\lambda}\). The bounds on \(T_{m_{\rm int}}\) and the shell norms give an operator-norm limit of \(T_{m_{\rm int}}+\sum_{j=j_0}^J T_{m_j}\) on \(L^p(\N_{\R^3})\) with norm at most \(C_{p,\lambda,\chi}\). The finite symbol sums converge uniformly to \(M_\lambda\), since \(\|m_j\|_{L^\infty(\R^3)}\lesssim_{\lambda,\chi}2^{-j\lambda}\). For \(F\in L^2\cap L^p(\N_{\R^3})\), Plancherel's identity gives convergence of the same partial sums applied to \(F\) in \(L^2\) to \(T_{M_\lambda}F\). The \(L^2\) and \(L^p\) limits agree by uniqueness in measure. Bounded elements whose support projections have finite trace are dense in \(L^p\) by Section~\ref{sec:prelim} and belong to \(L^2\cap L^p\). This gives the unique extension.
\end{proof}

\begin{corollary}\label{cone:small-order}
For \(0<\lambda\le1\) one has \(\|T_{M_\lambda}\|_{L^4(\N_{\R^3})\to L^4(\N_{\R^3})}\le C_\chi\lambda^{-3/2}\), uniformly over the coefficient algebras.
\end{corollary}

\begin{proof}
Write \(M_\lambda=m_{\rm int}+\sum_{j\ge j_0}m_j\) as in the proof of Theorem~\ref{cone:theorem}. The index \(j_0\), the partition function \(\varphi\), and the cutoffs \(b,c,\widetilde b,\widetilde c\) used there are fixed independently of \(\lambda\). By the triangle inequality, for \(0<\lambda\le1\) it suffices to show that
\[
 \|\check m_{\rm int}\|_{L^1(\R^3)}
 +\sum_{j\ge j_0}\|T_{m_j}\|_{L^4(\N_{\R^3})\to L^4(\N_{\R^3})}
 \le C_\chi\lambda^{-3/2}.
\]

On \(\operatorname{supp}m_{\rm int}\) we have \(u\ge2^{-j_0}\) and \(\rho\in[\rho_-,\rho_+]\), so
\begin{equation*}
 1-\frac{|\xi|^2}{\rho^2}
 =\frac{u(\rho+|\xi|)}{\rho^2}
 \ge\frac{2^{-j_0}}{\rho_+}\,.
\end{equation*}
There \(m_{\rm int}\) is the product of this base to the power \(\lambda\) and the fixed cutoff \(\chi(\rho)[1-\sum_{j\ge j_0}\varphi(2^j u)]\). The cutoff vanishes on the interior near \(u=0\) and equals \(\chi(\rho)\) near \(\xi=0\). Writing the positive power as the exponential of \(\lambda\) times the logarithm of the base gives \(\|m_{\rm int}\|_{C^4}\le C_\chi\) for \(0<\lambda\le1\). Applying \((1-\Delta_{\xi,\rho})^2\) in Fourier inversion bounds \(|\check m_{\rm int}(X)|\) by \(C_\chi(1+|X|^2)^{-2}\), hence \(\|\check m_{\rm int}\|_{L^1(\R^3)}\le C_\chi\).

For the amplitude \(A\) of \(2^{j\lambda}m_j\), the same exponential formula for \(v^\lambda\) and \(((\rho+|\xi|)/\rho^2)^\lambda\) gives a uniform \(C^8\) bound for \(0<\lambda\le1\). Taking \(\delta=2^{-j}\) in~\eqref{cone:shell-estimates} gives \(\|T_{m_j}\|_{L^4(\N_{\R^3})\to L^4(\N_{\R^3})}\le C_\chi 2^{-j\lambda}(1+j)^{1/2}\). Comparison on unit intervals and the substitution \(v=\lambda s\) yield
\begin{equation*}
 \sum_{j\ge j_0}2^{-j\lambda}(1+j)^{1/2}
 \lesssim\int_0^\infty 2^{-\lambda s}(1+s)^{1/2}\,ds
 =\lambda^{-3/2}\int_0^\infty 2^{-v}(\lambda+v)^{1/2}\,dv\,.
\end{equation*}
Since \(0<\lambda\le1\), the integral in \(v\) is at most \(\int_0^\infty 2^{-v}(1+v)^{1/2}\,dv\). Adding the bound \(\|\check m_{\rm int}\|_{L^1(\R^3)}\le C_\chi\) gives \(\|T_{M_\lambda}\|_{L^4(\N_{\R^3})\to L^4(\N_{\R^3})}\le C_\chi\lambda^{-3/2}\).
\end{proof}

\medskip
\phantomsection\label{cone:necessary}
We derive the necessary conditions for scalar boundedness by taking a slice of \(M_\lambda\) at fixed temporal frequency. Fix \(1<p<\infty\) and suppose \(\chi(\rho_0)\ne0\) for some \(\rho_0>0\). Scalar \(L^p(\R^3)\) boundedness of \(T_{M_\lambda}\) requires \(\lambda>\kappa(p)\) for \(\lambda>0\), and requires \(p=2\) at \(\lambda=0\). For the zeroth power we use the indicator of the open cone.

Let \(\lambda\ge0\) and suppose \(T_{M_\lambda}\) is bounded on scalar \(L^p(\R^3)\). To take a slice at temporal frequency \(\rho_0\), choose \(f\in\mathscr S(\R^2)\), \(0\ne g\in\mathscr S(\R)\), and for \(\epsilon>0\) use \(F_\epsilon(x,t)=f(x)\epsilon^{1/p}g(\epsilon t)e^{i\rho_0t}\). Put
\begin{equation*}
 G_\epsilon(x,s)=\epsilon^{-1/p}e^{-i\rho_0s/\epsilon}
                   (T_{M_\lambda}F_\epsilon)(x,s/\epsilon)
 \qquad((x,s)\in\R^2\times\R)\,.
\end{equation*}
Dilation and modulation give \(\|G_\epsilon\|_{L^p(\R^3)}\le\|T_{M_\lambda}\|_{L^p(\R^3)\to L^p(\R^3)}\|f\|_{L^p(\R^2)}\|g\|_{L^p(\R)}\) and
\begin{equation*}
 \widehat{G_\epsilon}(\xi,\eta)
 =M_\lambda(\xi,\rho_0+\epsilon\eta)\widehat f(\xi)\widehat g(\eta)
 \qquad((\xi,\eta)\in\R^2\times\R)\,.
\end{equation*}
For \(|\xi|\ne\rho_0\) the symbol converges to \(M_\lambda(\xi,\rho_0)\). Since it is bounded and \(\widehat f\otimes\widehat g\in L^1(\R^3)\), dominated convergence in Fourier inversion gives \(G_\epsilon(x,s)\to g(s)M_\lambda(D_x,\rho_0)f(x)\), where \(M_\lambda(D_x,\rho_0)\) is the spatial multiplier with symbol \(M_\lambda(\cdot,\rho_0)\). Fatou's lemma and the product norm give
\begin{equation*}
 \|g\|_{L^p(\R)}\|M_\lambda(D_x,\rho_0)f\|_{L^p(\R^2)}
 \le\|T_{M_\lambda}\|_{L^p(\R^3)\to L^p(\R^3)}\|f\|_{L^p(\R^2)}\|g\|_{L^p(\R)}\,.
\end{equation*}
Cancelling \(\|g\|_{L^p(\R)}>0\), dividing by \(|\chi(\rho_0)|\), and dilating the spatial variables gives boundedness of the unit disk multiplier with symbol \((1-|\xi|^2)_+^\lambda\) on \(L^p(\R^2)\).

As the disk symbol is real, duality gives boundedness at \(p'=p/(p-1)\) as well. Applying the multiplier to a Schwartz function whose transform is one on the closed unit disk gives the disk kernel \(K_\lambda\), hence \(K_\lambda\in L^q(\R^2)\) for \(q\in\{p,p'\}\). For \(r=|x|>0\), polar integration gives
\begin{equation*}
 K_\lambda(x)=\frac1{2\pi}\int_0^1 J_0(ru)(1-u^2)^\lambda u\,du
 =c_\lambda r^{-1-\lambda}J_{1+\lambda}(r),
 \qquad c_\lambda=\frac{2^\lambda\Gamma(1+\lambda)}{2\pi}\,.
\end{equation*}
The last equality follows by termwise integration of the power series for \(J_0\), using
\(\int_0^1u^{2k+1}(1-u^2)^\lambda\,du
 =\Gamma(k+1)\Gamma(1+\lambda)/(2\Gamma(k+2+\lambda))\) for \(k\ge0\). The fixed-order expansion~\cite[(10.17.3)]{dlmf} yields
\begin{equation*}
 K_\lambda(x)=c_\lambda\sqrt{\frac2\pi}\,r^{-3/2-\lambda}
   \left[\cos\!\left(r-\frac{\pi(1+\lambda)}2-\frac\pi4\right)
         +O_\lambda(r^{-1})\right]\qquad(r\to\infty)\,.
\end{equation*}
In every sufficiently large period there is an interval of fixed positive length on which the cosine has absolute value at least \(1/2\) and the error is at most \(1/4\). On the interval in the \(n\)th period, \(r\asymp n\), so radial integration against \(r\,dr\) gives
\begin{equation*}
 \|K_\lambda\|_{L^q(\R^2)}^q
 \gtrsim_{\lambda,q}\sum_{n\ge n_0} n^{1-q(3/2+\lambda)}
 \qquad(q\in\{p,p'\})
\end{equation*}
for some integer \(n_0\). Finiteness forces \(3/2+\lambda>2/p\) and \(3/2+\lambda>2/p'\). Together with \(\lambda>0\), these inequalities are equivalent to \(\lambda>\kappa(p)\). At order zero, Fefferman's disk multiplier theorem~\cite{fefferman-ball} excludes every \(p\ne2\), whereas Plancherel's identity gives boundedness at \(p=2\).

\subsection{Cone multipliers on three-dimensional quantum Euclidean spaces}\label{cone:quantum-section}

To prove Corollary~\ref{cone:quantum}, we use translations on the quantum Euclidean space to transfer the cone multiplier bounds of Theorem~\ref{cone:theorem} and Corollary~\ref{cone:small-order}.

Let \(\Theta\) be a real skew-symmetric \(3\times3\) matrix. We use the normalization of the Moyal model in Section~\ref{q:models}. The twisted group von Neumann algebra \(\R^3_\Theta\) is generated by unitaries \(\lambda_\Theta(\zeta)\), \(\zeta\in\R^3\), satisfying
\[
 \lambda_\Theta(\zeta)\lambda_\Theta(\eta)
 =e^{i\langle\eta,\Theta\zeta\rangle/2}\lambda_\Theta(\zeta+\eta)
 \qquad(\zeta,\eta\in\R^3)\,.
\]
For \(f\in\mathscr S(\R^3)\), define \(\lambda_\Theta(f)\in\R^3_\Theta\) by the weak operator integral \(\int_{\R^3} f(\zeta)\lambda_\Theta(\zeta)\,d\zeta\). The Plancherel trace satisfies \(\operatorname{tr}_\Theta(\lambda_\Theta(f))=f(0)\) and \(\|\lambda_\Theta(f)\|_{L^2(\R^3_\Theta)}=\|f\|_{L^2(\R^3)}\).

For a semifinite tracial algebra \(\mathcal C\), put \(\mathcal D=\mathcal C\bar\otimes\R^3_\Theta\) with the product trace. The Plancherel identification tensored with the identity on \(L^2(\mathcal C)\) gives a unitary map
\[
 \mathcal F_\Theta:L^2(\mathcal D)\longrightarrow L^2(\R^3;L^2(\mathcal C)),
 \qquad x\longmapsto\widehat x_\Theta\,.
\]
The trace-preserving dual action \(\beta_z:\mathcal D\to\mathcal D\), \(z\in\R^3\), acts trivially on \(\mathcal C\) and sends \(\lambda_\Theta(\zeta)\) to \(e^{iz\cdot\zeta}\lambda_\Theta(\zeta)\). For a bounded measurable symbol \(m:\R^3\to\mathbb C\), define \(T_m^\Theta:L^2(\mathcal D)\to L^2(\mathcal D)\) by the second formula below. For \(x\in L^2(\mathcal D)\) and almost every \(\zeta\in\R^3\),
\[
 \widehat{\beta_zx}_\Theta(\zeta)=e^{iz\cdot\zeta}\widehat x_\Theta(\zeta),
 \qquad
 \widehat{T_m^\Theta x}_\Theta(\zeta)=m(\zeta)\widehat x_\Theta(\zeta)\,.
\]

Trace preservation gives \(\|\beta_zx\|_{L^p(\mathcal D)}=\|x\|_{L^p(\mathcal D)}\) for \(x\in L^p(\mathcal D)\) and \(1\le p\le\infty\). The character formula gives \(\|\beta_zx-x\|_{L^2(\mathcal D)}\to0\) as \(z\to0\) by dominated convergence. For \(x=ab\in L^1(\mathcal D)\) with \(a,b\in L^2(\mathcal D)\) as in Section~\ref{q:models}, the trace H\"older inequality gives
\[
 \|\beta_zx-x\|_{L^1(\mathcal D)}\le
 \|\beta_z a-a\|_{L^2(\mathcal D)}\|b\|_{L^2(\mathcal D)}+\|a\|_{L^2(\mathcal D)}\|\beta_z b-b\|_{L^2(\mathcal D)}\longrightarrow0\,.
\]
For bounded \(x\in\mathcal D\) whose support projections have finite trace, writing \(d_z=\beta_zx-x\) gives the interpolation estimates
\begin{align*}
 \|d_z\|_{L^p(\mathcal D)}&\le (2\|x\|_{L^\infty(\mathcal D)})^{1-2/p}\|d_z\|_{L^2(\mathcal D)}^{2/p}
                  &&(2<p<\infty)\,,\\
 \|d_z\|_{L^p(\mathcal D)}&\le \|d_z\|_{L^1(\mathcal D)}^{2/p-1}\|d_z\|_{L^2(\mathcal D)}^{2-2/p}
                  &&(1<p<2)\,.
\end{align*}
Density of these elements in \(L^p(\mathcal D)\) and preservation of the \(L^p\) norm give continuous orbits for every \(1\le p<\infty\).

\begin{corollary}\label{cone:quantum}
For \(1<p<\infty\) and \(\lambda>\kappa(p)\), \(T_{M_\lambda}^\Theta\) extends from \(L^2\cap L^p(\mathcal D)\) to \(L^p(\mathcal D)\) with norm bounded by \(C_{p,\lambda,\chi}\) independently of \(\Theta,\mathcal C\). For \(p=4\) and \(0<\lambda\le1\) its norm is at most \(C_\chi\lambda^{-3/2}\).
\end{corollary}

\begin{proof}
To bound \(T_{M_\lambda}^\Theta\), use the decomposition \(M_\lambda=m_{\rm int}+\sum_{j\ge j_0}m_j\) from the proof of Theorem~\ref{cone:theorem}. We first show that each smooth symbol transfers with its semicommutative \(L^p\) operator norm. We then sum the transferred operators and identify their limit with \(T_{M_\lambda}^\Theta\).

Fix \(1<p<\infty\), \(m\in C_c^\infty(\R^3)\), and \(x\in L^p(\mathcal D)\). Define \(Q_m:L^p(\mathcal D)\to L^p(\mathcal D)\) by
\[
 Q_mx=\int_{\R^3}\check m(y)\beta_{-y}x\,dy\,.
\]
The orbit is norm continuous in \(L^p(\mathcal D)\), and the integrand's norm is bounded by \(|\check m(y)|\|x\|_{L^p(\mathcal D)}\). Since \(\check m\in L^1(\R^3)\), this is a Bochner integral defining a bounded operator. For \(x\in L^2(\mathcal D)\cap L^p(\mathcal D)\), the orbit is continuous in the sum norm \(\|\cdot\|_{L^2(\mathcal D)}+\|\cdot\|_{L^p(\mathcal D)}\), in which the same integral converges. In particular its integrals in \(L^2\) and \(L^p\) agree.

To compute the Fourier coefficients of \(Q_mx\) for such \(x\), apply the unitary \(\mathcal F_\Theta\) under the integral in \(L^2\). Pair the result with any \(g\in L^2(\R^3;L^2(\mathcal C))\). Cauchy--Schwarz in \(\zeta\) gives
\[
 \int_{\R^3}|\check m(y)|\int_{\R^3}
     \big|\langle g(\zeta),\widehat x_\Theta(\zeta)\rangle_{L^2(\mathcal C)}\big|
     \,d\zeta\,dy
 \le \|\check m\|_{L^1(\R^3)}\|g\|_{L^2(\R^3;L^2(\mathcal C))}\|x\|_{L^2(\mathcal D)}\,.
\]
Fubini and the formula for \(\widehat{\beta_{-y}x}_\Theta\) therefore give, for almost every \(\zeta\in\R^3\),
\[
 \widehat{Q_mx}_\Theta(\zeta)
 =\left(\int_{\R^3}\check m(y)e^{-iy\cdot\zeta}\,dy\right)\widehat x_\Theta(\zeta)
 =m(\zeta)\widehat x_\Theta(\zeta)\,.
\]

Let \(C_m\) be the operator norm of \(T_m\) on \(L^p(\R^3;L^p(\mathcal D))\), the space in \eqref{pre:fubini} with coefficient algebra \(\mathcal D\). We prove \(\|Q_m\|_{L^p(\mathcal D)\to L^p(\mathcal D)}\le C_m\) by comparing the two operators on large cubes. For \(R>0\), put \(B_R=[-R,R]^3\) and \(w_R=|B_R|^{-1/p}1_{B_R}\). Define \(J_R:L^p(\mathcal D)\to L^p(\R^3;L^p(\mathcal D))\) by \(J_Rx(z)=w_R(z)\beta_zx\) for \(z\in\R^3\). Since \(\|w_R\|_{L^p(\R^3)}=1\) and \(\|\beta_zx\|_{L^p(\mathcal D)}=\|x\|_{L^p(\mathcal D)}\), we have \(\|J_Rx\|_{L^p(\R^3;L^p(\mathcal D))}=\|x\|_{L^p(\mathcal D)}\).

By the triangle inequality,
\[
 \|Q_mx\|_{L^p(\mathcal D)}=\|J_RQ_mx\|_{L^p(\R^3;L^p(\mathcal D))}
 \le\|J_RQ_mx-T_m J_Rx\|_{L^p(\R^3;L^p(\mathcal D))}
       +C_m\|x\|_{L^p(\mathcal D)}\,.
\]

The convolution formula for \(T_m(J_Rx)\) and the group law \(\beta_z\beta_{-y}=\beta_{z-y}\) give, for almost every \(z\in\R^3\),
\[
 (T_m J_Rx-J_RQ_mx)(z)
 =\int_{\R^3}\check m(y)[w_R(z-y)-w_R(z)]\beta_{z-y}x\,dy\,.
\]
Minkowski's integral inequality and preservation of the \(L^p(\mathcal D)\) norm give
\[
 \|T_m J_Rx-J_RQ_mx\|_{L^p(\R^3;L^p(\mathcal D))}
 \le\|x\|_{L^p(\mathcal D)}\int_{\R^3}|\check m(y)|\,
              \|w_R(\cdot-y)-w_R\|_{L^p(\R^3)}\,dy\,.
\]
For each fixed \(y=(y_1,y_2,y_3)\in\R^3\),
\[
 \|w_R(\cdot-y)-w_R\|_{L^p(\R^3)}^p
 =\frac{|(B_R+y)\mathbin{\triangle}B_R|}{|B_R|}
 =2\left[1-\prod_{k=1}^3\left(1-\frac{|y_k|}{2R}\right)_+\right]
 \longrightarrow0 \qquad (R\to\infty)\,.
\]
The ratio is at most two, so dominated convergence makes the integral against \(|\check m(y)|\,dy\) tend to zero. Letting \(R\to\infty\) gives \(\|Q_mx\|_{L^p(\mathcal D)}\le C_m\|x\|_{L^p(\mathcal D)}\).

The operator norms of \(Q_{m_{\rm int}}\) and \(Q_{m_j}\), \(j\ge j_0\), on \(L^p(\mathcal D)\) are thus bounded by the operator norms of \(T_{m_{\rm int}}\) and \(T_{m_j}\), respectively, on \(L^p(L^\infty(\R^3)\bar\otimes\mathcal D)\) from the proof of Theorem~\ref{cone:theorem}. In particular \(Q_{m_{\rm int}}+\sum_{j\ge j_0}Q_{m_j}\) converges in operator norm on \(L^p(\mathcal D)\), with norm at most \(C_{p,\lambda,\chi}\). For \(x\in L^2(\mathcal D)\cap L^p(\mathcal D)\), Plancherel's identity bounds the \(L^2\) error by
\begin{align*}
 &\left\|Q_{m_{\rm int}}x+\sum_{j=j_0}^N Q_{m_j}x-T_{M_\lambda}^\Theta x\right\|_{L^2(\mathcal D)}\\
 &\qquad\le
 \left\|m_{\rm int}+\sum_{j=j_0}^N m_j-M_\lambda\right\|_{L^\infty(\R^3)}
 \|x\|_{L^2(\mathcal D)}
 \longrightarrow0 \qquad (N\to\infty)\,.
\end{align*}
Uniqueness of the limit in measure identifies the \(L^p\) limit with \(T_{M_\lambda}^\Theta x\) on this dense intersection. To obtain the estimate for \(p=4\) and \(0<\lambda\le1\), apply the uniform bounds on the interior term and the shells from the proof of Corollary~\ref{cone:small-order}. This gives
\begin{align*}
 \|T_{M_\lambda}^\Theta\|_{L^4(\mathcal D)\to L^4(\mathcal D)}
 &\le \|Q_{m_{\rm int}}\|_{L^4(\mathcal D)\to L^4(\mathcal D)}
 +\sum_{j\ge j_0}\|Q_{m_j}\|_{L^4(\mathcal D)\to L^4(\mathcal D)}\\
 &\le C_\chi\left(1+\sum_{j\ge j_0}2^{-j\lambda}(1+j)^{1/2}\right)
 \le C'_\chi\lambda^{-3/2}\,.
\end{align*}
\end{proof}

For \(\Theta=0\) and \(\mathcal C=\mathbb C\), the unitaries \(\lambda_0(\zeta)\) are identified with the classical characters \(y\mapsto e^{iy\cdot\zeta}\) on \(\R^3\). The Plancherel trace corresponds to integration against \((2\pi)^{-3}\,dy\), and this measure normalization leaves the multiplier norms unchanged. For \(\chi\not\equiv0\), the scalar necessary conditions of Section~\ref{cone:order} therefore apply to bounds uniform over \(\Theta\) and \(\mathcal C\).

\section{Stein row and column square functions}\label{sec:stein}

For a scalar Schwartz function \(f:\R^2\to\mathbb C\), Stein's square function takes the square root of the integral of \(|r\partial_r B_r^\nu f|^2\) against \(dr/r\), where \(B_r^\nu\) is the Bochner--Riesz mean of order \(\nu\) and radius \(r\). This is the radius normalization in~\cite[(1.1)--(1.2)]{gow}. Carbery used the planar scalar \(L^4\) square-function estimate to bound the maximal Bochner--Riesz operator~\cite{carbery}. With tracial coefficients there are two products to integrate, \((D_r^\nu F)^*D_r^\nu F\) and \(D_r^\nu F(D_r^\nu F)^*\), where \(D_r^\nu=r\partial_r B_r^\nu\). Theorem~\ref{stein:theorem} bounds the square roots of both integrals in \(L^p\) for \(2\le p<\infty\) and \(\nu>\max\{1/2,1-2/p\}\).

The behavior of the derivative symbol near the circle \(|\xi|=r\) is governed by \((1-|\xi|^2/r^2)^{\nu-1}\) for \(|\xi|<r\). On one dyadic interval of radii, Section~\ref{stein:one-scale} uses local smoothing to bound the contribution from a thin annulus about the circle, as in the scalar reduction of~\cite[Section 11]{openai}. Summing over intervals at \(p=4\) introduces traces of products of square integrals from different frequency scales. For nearby scales we estimate each square integral by the single-annulus bound. When the scales are far apart, Lemma~\ref{stein:four} estimates the error in commuting right multiplication by the lower-frequency square integral with the higher-frequency half-wave. The error involves a spatial increment of the lower-frequency square integral. Section~\ref{stein:across-scales} also gives a positive kernel estimate at finite \(q>2\) for the interpolation and the smooth amplitude estimates in Section~\ref{stein:amplitude}.

\subsection{The radius derivative of Bochner--Riesz means}
For real parameters \(\nu>0\) and \(r>0\), let \(B_r^\nu:L^2(\N_{\R^2})\to L^2(\N_{\R^2})\) be the spatial multiplier with symbol
\((1-|\xi|^2/r^2)_+^\nu\). On the linear span of tensors \(F=f\otimes a\), with \(f\in\mathscr S(\R^2)\) and \(a\in\M\) bounded with support projections of finite trace, define \(D_r^\nu\) by Fourier inversion with the scalar symbol
\begin{equation}\label{stein:derivative}
 d_r^\nu(\xi)=2\nu\,\frac{|\xi|^2}{r^2}
        (1-|\xi|^2/r^2)^{\nu-1}\mathbf1_{\{|\xi|<r\}},
        \qquad \xi\in\R^2 .
\end{equation}
The symbol is zero for \(|\xi|\ge r\). The power \(\nu-1>-1\) gives local integrability at the boundary, so \(d_r^\nu\widehat f\in L^1(\R^2)\). Thus \(D_r^\nu F(x)=(d_r^\nu\widehat f)^\vee(x)a\) for \(x\in\R^2\), and the linear map on the tensor span takes values in \(C_0(\R^2;\M)\).

For each \(\xi\in\R^2\), the symbol of \(B_r^\nu\) is locally absolutely continuous in \(r>0\), with derivative \(d_r^\nu(\xi)/r\) almost everywhere. For \(0<r_-<r_+<\infty\) we have
\begin{equation*}
 0\le \int_{r_-}^{r_+} d_r^\nu(\xi)\,\frac{dr}{r}
 = (1-|\xi|^2/r_+^2)_+^\nu-(1-|\xi|^2/r_-^2)_+^\nu\le1\,.
\end{equation*}
Together with \(\widehat f\in L^1(\R^2)\), this permits Fubini in the Fourier formula after testing in \(r\). For \(\varphi\in C_c^\infty((0,\infty))\), integration by parts in \(r\) gives
\begin{equation*}
 -\int_0^\infty B_r^\nu F(x)\varphi'(r)\,dr
 =\int_0^\infty D_r^\nu F(x)\varphi(r)\,\frac{dr}{r}\,.
\end{equation*}
Hence \(D_r^\nu F=r\partial_r B_r^\nu F\) as distributions in \(r\).

Let \(\mathcal A\) be a semifinite tracial algebra and \(2\le p<\infty\). Take a measure space \((S,\mu)\) with separable Hilbert space \(H=L^2(S,\mu)\). For a simple family \(X=(X_s)_{s\in S}\) of elements of \(L^p(\mathcal A)\) supported on a set of finite \(\mu\)-measure, set
\begin{equation*}
 \|X\|_{p,c}=\left\|\left(\int_S X_s^*X_s\,d\mu(s)\right)^{1/2}\right\|_{L^p(\mathcal A)},\qquad
 \|X\|_{p,r}=\left\|\left(\int_S X_sX_s^*\,d\mu(s)\right)^{1/2}\right\|_{L^p(\mathcal A)}\,.
\end{equation*}
The completions are denoted \(L^p(\mathcal A;H_c)\) and \(L^p(\mathcal A;H_r)\). Their intersection has norm \(\|X\|_{p,rc}=\max\{\|X\|_{p,c},\|X\|_{p,r}\}\), with the families identified by their coefficients in the same orthonormal basis of \(H\).

For the matrix realization, use \(B(\ell^2(\mathbb N))\bar\otimes\mathcal A\) with the product of the usual operator trace and the trace of \(\mathcal A\). Let \(e_{ij}\) be the matrix units on \(\ell^2(\mathbb N)\). For a finite rank family \(X_s=\sum_{i=1}^n h_i(s)x_i\), where \(h_1,\ldots,h_n\) are real orthonormal functions in \(H\) and \(x_i\in L^p(\mathcal A)\), write \(C_X=\sum_i e_{i1}\otimes x_i\) and \(R_X=\sum_i e_{1i}\otimes x_i\). Orthogonality gives
\begin{align*}
 C_X^*C_X &= e_{11}\otimes\sum_{i=1}^n x_i^*x_i
           = e_{11}\otimes\int_S X_s^*X_s\,d\mu(s)\,,\\
 R_XR_X^* &= e_{11}\otimes\sum_{i=1}^n x_i x_i^*
           = e_{11}\otimes\int_S X_sX_s^*\,d\mu(s)\,.
\end{align*}
Thus the norms of \(C_X\) and \(R_X\) are \(\|X\|_{p,c}\) and \(\|X\|_{p,r}\), respectively. Completion in a fixed real orthonormal basis identifies the two spaces with the corresponding first-column and first-row subspaces~\cite[Section 2.1]{hlw}.

We use the parameter norms to sum the derivative family. The positive square integrals extend to \(L^{p/2}(\mathcal A)\) by the following continuity estimates, stated first for simple families \(X,Y\).
\begin{equation}\label{stein:square-continuity}
\begin{aligned}
 \left\|\int_S (X_s^*X_s-Y_s^*Y_s)\,d\mu(s)\right\|_{L^{p/2}(\mathcal A)}
   &\le \|X-Y\|_{p,c}\bigl(\|X\|_{p,c}+\|Y\|_{p,c}\bigr),\\
 \left\|\int_S (X_sX_s^*-Y_sY_s^*)\,d\mu(s)\right\|_{L^{p/2}(\mathcal A)}
   &\le \|X-Y\|_{p,r}\bigl(\|X\|_{p,r}+\|Y\|_{p,r}\bigr).
\end{aligned}
\end{equation}
For the column estimate write \(X_s^*X_s-Y_s^*Y_s=(X_s-Y_s)^*X_s+Y_s^*(X_s-Y_s)\). To bound each mixed product, take two simple families \(U,V\) and use the normalized indicators of their common finite partition on their supports as parameter functions. If \(C_U,C_V\) are their first-column matrices in this basis, then
\[
 C_U^*C_V=e_{11}\otimes\int_S U_s^*V_s\,d\mu(s).
\]
H\"older's inequality for the trace in the matrix algebra tensored with \(\mathcal A\) bounds the \(L^{p/2}(\mathcal A)\) norm of the integral by \(\|C_U\|_{L^p(B(\ell^2(\mathbb N))\bar\otimes\mathcal A)}\|C_V\|_{L^p(B(\ell^2(\mathbb N))\bar\otimes\mathcal A)}=\|U\|_{p,c}\|V\|_{p,c}\).
Taking \((U,V)=(X-Y,X)\) and \((Y,X-Y)\) bounds the two terms by \(\|X-Y\|_{p,c}\|X\|_{p,c}\) and \(\|Y\|_{p,c}\|X-Y\|_{p,c}\), respectively. Their sum gives the column inequality. Applying it to the adjoint families \((X_s^*)_s,(Y_s^*)_s\) gives the row inequality, since taking adjoints sends the row norm to the column norm. Thus a Cauchy sequence in either family norm has a Cauchy sequence of positive integrals in \(L^{p/2}(\mathcal A)\). By the same estimate its limit is independent of the approximating simple families, and \eqref{stein:square-continuity} passes to each completion.

For the radius parameter take \(S=(0,\infty)\), \(d\mu(r)=dr/r\), and set
\[
 H=L^2((0,\infty),dr/r),\qquad
 \mathcal T_\nu F=(D_r^\nu F)_{r>0}=(r\partial_r B_r^\nu F)_{r>0}
\]
on the tensor span above, with the derivative understood distributionally in \(r\). The column and row square functions are the positive square roots
\[
 G_c^\nu F=\left(\int_0^\infty (D_r^\nu F)^*D_r^\nu F\,\frac{dr}{r}\right)^{1/2},
 \qquad
 G_r^\nu F=\left(\int_0^\infty D_r^\nu F(D_r^\nu F)^*\,\frac{dr}{r}\right)^{1/2}
\]
whenever the family belongs to the corresponding parameter spaces. In that case \(\|G_c^\nu F\|_{L^p(\N_{\R^2})}=\|\mathcal T_\nu F\|_{p,c}\) and \(\|G_r^\nu F\|_{L^p(\N_{\R^2})}=\|\mathcal T_\nu F\|_{p,r}\), respectively.

We first prove Theorem~\ref{thm:stein} at \(p=2\) by computing \(\int_0^\infty |d_r^\nu(\xi)|^2\,dr/r\) at each nonzero frequency \(\xi\).
For \(0<\nu<1\) and fixed \(r>0\), the symbol \(d_r^\nu(\xi)\) is unbounded as \(|\xi|\uparrow r\), so the multiplier at a single radius is unbounded on \(L^2(\N_{\R^2})\). For \(\nu>1/2\) we define the family jointly in the radius, as a map \(\mathcal T_\nu:L^2(\N_{\R^2})\to L^2((0,\infty),dr/r;L^2(\N_{\R^2}))\). For \(F\) in its domain, take the inverse spatial Fourier transform of the joint function \((r,\xi)\mapsto d_r^\nu(\xi)\widehat F(\xi)\) in the Hilbert space with measure \(dr/r\,d\xi\) and coefficients in \(L^2(\M)\). To check square integrability, for \(\xi\ne0\) the substitution \(u=|\xi|^2/r^2\) gives \(dr/r=-du/(2u)\) and
\begin{equation}\label{stein:mass}
 \int_0^\infty |d_r^\nu(\xi)|^2\,\frac{dr}{r}
 =2\nu^2\int_0^1 u(1-u)^{2\nu-2}\,du=\frac{\nu}{2\nu-1}.
\end{equation}
To evaluate the scalar integral in \eqref{stein:mass}, put \(w=1-u\). The integral of \(u(1-u)^{2\nu-2}\) is \(1/(2\nu-1)-1/(2\nu)\). The frequency zero has Lebesgue measure zero. In either product order the trace of the square integral is \(\int_0^\infty\|D_r^\nu F\|_{L^2(\N_{\R^2})}^2\,dr/r\) by the \(p=2\) matrix identification. Tonelli's theorem and Plancherel's identity in the spatial variables therefore give
\begin{align*}
 \|G_c^\nu F\|_{L^2(\N_{\R^2})}^2
 &=\|G_r^\nu F\|_{L^2(\N_{\R^2})}^2
  =\|\mathcal T_\nu F\|_{L^2((0,\infty),dr/r;\,L^2(\N_{\R^2}))}^2\\
 &=\frac{\nu}{2\nu-1}\|F\|_{L^2(\N_{\R^2})}^2 .
\end{align*}
Thus the Hilbert space definition agrees with Fourier inversion on these elementary tensors and determines the derivative family up to equality almost everywhere in the radius.
For \(0<\nu\le1/2\) the integral diverges at \(u=1\) at every nonzero frequency. For \(p>2\) we decompose the derivative by distance to the boundary \(|\xi|=r\). Proposition~\ref{stein:annuli} bounds each smooth boundary piece simultaneously on all dyadic intervals of radii. Its dependence on the thickness will determine which orders \(\nu\) can be summed in Section~\ref{stein:summation}.

\subsection{One frequency annulus}\label{stein:one-scale}
Write $\mathcal A=\N_{\R^2}$ and let $\Phi$ be its product trace. On a dyadic interval of radii $r=2^ks$, $k\in\mathbb Z$ and $1\le s<2$, the normalized signed distance to the boundary is $s-2^{-k}|\xi|$. Fix $0<c_0<C_0<\infty$ and a profile $\psi\in C_c^\infty(\R)$. For $0<\delta\le1/2$, $k\in\mathbb Z$ and $s\in\R$, define the spatial Fourier multiplier
\begin{equation*}
 V_{k,s}=\psi\big((s-2^{-k}|D|)/\delta\big),\qquad L_\delta=1+\log(1/\delta).
\end{equation*}
Each $V_{k,s}:L^2(\mathcal A)\to L^2(\mathcal A)$ is initially defined by its bounded symbol. We first bound the row and column square functions of $V_{0,s}$ with parameter measure $ds$. We then sum over $k$ with counting measure. The dependence on $\psi$ enters through $B_\psi=\sup_{v\in\R}(1+|v|)^{12}|\widehat\psi(v)|$, which is bounded by a constant times $\|\psi\|_{L^1(\R)}+\|\psi^{(12)}\|_{L^1(\R)}$. By homogeneity we may take $B_\psi\le1$ throughout the estimates for $V_{k,s}$.

\begin{lemma}\label{stein:one}
If $2\le p<\infty$ and $h\in L^p(\mathcal A)$ has Fourier support in $c_0\le|\xi|\le C_0$, then
\begin{equation}\label{stein:one-bound}
 \|(V_{0,s}h)_{s\in\R}\|_{p,rc}
 \lesssim_{p,c_0,C_0} B_\psi\,\delta^{1/2-\sigma(p)}L_\delta^{\beta(p)}\|h\|_{L^p(\mathcal A)}
\end{equation}
with Lebesgue measure $ds$.
\end{lemma}
\begin{proof}
We use Plancherel's identity in $s$ to express the norm in~\eqref{stein:one-bound} in terms of the time family $X_t=w_\delta(t)^{1/2}U(-t)h$, where
\begin{equation*}
 w_\delta(t)=\frac{\delta^2}{2\pi}|\widehat\psi(\delta t)|^2
 \qquad(t\in\R)\,.
\end{equation*}
First take $h$ to be a finite sum of scalar Schwartz functions with annular Fourier support and bounded coefficients whose support projections have finite trace. Choose $\chi\in C_c^\infty((0,\infty))$ equal to one on a neighborhood of $[c_0,C_0]$. By~\eqref{pre:kern}, the kernel of the symbol $e^{-it|\xi|}\chi(|\xi|)$ has $L^1(\R^2)$ norm at most $C_\chi(1+|t|)^4$ and depends continuously on $t$ in $L^1(\R^2)$. Since $\chi(|D|)h=h$, convolution with these kernels makes $t\mapsto U(-t)h$ norm continuous in $L^p(\mathcal A)$, hence strongly measurable. The rapid decay of $\widehat\psi(\delta t)$ gives absolute integrability. Fourier inversion in the profile variable then yields the Bochner integral
\begin{equation*}
 V_{0,s}h=\frac{\delta}{2\pi}\int_\R e^{ist}\widehat\psi(\delta t)U(-t)h\,dt\,.
\end{equation*}
Represent $\M$ faithfully on a Hilbert space. For each $x\in\R^2$ and vector $\upsilon$ in this space, the function $t\mapsto\widehat\psi(\delta t)(U(-t)h)(x)\upsilon$ is continuous and square integrable by the same kernel bound applied in the operator norm. Hilbert-space-valued Plancherel and equality of the quadratic forms give
\begin{equation}\label{stein:planch}
 \int_\R |V_{0,s}h|^2\,ds
 =\int_\R w_\delta(t)|U(-t)h|^2\,dt\,.
\end{equation}
Applying Plancherel's identity to $(V_{0,s}h)(x)^*\upsilon$ reverses the Fourier phase and gives the corresponding row identity with weight $w_\delta(t)$. Thus the left side of~\eqref{stein:one-bound} equals $\|(X_t)_{t\in\R}\|_{p,rc}$.

To apply the spacetime estimate~\eqref{cone:long-time}, we bound $\|(X_t)_{t\in J}\|_{p,rc}$ on dyadic time intervals $J$. For $u\in L^p(J;L^p(\mathcal A))$ on a finite time interval $J\subset\R$, the triangle inequality in $L^{p/2}(\mathcal A)$ and scalar H\"older give
\begin{equation*}
 \|(u(t))_{t\in J}\|_{p,rc}
 \le\bigg(\int_J\|u(t)\|_{L^p(\mathcal A)}^2\,dt\bigg)^{1/2}
 \le |J|^{1/2-1/p}\bigg(\int_J\|u(t)\|_{L^p(\mathcal A)}^p\,dt\bigg)^{1/p}.
\end{equation*}
For $R=2^n$ with $n\ge0$, let $J_R^+=[R,2R]$ and $J_R^-=[-2R,-R]$. The preceding inequality and~\eqref{cone:long-time} give, for $J=J_R^+$ or $J=J_R^-$,
\begin{align*}
 \|(X_t)_{t\in J}\|_{p,rc}
 &\le \sup_{t\in J}w_\delta(t)^{1/2}\,R^{1/2-1/p}
            \|U(-t)h\|_{L^p(J;L^p(\mathcal A))}\\
 &\lesssim \delta B_\psi(1+\delta R)^{-12}
          R^{1/2+\sigma(p)}[\log(2+R)]^{\beta(p)}\|h\|_{L^p(\mathcal A)}\,.
\end{align*}
On $J_R^-$, the substitution $s=-t$ gives the spacetime norm of $U(s)h$ on $[R,2R]$. On $|t|\le1$, Lemma~\ref{cone:long-wave} gives $\|(X_t)_{|t|\le1}\|_{p,rc}\lesssim\delta B_\psi\|h\|_{L^p(\mathcal A)}$.

To sum the time norms, choose $n_0\in\mathbb Z$ with $2^{n_0}\le\delta^{-1}<2^{n_0+1}$ and set $a=1/2+\sigma(p)$. For $0\le n\le n_0$ the logarithm is bounded by $CL_\delta$ and $\sum_{n=0}^{n_0}2^{na}\lesssim\delta^{-a}$. For $n=n_0+j$, $j\ge1$, we have $2^{na}\le\delta^{-a}2^{ja}$, $(1+\delta2^n)^{-12}\lesssim2^{-12j}$ and $\log(2+2^n)\lesssim L_\delta(1+j)$. Extending each time interval by zero and applying the triangle inequality therefore gives
\begin{equation*}
 \|(X_t)_{t\in\R}\|_{p,rc}
 \lesssim B_\psi\|h\|_{L^p(\mathcal A)}\,\delta^{1-a}L_\delta^{\beta(p)}
 \left(1+\sum_{j\ge1}2^{-j(12-a)}(1+j)^{\beta(p)}\right)\,.
\end{equation*}
Here $1/2\le a\le1$, so the series converges and $1-a=1/2-\sigma(p)$ gives~\eqref{stein:one-bound} on the tensor core.

For general $h$, take tensor approximants $h_n\to h$ in $L^p(\mathcal A)$ as in~\eqref{pre:localized-core}, with Fourier supports in a fixed slightly larger annulus. Choose $\chi$ to reproduce this annulus. For fixed $\delta$, the kernels of the symbols $\chi(|\xi|)\psi((s-|\xi|)/\delta)$ are continuous in $s$ in $L^1(\R^2)$ by~\eqref{pre:kern} and vanish outside a common compact $s$-interval. Thus $(V_{0,s}h_n)_s$ converges to $(V_{0,s}h)_s$ in $L^2(ds;L^p(\mathcal A))$, hence in both parameter norms. For the time variable, apply~\eqref{pre:kern} to $e^{-it|\xi|}\chi(|\xi|)$ and convolve with $h_n-h$. This gives
\begin{equation*}
 \|w_\delta^{1/2}U(-\cdot)(h_n-h)\|_{L^2(dt;L^p(\mathcal A))}
 \le C_\chi\left(\int_\R w_\delta(t)(1+|t|)^8\,dt\right)^{1/2}
      \|h_n-h\|_{L^p(\mathcal A)}\,.
\end{equation*}
The integral is finite for each fixed $\delta$ by the rapid decay of $\widehat\psi$, so the weighted families also converge in both parameter norms. Applying~\eqref{stein:square-continuity} to these limits gives both square identities in $L^{p/2}(\mathcal A)$ and extends~\eqref{stein:one-bound} to $h$.
\end{proof}

\subsection{Square-function estimates across frequency scales}\label{stein:across-scales}

On one scale Lemma~\ref{stein:one} gives the factor $\delta^{1/2}L_\delta^{1/2}$ at $p=4$. On a finite set $E\subset\mathbb Z$ of dyadic scales, the estimate below has the factor $\delta^{1/2}L_\delta^{3/4}$, uniformly in $E$.

\begin{lemma}\label{stein:four}
Let $E\subset\mathbb Z$ be finite, and let $F_k\in L^4(\mathcal A)$ have Fourier support in $c_02^k\le|\xi|\le C_02^k$ for each $k\in E$. Then for $d\in\{\mathrm c,\mathrm r\}$,
\begin{equation}\label{stein:glue4}
 \|(V_{k,s}F_k)_{k\in E,s\in\R}\|_{4,d}
 \lesssim_{c_0,C_0} B_\psi\,\delta^{1/2}L_\delta^{3/4}\|(F_k)_{k\in E}\|_{4,d}.
\end{equation}
\end{lemma}
\begin{proof}
It suffices to estimate the column norm on the annular tensor core and then apply the estimate to $(F_k^*)_{k\in E}$ with profile $\overline\psi$. With $B_\psi\le1$, put
\begin{equation*}
 H_k=\int_\R |V_{k,s}F_k|^2\,ds,\quad a_k=\|F_k\|_{L^4(\mathcal A)}^2,\quad
 Y=\Big\|\sum_{k\in E}H_k\Big\|_{L^2(\mathcal A)},\quad W=\|(F_k)_{k\in E}\|_{4,c}.
\end{equation*}
The target is $Y^2\lesssim\delta^2 L_\delta^3 W^4$. For positive $A,B\in L^2(\mathcal A)$, $\Phi(AB)=\Phi(B^{1/2}AB^{1/2})\ge0$. Write $A_k=|F_k|^2$. Expanding the two squared $L^2(\mathcal A)$ norms using the trace gives
\begin{equation*}
 Y^2=\sum_{k,\ell\in E}\Phi(H_kH_\ell),\qquad
 W^4=\sum_{k,\ell\in E}\Phi(A_k A_\ell)\ge\sum_{k\in E}\Phi(A_k^2)=\sum_{k\in E}a_k^2.
\end{equation*}
Each trace in these sums is nonnegative. By trace cyclicity, for any positive integer $L$ the expansion of $Y^2$ becomes
\[
 Y^2=\sum_{\substack{k,\ell\in E\\|k-\ell|\le L}}\Phi(H_kH_\ell)
     +2\sum_{\substack{k,\ell\in E\\k<\ell-L}}\Phi(H_kH_\ell).
\]
We choose $L$ only after estimating both sums. In the first sum the single-annulus estimate suffices. In the second we compare $H_\ell$ with $\delta\|\psi\|_{L^2(\R)}^2 A_\ell$ inside the trace, using the frequency separation between $H_k$ and $F_\ell$.

\medskip
For pairs with $|k-\ell|\le L$, we use the bound on each $H_k$ from Lemma~\ref{stein:one}. For $h(x)=F_k(2^{-k}x)$, the radial Fourier support of $h$ lies in $[c_0,C_0]$ and $V_{0,s}h(x)=(V_{k,s}F_k)(2^{-k}x)$. Cancelling the common factor $2^{k/2}$ in $\|h\|_{L^4(\mathcal A)}$ and $\|(V_{0,s}h)_{s\in\R}\|_{4,c}$ in~\eqref{stein:one-bound} gives $\|H_k\|_{L^2(\mathcal A)}\lesssim\delta L_\delta a_k$ after squaring. Extend $a_k$ by zero outside $E$. Cauchy--Schwarz for the trace and then for the sum in $k$ gives
\begin{equation}\label{stein:near}
 \sum_{|k-\ell|\le L}\Phi(H_kH_\ell)
 \lesssim\delta^2L_\delta^2\sum_{|n|\le L}\sum_k a_ka_{k+n}
 \lesssim(2L+1)\delta^2L_\delta^2\sum_k a_k^2.
\end{equation}

\medskip
For more widely separated scales, we express $\Phi(H_kH_\ell)$ using the half-wave at the higher frequency. For $k,\ell\in E$ with $k<\ell-L$ and $t\in\R$, write $U_{\ell,t}=U(-2^{-\ell}t)$, a unitary operator on $L^2(\mathcal A)$. Rescaling~\eqref{stein:planch} gives $H_\ell=\int_\R w_\delta(t)|U_{\ell,t}F_\ell|^2\,dt$. The polynomial kernel bound at the normalized scale gives $\int_\R w_\delta(t)\|U_{\ell,t}F_\ell\|_{L^4(\mathcal A)}^2\,dt<\infty$, so the continuous functional $X\mapsto\Phi(H_kX)$ on $L^2(\mathcal A)$ passes under this integral. Using $\int_\R w_\delta(t)\,dt=\delta\|\psi\|_{L^2(\R)}^2$ by scalar Plancherel, we obtain
\begin{align*}
 \Phi(H_kH_\ell)&=\delta\|\psi\|_{L^2(\R)}^2\Phi(H_kA_\ell)+e_{k,\ell}\,,\\
 e_{k,\ell}&=\int_\R w_\delta(t)
   \big[\Phi(H_k|U_{\ell,t}F_\ell|^2)-\Phi(H_kA_\ell)\big]\,dt\,.
\end{align*}
We bound $e_{k,\ell}$ using spatial increments of $H_k$. The Fourier support of $|V_{k,s}F_k|^2$ lies in $|\xi|\le2C_02^k$ by convolution of the reflected annulus with the annulus, and integration in $s$ preserves this support. Hence $\|\nabla H_k\|_{L^2(\mathcal A)^2}\le2C_02^k\|H_k\|_{L^2(\mathcal A)}$ by Plancherel in the spatial variables. The product $F_\ell H_k$ has Fourier support in the radial interval
\begin{equation*}
 [c_02^\ell-2C_02^k,\,C_02^\ell+2C_02^k]
 \subset [(c_0/2)2^\ell,\,2C_02^\ell]\,,
\end{equation*}
provided $L\ge L_0$, where $L_0$ is a positive integer with $2C_02^{-L_0}<c_0/2$. A fixed real $\kappa\in C_c^\infty((0,\infty))$ equal to one on $[c_0/2,2C_0]$ therefore makes $\kappa(2^{-\ell}|D|)$ reproduce both $F_\ell$ and $F_\ell H_k$. Unitarity on $L^2(\mathcal A)$ gives
\begin{equation*}
 \Phi\big(H_k|U_{\ell,t}F_\ell|^2\big)-\Phi\big(H_k|F_\ell|^2\big)
 =\big\langle U_{\ell,t}F_\ell,\,
 (U_{\ell,t}F_\ell)H_k-U_{\ell,t}(F_\ell H_k)\big\rangle .
\end{equation*}
Indeed, for the trace inner product $\langle X,Y\rangle=\Phi(X^*Y)$, cyclicity gives $\langle U_{\ell,t}F_\ell,(U_{\ell,t}F_\ell)H_k\rangle=\Phi(H_k|U_{\ell,t}F_\ell|^2)$. The other inner product equals $\langle F_\ell,F_\ell H_k\rangle=\Phi(H_k|F_\ell|^2)$ by unitarity. Both $F_\ell$ and $F_\ell H_k$ belong to $L^2(\mathcal A)$ on the bounded tensor core.

We pair $U_{\ell,t}F_\ell$ in $L^4(\mathcal A)$ with $(U_{\ell,t}F_\ell)H_k-U_{\ell,t}(F_\ell H_k)$ in $L^{4/3}(\mathcal A)$. Put $K_t=(e^{-it|\cdot|}\kappa(|\cdot|))^\vee$ and $K_{\ell,t}(y)=2^{2\ell}K_t(2^\ell y)$. Three integrations with $1-\Delta_\xi$ on the fixed annulus give
$|K_t(y)|\lesssim (1+|t|)^6(1+|y|)^{-6}$. In particular
\begin{equation*}
 \|K_{\ell,t}\|_{L^1(\R^2)}+2^\ell\int_{\R^2}|y|\,|K_{\ell,t}(y)|\,dy\lesssim (1+|t|)^6.
\end{equation*}
The kernel formula gives
\begin{equation*}
 \big[(U_{\ell,t}F_\ell)H_k-U_{\ell,t}(F_\ell H_k)\big](x)
 =\int_{\R^2}K_{\ell,t}(y)F_\ell(x-y)\big[H_k(x)-H_k(x-y)\big]\,dy\,.
\end{equation*}
Integrating $\nabla H_k$ along the segment from $x-y$ to $x$ bounds $\|H_k-H_k(\cdot-y)\|_{L^2(\mathcal A)}$ by $|y|\|\nabla H_k\|_{L^2(\mathcal A)^2}$. H\"older's inequality for the product of $F_\ell(\cdot-y)$ and this increment now gives
\begin{align*}
 \|(U_{\ell,t}F_\ell)H_k-U_{\ell,t}(F_\ell H_k)\|_{L^{4/3}(\mathcal A)}
 &\le \|F_\ell\|_{L^4(\mathcal A)}\|\nabla H_k\|_{L^2(\mathcal A)^2}
       \int_{\R^2}|y|\,|K_{\ell,t}(y)|\,dy\\
 &\lesssim 2^{k-\ell}(1+|t|)^6\|F_\ell\|_{L^4(\mathcal A)}\|H_k\|_{L^2(\mathcal A)} .
\end{align*}
Since $\|K_{\ell,t}\|_{L^1(\R^2)}\lesssim(1+|t|)^6$, convolution gives
$\|U_{\ell,t}F_\ell\|_{L^4(\mathcal A)}\le C(1+|t|)^6\|F_\ell\|_{L^4(\mathcal A)}$.
Trace H\"older bounds $|\Phi(H_k|U_{\ell,t}F_\ell|^2)-\Phi(H_kA_\ell)|$ by $C2^{k-\ell}(1+|t|)^{12}\|H_k\|_{L^2(\mathcal A)}a_\ell$. To integrate this bound against $w_\delta(t)$, set $v=\delta t$. Since $1+|v|/\delta\le(1+|v|)/\delta$ and $B_\psi\le1$, we obtain
\begin{equation*}
 \int_\R w_\delta(t)(1+|t|)^{12}\,dt
 \le \frac{\delta^{-11} B_\psi^2}{2\pi}\int_\R(1+|v|)^{-12}\,dv
 \lesssim \delta^{-11}\,.
\end{equation*}
It follows that
\begin{equation}\label{stein:far}
 |e_{k,\ell}|\lesssim 2^{k-\ell}\delta^{-11}\|H_k\|_{L^2(\mathcal A)} a_\ell\,.
\end{equation}

\medskip
To bound the sum of $\Phi(H_kA_\ell)$ over $k<\ell-L$, we use positivity to include all pairs in $E\times E$. H\"older's inequality for the trace then gives
\begin{equation*}
 \sum_{k<\ell-L}\Phi(H_kA_\ell)
 \le \Phi\bigg(\Big(\sum_{k\in E}H_k\Big)\Big(\sum_{\ell\in E}A_\ell\Big)\bigg)
 \le Y\Big\|\sum_{\ell\in E}A_\ell\Big\|_{L^2(\mathcal A)}=YW^2.
\end{equation*}
For the error terms sum first at each displacement $n=\ell-k$. Using $\|H_k\|_{L^2(\mathcal A)}\lesssim\delta L_\delta a_k$ in~\eqref{stein:far} and $\sum_k a_k a_{k+n}\le\sum_k a_k^2\le W^4$, we obtain
\begin{equation*}
 \sum_{k<\ell-L}|e_{k,\ell}|
 \lesssim\delta^{-10} L_\delta\sum_{n>L}2^{-n}\sum_k a_k a_{k+n}
 \lesssim \delta^{-10}L_\delta 2^{-L} W^4.
\end{equation*}
Choose $L=L_0+\lceil14\log_2(2/\delta)\rceil$. Then $2L+1\lesssim_{c_0,C_0}L_\delta$ in~\eqref{stein:near}, and $2^{-L}\le(\delta/2)^{14}$ gives $\delta^{-10}L_\delta2^{-L}\lesssim\delta^2$. Substituting the nearby bound and the sum of the separated terms into the expansion of $Y^2$ yields
\begin{equation*}
 Y^2\le C\delta YW^2+C\delta^2L_\delta^3W^4
 \le \tfrac12Y^2+\tfrac12 C^2\delta^2W^4+C\delta^2L_\delta^3W^4 .
\end{equation*}
It follows first that $Y\lesssim\delta L_\delta^{3/2} W^2$ and then that $\|(V_{k,s}F_k)_{k,s}\|_{4,c}=Y^{1/2}\lesssim\delta^{1/2}L_\delta^{3/4}W$.

For the row norm use $\widehat{F_k^*}(\xi)=\widehat F_k(-\xi)^*$ and $|-\xi|=|\xi|$. These give
\begin{equation*}
 (V_{k,s}F_k)^*
 =\overline\psi\big((s-2^{-k}|D|)/\delta\big)F_k^* .
\end{equation*}
The profile $\overline\psi$ has the same $B_\psi$ because $\widehat{\overline\psi}(v)=\overline{\widehat\psi(-v)}$. Applying the column estimate with $\overline\psi$ to $(F_k^*)_{k\in E}$ gives the row estimate by $\|X^*\|_{4,c}=\|X\|_{4,r}$ for $X=(F_k)_{k\in E}$ and $X=(V_{k,s}F_k)_{k\in E,s\in\R}$.

To remove the tensor assumption, use~\eqref{pre:localized-core} on each $F_k$ with the same relative enlargement of its annulus. For each $d\in\{\mathrm c,\mathrm r\}$, the triangle inequality for the sum of squares bounds $\|(F_k^{(n)}-F_k)_{k\in E}\|_{4,d}$ by $(\sum_{k\in E}\|F_k^{(n)}-F_k\|_{L^4(\mathcal A)}^2)^{1/2}$, which tends to zero since $E$ is finite. The family $(V_{k,s}F_k^{(n)})_{k\in E,s\in\R}$ converges in $\|\cdot\|_{4,d}$ to $(V_{k,s}F_k)_{k\in E,s\in\R}$ by the fixed-scale approximation in Lemma~\ref{stein:one}. This proves~\eqref{stein:glue4}.
\end{proof}

\medskip
We also need a bound uniform over finite sets of scales for each $2<q<\infty$. For $b>0$, define the positive scalar kernel $R_b:\R^2\to(0,\infty)$ by $R_b(y)=b^2(1+b|y|)^{-4}$. We bound the column square integral $\int_\R |V_{k,s}F_k|^2\,ds$ by a multiple of $R_{\delta2^k}*|F_k|^2$ and sum these majorants by duality.

\begin{lemma}\label{stein:coarse}
For $2<q<\infty$, $d\in\{\mathrm c,\mathrm r\}$, and finite annular families $F_k\in L^q(\mathcal A)$ with the supports of Lemma~\ref{stein:four},
\begin{equation}\label{stein:coarse-bound}
 \|(V_{k,s}F_k)_{k\in E,s\in\R}\|_{q,d}\lesssim_{q,c_0,C_0}B_\psi\,\|(F_k)_{k\in E}\|_{q,d}.
\end{equation}
\end{lemma}
\begin{proof}
With $H_k=\int_\R |V_{k,s}F_k|^2\,ds$ and $B_\psi\le1$, the column estimate follows from the two inequalities
\[
 H_k\lesssim R_{\delta2^k}*|F_k|^2,\qquad
 \Big\|\sum_{k\in E} R_{\delta2^k}*|F_k|^2\Big\|_{L^{q/2}(\mathcal A)}
 \lesssim_q\Big\|\sum_{k\in E}|F_k|^2\Big\|_{L^{q/2}(\mathcal A)}.
\]
We first establish the pointwise operator inequality on the tensor core. For the annular kernel $K_t=(e^{-it|\cdot|}\kappa(|\cdot|))^\vee$ from the proof of Lemma~\ref{stein:four} we need, for each integer $m\ge1$,
\begin{equation}\label{stein:wavefront}
 |K_t(y)|\lesssim_m (1+|y|)^{-1/2}(1+\big||t|-|y|\big|)^{-m},
 \qquad \|K_t\|_{L^1(\R^2)}\lesssim (1+|t|)^{1/2}.
\end{equation}
Put $u=|y|$. The angular integral in polar coordinates is $2\pi J_0(u\rho)$, so our Fourier convention gives
\begin{equation*}
 K_t(y)=\frac1{2\pi}\int_0^\infty e^{-it\rho}J_0(u\rho)\kappa(\rho)\rho\,d\rho .
\end{equation*}
For bounded $u$, integration by parts in $\rho$ gives $O_m((1+|t|)^{-m})$, since all the radial derivatives of $J_0(u\rho)$ are uniformly bounded on $\operatorname{supp}\kappa$. For $u$ sufficiently large, $u\rho\ge1$ on the support of $\kappa$. By~\eqref{max:hankel} at order zero, $J_0(u\rho)$ is the sum of $e^{\pm iu\rho}$ times amplitudes whose $\rho$-derivatives through order $m$ are $O_m(u^{-1/2})$ on this support. Integration by parts now yields the decays $(1+|t\mp u|)^{-m}$.

To bound $\|K_t\|_{L^1(\R^2)}$ in~\eqref{stein:wavefront}, take $m=2$. In the radial integral use $u(1+u)^{-1/2}\le(1+u)^{1/2}$ and $1+u\le(1+|t|)(1+|u-|t||)$. This gives
\begin{equation*}
 \|K_t\|_{L^1(\R^2)}\lesssim (1+|t|)^{1/2}\int_0^\infty(1+|u-|t||)^{-3/2}\,du
 \lesssim (1+|t|)^{1/2}.
\end{equation*}
To bound $H_k$ using~\eqref{stein:planch}, we prove the weighted kernel estimate
\begin{equation}\label{stein:kernel-major}
 \delta^2\int_\R |\widehat\psi(\delta t)|^2\|K_t\|_{L^1(\R^2)} |K_t(y)|\,dt
 \lesssim B_\psi^2 \delta^2(1+\delta |y|)^{-4}.
\end{equation}
To prove~\eqref{stein:kernel-major}, write
$v=\big||t|-|y|\big|$. Then $1+\delta|y|\le(1+\delta|t|)(1+v)$ and $(1+|t|)/(1+|y|)\le1+v$. Taking $m=8$ in~\eqref{stein:wavefront} gives
\[
 |\widehat\psi(\delta t)|^2\|K_t\|_{L^1(\R^2)}|K_t(y)|
 \lesssim B_\psi^2(1+\delta|t|)^{-24}(1+v)^{-8+1/2}
 \le B_\psi^2(1+\delta|y|)^{-4}(1+v)^{-8+4+1/2}.
\]
The integral of $(1+\big||t|-|y|\big|)^{-7/2}$ in $t\in\R$ is bounded independently of $y$. Multiplying by $\delta^2$ and integrating proves~\eqref{stein:kernel-major}.

For a scalar kernel $K\in L^1(\R^2)$ and a bounded test tensor sum $G$, the operator inequality $|K*G|^2\le\|K\|_{L^1(\R^2)} (|K|*|G|^2)$ follows by applying the Cauchy--Schwarz inequality to $(K*G)(x)$ acting on a Hilbert space vector. Apply it with $K=K_{k,t}$, multiply by the time weight in~\eqref{stein:planch}, and use~\eqref{stein:kernel-major} with $y$ replaced by $2^ky$. The spatial rescaling contributes a factor $2^{2k}$, giving $H_k\lesssim (\delta2^k)^2(1+\delta2^k|\cdot|)^{-4}*|F_k|^2$.

To sum the bounds for $H_k$, we allow arbitrary scales $b_k>0$, $k\in E$. This form of the estimate will also handle the smooth amplitude in Proposition~\ref{stein:annuli}. For $1<v<\infty$, nonnegative $A_k\in L^v(\mathcal A)$, and any $b_k>0$,
\begin{equation}\label{stein:positive-sum}
 \Big\|\sum_{k\in E}R_{b_k}*A_k\Big\|_{L^v(\mathcal A)}\lesssim_v \Big\|\sum_{k\in E} A_k\Big\|_{L^v(\mathcal A)} .
\end{equation}
Test against positive $Z\in L^{v'}(\mathcal A)$ of unit norm, where $v'=v/(v-1)>1$. For $R>0$, write $B_R Z=|B(0,R)|^{-1}\int_{|y|\le R} Z(\cdot-y)\,dy$ for the centered ball average. On $|y|\le b^{-1}$ the kernel is bounded by $b^2$. On $2^{n-1}<b|y|\le2^n$, $n\ge1$, it is bounded by $Cb^2 2^{-4n}$. Multiplying by the area of $B(0,2^n/b)$ gives the weight $C2^{-2n}$ for that average. The noncommutative Hardy--Littlewood maximal inequality for balls in $\R^2$, used for the low-frequency means in~\eqref{max:low} and in~\cite[proof of Lemma~3.2]{hlw}, supplies a single positive $Z_0\in L^{v'}(\mathcal A)$ with $\|Z_0\|_{L^{v'}(\mathcal A)}\lesssim_{v'}1$ and $B_R Z\le Z_0$ for all $R>0$. Thus
\begin{equation*}
 R_b*Z\lesssim\sum_{n\ge0}2^{-2n}B_{2^n/b}Z\lesssim Z_0
 \qquad(b>0).
\end{equation*}
By evenness of $R_b$ and H\"older's inequality for the trace,
\begin{equation*}
 \sum_k\Phi((R_{b_k}*A_k)Z)=\sum_k\Phi(A_k(R_{b_k}*Z))
 \lesssim \Phi\big((\textstyle\sum_k A_k)Z_0\big)\lesssim_{v'}\Big\|\sum_k A_k\Big\|_{L^v(\mathcal A)} .
\end{equation*}
Taking the supremum proves~\eqref{stein:positive-sum}. More generally, if $0\le X\le Y$ are positive $L^v$ elements, then $\Phi(XZ)\le\Phi(YZ)$ for every positive dual test $Z$, since $\Phi((Y-X)Z)\ge0$. Thus $\|X\|_{L^v(\mathcal A)}\le\|Y\|_{L^v(\mathcal A)}$ by the same positive duality. For $2<q<\infty$, the choice $v=q/2$ gives $v'=q/(q-2)>1$. Applying the norm inequality to $\sum_k H_k\lesssim\sum_k R_{\delta2^k}*|F_k|^2$ and using~\eqref{stein:positive-sum} yields the column estimate in~\eqref{stein:coarse-bound} after taking square roots. Applying the same estimate to $F_k^*$ and $\overline\psi$ proves the row estimate. The annular approximation~\eqref{pre:localized-core} extends both estimates to every family $(F_k)_{k\in E}$ satisfying the hypotheses of Lemma~\ref{stein:coarse}.
\end{proof}

\subsection{Square-function estimates with smooth amplitudes}\label{stein:amplitude}
The derivative symbol near the disc boundary is a product of the thin-shell profile and a smooth amplitude in the normalized radius and frequency. For the shell family we interpolate the $L^2$ bound from~\eqref{stein:planch} with the $L^4$ estimate~\eqref{stein:glue4} when $2<p<4$. When $p>4$, we interpolate~\eqref{stein:glue4} with the $L^q$ estimate~\eqref{stein:coarse-bound} at a sufficiently large finite $q>p$. We then expand the amplitude in a Fourier series. The kernel bound~\eqref{stein:positive-sum} for $p>2$ and Plancherel's identity in the spatial variables for $p=2$ control the frequency factors of the series uniformly in the dyadic scales.

Let $I=[1,2]$ and take an amplitude $A\in C_c^\infty(\R\times(0,\infty))$ supported in a fixed compact rectangle.
\begin{proposition}\label{stein:annuli}
For $2\le p<\infty$, $\varepsilon>0$, and $0<\delta\le1/2$, put
\begin{equation*}
 b_{p,\varepsilon}(\delta)=
 \begin{cases}
 \delta^{1/2}L_\delta^{3/2-3/p},& 2\le p\le4,\\
 \delta^{2/p-\varepsilon},&4<p<\infty.
 \end{cases}
\end{equation*}
For every $F\in L^p(\mathcal A)$ and every finite $E\subset\mathbb Z$,
\begin{align}
 \big\|\big(A(s,2^{-k}|D|)V_{k,s} F\big)_{k\in E,s\in I}\big\|_{p,rc}
 &\lesssim_{p,\varepsilon} \|A\|_{C^{10}}B_\psi\, b_{p,\varepsilon}(\delta)\|F\|_{L^p(\mathcal A)},\label{stein:amplified}\\
 \big\|\big(A(s,2^{-k}|D|)F\big)_{k\in E,s\in I}\big\|_{p,rc}
 &\lesssim_p \|A\|_{C^{10}}\|F\|_{L^p(\mathcal A)}.\label{stein:smooth-family}
\end{align}
The constants may depend on the fixed rectangle. Both bounds hold for counting measure on $E$ and either $ds$ or $ds/s$ on $I$.
\end{proposition}
\begin{proof}
By homogeneity take $B_\psi\le1$. Choose $P_k=P(2^{-k}|D|)$ with $P\in C_c^\infty((0,\infty))$ real and equal to one throughout the radial support of $A$. Choose the annulus $[c_0,C_0]$ in Lemmas~\ref{stein:four} and~\ref{stein:coarse} to contain $\operatorname{supp}P$. We first bound the map
\begin{equation*}
 T:F\longmapsto(V_{k,s}P_kF)_{k\in E,s\in\R}
\end{equation*}
from $L^p(\mathcal A)$ to $L^p(\mathcal A;L^2(E\times\R)_d)$ for $d\in\{\mathrm c,\mathrm r\}$, with counting measure on $E$ and Lebesgue measure on $\R$. At every exponent we start from the same formula on the finite sums of Schwartz test tensors from Section~\ref{sec:prelim}.
To apply~\eqref{stein:glue4} and~\eqref{stein:coarse-bound} with $F_k=P_kF$, we first bound the row and column norms of $(P_kF)_{k\in E}$. For independent symmetric signs $\epsilon_k$, write $S_\epsilon=\sum_{k\in E}\epsilon_kP_k$. The annular overlap bounds the symbol of $S_\epsilon$ uniformly. Off the origin, summing the rescaled Schwartz kernel and its gradient gives bounds $C|x|^{-2}$ and $C|x|^{-3}$, respectively, also when $k$ ranges over all integers. Parcet's scalar-kernel Calder\'on--Zygmund theorem, applied as in Section~\ref{lp:sum-annuli}, gives a uniform bound for $S_\epsilon$ on $L^p(\mathcal A)$. Since $\mathbb E|S_\epsilon F|^2=\sum_{k\in E}|P_kF|^2$, the triangle inequality in $L^{p/2}(\mathcal A)$ gives
\begin{equation}\label{stein:blocks}
 \|(P_kF)_k\|_{p,c}^2\le \mathbb E\|S_\epsilon F\|_{L^p(\mathcal A)}^2\lesssim_p \|F\|_{L^p(\mathcal A)}^2.
\end{equation}
The same inequality applied to $F^*$ gives the row bound because $P$ is real and radial, so $P_k(F^*)=(P_kF)^*$. This is the annular Littlewood--Paley upper estimate of~\cite[Lemma~2.6]{hlw} for the chosen cutoff $P$.

For $p=2$,~\eqref{stein:planch} and unitarity give $\delta^{1/2}\|\psi\|_{L^2(\R)}(\sum_{k\in E}\|P_kF\|_{L^2(\mathcal A)}^2)^{1/2}$ for $\|TF\|_{2,d}$, $d\in\{\mathrm c,\mathrm r\}$. Scalar Plancherel and the definition of $B_\psi$ give $\|\psi\|_{L^2(\R)}\lesssim B_\psi\le1$. At $p=4$, applying~\eqref{stein:glue4} to $F_k=P_kF$ and using~\eqref{stein:blocks} gives the operator norm bound $C\delta^{1/2}L_\delta^{3/4}$. For finite $q>4$, applying~\eqref{stein:coarse-bound} and~\eqref{stein:blocks} to $(P_kF)_{k\in E}$ gives the bound $C_q$. Fix one real orthonormal basis of $L^2(E\times\R)$ for all exponents. In that basis the spaces $L^p(\mathcal A;L^2(E\times\R)_d)$ for $d=\mathrm c,\mathrm r$ are realized as the first-column and first-row subspaces, respectively, of the same tracial matrix algebra, as in Section~\ref{sec:stein}. The interpolation identity in~\cite[Lemma~2.4]{hlw} therefore applies for each $d\in\{\mathrm c,\mathrm r\}$. For $2<p<4$, solve $1/p=(1-\theta)/2+\theta/4$, giving $\theta=2-4/p$. Interpolation gives
\begin{equation*}
 \|T\|_{p\to(p,d)}\lesssim_p
 (\delta^{1/2})^{1-\theta}(\delta^{1/2}L_\delta^{3/4})^\theta
 =\delta^{1/2}L_\delta^{3/2-3/p}.
\end{equation*}
Here $\|T\|_{p\to(p,d)}$ denotes the operator norm with the parameter norm $\|\cdot\|_{p,d}$ on the range. At $p=2$ the bound follows from~\eqref{stein:planch}, and at $p=4$ it follows from~\eqref{stein:glue4} and~\eqref{stein:blocks}. At $4<p<q<\infty$ write $1/p=\theta/4+(1-\theta)/q$. The same map now satisfies
\begin{equation*}
 \theta=\frac{1/p-1/q}{1/4-1/q},\qquad
 \|T\|_{p\to(p,d)}\lesssim_{p,q}
 (\delta^{1/2}L_\delta^{3/4})^\theta\cdot 1^{1-\theta}.
\end{equation*}
Since $2/p-\theta/2=2(p-4)/(p(q-4))$, choose a finite $q>p$ large enough to make this difference less than $\varepsilon/2$. Then $\delta^{\theta/2}\le\delta^{2/p-\varepsilon/2}$ for $0<\delta\le1/2$, and $L_\delta^{3\theta/4}\lesssim_\varepsilon\delta^{-\varepsilon/2}$. These bounds give $\|T\|_{p\to(p,d)}\lesssim_{p,\varepsilon}b_{p,\varepsilon}(\delta)$, uniformly in $E$.

The Fourier series for $A$ introduces the frequency factors $M_{k,b}=e^{ib2^{-k}|D|}$, $b\in\R$, on the annuli of Lemma~\ref{stein:four}. Insert the fixed reproducing cutoff $\kappa(2^{-k}|D|)$. Applying $(1-\Delta_\xi)^2$ under Fourier inversion on the fixed annulus gives
\begin{equation*}
 |K_{k,-b}(y)|\lesssim(1+|b|)^4 R_{2^k}(y),
 \qquad \|K_{k,-b}\|_{L^1(\R^2)}\lesssim(1+|b|)^4,
\end{equation*}
since $\|R_{2^k}\|_{L^1(\R^2)}=\|R_1\|_{L^1(\R^2)}$. Apply the convolution inequality $|K*G|^2\le\|K\|_{L^1}(|K|*|G|^2)$ from Lemma~\ref{stein:coarse} with $K=K_{k,-b}$ and $G=F_k$. Since $\kappa(2^{-k}|D|)F_k=F_k$, this gives
\begin{equation*}
 |M_{k,b}F_k|^2\lesssim (1+|b|)^8R_{2^k}*|F_k|^2.
\end{equation*}
Use~\eqref{stein:positive-sum} at $v=p/2$ when $p>2$, and Plancherel's identity in the spatial variables separately for each $k$ when $p=2$. For rows one applies the column bound to $F_k^*$ with phase $-b$, since $(M_{k,b}F_k)^*=M_{k,-b}(F_k^*)$. Thus
$\|(M_{k,b}F_k)_k\|_{p,d}\lesssim_p(1+|b|)^4\|(F_k)_k\|_{p,d}$ for $d\in\{\mathrm c,\mathrm r\}$.

Expand $A$ periodically on a rectangle whose interior contains $I\times\operatorname{supp}P$ and $\operatorname{supp}A$. Its Fourier coefficients $c_{mn}$, $m,n\in\mathbb Z$, satisfy $|c_{mn}|\lesssim\|A\|_{C^{10}}(1+|m|+|n|)^{-10}$ by integration by parts, with frequencies $a_m,b_n$ fixed multiples of $m,n$, respectively. Since $P(\rho)=1$ on the radial support of $A$, Fourier inversion on the rectangle gives, for $s\in I$ and all $\rho>0$,
\begin{equation*}
 A(s,\rho)=P(\rho)\sum_{m,n\in\mathbb Z}c_{mn}e^{ia_ms}e^{ib_n\rho}.
\end{equation*}
Outside $\operatorname{supp}P$ both sides vanish. For each $m,n\in\mathbb Z$ the corresponding family in~\eqref{stein:amplified}, before multiplying by $c_{mn}$, is
\begin{equation*}
 (e^{ia_ms}V_{k,s}M_{k,b_n}P_kF)_{k\in E,s\in I}.
\end{equation*}
Multiplication by $e^{ia_ms}$ is an isometry in either parameter norm. Since $M_{k,b_n}$ preserves Fourier support, the endpoint estimates apply with $F_k=M_{k,b_n}P_kF$ for each $k\in E$. The sequence bound for $M_{k,b_n}$ and~\eqref{stein:blocks} give the factor $(1+|b_n|)^4$ at each endpoint. Interpolating $F\mapsto(V_{k,s}M_{k,b_n}P_kF)_{k,s}$ gives the same factor at $p$, so this family has norm at most $C_{p,\varepsilon}(1+|b_n|)^4 b_{p,\varepsilon}(\delta)\|F\|_{L^p(\mathcal A)}$. With $V_{k,s}$ omitted, the square integral in $s$ gives the factor $\|e^{ia_ms}\|_{L^2(I)}=1$ and the sequence bound gives $C_p(1+|b_n|)^4\|F\|_{L^p(\mathcal A)}$. The Fourier expansions of the families in~\eqref{stein:amplified} and~\eqref{stein:smooth-family} therefore converge absolutely in both the row and column norms, since $\sum_{m,n\in\mathbb Z}|c_{mn}|(1+|b_n|)^4\lesssim\|A\|_{C^{10}}$.

To identify the sums with the stated multiplier families, take $F\in L^2(\mathcal A)\cap L^p(\mathcal A)$. Write $A_N(s,\rho)$ for the finite Fourier sum with $|m|,|n|\le N$, $N\in\mathbb N$, multiplied by $P(\rho)$. It converges uniformly to $A(s,\rho)$ for $s\in I$ and $\rho\ge0$, using the coefficient decay on $\operatorname{supp}P$ and the zero extension elsewhere. For fixed finite $E$, Plancherel gives, for example,
\begin{align*}
 &\big\|\big((A_N-A)(s,2^{-k}|D|)V_{k,s}F\big)_{k\in E,s\in I}\big\|_{2,c}\\
 &\qquad\le (|E||I|)^{1/2}\|A_N-A\|_{L^\infty(I\times[0,\infty))}
       \|\psi\|_{L^\infty(\R)}\|F\|_{L^2(\mathcal A)}
       \longrightarrow0 .
\end{align*}
The $L^2$ column and row norms coincide with the joint Hilbert space norm. Omitting $V_{k,s}$ gives the same convergence with the factor $\|\psi\|_{L^\infty(\R)}$ replaced by $1$. In each matrix realization, uniqueness of the limit in measure identifies the $L^p$ sums with the multiplier families in~\eqref{stein:amplified} and~\eqref{stein:smooth-family}, respectively. Density of $L^2(\mathcal A)\cap L^p(\mathcal A)$ and restoration of $B_\psi$ give~\eqref{stein:amplified} and~\eqref{stein:smooth-family}. Finally, $ds/s\le ds$ on $I$, so each positive square integral against $ds/s$ is dominated by the same ordered integral against $ds$. Monotonicity of the $L^{p/2}$ norm gives~\eqref{stein:amplified} and~\eqref{stein:smooth-family} with $ds/s$.
\end{proof}

\subsection{The full-radius square-function estimate}\label{stein:summation}

We now complete the proof of Theorem~\ref{thm:stein} by summing the estimates for the interior and boundary pieces of \(D_r^\nu\) in the row and column norms.

The factor \((1-|\xi|^2/r^2)^{\nu-1}\) in
\eqref{stein:derivative} has size \(\delta^{\nu-1}\) on a boundary
piece of relative thickness \(\delta\). Multiplying this weight by
\(b_{p,\varepsilon}(\delta)\) from Proposition~\ref{stein:annuli}
gives the summability condition in the following theorem.

\begin{theorem}\label{stein:theorem}
Let \(2\le p<\infty\), let \(\nu>\nu_*(p):=\max\{1/2,1-2/p\}\) be real, and take \(H=L^2((0,\infty),dr/r)\). The restriction of the Hilbert space operator \(\mathcal T_\nu\) to \(L^2(\N_{\R^2})\cap L^p(\N_{\R^2})\) extends uniquely to a bounded linear map
\[
 \mathcal T_\nu:L^p(\N_{\R^2})\longrightarrow
       L^p(\N_{\R^2};H_c)\cap L^p(\N_{\R^2};H_r).
\]
With the maximum norm on the intersection, the estimate is equivalently
\[
 \|\mathcal T_\nu F\|_{p,rc}
   =\max\{\|G_c^\nu F\|_{L^p(\N_{\R^2})},\|G_r^\nu F\|_{L^p(\N_{\R^2})}\}
   \le C_{p,\nu}\|F\|_{L^p(\N_{\R^2})} \qquad (F\in L^p(\N_{\R^2})).
\]
The constant is uniform over all semifinite tracial coefficient algebras \(\M\).
\end{theorem}

\begin{proof}
Fix \(F\in L^2(\N_{\R^2})\cap L^p(\N_{\R^2})\). To bound \(\|\mathcal T_\nu F\|_{p,rc}\), we first estimate finite multiplier sums on finitely many dyadic intervals of radii. The interior symbol vanishes quadratically at the origin, which lets us sum its pieces on annuli shrinking to zero. The boundary terms are indexed by their distance to \(|\xi|=r\) and are estimated by Proposition~\ref{stein:annuli}. We then recover \(\mathcal T_\nu F\) in \(L^2\) and use the uniform bounds to pass to all radii.

Write \(r=2^k s\), \(k\in\mathbb Z\), \(1\le s<2\), and put \(\rho=2^{-k}|\xi|\). On the interval \([2^k,2^{k+1})\) the change \(r=2^ks\) gives \(dr=2^k ds\), so for any nonnegative measurable function \(h:(0,\infty)\to[0,\infty]\)
\[
 \int_0^\infty h(r)\,\frac{dr}{r}
  =\sum_{k\in\mathbb Z}\int_1^2 h(2^ks)\,\frac{ds}{s}.
\]
Taking \(h=|g|^2\) for \(g\in H\) identifies the parameter Hilbert space isometrically with the square sum of copies of \(L^2([1,2),ds/s)\). Fix a finite set \(E\subset\mathbb Z\) and put \(S_E=\bigcup_{k\in E}[2^k,2^{k+1})\). Write \(v=s-\rho\). On \(v>0\), substituting \(|\xi|/r=\rho/s\) into~\eqref{stein:derivative} and factoring \(1-\rho^2/s^2=(s+\rho)v/s^2\) gives the normalized symbol
\[
 d(s,\rho)=A_\nu(s,\rho)v^{\nu-1},\qquad
 A_\nu(s,\rho)=2\nu\rho^2s^{-2\nu}(s+\rho)^{\nu-1}.
\]
Set \(d(s,\rho)=0\) for \(v\le0\), so \(d_r^\nu(\xi)=d(s,\rho)\) at \(r=2^k s\), \(\rho=2^{-k}|\xi|\). Choose a nonnegative decreasing smooth cutoff \(\chi\) on \([0,\infty)\), equal to one on \([0,1]\) and zero on \([2,\infty)\). Put \(\eta(w)=\chi(w)-\chi(2w)\) for \(w>0\) and extend \(\eta\) by zero for \(w\le0\). Its support is contained in \([1/2,2]\). For \(v>0\) and any integer \(J\ge3\), the sum \(\sum_{j=3}^J\eta(2^jv)=\chi(8v)-\chi(2^{J+1}v)\) increases to \(\chi(8v)\). The interior symbol is \(b_\nu(s,\rho)=(1-\chi(8v))d(s,\rho)\) for \(v>0\) and zero for \(v\le0\).

The support of \(b_\nu\) can reach \(\rho=0\), so we use a dyadic partition to keep each interior amplitude away from zero in its rescaled radial variable. Choose a nonnegative \(\zeta\in C_c^\infty((1/2,2))\) with \(\sum_{\ell\in\mathbb Z}\zeta(2^{-\ell}\rho)=1\) for \(\rho>0\). Since \(b_\nu(s,\rho)\) vanishes for \(\rho\ge2\) when \(s\in[1,2)\), only \(\ell\le1\) occur. For integers \(L\le1\) and \(J\ge3\), and \(r=2^ks\in S_E\), put
\begin{equation*}
 m_{r,L,J}(\xi)
   = b_\nu(s,\rho)\sum_{\ell=L}^{1}\zeta(2^{-\ell}\rho)
      +d(s,\rho)\sum_{j=3}^{J}\eta(2^jv).
\end{equation*}
Write \(I_{\ell,E}F\) for the family with symbol \(b_\nu(s,\rho)\zeta(2^{-\ell}\rho)\) and \(S_{j,E}F\) for the family with symbol \(d(s,\rho)\eta(2^jv)\), indexed by \((k,s)\in E\times[1,2)\). The triangle inequality for the finite sum gives
\begin{equation}\label{stein:split}
 \|(m_{r,L,J}(D)F)_{r\in S_E}\|_{p,rc}
 \le \sum_{\ell=L}^{1}\|I_{\ell,E}F\|_{p,rc}
       +\sum_{j=3}^{J}\|S_{j,E}F\|_{p,rc}.
\end{equation}
Here both descriptions of the family use the measure \(dr/r=ds/s\) on each dyadic interval. We will show that the two series of norms on the right are summable uniformly in \(E\).

\medskip
We first bound the interior terms in~\eqref{stein:split}. For \(\ell\le1\) and \(w>0\), extract the factor \(2^{2\ell}\) by defining
\[
 b_{\nu,\ell}(s,w)=2^{-2\ell} b_\nu(s,2^\ell w)\zeta(w).
\]
To check the bound after extracting \(2^{2\ell}\), put \(q_\nu(s,\rho)=b_\nu(s,\rho)/\rho^2\) for \(\rho>0\). On \(s-\rho>0\) this quotient equals
\[
 [1-\chi(8(s-\rho))]\,2\nu s^{-2}(1-\rho^2/s^2)^{\nu-1}.
\]
It extends smoothly to \(\rho=0\). The factor \(1-\chi(8(s-\rho))\) vanishes for \(0<s-\rho\le1/8\). On the remaining region \(s-\rho>1/8\), the base \(1-\rho^2/s^2\) is bounded away from zero for \(s\) in a fixed neighborhood of \([1,2]\). The zero extension across \(s=\rho\) is smooth. Thus \(b_{\nu,\ell}(s,w)=w^2 q_\nu(s,2^\ell w)\zeta(w)\). Each derivative on the second argument of \(q_\nu\) contributes a factor \(2^\ell\le2\), and \(w\) stays in the fixed support of \(\zeta\). Extend the \(s\)-amplitudes a short distance beyond \([1,2]\) by the same formulas and multiply by a fixed smooth cutoff equal to one on \([1,2]\). This gives uniform \(C^{10}\) bounds and support in a fixed compact rectangle in \(\R\times(0,\infty)\).

Since \(2^{-\ell}\rho=2^{-(k+\ell)}|\xi|\), apply~\eqref{stein:smooth-family} to the finite set \(E+\ell\) of frequency scales. Relabeling this set by \(E\) preserves the parameter norms and gives
\[
 \|I_{\ell,E}F\|_{p,rc}
 = 2^{2\ell}\|(b_{\nu,\ell}(s,2^{-k-\ell}|D|)F)_{k\in E,s}\|_{p,rc}
 \le C_{p,\nu}2^{2\ell}\|F\|_{L^p(\N_{\R^2})} .
\]
Summing over \(\ell\le1\) controls the entire interior series in~\eqref{stein:split}.

\medskip
For the boundary terms in~\eqref{stein:split}, we reduce to the thin annuli of Proposition~\ref{stein:annuli}. The support of every boundary piece for \(s\in[1,2]\) is contained in the fixed closed rectangle \([1,2]\times[3/4,2]\), since \(\eta(2^jv)\ne0\) with \(j\ge3\) implies \(0<v\le1/4\). Choose \(\theta\in C_c^\infty((0,\infty)^2)\) equal to one on this rectangle and replace \(A_\nu\) by \(\theta A_\nu\) in the boundary symbol. For \(j\ge3\) put \(\delta=2^{-j}\) and set \(\psi_\nu(w)=w^{\nu-1}\eta(w)\) for \(w>0\), extended by zero. The support condition \(\operatorname{supp}\eta\subset[1/2,2]\) makes \(\psi_\nu\) a smooth compactly supported function with \(B_{\psi_\nu}\) independent of \(j\). The symbol of \(S_{j,E}\) is
\[
 \delta^{\nu-1} A_\nu(s,\rho)\psi_\nu((s-\rho)/\delta).
\]
Apply~\eqref{stein:amplified} with \(\psi=\psi_\nu\) and \(\delta=2^{-j}\) to the family indexed by \(k\in E\). For \(2\le p\le4\), summing in \(j\) bounds \(\sum_{j\ge3}\|S_{j,E}F\|_{p,rc}\) by
\[
 C_{p,\nu}\|F\|_{L^p(\N_{\R^2})}\sum_{j\ge3}2^{-j(\nu-1/2)}(1+j)^{3/2-3/p}.
\]
For \(p>4\), take \(0<\varepsilon<\nu-1+2/p\). In this case the two factors in the \(j\)-th estimate are the boundary weight \(\delta^{\nu-1}\) and the square function gain \(\delta^{2/p-\varepsilon}\), so
\[
 \sum_{j\ge3}\|S_{j,E}F\|_{p,rc}
 \le C_{p,\nu,\varepsilon}\|F\|_{L^p(\N_{\R^2})}
           \sum_{j\ge3}2^{-j(\nu-1+2/p-\varepsilon)} .
\]
The exponent \(\nu-1+2/p-\varepsilon\) is positive. For \(p\le4\), the positive exponent \(\nu-1/2\) likewise gives geometric decay, even after multiplication by \((1+j)^{3/2-3/p}\). Together the interior and boundary estimates bound the right side of~\eqref{stein:split} by \(C_{p,\nu}\|F\|_{L^p(\N_{\R^2})}\), independently of \(E,L,J\). For each finite \(E\subset\mathbb Z\), both \(\sum_{\ell\le1}I_{\ell,E}F\) and \(\sum_{j\ge3}S_{j,E}F\) converge absolutely in the column and row norms.

\medskip
The summable bounds for the two series in~\eqref{stein:split} let us recover the family $(D_r^\nu F)_{r>0}$. Keep \(E\) fixed. For \(v=s-\rho>0\), the radial partition and the telescoping identity give, for integers \(L\le1\) and \(J\ge3\),
\begin{equation*}
 \begin{aligned}
 0&\le (1-\chi(8v))\sum_{\ell=L}^1\zeta(2^{-\ell}\rho)\le 1-\chi(8v),\\
 0&\le \sum_{j=3}^J\eta(2^jv)=\chi(8v)-\chi(2^{J+1}v)\le\chi(8v)\,.
 \end{aligned}
\end{equation*} Hence \(0\le m_{r,L,J}(\xi)\le d_r^\nu(\xi)\). Both sides are zero when \(v\le0\) or \(\rho=0\). For \(\rho>0\) and \(v>0\), the radial partition gives \(\sum_{\ell=L}^{1}\zeta(2^{-\ell}\rho)\to1\) on the support of \(b_\nu\) as \(L\to-\infty\), and telescoping gives \(\sum_{j=3}^J\eta(2^jv)\to\chi(8v)\) as \(J\to\infty\). It follows that \(m_{r,L,J}(\xi)\to d_r^\nu(\xi)\). In the joint Hilbert space on \(S_E\), the squared error is the integral of
\[
 (2\pi)^{-2}|m_{r,L,J}(\xi)-d_r^\nu(\xi)|^2\,
            \|\widehat F(\xi)\|_{L^2(\M)}^2
 \quad\text{over }\R^2\times S_E\text{ against }d\xi\,dr/r .
\]
Since \(|m_{r,L,J}(\xi)-d_r^\nu(\xi)|\le d_r^\nu(\xi)\), equation~\eqref{stein:mass} supplies an integrable majorant. Dominated convergence therefore shows that the squared error tends to zero. In either matrix realization of the parameter family, the families \((m_{r,L,J}(D)F)_{r\in S_E}\) have both an \(L^2\) and an \(L^p\) limit as \(L\to-\infty\) and \(J\to\infty\). Convergence in either trace norm implies convergence in measure, which has unique limits. Thus the row and column sums on \(S_E\) both represent \((D_r^\nu F)_{r\in S_E}\), with norm bounded by \(C_{p,\nu}\|F\|_{L^p(\N_{\R^2})}\).

To remove the restriction on the radii, let \(E_n=\{-n,\ldots,n\}\) for \(n\ge1\) and regard each family on \(S_{E_n}\) as a family on \((0,\infty)\) by zero extension. By~\eqref{stein:mass} these families converge in \(L^2\) to the full derivative family. They have \(L^p\) norm at most \(C_{p,\nu}\|F\|_{L^p(\N_{\R^2})}\) in both matrix realizations, uniformly in \(n\). In either realization, convergence in trace \(L^2\) implies convergence in measure, so Fatou's inequality gives the same \(L^p\) bound for the limit~\cite[Lemma~3.4(i)]{fk}. The estimate therefore holds for the Hilbert space derivative of every \(F\in L^2\cap L^p\). Density of this intersection in \(L^p\) gives the unique extension in the maximum of the column and row norms.
\end{proof}

\begin{remark}\label{stein:constant}
At \(p=4\), one can take \(C_{4,1/2+\varepsilon}\lesssim\varepsilon^{-7/4}\) for \(0<\varepsilon\le1/2\). Indeed, for \(\nu\in[1/2,1]\) the differentiated factors \(w^{\nu-1}\) on \(\operatorname{supp}\eta\), \(s^{-2\nu}(s+\rho)^{\nu-1}\) on \(\operatorname{supp}\theta\), and \((1-\rho^2/s^2)^{\nu-1}\) on the support of \(b_\nu\) for \(1\le s\le2\) have uniformly bounded derivatives of every fixed order. This follows by differentiating the powers on compact sets where their bases have a positive lower bound. Thus the constants in the estimates for \(I_{\ell,E}F\) and \(S_{j,E}F\) at \(p=4\) are uniform in \(\nu\). For the remaining shell sum, comparison on unit intervals and the substitution \(u=\varepsilon t\) give
\begin{equation*}
 \sum_{j\ge3}2^{-j\varepsilon}(1+j)^{3/4}
 \lesssim\int_0^\infty 2^{-\varepsilon t}(1+t)^{3/4}\,dt
 =\varepsilon^{-7/4}\int_0^\infty 2^{-u}(\varepsilon+u)^{3/4}\,du .
\end{equation*}
The last integral is bounded uniformly for \(0<\varepsilon\le1/2\).
\end{remark}

\subsection{Quantum tori and Moyal planes}\label{stein:qsection}

The proof of Theorem~\ref{stein:theorem} bounds each smooth part of
the radius derivative uniformly over the tracial coefficient algebra.
We use these bounds on the orbits of the translations from
Section~\ref{q:models} to estimate the quantum derivative. Recall that
\(\mathcal D=\mathcal C\bar\otimes\mathcal Q_\Theta\) has the product
trace and \(\beta_z=\mathrm{id}_{\mathcal C}\otimes\gamma_z\) for
\(z\in\R^2\), where \(\mathcal C\) is a semifinite tracial von Neumann
algebra.

Let \((\Lambda_\Theta,\mu_\Theta)\) denote \(\mathbb Z^2\) with counting
measure on the torus and \(\R^2\) with Lebesgue measure on the plane.
For real \(\nu>1/2\) and \(x\in L^2(\mathcal D)\), let
\(\mathcal T_{\nu,\Theta}x=(D_{r,\Theta}^\nu x)_{r>0}\) be the inverse
quantum Fourier transform in \(\xi\) of
\((r,\xi)\mapsto d_r^\nu(\xi)\widehat x_\Theta(\xi)\), where
\((r,\xi)\in(0,\infty)\times\Lambda_\Theta\) and \(d_r^\nu\) is
given by~\eqref{stein:derivative}. By~\eqref{stein:mass} and
\(d_r^\nu(0)=0\), this defines a linear map from \(L^2(\mathcal D)\)
to \(L^2((0,\infty),dr/r;L^2(\mathcal D))\).

\begin{corollary}\label{stein:quantum}
For either quantum model, let \(2\le p<\infty\) and let \(\nu>\nu_*(p)\)
be real. With \(H=L^2((0,\infty),dr/r)\), the restriction of
\(\mathcal T_{\nu,\Theta}\) to \(L^2(\mathcal D)\cap L^p(\mathcal D)\)
extends uniquely to a bounded linear map
\[
 \mathcal T_{\nu,\Theta}:L^p(\mathcal D)\longrightarrow
     L^p(\mathcal D;H_c)\cap L^p(\mathcal D;H_r).
\]
Write \(G_{c,\Theta}^\nu x\) and \(G_{r,\Theta}^\nu x\) for the positive
square roots of the column and row integrals of this family,
respectively. For every \(x\in L^p(\mathcal D)\),
\[
 \|\mathcal T_{\nu,\Theta}x\|_{p,rc}
 =\max\{\|G_{c,\Theta}^\nu x\|_{L^p(\mathcal D)},\|G_{r,\Theta}^\nu x\|_{L^p(\mathcal D)}\}
 \le C_{p,\nu}\|x\|_{L^p(\mathcal D)},
\]
where \(C_{p,\nu}\) is independent of \(\Theta\) and \(\mathcal C\) and
satisfies the bound in Remark~\ref{stein:constant} at \(p=4\).

For every \(x\in L^2(\mathcal D)\), the squared norms are
\[
 \|G_{c,\Theta}^\nu x\|_{L^2(\mathcal D)}^2=\|G_{r,\Theta}^\nu x\|_{L^2(\mathcal D)}^2
 =\frac{\nu}{2\nu-1}
 \begin{cases}
   \|x-\mathbb E_0x\|_{L^2(\mathcal D)}^2,&\mathcal Q_\Theta=\T^2_\Theta,\\
   \|x\|_{L^2(\mathcal D)}^2,&\mathcal Q_\Theta=\R^2_\Theta .
 \end{cases}
\]
Here \(\mathbb E_0=\mathrm{id}_{\mathcal C}\otimes\operatorname{tr}_\Theta\)
is the torus expectation onto \(\mathcal C\otimes 1\) that extracts the
Fourier coefficient at \(n=0\).
\end{corollary}

\begin{proof}
We first compute the \(L^2\) norms. We then transfer the bounds for
the smooth symbols in \(m_{r,L,J}\) on finitely many dyadic intervals
of radii and identify their sum with \(\mathcal T_{\nu,\Theta}x\) on
\(L^2(\mathcal D)\cap L^p(\mathcal D)\).

The trace of \(a^*a\) equals that of \(aa^*\) for \(a\in L^2(\mathcal D)\),
so the two squared norms at \(p=2\) both equal the joint squared norm
\[
 \int_0^\infty\|D_{r,\Theta}^\nu x\|_{L^2(\mathcal D)}^2\,\frac{dr}{r}
 =\int_{\Lambda_\Theta}\|\widehat x_\Theta(\xi)\|_{L^2(\mathcal C)}^2
       \left(\int_0^\infty|d_r^\nu(\xi)|^2\,\frac{dr}{r}\right)
       d\mu_\Theta(\xi)
\]
by Plancherel's identity and Tonelli's theorem. The inner integral is
\(\nu/(2\nu-1)\) at every \(\xi\ne0\) by~\eqref{stein:mass}, and is zero
at \(\xi=0\). On the torus, the sum of the squared Fourier coefficients
over \(n\in\mathbb Z^2\setminus\{0\}\) is \(\|x-\mathbb E_0x\|_{L^2(\mathcal D)}^2\).
On the plane, \(\{0\}\) has Lebesgue measure zero and the integral of
the squared coefficients is \(\|x\|_{L^2(\mathcal D)}^2\).

Fix \(2<p<\infty\), a finite set \(E\subset\mathbb Z\), and
\(S_E=\bigcup_{k\in E}[2^k,2^{k+1})\). Let \(d\mu(r)=dr/r\) and
\(H_E=L^2(S_E,d\mu)\). Fix one interior index \(\ell\le1\) or one
boundary index \(j\ge3\), and for \(r\in S_E\) let \(m_r:\R^2\to\mathbb C\)
be the symbol of \(I_{\ell,E}\) or \(S_{j,E}\) from~\eqref{stein:split},
respectively. Put
\(K_r=m_r^\vee\) and define \(Q_r:L^p(\mathcal D)\to L^p(\mathcal D)\)
by~\eqref{q:orbit} with kernel \(K_r\).
Let \(\mathcal S_{\mathcal D}=L^\infty(\R^2)\bar\otimes\mathcal D\)
with its product trace. For a semifinite tracial algebra
\(\mathcal A\), write \(\|\cdot\|_{p,d;\mathcal A}\) for the radius parameter norm of
type \(d\in\{\mathrm c,\mathrm r\}\) over \(\mathcal A\).
Let \(C_m\) be the larger of the Euclidean column and row operator
norms for the family with symbol \(m_r\), \(r\in S_E\). Thus
\[
 \|(K_r*F)_{r\in S_E}\|_{p,d;\mathcal S_{\mathcal D}}
 \le C_m\|F\|_{L^p(\mathcal S_{\mathcal D})}
 \qquad(F\in L^p(\mathcal S_{\mathcal D})).
\]
The estimates for \(I_{\ell,E}\) and \(S_{j,E}\) in the proof of Theorem~\ref{stein:theorem} give
\begin{equation}\label{stein:block-constants}
 C_m\le
 \begin{cases}
 C_{p,\nu}2^{2\ell}, & m_r \text{ the symbol of } I_{\ell,E},\quad \ell\le1,\\
 C_{p,\nu,\varepsilon}2^{-j(\nu-1)}b_{p,\varepsilon}(2^{-j}),
    & m_r \text{ the symbol of } S_{j,E},\quad j\ge3 .
 \end{cases}
\end{equation}
For \(p>4\) choose
\(0<\varepsilon<\nu-1+2/p\). For \(p\le4\), \(b_{p,\varepsilon}\)
is independent of \(\varepsilon\). Both bounds are independent of \(E\) and of the
coefficient algebra.

The symbols \(m_r\) have a common compact support and uniformly bounded
frequency derivatives through order four for this fixed choice of
\(E\) and term. Integration by parts gives
\(|K_r(y)|\le C(1+|y|)^{-4}\) for \(y\in\R^2\) and \(r\in S_E\).
The dependence of \(m_r\) on \(r\) is smooth on each dyadic interval,
so the kernels are jointly measurable. Because \(\mu(S_E)<\infty\),
the map \(y\mapsto K_\cdot(y)\) is continuous into \(H_E\) by dominated
convergence, and
\[
 \int_{\R^2}\|K_\cdot(y)\|_{H_E}\,dy<\infty .
\]
For \(h\in H_E\) and \(a\in L^p(\mathcal A)\), the elementary family
\((h(r)a)_{r\in S_E}\) has norm \(\|h\|_{H_E}\|a\|_{L^p(\mathcal A)}\) for both
product orders, first by the simple-family formulas and then by
density in \(H_E\). More generally,
\(L^2(S_E,d\mu;L^p(\mathcal A))\) maps continuously into both
parameter completions, since the triangle inequality in
\(L^{p/2}(\mathcal A)\) bounds either squared norm of a simple family
\(U\) by \(\int_{S_E}\|U_r\|_{L^p(\mathcal A)}^2\,d\mu(r)\). Fix
\(x\in L^p(\mathcal D)\). The family
\((K_r(y)\beta_{-y}x)_{r\in S_E}\) is continuous in \(y\)
in the Bochner space and in both parameter completions, by
continuity of \(K_\cdot(y)\) and of the \(L^p(\mathcal D)\) orbit.
Its integral in \(y\) exists in the Bochner space by the bound
\(\|K_\cdot(y)\|_{H_E}\|x\|_{L^p(\mathcal D)}\). Bochner Fubini identifies it with
\((Q_rx)_{r\in S_E}\), and the continuous maps into the row and column
spaces take it to the integrals of the same elementary families
in those spaces.

To compare \(Q_r\) with \(K_r*\), choose
\(\varphi\in C_c^\infty(\R^2)\) with \(\|\varphi\|_{L^p(\R^2)}=1\).
For \(T>0\) put \(\varphi_T(z)=T^{-2/p}\varphi(z/T)\) and define
\(J_T:L^p(\mathcal D)\to L^p(\mathcal S_{\mathcal D})\) by
\((J_Tx)(z)=\varphi_T(z)\beta_zx\). It is an isometry since \(\beta_z\)
preserves \(\|x\|_{L^p(\mathcal D)}\) and \(\int_{\R^2}|\varphi_T(z)|^p\,dz=1\).
For a simple family \(X=(X_r)_{r\in S_E}\),
\begin{align*}
 \int_{S_E}(J_TX_r)(z)^*(J_TX_r)(z)\,d\mu(r)
 &=|\varphi_T(z)|^2\beta_z\left(\int_{S_E}X_r^*X_r\,d\mu(r)\right),\\
 \int_{S_E}(J_TX_r)(z)(J_TX_r)(z)^*\,d\mu(r)
 &=|\varphi_T(z)|^2\beta_z\left(\int_{S_E}X_rX_r^*\,d\mu(r)\right).
\end{align*}
Taking square roots commutes with \(\beta_z\). Thus \(J_T\) also
preserves each of the two parameter norms and extends isometrically
to the two completions. On Bochner families in
\(L^2(S_E,d\mu;L^p(\mathcal D))\) the extension agrees with
\(r\mapsto J_TX_r\) by simple-family approximation.

Set \(\mathcal E_{T,r}=K_r*J_Tx-J_T(Q_rx)\) for \(r\in S_E\).
By the triangle inequality and the Euclidean estimate,
\[
 \|(Q_rx)_{r\in S_E}\|_{p,d;\mathcal D}
 =\|(J_T(Q_rx))_{r\in S_E}\|_{p,d;\mathcal S_{\mathcal D}}
 \le C_m\|x\|_{L^p(\mathcal D)}+
      \|(\mathcal E_{T,r})_{r\in S_E}\|_{p,d;\mathcal S_{\mathcal D}} .
\]
The action law gives
\[
 \mathcal E_{T,r}(z)=\int_{\R^2} K_r(y)a_{T,y}(z)\,dy,\qquad
 a_{T,y}(z)=[\varphi_T(z-y)-\varphi_T(z)]\beta_{z-y}x .
\]
The map \(y\mapsto a_{T,y}\) is continuous in
\(L^p(\mathcal S_{\mathcal D})\), since
\(a_{T,y}=(J_Tx)(\cdot-y)-J_T(\beta_{-y}x)\) and both
translations and the \(\beta\)-orbit are norm continuous.
Since \(\beta_{z-y}\) preserves the trace,
\[
 \|a_{T,y}\|_{L^p(\mathcal S_{\mathcal D})}
   =\|\varphi(\cdot-y/T)-\varphi\|_{L^p(\R^2)}\|x\|_{L^p(\mathcal D)} .
\]
The elementary families \((K_r(y)a_{T,y})_{r\in S_E}\) are
therefore continuous in \(y\) in the row and column spaces and in
\(L^2(S_E,d\mu;L^p(\mathcal S_{\mathcal D}))\). As for \(Q_rx\),
integrating the families in the Bochner space yields
\((\mathcal E_{T,r})_{r\in S_E}\). The triangle inequality in each
parameter space yields
\[
 \|(\mathcal E_{T,r})_{r\in S_E}\|_{p,d;\mathcal S_{\mathcal D}}
 \le \|x\|_{L^p(\mathcal D)}\int_{\R^2}\|K_\cdot(y)\|_{H_E}
          \|\varphi(\cdot-y/T)-\varphi\|_{L^p(\R^2)}\,dy.
\]
For each \(y\) the last factor tends to zero as \(T\to\infty\), and
it is bounded by \(2\) for all \(y,T\). Dominated convergence now proves
\(\|(Q_rx)_{r\in S_E}\|_{p,d;\mathcal D}\le C_m\|x\|_{L^p(\mathcal D)}\) for each
\(d\in\{\mathrm c,\mathrm r\}\).

Let \(Q_{r,L,J}\) denote the operator in~\eqref{q:orbit} with
kernel \(m_{r,L,J}^\vee\). Apply the estimate for \((Q_rx)_{r\in S_E}\) separately to the symbols of
\(I_{\ell,E}\) for \(L\le\ell\le1\) and \(S_{j,E}\) for \(3\le j\le J\).
Since the constants in~\eqref{stein:block-constants} are summable over
\(\ell\le1\) and \(j\ge3\), the families \((Q_{r,L,J}x)_{r\in S_E}\)
converge in the column and row norms as \(L\to-\infty\) and
\(J\to\infty\). The bounds on the two limits are uniform in \(E\)
and have the dependence on \(\nu\) from Remark~\ref{stein:constant}
at \(p=4\).

For \(x\in L^2(\mathcal D)\cap L^p(\mathcal D)\), the formula for
\(A_K\) in Section~\ref{q:secspacetime} gives
\(\widehat{Q_{r,L,J}x}_\Theta(\xi)=m_{r,L,J}(\xi)\widehat x_\Theta(\xi)\).
The integrals defining \(Q_{r,L,J}x\) in \(L^2\) and \(L^p\)
agree by taking the integral in the intersection with norm
\(\|\cdot\|_{L^2(\mathcal D)}+\|\cdot\|_{L^p(\mathcal D)}\). This also gives common simple-family
approximants in \(L^2(S_E,d\mu;L^2(\mathcal D)\cap L^p(\mathcal D))\).
In the joint Hilbert space the squared error of
\((Q_{r,L,J}x)_{r\in S_E}\) as an approximation to
\((D_{r,\Theta}^\nu x)_{r\in S_E}\) is
\[
 \int_{S_E}\int_{\Lambda_\Theta}
   |m_{r,L,J}(\xi)-d_r^\nu(\xi)|^2
   \|\widehat x_\Theta(\xi)\|_{L^2(\mathcal C)}^2
     \,d\mu_\Theta(\xi)\,\frac{dr}{r}.
\]
Here \(0\le m_{r,L,J}(\xi)\le d_r^\nu(\xi)\) and the symbols converge
pointwise to \(d_r^\nu(\xi)\) as \(L\to-\infty\) and \(J\to\infty\).
Equation~\eqref{stein:mass} bounds the integral of the majorant
\(|d_r^\nu(\xi)|^2\|\widehat x_\Theta(\xi)\|_{L^2(\mathcal C)}^2\),
so the squared error tends to zero. In each matrix realization of the
parameter space, the same partial sums converge in \(L^2\) and in
\(L^p\). Both imply convergence in measure. By uniqueness of the
limit in measure, the \(L^p\) limits are the row and column realizations of
\((D_{r,\Theta}^\nu x)_{r\in S_E}\).

Finally, for each \(n\ge1\) take \(E_n=\{-n,\ldots,n\}\) and extend
\((D_{r,\Theta}^\nu x)_{r\in S_{E_n}}\) by zero to all radii.
These families converge in the joint Hilbert space to
\(\mathcal T_{\nu,\Theta}x\) because its squared norm is integrable in
\(r\). The uniform bounds on \(S_{E_n}\) and Fatou's lemma for the
trace under convergence in measure give the claimed \(L^p\) bound
in both matrix realizations \cite{fk}. Density of
\(L^2(\mathcal D)\cap L^p(\mathcal D)\) in \(L^p(\mathcal D)\) gives
the unique extension in the maximum of the row and column norms.
\end{proof}

\medskip 

\begingroup
\renewcommand{\addcontentsline}[3]{}
\subsection*{AI statement}
The authors use ChatGPT for language polishing, LaTeX editing and exploratory mathematical discussions. 
GPT-6-Astra provides some new methods, proofs and improves our original results. 
The authors checked the suggested ideas, and take full responsibility for the conclusions of this work. 
\par
\endgroup

\bibliographystyle{plain}
\bibliography{references}

\end{document}